\documentclass [twoside,12pt] {report}
\usepackage{amssymb}   
\usepackage{mathrsfs}
\usepackage{amsmath}  
\usepackage{delarray}
\usepackage{latexsym} 
\usepackage{epsfig,psfrag} 
\usepackage{amsthm} 
\usepackage[all]{xy}
\usepackage{bbm}
\newenvironment {BEWEIS} {\noindent{\textbf{Proof.}}}{{\hspace*{\fill}
     $\Box$}}

\newtheorem {LEMMA} {Lemma} [section]

\newtheorem {DEFINITION} [LEMMA] {Definition}
\newtheorem {SATZ} [LEMMA] {Theorem}
\newtheorem {PROPOSITION} [LEMMA] {Proposition}

\newtheorem {KOROLLAR} [LEMMA] {Corollary}

\theoremstyle {definition} 
\newtheorem {example}{Example}
\newcommand {\Tor} {\operatorname{Tor}}
\newcommand {\Ext} {\operatorname{Ext}}
\newcommand {\Hom} {\operatorname{Hom}}
\newcommand {\End} {\operatorname{End}}
\newcommand {\Aut} {\operatorname{Aut}}
\newcommand {\poly} {\mathfrak X}
\newcommand {\sign} {\operatorname{sign}}
\newcommand {\qhop} {\mbox{\boldmath $h$}}
\newcommand {\qHop} {\mbox{\boldmath $H$}}
\newcommand {\qdz} {\mbox{\boldmath $d$}}
\newcommand {\brsop} {\mbox{$\mathcal D$}} 
\newcommand {\qbrst} {\mbox{\boldmath $\mathcal D$}} 
\newcommand {\qch} {\mbox{\boldmath $\theta$}}
 \newcommand {\qkos} {\mbox{\boldmath${\partial}$}}
\newcommand {\qce} {\mbox{\boldmath$\delta$}}
\newcommand {\qrep} {\mbox{\boldmath{$\operatorname{L}$}}}
\newcommand {\qPhi} {\mbox{\boldmath{$\Phi$}}}
\newcommand {\q} {\mbox{\boldmath $q$}}
\newcommand {\p} {\mbox{\boldmath $p$}}
\newcommand {\id} {\operatorname{id}}
\newcommand {\kernel} {\operatorname{ker}}
\newcommand {\ad} {\operatorname{ad}}
\newcommand {\degree} {\operatorname{deg}}
\newcommand {\Der} {\operatorname{Der}}
\newcommand {\Coder} {\operatorname{Coder}}
\newcommand {\Sym} {\operatorname{Sym}}
\newcommand {\Alt} {\operatorname{Alt}}
\newcommand {\res} {\operatorname{res}}
\newcommand {\qres} {\mbox{\boldmath{$\operatorname{res}$}}}
\newcommand {\prol} {\operatorname{prol}}
\newcommand {\ddxi}[1]{\frac{\partial }{\partial \xi^{#1}}}
\newcommand {\uddxi}[1]{\frac{\partial }{\partial \xi_{#1}}}
\newcommand {\uddix}[1]{\frac{\partial }{\partial x_{#1}}}
\newcommand {\ddAG}[2]{\frac{\partial}{\partial \xi^{(#1)}_{#2}}}
\newcommand {\ddG}[2]{\frac{\partial}{\partial \xi_{(#1)}^{#2}}}
\newcommand {\xAG}[2]{\xi^{(#1)}_{#2}}
\newcommand {\xG}[2]{\xi_{(#1)}^{#2}}
\newcommand {\half} {\frac{1}{2}}

\newcommand {\qimp} {\mbox{\boldmath{$J$}}}
\newcommand {\launu} {((\nu))}
\newcommand {\fnu} {[[\nu]]}
\newcommand {\clopp} {\cdot_{\operatorname{opp}}}
\newcommand {\stopp} {*_{\operatorname{opp}}}
\newcommand {\redind} {\operatorname{red}}
\newcommand {\annind} {{\operatorname{ann}}}
\newcommand {\diffind} {\operatorname{diff}}
\newcommand {\contind} {\operatorname{cont}}
\newcommand {\oppind} {\operatorname{opp}}
\newcommand {\ii} {\sqrt{-1}}
\newcommand {\ee} {\operatorname{e}}
\newcommand {\dext} {\mathrm d}
\newcommand {\Gmod}[1] {{#1}\mathbf{-Mod}^{\Gamma}}
\newcommand {\Zmod}[1] {{#1}\mathbf{-Mod}^{\mathbbm Z}}
\newcommand {\chain}[1] {{#1}\mathbf{\mathbf{-Compl}}}
\newcommand {\Zvect} {\Zmod{\mathbbm K}}
\newcommand {\Zzmod}[1] {{#1}\mathbf{-Mod}^{\mathbbm Z_2}}
\newcommand {\Zzvect} {\Zzmod{\mathbbm K}}
\newcommand {\brsalgo}{\mathscr A }
\newcommand {\brsalg} {\mathscr A_{\mathcal V}}
\newcommand {\ibrsalg}[1] {\mathscr A_{\mathcal V^{\le{#1}}}}

\begin{document}

\begin {titlepage}
\vspace*{1cm}
\begin{center}
\Large{\bf{Variations on Homological Reduction}}\\
\vspace{2cm}

\Large {Dissertation}\\
\large{zur Erlangung des Doktorgrades}\\
\large{der Naturwissenschaften}
\vspace{2cm}
\\\large{vorgelegt von}\\
\Large{Hans-Christian Herbig}\\
\large{aus Dresden}\\
\vspace{2cm}

\large{beim Fachbereich Mathematik\\}
  \large{der Johann Wolfgang Goethe-Universit\"at\\
in Frankfurt am Main}\\
\vspace{2cm}
\large{Frankfurt 2006}\\
(D 30)
\end{center}
\newpage 
\thispagestyle{empty}
\vspace*{170mm}
\noindent Vom Fachbereich Mathematik \\
der Johann Wolfgang Goethe -- Universit\"at \\
als Dissertation angenommen.\\

\noindent Prodekan: Prof. Dr. Klaus Johannson\\

\noindent Gutachter: Prof. Dr. Markus Pflaum (Frankfurt) und \\
\hspace*{20mm}        Prof. Dr. Martin Bordemann (Mulhouse, Frankreich)\\

\noindent Datum der Disputation: 
\end{titlepage}

\thispagestyle{empty}
\newenvironment{dedication}
   {\cleardoublepage \thispagestyle{empty} \vspace*{\stretch{1}}  
\begin{center} \em}
   {\end{center} \vspace*{\stretch{3}} \clearpage}
\begin{dedication}
\large{Meinen lieben Eltern gewidmet.}
\end{dedication}
\thispagestyle{empty} \cleardoublepage

\setcounter{page}{0}
\begin{abstract}
Die vorliegende Arbeit besch\"aftigt sich mit der BFV-Reduktion von Hamiltonschen Systemen mit erstklassigen Zwangsbedingungen im Rahmen der klassischen Hamiltonschen Mechanik und im Rahmen der Deformationsquantisierung. Besondere Aufmerksamkeit wird dabei Zwangsbedingungen zuteil, die als Nullfaser singul\"arer \"aquivarianter Impulsabbildungen entstehen. Es ist schon l\"anger bekannt, da{\ss} f\"ur Nullfasern  regul\"arer \"aquivarianter Impulsabbildungen die in der theoretischen Physik gebr\"auchliche Methode der BFV-Reduktion  zur Phasenraumreduktion nach Marsden/Weinstein \"aquivalent ist. In \cite{BHW} konnte gezeigt werden, da{\ss} in dieser Situation  die BFV-Reduktion sich auch im Rahmen der Deformationsquantisierung nat\"urlich formulieren l\"a{\ss}t und erfolgreich zur Konstruktion von Sternprodukten auf Marsden/Weinstein-Quotienten verwendet werden kann. Ein Hauptergebnis der vorliegenden Arbeit besteht  in der Verallgemeinerung der Ergebnisse aus \cite{BHW} auf den Fall singul\"arer Impulsabbildungen, deren Komponenten 1.) das Verschwindungsideal der Zwangsfl\"ache erzeugen und 2.) einen vollst\"andigen Durchschnitt bilden. Die Argumentation von \cite{BHW} wird durch Gebrauch der St\"orungslemmata aus dem Anhang \ref{HPT} systematisiert und vereinfacht. Zum Existenzbeweis von stetigen Homotopien und stetiger Fortsetzungsabbildung f\"ur die Koszulaufl\"osung werden der Zerf\"allungssatz und der Fortsetzungssatz von Bierstone und Schwarz \cite{Schwarzbier} benutzt. Au{\ss}erdem wird ein 'Jacobisches Kriterium' f\"ur die \"Uberpr\"ufung von Bedingung 2.) angegeben. Basierend auf diesem Kriterium und Techniken aus \cite{AGJ} werden die Bedingungen 1.) und 2.) an einer Reihe von Beispielen getestet. Als Korollar erh\"alt man den Beweis daf\"ur, da{\ss} es symplektisch stratifizierte R\"aume gibt, die keine Orbifaltigkeiten sind und dennoch  eine stetige Deformationsquantisierung zulassen. Ferner wird (\"ahnlich zu  \cite{Sevost}) eine konzeptionielle Erkl\"arung daf\"ur gegeben, warum im Fall vollst\"andiger Durchschnitte das Problem der Quantisierung der BRST-Ladung eine so einfache L\"osung hat. 

Bildet die Impulsabbildung eine erstklassige Zwangsbedingung, ist aber \emph{kein} vollst\"andiger Durchschnitt, dann ist es im allgemeinen nicht bekannt, wie entsprechende Quantenreduktionsresultate zu erzielen sind. Ein Hauptaugenmerk der Untersuchung wird es deshalb sein, in dieser Situation die klassische BFV-Reduktion besser zu verstehen -- nat\"urlich in der Hoffnung, Grundlagen f\"ur eine etwaige (Deformations-)Quantisierung zu liefern. Wir werden feststellen, da{\ss} es zwei Gr\"unde gibt, die Tate-Erzeuger (alias: Antigeister h\"oheren Niveaus) notwendig machen: die Topologie der Zwangsfl\"ache  und die Singularit\"atentheorie der Impulsabbildung. Die Zahl der Tate-Erzeuger kann durch \"Ubergang zu projektiven  Tate-Erzeugern, also Vektorb\"undeln, verringert werden. Allerdings sorgt Halperins Starrheitssatz \cite{Halperin} daf\"ur, da{\ss} im wesentlichen alle F\"alle, f\"ur die die Zwangsfl\"ache kein lokal vollst\"andiger Durchschnitt ist, zu unendlich vielen  Tate-Erzeugern  f\"uhren. Erzeugen die Komponenten einer Impulsabbildung einer linearen symplektischen Gruppenwirkung das Verschwindungsideal der Zwangsfl\"ache, so kann man eine lokal endliche Tate-Aufl\"osung  finden. Diese besitzt nach dem Fortsetzungssatz und dem Zerf\"allungssatz von Bierstone und Schwarz stetige, kontrahierende Homotopien. Ausgehend von einer solchen Tate-Aufl\"osung konstruieren wir, die klassische BFV-Konstruktion f\"ur vollst\"andige Durchschnitte verallgemeinernd, eine graduierte superkommutative Algebra. Wir k\"onnen zeigen, da{\ss} diese graduierte Algebra auch im Vektorb\"undelfall eine graduierte Poissonklammer besitzt, die sogenannte Rothstein-Poissonklammer. Die Existenz einer solchen Poissonklammer war bereits von Rothstein \cite{Rothstein} f\"ur die einfachere Situation einer symplektischen Supermannigfaltigkeit bewiesen worden. Dar\"uberhinaus werden wir sehen, da{\ss} es auch im Vektorb\"undelfall eine BRST-Ladung gibt. Diese sieht im Fall von Impulsabbildungen etwas einfacher aus als f\"ur allgemeine erstklassige Zwangsbedingungen. Insgesamt wird also die klassische BFV-Konstruktion \cite{Stashbull} auf den Fall projektiver Tate-Erzeuger verallgemeinert, und als eine Homotopie\"aquivalenz in der additiven Kategorie der Fr\'echet-R\"aume interpretiert.
\end{abstract}
\tableofcontents

\chapter{Introduction}

If we would adopt the terminology of the physicists the topic of this work
would probably most succinctly be described as `BRST-deformation quantization'.
Let us read this term backwards, and say a few words about `quantization' first. 
In physics, by quantization one loosely means a rule how to assign to a classical mechanical system a quantum mechanical system. For example, the quantum mechanical analog of the two body gravitational system (the Kepler system) is the hydrogen atom, i.e., the system of a proton and an electron subjected to electromagnetic force. In classical mechanics the states of a system are the points of the phase space (or, more generally, probability distributions on the phase space) and the observable quantities (for short observables), are the (smooth) functions on the phase space. Every observable makes the phase space into a dynamical system, the dynamics is governed by the Hamiltonian equations of motion, which is a nonlinear first order ordinary differential equation.
In quantum mechanics the states of a system are (ray equivalence classes of) vectors in a separable Hilbert space (or, more generally, density matrices). An observable is a self adjoint operator acting on this Hilbert space. The dynamics in quantum mechanics is determined by the Schr\"odinger equation, which is a linear partial differential equation. One wishes that the quantization procedure should be structurally clear, e.g., symmetry properties of the system should be preserved. This is already an interesting issue for the hydrogen atom (see \cite{Kepler}).

At this point the consensus ends and the dissent begins. About thirty years ago somebody coined the sentence `First quantization is a mystery, but second quantization is a functor!' \cite{ReedSimon}. It is the opinion of the author that the situation did not change too dramatically ever since. Up to now there is no completely satisfactory,  mathematical well-defined and universally applicable  theory 
of (`first') quantization. This difficulty is already notorious if one restricts to systems with finitely many degrees of freedom. The quantum mechanics textbook approach, which is usually called `canonical quantization', is clearly satisfactory for the working physicist, but should be viewed merely as heuristics with a fairly limited domain of applicability rather than a true theory. More than fifty years after its invention by R. Feynman, the path integral approach, which is highly used in theoretical physics and which is certainly pretty universal, still deserves full mathematical justification. Due to the tremendous success of applying path intergral techniques to deep problems in pure mathematics in the last two decades the attitude of the mathematical world to the path integral changed from brusque rejection to some sort of (neurotic) admiration.  An extensively studied, rigorous approach to the quantization problem is \emph{geometric quantization}\footnote{The number of publications related to geometric quantization is presumably already of order $10^3$.}, which goes back to works of B. Kostant and J.-M. Soriau (see for example the monograph \cite{Woodhouse}). A serious drawback of this method is that the set of observables which can be quantized is, in general, too small. One should also mention the operator algebraic approach of \emph{strict quantization} which goes back to M. Rieffel (being advocated  in the monograph \cite{Landsman}).

The approach to quantization that will be pursued in this work is that of \emph{formal deformation quantization} (the `formal' will be dropped for convenience). In deformation quantization one deliberately neglects all functional-analytic and convergence questions, and uses the gained freedom to focus on the algebraic content of the quantization problem. In this way one obtains a (not completely satisfactory) mathematically well-defined, universally applicable theory of quantization of mechanical systems with finitely many degrees of freedom. Inspired by the symbol calculus of differential operators and the deformation theory of associative algebras \cite{Gerstenhaber} the founding fathers \cite{BFFLS77,BFFLSI78,BFFLSII78} of the theory proposed to view the quantization problem as a deformation problem for the algebra of smooth functions on the phase space seen as an \emph{associative} algebra. Accordingly, the basic objects of study are so-called \emph{star products} which are associative formal deformations of the algebra of smooth functions on a Poisson manifold given by formal series of bidifferential operators and which reproduce the original Poisson structure as a semiclassical limit. The classification of star products on symplectic manifolds has been achieved by De Wilde and Lecomte \cite{dWL} and Deligne \cite{Deligne} using sheaf theoretic methods and by Fedosov \cite{Fed} using global, geometric methods. In 1997 in a preprint (which is meanwhile published \cite{Kontsevich}) Kontsevich obtained a proof of his formality theorem, which says that the differential Hochschild cochain complex of the algebra of smooth functions on a smooth manifold is $L_\infty$-quasiisomorphic to its cohomology, i.e., the $\mathbbm Z$-graded Lie algebra of polyvector fields. For the flat space the  $L_\infty$-quasiisomorphism is given by a remarkable,  explicit formula. The formality theorem  entails the classification of star products for Poisson manifolds. Ironically, the `magic' of Kontsevich's universal quantization formula appears to be natural in the light of topological quantum field theory. As it has been explained by Cattaneo and Felder \cite{CF1} the associativity of the Kontsevich star product can be understood as a Ward identity in the perturbative expansion of a certain topological sigma model. A more refined version of the formality theorem with a more conceptual method of proof has been found by Tamarkin \cite{Tamarkin} (see also the recent preprint \cite{DolgTamTsy}). Since the appearance of the Kontsevich formality theorem the theory underwent a noteworthy expansion/metamorphosis\footnote{According to the citation index of the AMS there are over 500 publications related to deformation quantization (MSC 53D55) since 1998.}, which makes it somehow difficult to give a clear picture of the current status. For a more detailed exposition and references we refer to \cite{Dito}. At this point we would like to emphasize that (even though some progress has been made in the case of orbifolds \cite{PflDQSO}), in general, the above mentioned methods do not work if the phase space has singularities. 

It is a delicate task to find out the precise number of publications\footnote{The guess $10^4$ is not too far-fetched.} which employ BRST-like methods, since they are part of the collective subconsciousness of particle physics. Roughly speaking, the idea attributed to Becchi, Rouet, Stora \cite{BRS1,BRS2,BRS3} and Tyutin \cite{Tyutin}  was to tackle the problem of gauge invariance, which makes the scattering amplitudes in the perturbative expansion of a gauge theory formally infinite, by introducing new artifical fermionic field variables, the so-called `ghosts', and to exploit a certain transformation (the BRST-symmetry) on the field variables, which leaves the action invariant and is of square zero. Still in a field theoretic spirit, in a series of papers Batalin, Fradkin and Vilkovisky \cite{BF,BV1,BV2,BV3} formulated a quite general Hamiltonian version of the theory. In the second half of the eighties it has been realized by the workers in the field that for mechanical systems of finitely many degrees of freedom, this BFV-reduction method is an incarnation of phase space reduction. For example, in the seminal article of Kostant and Sternberg \cite{KostStern} it has been explained that the reduced algebra of a regular Marsden-Weinstein reduction at zero level of the moment map can be interpreted as the zeroth cohomology of an irreducible first class  BRST-algebra. Moreover, it is shown in \cite{KostStern} that, in this case, the important problem of the quantization of the BRST-charge has a straight forward solution.  As a consequence, in the regular case the techniques of BFV-reduction have been successfully applied to quantum phase space reduction in the context of geometric quantization program, see e.g., \cite{DET}. 

Inspired by these works in \cite{BHW} the authors have been able to  show that the BFV-technique can be successfully employed to construct differential star products on phase spaces which are obtained by regular Marsden-Weinstein reduction with respect to a proper (locally) free  Hamiltonian Lie group action (a similar result has been proved before by Fedosov \cite{Fed98} using other methods). Based on techniques akin to standard  homological perturbation theory (for HPT see e.g., \cite{HuebKad} and the references therein) the reduced star product was given by a formula involving a series of differential operators, which are recursively determined. The key ingredients for the proof that this star product is in fact differential have been that 1.) in the regular case, there are explicit integral formulas for contracting homotopies of the Koszul-resolution, 2.) there is a \emph{multiplicative} prolongation map. Moreover, in certain examples the above recursion can be solved, and the reduced star product can be given by an explicit formula.

In this paper we will see that the techniques of  \cite{BHW}, suitably modified, apply also to certain cases of singular phase space reduction.
More specifically, we will be concerned with a symplectic manifold $M$ with a  Hamiltonian action of a compact, connected Lie group $G$ (with Lie algebra $\mathfrak g$), such that  the moment map $J:M\to\mathfrak g^*$ satisfies the following two conditions
\begin{enumerate}
\item \label{gh1}\emph{generating hypothesis:} the components of the moment map generate the ideal of the zero fibre $Z=J^{-1}(0)$,
\item \label{ch}\emph{complete intersection hypothesis:} the Koszul complex on $J$ is acyclic. 
\end{enumerate}
The generating hypothesis entails that the constraint set $Z$ is first class. In case the moment map satisfies condition \ref{ch}., according to the physicist's terminology, one also says  that $J$ is an \emph{irreducible constraint}. Hence, the above setup is a special case of what is called in the physics literature an \emph{irreducible first class constraint}. 
It is well-known (and not difficult to see) that these conditions are fulfilled if $0\in\mathfrak g^*$ is a regular value of $J$. However, if $0\in\mathfrak g^*$ is a \emph{singular} value of $J$ it is not at all a straight forward to check the above conditions (in particular condition 1.). In order to decide wether the generating hypothesis holds we will make use of the techniques developed in the seminal paper of Arms, Gotay and Jennings \cite{AGJ}. Accordingly, the generating hypothesis can be reduced to an elementary problem in algebraic geometry. In order to check the complete intersection hypothesis we have found a nice `Jacobian criterion' (Theorem \ref{acyckrit}), which applies if the generating hypothesis is true. In contrast to the regular case, in the singular case the continuous prolongation map and the continuous contracting homotopies of the Koszul resolution are not explicitly given, but appear as \emph{Deus ex machina} as a consequence of the extension and the division theorem of \cite{Schwarzbier}. Moreover, in the singular case the prolongation map is not multiplicative. As a result, if condition \ref{gh1}.) and \ref{ch}.) are fulfilled, we observe that the Koszul resolution on the moment map is actually a contraction in the additive category of Fr\'echet spaces (cf. Appendix \ref{HPT}). This contraction is the main ingredient for the construction of the classical BRST-algebra (cf. Section \ref{classred}). This BRST-algebra is a differential graded Poisson algebra. At the same time it contracts to the Lie algebra cohomology of the $\mathfrak g$-module $\mathcal C^\infty(Z)$. Therefore, $\mathrm H^\bullet(\mathfrak g,\mathcal C^\infty(Z))$ aquires a $\mathbbm Z$-graded Poisson algebra structure. If the generating hypothesis is true it will be explained (cf. Proposition \ref{redcomparison}) that the Poisson subalgebra $\mathrm H^0(\mathfrak g,\mathcal C^\infty(Z))=C^\infty(Z)^{\mathfrak g}$ is isomorphic to the Poisson algebra of functions of the stratified symplectic space $Z/G$ (cf. \cite{SjLerm}). Along the lines of Kostant/Sternberg \cite{KostStern} and \cite{BHW} we construct in Section \ref{brstalgconstr} a deformation quantized version of the classical BRST-algebra, the main ingredients being a star product on $M$ and a quantum moment map. To this end we follow the observation made in \cite{BHW} that it is most comfortable
to use the standard ordered (aka normal ordered) Clifford multiplication. Afterwards (cf. Section \ref{extalg}) we give a conceptual explanation (similar to \cite{Sevost}) why this construction of a quantum BRST algebra works, i.e., why the `miracle' of quantization of the BRST-charge happens. In the final Section \ref{brstcomp} we show, using the homological perturbation theory techniques of appendix \ref{HPT}, that the quantum BRST-complex contracts to a Lie algebra cohomology complex of a certain deformed respresentation. Unfortunately, the space of invariants of this representation does, in general, not coincide with the topological free module generated by the space of invariants of the classical moment map. However, we are able to show that this problem can be circumvented if either we choose a strongly invariant star product, or the Lie group $G$ is a compact, connected and \emph{semisimple}. In these cases we obtain continuous star products which deform the Poisson algebra of smooth functions on the singular reduced space $Z/G$.

In order to illustrate our reduction methods we spend considerable effort in checking the generating and complete intersection hypothesis on a number of examples (all of them appeared in the literature before), see the list of examples at the end of Subsection \ref{hamgract}. As a matter of principle, the generating hypothesis is much more easy to verify for torus actions than for nonabelian group actions (cf. Proposition \ref{nonpcond} and Theorem \ref{torus}). In particular, if a moment map of an $S^1$-action changes sign in a open neighborhood of $z$ in $M$ for every $z\in Z$, then the generating hypothesis and the complete intersection hypothesis are true. An example of such a moment map is provided by the so-called $(1,1,-1,-1)$-resonance of \cite{CushSjam}, for which it is known that the reduced space is not an orbifold, but a genuine stratified symplectic space. Based on results of  \cite{BrenPiVas} we have also found a nonabelian example, a so-called \emph{commuting variety}, for which the generating and the complete intersection hypothesis hold. Note that the results of \cite{BrenPiVas} have been generalized to moment maps of the isotropy representations of symmetric spaces of maximal rank \cite{Pany}. 

This closes the first circle of ideas which will be discussed here. The second topic of this work addresses the question what happens if the complete intersection hypothesis is dropped, i.e., if, in the physicist's terminology, $J$ is a \emph{first class reducible} constraint. A prominent example, where this happens is given by the system of one particle of zero angular momentum in dimension $\ge 3$. The original suggestion of Batalin, Fradkin and Vilkovisky \cite{BF,BV1,BV2,BV3} was to adjoin successively higher antighost variables, in order to kill the homology degree by degree. In this way one acquires a resolution of the $\mathcal C^\infty(M)$-module $\mathcal C^\infty(Z)$, which shares important structural properties with the Koszul complex. It was noted by J. Stasheff \cite{Stashbull} that this adjunction process already appeared in a work of J. Tate from 1957 \cite{Tate}. Following the customs of mathematical physics, we shall call these infinite resolutions \emph{Koszul-Tate resolutions}. The essential structural properties of Koszul-Tate resolution seem to be that 1.) it is a (semifree) super-commutative differential graded algebra resolution of $\mathcal C^\infty(Z)$ and 2.) this algebra structure is part of a super-commutative, super-cocommutative $\mathcal C^\infty(M)$-bialgebra structure. The Koszul complex over a moment map is, due to the equivariance, naturally a $\mathfrak g$-module. It is not known to the author, whether a Koszul-Tate resolution may have an analogous feature. The next step in the BFV-construction is to adjoin ghost variables, which are dual to the antighosts, and to extend, by using dual pairing, the original Poisson bracket on the base manifold to acquire a $\mathbbm Z$-graded super Poisson bracket. The main theorem of classical BFV-reduction \cite{Stashbull} is that this Poisson algebra possesses a homological Hamiltonian vector field, the so-called \emph{BRST-differential}, which is a perturbation of the Koszul-Tate differential. This differential depends on the Poisson structure of the base manifold and is not $\mathcal C^\infty(M)$-linear. The associated Hamilton function, is called the \emph{BRST-charge}. As a consequence, one obtains a differential graded Poisson algebra, the so-called \emph{BRST}-algebra, whose zeroth cohomology is isomorphic to the singular reduced Poisson algebra. The whole construction is done in a purely formal manner, i.e., it does not make use of any particular feature of the Koszul-Tate resolution. At this point we would like to mention that the Koszul-Tate-resolution is linked to the singularity theory of constraint surface (e.g. to the homotopy Lie algebra of the constraint). In particular, the number of antighosts (aka Tate generators) is bounded from below by certain homological invariants of the ring, the so-called \emph{deviations}. There is an important theorem of Halperin \cite{Halperin} (which appears to be unmentioned in the physics literature), which essentially says that if any of the deviations vanish the variety in question has to be a complete intersection. This means that in essentially all cases which force us to introduce anitghosts of higher level the adjunction procedure will not terminate. 

In the following we  will see that the  BFV-construction can also be done in the vector bundle setting. Actually, this appears to be a new result already for the case of projective Koszul resolutions. In fact, the passage from free to projective Koszul-Tate resolutions (being defined in Section \ref{projkostate}) can be used to reduce the number of generators. This is illustrated by the special case of closed submanifolds, see Section  \ref{projkos}. However, the lower bounds given by the deviations, being of local nature, govern the construction in the projective case as well. In order to construct the BRST-algebra the first nontrivial task is to find a meaningful way how to adjoin to the Tate generators (=antighosts) momenta (=ghosts). This will be explained in Section \ref{ggmf}. Besides, even though the necessity of a `completion' is noticed in \cite{Kimura} we have not been able to find in the literature a clear definition of the BRST-algebra as a space. The next step is to make the BRST-algebra into a $\mathbbm Z$-graded super-Poisson algebra. It has been shown by M.Rothstein \cite{Rothstein} that one can lift a symplectic Poisson structure from a manifold $M$ to a super-Poisson structure on a Gra{\ss}mann-algebra bundle over $M$. From this explicit formula it is easy to guess a formula for the   $\mathbbm Z$-graded super-Poisson bracket we are looking for. In order to prove that this \emph{Rothstein-bracket} is actually a $\mathbbm Z$-graded super-Poisson bracket it is most comfortable to work with the Schouten-algebra over the BRST-algebra and to use the derived bracket construction of Koszul and Kosmann-Schwarzbach (this construction will be recalled in Section \ref{derbrsection}). Actually, this computation, which is done in Section \ref{RothPBfin}, does also apply if the base manifold is a genuine Poisson manifold. Unfortunately, it is well-defined merely in the case of finitely many Tate generators. The point is that neither geometric nor the algebraic part of the bracket is a genuine bivector field (this problem is also present in the free version of the theory). The way out is to consider a slightly bigger Gerstenhaber algebra which we call the algebra of multiderivations, which contains the Schouten-algebra of polyvector fields as a subalgebra (see Section \ref{derbrsection}) and to view the Rothstein-Poisson bracket  as a derived bracket for this bigger algebra. In Section \ref{RothPBinfin} we try to convince the reader that the Rothstein bracket is a super-Poisson bracket also in the infinitely generated case. Admittedly, the argument is still merely heuristics. Next, we show that the 'Existence of Charge'-theorem holds also in the projective setup (see Theorem \ref{chargeexists}). This will be done by refining the argument of Kimura \cite{Kimura} (which has been a refinement of the argument of Stasheff \cite{Stashbull} by itself). In fact, we need not only approximations to the Tate-differential, but also approximations to the Rothstein bracket. Next we show that for linear Hamiltonian group actions we find a locally finite Koszul-Tate resolution together with continuous contracting homotopies (see Theorem \ref{existsfiniteTate}). Using the homological perturbation theory techniques of appendix \ref{HPT} we show that the BRST-algebra contracts to a complex which we call the \emph{vertical complex}, which generalizes the Lie algebra cohomology complex of the $\mathfrak{g}$-module $\mathcal C^\infty(Z)$. In this way, the vertical cohomology aquires the structure of a $\mathbbm Z$-graded super-Poisson algebra. 
 
The question of uniqueness of the BRST-algebra (the vertical complex, respectively) has yet to be settled (we conjecture that it is unique up to $P_\infty$-quasiisomorphism in the sense of \cite{relform}.) Moreover, in the special case of a projective Koszul resolution there are alternative constructions of a $\mathbbm Z$-graded Poisson structure on the vertical complex, one being the derived bracket of Theorem \ref{redderbr} the other the $P_\infty$-structure of \cite{relform}, the relations still have to be clarified. To the authors knowledge, up to now, there have been no notable attempts to find a quantization of the BRST-charge in the reducible (aka non-complete intersection) situation (a possible way could be to proceed along the lines of \cite{relform}). The main reason seems to be that the nature of the ghost and antighost variables is not well enough understood. In particular, in the reducible(=noncomplete intersection) case, the author does not know of any interesting example for which the vertical cohomology, has been computed. Another important open question is, whether the formality theorem for super-manifolds \cite{relform} generalizes to the above setup.  
\paragraph{How to read this paper.} We suggest that the reader takes first of all a look on the material of appendix \ref{HPT} since the techniques explained there will be used at several places. If the reader is primarily interested in the complete intersection case, then he or she may skip all sections of chapter \ref{BFVchapter}, except the Sections \ref{kc} and \ref{classred}. This reading strategy might also be useful for the first time reading, if the reader is not familiar with the concepts discussed here.  In Section \ref{classred} and Chapter \ref{qbrst} we use a slightly different super-Poisson structure than in the sections before.  This causes the traditional (but somehow unaesthetic) factor of $2$ in the decomposition of the BRST differential $\brsop=2\partial+\delta$. For readers who dislike this factor of $2$ there is an easy way to get rid of it: halve the odd part of the BRST-Poisson structure (equation (\ref{oddpoiss})) and replace $\theta=-\frac{1}{4}\sum_{abc} f_{ab}^c\xi^a\xi^b\xi_c+\sum_a J_a\xi^a$ with $-\frac{1}{2}\sum_{abc} f_{ab}^c\xi^a\xi^b\xi_c+\sum_a J_a\xi^a$. Analogously, cancel the $2$ in the exponent of formula (\ref{clmult}) and replace $\qch:= -\frac{1}{4}\sum_{a,b,c}\:f_{ab}^c\:\xi^a\xi^b\xi_c+ \sum_a\qimp_a\:\xi^a+\frac{\nu}{2}\sum_{a,b} f^b_{ab}\:\xi^a$ with  $-\frac{1}{2}\sum_{a,b,c}\:f_{ab}^c\:\xi^a\xi^b\xi_c+ \sum_a\qimp_a\:\xi^a+\frac{\nu}{2}\sum_{a,b} f^b_{ab}\:\xi^a$. We assume that the reader is familiar with the basics of differential geometry and homological algebra. In particular, we tacitely assume some familiarity with the cohomology of Lie algebras and the Theorem of Serre and Swan. The letter $\mathbb K$ will stand for the field of real numbers $\mathbbm R$ or the field of complex numbers $\mathbbm C$.  

\paragraph{Acknowledgments.} First of all, I would like to express my gratitude to Markus Pflaum who gave me the opportunity (and freedom) to rethink these seemingly oldfashioned concepts from a different, and in my eyes very interesting perspective. Not  least, he has been able to realize excellent working conditions and create a productive and open-minded athmosphere. Second, I would like to thank Martin Bordemann for maintaining the dialogue on mathematics (and somehow physics) also through times of stagnation and frustration. Many of the threads of this work originate from our ongoing discourse. I would like to thank Lucho Avramov and Srikanth Iyengar for encouraging my attempts to earth the BFV-construction to the solid grounds of commutative algebra. A good portion of the project has been financed by the Deutsche Forschungsgemeinschaft. I would also like to thank the Herrmann-Willkomm-Stiftung for financial support. It is taken for granted that the proof-reader is responsible for every error that has survived.

\chapter{Preparatory material}
\section{Hamiltonian reduction}
In this section we recall basic notions of classical Hamiltonian mechanics. We will exhibit a list of examples of (mostly singular) moment maps, which will serve as an illustration to the methods presented in the course this work. We will recall and compare two different notions of singular reduction: universal reduction and Dirac reduction. We will present the toolbox of Arms, Gotay and Jennings \cite{AGJ} in order to investigate the $\mathcal C^\infty$-algebraic geometry of the singular moment maps from our list.
\subsection{Hamiltonian group actions}\label{hamgract}
In order to fix notation and sign conventions, let us recall the some basic notions from Hamiltonian mechanics. Even though a considerable part of this work applies only for symplectic manifolds, let us talk for the moment, more generally, about Poisson manifolds. A Poisson manifold is a manifold $M$, which carries a \emph{Poisson tensor} $\Pi\in\Gamma^\infty(M,\wedge^2 TM)$. This Poisson tensor has to satisfy the following first order differential equation
\begin{align}\label{master}
[\Pi,\Pi]=0,
\end{align}
where the bracket $[\:,\:]:\Gamma^\infty(M,\wedge^i TM)\times\Gamma^\infty(M,\wedge^j TM)\to \Gamma^\infty(M,\wedge^{i+j-1} TM)$ is the Schouten-Nijenhuis bracket. Using local coordinates $x^1,\dots, x^n$ we write $\Pi=\half\sum_{i,j}\Pi^{ij}\frac{\partial}{\partial x^i}\wedge\frac{\partial}{\partial x^j}$ and rewrite equation (\ref{master}) as
\begin{align}
\sum_{m}\Big(\Pi^{im}\frac{\partial}{\partial x^m}\Pi^{jk}+\Pi^{jm}\frac{\partial}{\partial x^m}\Pi^{ki}+\Pi^{km}\frac{\partial}{\partial x^m}\Pi^{ij}\Big)=0.
\end{align}
From the Poisson tensor we derive a Poisson bracket $\{\:,\:\}$. This is a bilinear antisymmetric map $\mathcal C^\infty(M)\times \mathcal C^\infty(M)\to \mathcal C^\infty(M)$. For the definition we again use the Schouten-Nijenhuis bracket
\begin{align}
\{f,g\}:=-[[\Pi,f],g],
\end{align}
where $f,g \in \mathcal C^\infty(M)$. In local cordinates the Poisson bracket is given by
\begin{align}
\{f,g\}=\sum_{i<j}\Pi^{ij}\Big(\frac{\partial f}{\partial x^i}\frac{\partial g}{\partial x^j}-\frac{\partial g}{\partial x^i}\frac{\partial f}{\partial x^j}\Big).
\end{align}
As an easy consequence of the Gerstenhaber algebra identities for $[\:,\:]$ the bracket $\{\:,\:\}$ satisfies Leibniz rule in every argument. Equation (\ref{master}) implies that $\{\:,\:\}$ satisfies Jacobi identity and vice versa. Summerizing, we say that $(\mathcal C^\infty(M),\cdot,\{\:,\:\})$  satisfies the axioms of a \emph{Poisson algebra}, i.e., it is a commutative algebra with a bilinear composition $\{\:,\:\}$ such that 
\begin{enumerate}
\item $\{f,g\}=-\{g,f\}$,
\item $\{f,gh\}=\{f,g\}h+g\{f,h\}$,
\item $\{f,\{g,h\}\}+\{g,\{h,f\}\}+\{h,\{f,g\}\}=0$
\end{enumerate}
for all $f,g,h\in \mathcal C^\infty(M)$. Given such a Poisson algebra the center of the Lie algebra $(\mathcal C^\infty(M),\{\:\})$ is called the space of \emph{Casimir functions}, or, for short, of Casimirs. The space of Casimir functions will be denoted by $\mathrm H^0_\Pi(M)$. Using the Poisson structure  one may associate to every function $f\in\mathcal C^\infty(M)$ a vector field
\begin{align}
X_f:=-[\Pi,f]=\{f,\:\},
\end{align}
which is called the \emph{Hamiltonian vector field} associated to $f$. If $X=X_f$ then $f$, which is unique up to a Casimir, is sometimes called a Hamiltonian function for $X$. An easy computation yields, that this assignment is actually  a morphism of Lie algebras, i.e., we have
\begin{align}
[X_f,X_g]=X_{\{f,g\}}
\end{align}
for all $f,g \in \mathcal C^\infty(M)$. Hamiltonian vector fields commute with the Poisson bivector field. Vector fields $X$ with this property, i.e., $[X,\Pi]=0$, are called \emph{Poisson} vector fields. They are infinitisemal versions of Poisson diffeomorphism. A Poisson diffeomorphism $\Phi$ of a Poisson manifold is a diffeomorphism, whose pullback $\Phi^*$ is an automorphism of the Lie algebra $(\mathcal C^\infty(M),\{\:,\:\})$.

\emph{Symplectic manifolds} constitute an important subclass of the class of Poisson manifolds. They arise if the Poisson tensor is everywhere nondegenerate. The inverse of the Poisson tensor is the symplectic form $\omega\in\Omega^2(M)$, which is nondegenerate and is uniquely determined by the requirement
\begin{align}
i(X_f)\:\omega=df \quad\forall f\in\mathcal C^\infty(M).
\end{align} 
Equation (\ref{master}) is equivalent to $\omega$ being closed. Since $\omega$ is nondegenerate the dimension of $M$ has to be even. 

Examples of noncompact symplectic manifolds are provided by the \emph{cotangent bundle} $T^*N$ of an arbitrary smooth manifold $N$. Using coordinates $q^1,\dots,q^n$ for $N$ there is a canonical frame $dq^1,\dots,dq^n$ for $T^*N$. The bundle coordinates with respect to this frame provide the \emph{canonical coordinates} $q^1,\dots,q^n,p_1,\dots,p_n$ for $T^*N$. The symplectic form $\omega=-d\theta_0$ on $T^*N$ is given (up to a sign) by the exterior differential of the canonical one form, which is given in canonical coordinates by $\theta_0=\sum_i p_i\:dq^i$.

In many examples, the symplectic manifold $(M,\omega)$ carries in addition  an almost complex structure $I\in\Gamma ^\infty(M,\End TM)$, $I^2=-\id$, which is compatible with $\omega$, i.e., 
\begin{align}
g(X,Y):=\omega(X,IY)
\end{align}
defines a Riemannian metric on $M$, where  $X,Y\in\Gamma^\infty(M,TM)$.
In this case $M$ is called an \emph{almost K\"ahler manifold}. If the complex structure is integrable, i.e. $[I,I]_{RN}=0$, where $[\:,\:]_{RN}$ denotes the Richardson-Nijenhuis bracket, $M$ is called a \emph{K\"ahler manifold}. In particular, if $M$ is K\"ahler, it follows from the Newlander-Nirenberg theorem, that $M$ is a complex manifold and $I$ coincides with the standard complex structure. Typical examples of K\"ahler manifolds are the affine space $(\mathbbm  C^n,\omega)$ and the complex projective  space $(\mathbbm  C P^n,\omega_{FS})$. More precisely, the standard K\"ahler structure on $\mathbbm  C^n$ is given by $\omega=\frac{\ii}{2}\sum_i dz_i\wedge d\bar{z}_i$. The \emph{Fubini-Study form} $\omega_{FS}$ on  $\mathbbm  C P^n$ is given as follows. Let $\pi:\mathbbm  C^{n+1}-\{0\}\to  \mathbbm  C P^n$ be the standard projection and let $Z=(Z_0,\dots,Z_n): \mathbbm  C P^n\to\mathbbm  C^{n+1}-\{0\}$ be a holomorphic lift for $\pi$, i.e. $\pi\:Z=\id$. Then 
\begin{align*}
\omega_{FS}:=\frac{\ii}{2}\partial\bar{\partial}\ln\Big(\sum_{i=0}^n\ Z_i^2\Big)
\end{align*}
does not depend on the choice of $Z$ and defines an ($U(n+1)$-invariant) K\"ahler structure on  $\mathbbm  C P^n$.

The most basic examples of  Poisson manifolds which are not symplectic are the \emph{constant} and the \emph{linear Poisson structures}. Linear Poisson structures arise as follows. The phase space $M=\mathfrak h^*$ is the linear dual space of a finite dimensional real Lie algebra $\mathfrak h$. The Lie bracket $[\:,\:]$ can be naturally interpreted as a linear bivector field $\Pi\in\Gamma^\infty(M,\wedge^2 TM)$. Equation (\ref{master}) is equivalent to the Jacobi identity for $[\:,\:]$. More specifically, if we use  a basis $e^1,\dots,e^n$ of $\mathfrak h^*$ in order to write $x=\sum_{i=1}^n x_i e^i\in \mathfrak h^*$, then $\Pi=\half\sum_{i,j,k=1}^n f_{ij}^k\:x_k\:\frac{\partial}{\partial x_i}\wedge\frac{\partial}{\partial x_j}$. Here, $f_{ij}^k$ denote the structure constants of the Lie algebra $\mathfrak h$. It is well known that the symplectic leaves of $\mathfrak h^*$ are precisely the orbits under the coadjoint action of $H$ (a connected Lie goup with Lie algebra $\mathfrak h$) on $\mathfrak h^*$ (see e.g. \cite[section 3.1]{Vaisman}).  

A \emph{(left) action} $\Phi$ of a Lie group $G$ on a manifold $M$ is a group homomorphism from $G$ to the diffeomorphism group of $M$. We will assume that this group action is \emph{effective}, i.e. the aforementioned group homomorphism is injective. The action of an element $g$ of $G$ on a point $m\in M$ will be written $m\mapsto g.m:=\Phi_g(m)$. An action of the Lie group $G$ on $M$ induces a morphism of Lie algebras from the Lie algebra $\mathfrak g$ of $G$ into the Lie algebra of $\Gamma^\infty(M,TM)$ vector fields of $M$
\begin{align*}
 \mathfrak g\to\Gamma^\infty(M,TM),\quad X\mapsto  X_M(m):=\frac{d}{dt}_{|t=0}\Phi_{\exp (-tX)}(m).
\end{align*} 
The image $X_M$ of a vector $X\in\mathfrak g$ is called the \emph{fundamental vector field} associated to X. If $M$ is a Poisson manifold, then we are interested in group actions preserving the Poisson structure, i.e. acting by Poisson diffeomorphisms. Here, most useful is the situation where the fundamental vector fields of the action are given by Hamiltonian vector fields, such that the Hamiltonian functions can be chosen in a coherent way.

\begin{DEFINITION} An action of a Lie group $G$ on a Poisson manifold $M$ is called a \emph{Hamiltonian  action with an moment map $J: M\to \mathfrak g^*$} if the following conditions are true:
\begin{enumerate}
\item $\xi_M=X_{J(\xi)}$ for all $\xi\in\mathfrak g$.
\item $J$ is equivariant (here $G$ acts on $\mathfrak g^*$ via the coadjoint action).
\end{enumerate}
Here $J(\xi)$ denotes the function, which is obtained by evaluating $J$ on $\xi\in \mathfrak g$.   
\end{DEFINITION}
In the same way, if we are given merely a $\mathfrak g$-action on $M$, i.e. a morphism of Lie algebras $\mathfrak g\to\Gamma^\infty(M,TM)$, we say that this action is \emph{Hamiltonian}, if the two conditions above are true (clearly, in 2. we have to replace $G$-equivariance by $\mathfrak g$-equivariance). At this point, let us stipulate that, unless otherwise stated,  for a Hamiltonian action as above
\begin{center}
 the Lie group $G$ will be assumed to be \emph{connected}
\end{center} 
and, hence, there is no need to distinguish between $\mathfrak g$- and $G$-equivariance. For our purposes the most useful form of the equivariance property is
\begin{align}\label{JLiemorph}
\{J(\xi),J(\eta)\}=J([\xi,\eta])\quad\forall \xi,\eta\in \mathfrak g.
\end{align}
This means, that $J$ may equally well be considered as a Lie algebra morphism from $\mathfrak g$ to the Lie algebra $(\mathcal C^\infty(M),\{\:,\:\})$. Note that a moment map is a Poisson morphism from $M$ to $\mathfrak g^*$, where $\mathfrak g^*$ is endowed with the linear Poisson structure. For a proof of this fact and of equation (\ref{JLiemorph}) see e.g. \cite[Proposition 7.30]{Vaisman}. Obviously, we have the freedom to add to $J$ a Casimir function, which vanishes on $[\mathfrak g,\mathfrak g]$, or, in other words, a Lie algebra 1-cocyle from $Z^1\big(\mathfrak g,\mathrm H^0_{\Pi}(M)\big)$. An example of such a cocycle is given by the trace form $\chi(\xi):=\half \operatorname{trace}(\ad(\xi))$.
 
Given a Poisson action of a Lie group $G$ on a Poisson manifold $M$, there are some obstructions for the existence of a moment map (see e.g. \cite[Part III, section 7]{CdSWein}). For instance, if $M$ is symplectic it is sufficient that $\mathrm H^1(\mathfrak g,\mathbbm  R)=0=\mathrm H^2(\mathfrak g,\mathbbm  R)$, or that $M$ is compact and $\mathrm H^1_{dR}(M,\mathbbm  R)=0$, for a moment map to exist. Another class of examples of Hamiltonian actions is provided by cotangent lifted actions, which arise as follows. Every  diffeomorphism $\varphi:N\to N$  of the base manifold $N$ can be lifted to a diffeomorphism of the cotangent  bundle $T^*\varphi:T^*N\to T^*N$. In local coordinates coordinates this \emph{cotangent lift} of $\varphi$ is given by
$T^*\varphi:(q^1,\dots ,q^n,p_1,\dots p_n)\mapsto (Q^1,\dots ,Q^n,P_1,\dots ,P_n)$, where $Q^i$ is the $i$th component of $\varphi(q^1,\dots ,q^n,p_1,\dots p_n)$ and $p_i=\sum_{j=1}^n \frac{\partial Q^j}{\partial q^i} P_j$. Since such  cotangent lifts (which are also known as point transformations) preserve the canonical one-form $(T^*\varphi)^*\theta_0= \theta_0$, they are in fact Poisson diffeomorphism. They obey the composition rule $T^*(\varphi\circ\psi)=T^*\psi\circ T^*\varphi$. Accordingly, a left action $\Phi:G\times N\to N$ of the Lie group $G$ induces  $g\mapsto T^*\Phi_{g^{-1}}$ a left action of $G$ on $T^*N$ by Poisson diffeomorphisms, the so-called \emph{cotangent lifted action}. Note that according to our definition the Poisson bracket on the cotangent bundle writes $\{f,g\}=\sum_i \frac{\partial f}{\partial p_i} \frac{\partial g}{\partial q^i}- \frac{\partial f}{\partial q^i} \frac{\partial g}{\partial p_i}$, this is $-1$ times the physicist's convention.
\begin{PROPOSITION} Any cotangent lifted $G$-action is Hamiltonian. In the notation as above, a moment map $J:T^*N\to\mathfrak g^*$ is provided in canonical coordinates $(\q,\p)=(q^1,\dots,q^n,p_1\dots,p_n)$ by the formula
\begin{align}
J(\xi)(\q,\p)=-\sum_i p_i\:\xi^i_N(\q),
\end{align}
where $\xi_N=\sum_i \xi^i_N\frac{\partial}{\partial q^i}$ is the fundamental vector field for the action $\xi\in\mathfrak g$ on $N$.
\end{PROPOSITION} 
\begin{BEWEIS} See \cite[p.282--283]{AM78}.
\end{BEWEIS}
\\

A symmetry can be used to reduce the number of degrees of freedom of the phase space by taking only ``gauge invariant quantities'' into account. In Hamiltonian mechanics this idea involves a two step procedure, which is illustrated by the following reduction theorem, attributed to Marsden and Weinstein \cite{MW} and Meyer \cite{Meyer}.

\begin{SATZ}[Regular point reduction] Let $M$ be a symplectic manifold with a proper free Hamiltonian action of a Lie group $G$ with moment map $J:M\to \mathfrak g^*$. Let $\mu\in \mathfrak g^*$ be a regular value of $J$ and $G_\mu$ the isotropy group of $\mu$. Then the \emph{reduced space} $M_\mu:=J^{-1}(\mu)/G_\mu$ is a symplectic manifold with symplectic form $\omega_\mu$, which is uniquely determined by the requirement $\pi_\mu^*\omega_\mu=i_\mu^*\omega$. Here $\pi_\mu$ and $i_\mu$
\begin{align*}
\xymatrix{      &\qquad J^{-1}(\mu)\ar[dl]_{\pi_\mu}\ar@{^{(}->}[dr]^{i_\mu}\qquad &\\
          M_\mu &                                                    &M}
\end{align*} 
are the obvious projection and injection, repsectively.
\end{SATZ}
\begin{BEWEIS} See e.g. \cite[p.299--300]{AM78}.
\end{BEWEIS}
\\

There are numerous generalizations and versions of this theorem. For instance,  there is a generalization to Poisson manifolds \cite[Theorem 7.31]{Vaisman}. Instead of taking $\mu$ one may take the preimage of a whole coadjoint orbit $J^{-1}(\mathcal O_\mu)$ and divide out by the action of $G$. This approach is called \emph{orbit reduction}, which is a special case of coisotropic reduction. Regular orbit and point reduction essentially coincide \cite[chapter II, section 26]{symptech}. If one drops the freeness assumption, the reduced space will be a symplectic orbifold. For an intelligible treatment of this issue the reader may consult \cite{BatesCush}.

In the following we shall address  the case, when the regularity assumption is dropped. Even though the set of regular values is generic, the singular values $\mu \in \mathfrak g^*$ are particularly interesting. One of the reasons is that the points in the $\mu$-fibre $J^{-1}(\mu)$ tend to have larger isotropy groups. Since we will use the BRST-method, we are forced to treat merely the case of  reduction at $0\in\mathfrak g^*$. If $G$ is abelian this does not cause a restriction at all, since we are free to add a constant $\mu\in\mathfrak g^*$ to $J$. For nonabelian $G$ one has to use the \emph{shifting trick}, i.e. adjoin a coadjoint orbit $\mathcal O_\mu$. For the shifting trick in singular reduction we refer to \cite{CushSjam}.

We close this section by giving a list of examples of moment maps, which we will frequently refer to in the course of this work. Note that for all of these examples, except Examples \ref{irrflow} and \ref{linPoiss}, zero is a singular value. Furthermore, in all examples, except Example \ref{standardex}, the group $G$ is connected. The elementary Examples \ref{harmosc},\ref{allwrong} and \ref{irrflow}  and the physically interesting Example \ref{drehimpuls} will serve to illustrate the limitations of the methods presented in this article. Example \ref{standardex} will be important in connection with the normal form theorem, Theorem \ref{normcoord}. The remaining examples will turn out as instances, where our ultimate goal, i.e., a quantum phase space reduction will be achieved.
 
\begin{example}[Harmonic oscillator]\label{harmosc} Consider the $S^1=\mathbbm  R/2\pi \mathbbm  Z$-action on $\mathbbm  C$ given by $(\vartheta,z)\mapsto \ee^{\ii\vartheta} z$. The moment map for this action is $J(z)=\half |z|^2$, which has $0$ as the only critical value. 
\end{example}
\begin{example}[Free particle on the line]\label{allwrong} Here the phase space is $T^*\mathbbm  R=\mathbbm  R^2$. The $\mathbbm  R$-action which is generated by the kinetic energy of the free particle $J(q,p)=\half p^2$ is given by $(q,p)\mapsto (q+tp,p)$.  
\end{example}
\begin{example}[Cotangent lift of an irrational flow on $\mathbbm  T^2$]\label{irrflow} Let $M:=T^*\mathbbm  T^2$ be the cotangent bundle of the 2-torus $\mathbbm  T^2=(\mathbbm   R/{2\pi \mathbbm  Z})^2$. An element $t\in\mathbbm  R$ acts on $\mathbbm  T^2$ by the formula $(\vartheta_1,\vartheta_2)\mapsto (\vartheta_1+2\pi t,\vartheta_1+\alpha 2\pi t)$, where the slope parameter $\alpha\in \mathbbm  R\backslash\mathbbm Q$ is irrational. Every orbit of this $\mathbbm  R$-action is dense. A moment map for the cotangent lifted $\mathbbm  R$-action on $T^*\mathbbm  T^2$ is given by $J( \vartheta_1\vartheta_2,p_1,p_2)=2\pi(p_1+\alpha p_2)$.
\end{example}
\begin{example}[Standard example]\label{standardex}
Let $M:=\mathbbm  C^n$ with its standard K\"ahler structure. Let $G$ be any subgroup of the unitary group $U(n)$. We identify the Lie algebra $\mathfrak u(n)$ with the space of skew hermitian matrices $\xi=(\xi_{ij})$, $\bar{\xi}_{ij}=-\xi_{ji}$. The moment map for the action is given by the formula 
\begin{align}
J(\xi)=\frac{\ii}{2}\sum_{ij}\xi_{ij}\bar z_i z_j,
\end{align}
for $\xi=(\xi_{ij})\in \mathfrak g\subset \mathfrak u(n)$. Using real coordinates $z_i=x_i+\ii y_i$ we obtain
\begin{align}
J(\xi)=-\frac{1}{2}\sum_{i,j} \big(A_{ij}(x_iy_j-x_jy_i)+ S_{ij}(x_iy_j+x_jy_i)\big),
\end{align} 
where $A_{ij}=-A_{ji}$ and $S_{ij}=S_{ji}$ are the real and imaginary part of $\xi_{ij}$, respectively. In all cases, except the trivial case when $G$ is discrete, zero is a singular value of $J$.
\end{example}
\begin{example}[One particle in dimension $n$ with zero angular momentum]\label{drehimpuls} Let $M:=T^*\mathbbm  R^n$ be the cotangent bundle of the euclidean space $\mathbbm  R^n$ together with the cotangent lift of the obvious $SO(n,\mathbbm  R)$-action on $\mathbbm  R^n$. Let us write the canonical ccordinates $\q=(q^1,\dots,q^n)^t$ and $\p=(p_1,\dots,p_n)^t$ in column vector form. The cotangent lifted action is just the diagonal $SO(n)$-action on $T^*\mathbbm  R^n=\mathbbm  R^n\times\mathbbm  R^n$ ,i.e., an orthogonal matrix $O\in SO(n)$ acts according to $(\q,\p)\mapsto (O\q,O\p)$. Using the euclidian scalar product $<,>$ on $\mathbbm  R^n$ we identify $\wedge^2\mathbbm  R^n$ and $\mathfrak{so}(n)$ by letting $u\wedge v\in\wedge^2\mathbbm  R^n$ act on $w\in \mathbbm  R^n$ according to  $(u\wedge v)w:=<u,w>v-<v,w>u$. Using the invariant, positive definite scalar product $(\:,\:)$ on $\mathfrak{so}(n)$, $(A,B):=-\half\operatorname{tr}(AB)$, we identify $\mathfrak{so}(n)$ and $\mathfrak{so}(n)^*$. With these identifications understood the moment map, the so-called `angular momentum', writes 
\begin{align}
J:T^*\mathbbm  R^n\to \mathfrak{so}(n)^*\cong \wedge^2\mathbbm  R^n,\quad(\q,\p)\mapsto\q\wedge\p.
\end{align}
Of course, by identifying the cotangent lifted action on $T^*\mathbbm  R^n$ with the compexified $SO(n)$-action on $\mathbbm  C^n$ this example can be viewed as a special case of the preceding Example \ref{standardex}. For completeness, let us mention the physically important special case $n=3$. Here, accidentally $\mathbbm  R^3$ and $\mathfrak{so}(3)$ are isomorphic vector spaces. More specifically, the linear isomorphism 
\begin{align*}
(v_1,v_2,v_3)\mapsto
\left(\begin{array}{ccc}
0&-v_3&v_2\\
v_3&0&-v_1\\
-v_2&v_1&0
\end{array}
\right)
\end{align*}
is in fact an isomorphism of metric Lie algebras $\left(\mathfrak{so}(3),[\:,\:],(\:,\:)\right)$ and $\left(\mathbbm  R^3,\times,<,>\right)$, where $\times$ is the well known vector product. The angular momentum now is the $\mathbbm  R^3$-valued function $J(\q,\p)=\q\times\p$.
\end{example}
\begin{example}[Zero angular momentum for $m$ particles in the plane]\label{mpart} We consider 
the system of $m$ particles in $\mathbbm  R^2$ with zero total angular momentum. More 
precisely, the phase space is $M:=(T^*\mathbbm  R^2)^m$ and we let 
$SO(2,\mathbbm  R)\cong S^1$ act on it by lifting the diagonal action, i.e.,
\begin{eqnarray*}
  SO(2)\times M&\to& M\\
  (g,(\q_1,\p^1,\dots,\q_m,\p^m))&\mapsto& (g\q_1,g\p^1,\dots,g\q_m,g\p^m),
\end{eqnarray*} 
where $\q_i=(q^1_i,q^2_i)^t$ and $\p^i=(p_1^i,p_2^i)^t$ for $i=1,\dots,m$.
The moment map $J:M\to\mathfrak{so}(2)^*=\mathbbm  R$ is given by 
$J(\q,\p)=\sum_{i=1}^m q^1_ip_2^i-q^2_ip_1^i$.
\end{example}
\begin{example}[The `lemon']\label{lemon}
Let $S^1=\mathbbm  R/2\pi\mathbbm  Z$ act on the product $M:=\mathbbm  C P^1\times\mathbbm  C P^1$ according to the formula 
\begin{align*}
\left((z_1:z_2),(z_3:z_4)\right)\mapsto\left((\ee^{\ii\vartheta}z_1:\ee^{-\ii\vartheta}z_2),(\ee^{\ii\vartheta}z_3:\ee^{-\ii\vartheta}z_4)\right).
\end{align*}
The fix points of this action are $((1:0),(1:0))$, $((0:1),(0:1))$, $((1:0),(0:1))$  and $((0:1),(1:0))$.
A moment map for this action is
\begin{align*}
J\big((z_1:z_2),(z_3:z_4)\big)=\half\left(\frac{|z_1|^2-|z_2|^2}{|z_1|^2+|z_2|^2}+\frac{|z_3|^2-|z_4|^2}{|z_3|^2+|z_4|^2}\right).
\end{align*}
The critical set of $J$ constists of the points $-1,0$ and $1$.
\end{example}
\begin{example}[(1,1,-1,-1)-resonance] \label{stratres}Consider the $S^1$-action 
on $\mathbbm  C^4$, endowed with the standard K\"ahler structure, given by 
$(z_1,z_2,z_3,z_4)\mapsto
(\ee^{\ii\vartheta}z_1,\ee^{\ii\vartheta}z_2,
\ee^{-\ii\vartheta}z_3,\ee^{-\ii\vartheta}z_4)$. 
The moment map for the action is
\[J(z_1,z_2,z_3,z_4)=\frac{1}{2}(|z_1|^2+|z_2|^2-|z_3|^2-|z_4|^2).\]
\end{example}
\begin{example}[A $\mathbbm  T^2$-action on $\mathbbm  C^4$]\label{tzwo} The action of $\mathbbm  T^2=(\mathbbm  R/2\pi \mathbbm  Z)^2$ on $\mathbbm  C^4$ is given by,
\begin{align*} 
((\vartheta_1,\vartheta_2),(z_1,z_2,z_3,z_4))
 \mapsto (\ee^{\ii (\alpha  \vartheta_1+\beta\vartheta_2)}z_1, 
 \ee^{-\ii\vartheta_2}z_2,\ee^{\ii\vartheta_1}z_3,
 \ee^{-\ii\vartheta_2}z_4)
\end{align*} 
for $\alpha,\beta \in \mathbbm Z$. A moment map for the action is 
\begin{align*}
J:\mathbbm C^4\to \mathbbm R^2,\quad J(z_1,z_2,z_3,z_4):=
-\frac {1}{2}(\alpha|z_1|^2+|z_3|^2,\beta|z_1|^2-|z_2|^2+|z_4|^2). 
\end{align*}
\end{example}
\begin{example}[Commuting varieties]\label{commvar}
Let $S$ the space of symmetric $n\times n$-matrices with real entries. We let 
$SO(n)$ act on $S$ by conjugation and we lift this action to an action of $SO(n)$ on 
the cotangent bundle $T^*S=S\times S$. This action is Hamiltonian with the moment map
\begin{eqnarray*}
  J:S\times S&\to&\wedge^2\mathbbm R^n=\mathfrak{so}(n)^*\\
  (Q,P)&\mapsto&[Q,P].
\end{eqnarray*}
\end{example}
\begin{example} [Linear Poisson structures]\label{linPoiss}Consider $M:=\mathfrak h^*$ the dual space of an $n$-dimensional $\mathbbm  R$-Lie algebra $\mathfrak h$ with Poisson linear structure, i.e., 
\begin{align*}\Pi(x_1,\dots,x_n)=\half\sum_{i,j,k} f_{ij}^k\:x_k\:\uddix{i}\wedge\uddix{j},
\end{align*}
where $x_1,\dots,x_n$ are the linear coordinates with respect to a
chosen basis $e^1,\dots,e^n$ of $\mathfrak h^*$ and
$f_{ij}^k=<e^k,[e_i,e_j]>$ are the structure constants with respect to
the dual basis $e_1,\dots,e_n$ of $\mathfrak h$. Given a Lie
subalgebra $\mathfrak g$ of $\mathfrak h$ we obtain an $\mathfrak
h$-action on $M$ by restricting the coadjoint $\mathfrak h$-action on  $M$ to $\mathfrak g$. For any connected Lie group $G$ with Lie algebra $\mathfrak g$ this $\mathfrak g$-action integrates to a Hamiltonian $G$-action on $M$. The moment map for this action is given by restriction
$J:M=\mathfrak h^*\to \mathfrak g^*$, $\alpha\mapsto \alpha_{|\mathfrak g}$.
\end{example}
\subsection{Universal reduction}
Contrary to the regular case, in the case when zero is a singular value of the moment map, there are several possible approaches to phase space reduction (see e.g. \cite{Dirac,SnW, AGJ,Henneaux}). In the beginning of the nineties, the attempts to find the most `correct' reduction procedure  lead to the notion of universal reduction of Arms, Cushman and Gotay \cite{ACG}, who noticed that $\mathcal C^\infty(M)^G/I_Z^G$ is a good candidate for the reduced Poisson algebra. Here we shall present  a (slightly different) version $\mathcal C^\infty(M)^{\mathfrak g}/I_Z^{\mathfrak g}$ of this algebra, where we have replaced $G$-invariance by $\mathfrak g$-invariance. Due to our overall assumption that $G$ is connected, this will make no difference for us. Universal reduction does not make the other aproaches -- in particular Dirac reduction -- obsolete, but serves rather as a benchmark. 

\begin{SATZ} \label{geomRed} The space of invariant functions 
\[\mathcal C^\infty(M)^{\mathfrak g}\subset\mathcal C^\infty(M)\] is a Poisson subalgebra, which contains the ideal of invariant functions vanishing on $Z$
\[I^\mathfrak g_Z:=\{f\in\mathcal C^\infty(M)^{\mathfrak g}\:|\: f(z)=0\:\forall z\in Z\}\] as a Poisson ideal. Hence the quotient $\mathcal C^\infty(M)^{\mathfrak g}/I_Z^\mathfrak g$ is in a natural way a Poisson algebra. More precisely, if $f$ and $g\in \mathcal C^\infty(M)^{\mathfrak g}$ are representatives of $[f]$ and $[g] \in  \mathcal C^\infty(M)^{\mathfrak g}/I_Z^\mathfrak g$, respectively, then the Poisson bracket is defined as
\begin{align*}
\{[f],[g]\}:=\big[\{f,g\}\big].
\end{align*}  
If $M$ is a symplectic manifold and the action of $G$ on M is proper, then this Poisson structure is nondegenerate, i.e., there are no nontrivial Casimirs.
\end{SATZ}
\begin{BEWEIS} Since $\xi\in\mathfrak g$ acts on $\mathcal C^\infty(M)$ via $\{J(\xi),\:\}$, the claim that $\mathcal C^\infty(M)^\mathfrak g\subset \mathcal C^\infty(M)$ is a Poisson subalgebra follows immediately from the Jacobi identity and the Leibniz rule for $\{\:,\:\}$. Let $f$ be in $\mathcal C^\infty(M)^\mathfrak g$. This means, that 
\begin{align*} \{J(X),f\}=0\quad \forall X\in\mathfrak g.
\end{align*}
Hence, for all $X\in\mathfrak g$ the function $J(X)$ is constant along the integral curve of $X_f$. In particular, we conclude that the set $Z=J^{-1}(0)$ is stable under the local flow of $X_f$. Now, let $z\in Z$ and $\gamma:]-\epsilon,\epsilon[\to Z\subset M$ be an intergral curve for $X_f$ such that $\gamma(0)=z$. It follows that for every $g\in I_Z$ we have
\begin{align*}
0=\frac{d}{dt}_{|t=0}\:g(\gamma(t))=\{f,g\}(z).
\end{align*}
This implies that $I_Z^\mathfrak g\subset  \mathcal C^\infty(M)^\mathfrak g$ is a Poisson ideal. A proof for the nondegeneracy of the reduced Poisson structure can be found in \cite[section 3]{ACG}.
\end{BEWEIS}
\\

The universal reduction may be interpreted geometrically as follows. We have the commutative diagram
\begin{align*}
\xymatrix{Z=J^{-1}(0)\:\ar@{^{(}->}[r]\ar[d]&M\ar[d]\\
   M_{\redind}=Z/G\:\ar@{^{(}->}[r]^{\quad i_0}&M/G.}
\end{align*}
Since $0\in \mathfrak g^*$ may be a singular value of $J$, the zero fibre $Z$ is no longer a submanifold, but, let us say at least, a closed subset of $M$. The (not necessarily free) action of $G$ on $Z$ and  on $M$ yield quotient spaces which are singular. The map $i_0$ associates to every orbit in $Z$ the corresponding orbit in $M$. It may be shown to be a homoemorphism onto its image. If the action  of $G$ on $M$ is not too pathological (e.g. if it is proper), then  the space of invariant functions $\mathcal C^\infty(M)^G=\mathcal C^\infty(M)^\mathfrak g$  is a good substitute for the space of functions on $M/G$. The ideal of invariant functions $I^\mathfrak g_Z$ vanishing on $Z$ then may be thought of as the defining ideal of the ``subvariety'' $M_{\redind}$ of $M/G$. The above theorem says that $i_0$ is an embedding of the Poisson ``variety'' $M_{\redind}$ into the Poisson ``variety'' $M/G$. Note, that there are some subtleties in case the image of $i_0$ is not closed.

It is easy to give examples of nonproper Hamiltonian group actions, for which the space of invariant functions is too small to give a meaningful description of the quotient  spaces $M/G$ and $Z/G$, respectively. For instance, one may consider the cotangent lift of an irrational torus, Example \ref{irrflow}. Here for any value of the moment map the reduced space is $M_{\redind}=(\mathbbm  T^2/\mathbbm  R)\times \mathbbm  R$ and, since the orbit of the $\mathbbm  R$-action is dense, the reduced algebra is identified with $\mathcal C^\infty(\mathbbm  R)$ with the trivial Poisson structure. The great advantage of universal reduction is that it is always applicable and gives sensible results for proper Hamiltonian group actions.  

For symplectic manifolds with a Hamiltonian action of a compact Lie group the nature of the singular geometry of the reduced space $M_{\redind}$ has been clarified by the important singular reduction theorem of Sjamaar and Lerman \cite{SjLerm}. 

\begin{SATZ}[Sjamaar/Lerman]\label{SL}
Let $(M,\omega)$ be a symplectic manifold and $G$ be a compact Lie group acting on $M$ in a Hamiltonian fashion with moment map $J:M\to \mathfrak g^*$ and let $Z=J^{-1}(0)$. Then for every subgroup $H\subset G$ the intersection $M_{(H)}\cap Z$ and the quotient space  
\begin{align*}
(M_{\redind})_{(H)}:=(M_{(H)}\cap Z)/G,
\end{align*} 
are manifolds. Here $M_{(H)}$ is the set of points in $M$ whose isotropy group is conjugate to $H$. There exists a unique symplectic form $\omega_{(H)}$ on $(M_{\redind})_{(H)}$ such that the pullback of $\omega_{(H)}$ to  $M_{(H)}\cap Z$ coincides with the restriction of $\omega$.  The disjoint union over the conjugacy classes $(H)$ of subgroups $H\subset G$
\begin{align*}
M_{\redind}=\coprod (M_{\redind})_{(H)}
\end{align*}
is in fact a symplectic stratification. Here the set of conjugacy classes is understood  to be ordered by reverse subconjugacy. The Poisson algebra of smooth functions on $M_{\redind}$ is isomorphic  to the Poisson algebra obtained by universal reduction $\mathcal C^\infty(M)^G/I_Z^G$. 
\end{SATZ}

For a detailed discussion of symplectic stratifications we refer to \cite{SjLerm,LerMontSj} and the monographs \cite{PflaumHab,OrtRat}. Note that in the orbit type decompositions as above we allow the pieces to have components of different dimension.  Furthermore, let us mention that the symplectic pieces $(M_{\redind})_{(H)}$ can  also be obtained by regular Marsden-Weinstein reduction with respect to the Hamiltonian action of a, in general smaller, Lie group.  Theorem \ref{SL} has been generalized to singular orbit reduction of proper Hamiltonian Lie group actions \cite{BatesLerman}.

Let us look at our list of examples in the light of Thereom \ref{SL}. In the case of one particle of angular momentum zero in dimension $n\ge 2$, Example \ref{drehimpuls}, it is well known that the reduced phase space is the symplectic orbifold $\mathbbm  R^2/\mathbbm  Z_2$, where $\mathbbm  Z_2$ acts on $\mathbbm  R^2$ by $(x_1,x_2)\mapsto (-x_1,-x_2)$. The first systematic treatment of this example, appears to be \cite{BosGotay}.
For the commuting variety, Example \ref{commvar}, the reduced space has been identified in \cite{LerMontSj} as the symplectic orbifold $\mathbbm  R^n\times \mathbbm  R^n/S_n$, where the symmetric group $S_n$ is understood to act diagonally. Example \ref{lemon} (the `lemon') has been discussed at length in \cite{CushSjam}. Accordingly, the reduced spaces at the singular values $\pm 1$ are points. The reduced spaces for regular values of $J$ in the open intervals $]-1,0[$ and $]0,1[$ are $\mathbbm  C P^1$.  Finally, the reduced space at singular value $0$ (being homeomorphic to $\mathbbm  C P^1$) is in fact a symplectic orbifold, which can be pictured as two copies of a quadratic half-cone in $\mathbbm  R^3$ being glued together (whence the name). An example of a singular momentum map whose reduced space is not an orbifold is provided by the $(1,1,-1,-1)$-resonance, Example \ref{stratres}, which has been discussed in \cite[Example 2.4]{CushSjam}. As a result the reduced space is a real cone over the fibred product $S^3\times_{S^1}S^3$, where $S^1$ acts on the first half $S^3\subset \mathbbm  C^2$ by $(z_1,z_2)\mapsto (\ee^{\ii\vartheta}z_1,  \ee^{\ii\vartheta}z_2)$ and on the second half by $(z_3,z_4)\mapsto (\ee^{-\ii\vartheta}z_3,  \ee^{-\ii\vartheta}z_4)$. In \cite{CushSjam} it has been argued that this reduced space is not a rational homology manifold, and hence no orbifold.

Even though we will make no use of it, let us briefly explain how \emph{invariant theory} can be utilized to describe the geometry of singular reduced spaces. Here, we have to restrict to the case of a linear Hamiltonian action of the compact (for the moment, not necessarily connected) Lie group $G$ on a symplectic vector space $M=\mathbbm R^{2n}$. This is not a severe restriction since, due to a theorem of Gotay and Tuynman \cite{GotayTuynman}, every symplectic  Hamiltonian $G$-manifold with finitely generated homology can be equivariantly obtained by phase space reduction from flat space.  The algebra of $G$-invariant real polynomials on $M$ is finitely generated. A system $\sigma_1,\dots,\sigma_k$ generators for this algebra is called a real Hilbert basis. According to \cite{Schwarz} $\sigma_1,\dots,\sigma_k$ generate the algebra of smooth invariant functions $\mathcal C^\infty(M)^G$. The \emph{Hilbert map}
\begin{align*}
\sigma=(\sigma_1,\dots,\sigma_k):M=\mathbbm R^{2n}\to \mathbbm R^k,\quad m\mapsto  (\sigma_1(m),\dots,\sigma_k(m))
\end{align*}
seperates $G$-orbits, and hence gives rise to an injective map $\bar{\sigma}:M/G\to \mathbbm R^k$. Being the image of a real polynomial map the image of $\sigma$ is a semi-algebraic set as a consequence of the Tarski-Seidenberg theorem. For a comparison of the associated Whitney stratification and the orbit type stratification see \cite{Bierstone75,PflaumHab}. Furthermore, it can be shown that the map $\sigma$ is proper and the pullback $\sigma^*:\mathcal C^\infty(\mathbbm R^k)\to \mathcal C^\infty(M)$ is a split surjective map of Fr\'echet spaces \cite{Mather}. In many cases (but not all, see \cite{Egilsson}) there is a Poisson structure on $\mathbbm R^k$, such that $\sigma^*$ is a Poisson algebra morphism. Real Hilbert bases are known for essentially all linear examples from our list. In \cite[section 5]{Egilsson} we find a formula for the real Hilbert basis of a general linear Hamiltonian circle action. For a real Hilbert basis for Example \ref{tzwo} see \cite[Example 7.7]{AGJ}, and for a real Hilbert basis for Example \ref{commvar} see \cite[p.145]{LerMontSj}. Most favorable is the situation, when the algebra of invariant polynomials is generated by \emph{quadratic polynomials}. In this case the Hilbert map itself is a moment map of a Hamiltonian action of a certain Lie subgroup $H$ of $Sp(n,\mathbbm R)$. If, moreover, $G$ and $H$ form a \emph{reductive dual pair}. i.e., $G$ and $H$ are reductive Lie subgroups of $Sp(n,\mathbbm R)$ which centralize each other in $Sp(n,\mathbbm R)$, then the orbit-reduced spaces for the $G$-action are bijectively mapped via $\bar{\sigma}$ onto closures of coadjoint orbits in $\mathfrak h^*$. This map is compatible with the stratifications and Poisson structures in an appropriate sense (for details see \cite[Theorem 4.4]{LerMontSj}). In particular, if  $H$ is semisimple, then the reduced space at zero $M_{\redind}=J^{-1}(0)/G$ is bijectively mapped onto the closure of a nilpotent coadjoint orbit in $\mathfrak h^*$. An important example of a reductive dual pair is $O(d),Sp(m,\mathbbm R)\subset Sp(md,\mathbbm R)$, which corresponds to total angular momentum of $m$ particles in dimension $d$, and which generalizes Examples \ref{drehimpuls} and \ref{mpart} (for more details see \cite[section 5]{LerMontSj}). It is indicated in \cite{Kepler}, that the $S^1$-action of Example \ref{stratres} is linked to the Kepler problem via the reductive dual pair $U(1),SU(2,2)\subset Sp(4,\mathbbm R)$. 
\subsection{Dirac  reduction}\label{Diracsubsection}

There is a second notion of algebraic phase space reduction, which goes back to works of Dirac \cite{Dirac}, and considerably predates the universal reduction. In the following we will exclusively be concerned with the Dirac reduction in its most simple form, i.e., the first class formalism. Let us recall the requisite terminology.

Let $Z$ be  a closed subset of $M$. The \emph{vanishing ideal} of $Z$ is defined to be 
\[I_Z:=\{f\in\mathcal C^\infty(M)\:|\:f(z)=0\quad\forall z\in Z\}.\]
A function on $Z$ is called \emph{smooth} if it arises as a restriction of a smooth function on $M$. The space of smooth functions on $Z$ is denoted by $\mathcal C^\infty(Z)$. Sometimes $\mathcal C^\infty(Z)$ is also callled the space of Whitney smooth functions on $Z$ -- not to be confused with the Whitney functions. Clearly, the restriction map $\res:\mathcal C^\infty(M)\to \mathcal C^\infty(Z)$ is onto, and a morphism of $\mathbbm  K$-algebras.  The kernel of $\res$ being $I_Z$, we obtain a short exact sequence 
\[0\to I_Z\to\mathcal C^\infty (M)\stackrel{\res}{\to}\mathcal C^\infty (Z)\to 0\]
of commutative $\mathbbm  K$-algebras.
\begin{DEFINITION} 
Let $Z$ be a closed subset of a Poisson manifold $M$. A smooth function $f\in \mathcal C^\infty(M)$ is called a \emph{constraint} if its restriction to $Z$ vanishes,
 i.e.,  $f\in I_Z$. It is called \emph{first class} if $\{\:f,I_Z\:\}\subset I_Z$, otherwise it is called \emph{second class}. If $I_Z$ consists of first class functions, i.e., if it is a Poisson subalgebra of $\mathcal C^\infty(M)$, then we call the ideal $I_Z$ \emph{coisotropic} and $Z$ \emph{first class}.
\end{DEFINITION}

If $Z$ is first class, then there is a canonical action of the Lie algebra $I_Z$ on $\mathcal C^\infty (Z)=\mathcal C^\infty(M)/I_Z$. If $F\in \mathcal C^\infty (M)$ is a representative of $f\in \mathcal C^\infty(Z)$, i.e., $f=\res\:F$, then $h\in I_Z$ acts on $f$ via $h.f:=\res\:(\{h,F\})$. For the space of invariants of this action we write $\mathcal C^\infty(Z)^{I_Z}$. According to Dirac it is a Poisson algebra, and we call it the \emph{Dirac reduced algebra}.

\begin{SATZ} [Dirac reduction] Given a first class constraint set $Z\subset M$, then $\mathcal C^\infty(Z)^{I_Z}$ carries a canonical Poisson bracket, which is given as follows. Let $f,g\in \mathcal C^\infty(Z)^{I_Z}$ and let $F$ and $G$ be smooth functions on $M$ such that $f=\res\:F$ and $g=\res\:G$. Then the bracket is given by $\{f,g\}:=\res\{F,G\}$.
\end{SATZ}
\begin{BEWEIS} The only thing we have to check is, that the bracket is well defined. Note, that if $f=\res F$ is invariant, then $F$ has to be in the Lie normalizer $N(I_Z):=\{\:h\in\mathcal C^\infty(M)\:|\:\{h,I_Z,\}\subset I_Z\:\}$ of the ideal $I_Z$. By Jacobi's identity $N(I_Z)$ is a Lie subalgebra of $\mathcal C^\infty(M)$. Since $Z$ is first class, it is in fact a Poisson subalgebra, which contains $I_Z$ as a Poisson ideal. In fact, if $f,h\in\mathcal C^\infty(Z)^{I_Z}$ are represented by $\res F=f=\res F'$ and $h=\res H$, then $F-F'\in I_Z$ and it follows that $\res(\{\:F-F',H\:\})=0$ since $H$ is in $N(I_Z)$.
\end{BEWEIS}
\\
 
Of primary interest is, of course, the situation, when $Z=J^{-1}(0)$ is the zero fibre of an equivariant moment map $J$. If $0\in \mathfrak g^*$ is a regular value of $J$, then $Z$ is a coisotropic submanifold of $M$, that is a first class constraint set. Unfortunately, for many singular moment maps of interest $Z$ is \emph{not} first class. An easy example of a moment map with a second class constraint set is the harmonic oscillator (example \ref{harmosc}).  For this reason Dirac reduction is not called universal. At least for compact group actions it is quite obvious, that if Dirac reduction works both procedures give the same result. 

\begin{PROPOSITION}\label{redcomparison}Let $M$ be a Poisson manifold with a Hamiltonian action of a compact, connected Lie group $G$ with equivariant moment map $J:M\to \mathfrak g^*$, such that $Z:=J^{-1}(0)$ is first class. Then the Dirac reduced algebra is isomorphic to the Poisson algebra obtained by  universal reduction.
\end{PROPOSITION}
\begin{BEWEIS} It suffices to show that every smooth invariant function $f$ on $Z$ has an invariant representative $F\in\mathcal C^\infty(M)$. But that can be easily obtained by averaging.
\end{BEWEIS}
\subsection{Normal coordinates}
The basic requisite to make Dirac reduction a powerful tool is a good description of the vanishing ideal, e.g. by giving a set of generators. In \cite{AGJ} the authors developed techniques for deciding whether the components of a given moment map $J_1,\dots,J_\ell$ are a set of generators for the ideal of the zero fibre $Z=J^{-1}(0)$. They rely on the following notion of normal coordinates for the moment map.   
\begin{SATZ}[Existence of normal coordinates] \label{normcoord} Let $M$ be an almost K\"ahler manifold of dimension $2n$ and let the Lie group $G$ act properly on $M$ by automorphisms. Moreover, let $J$ be an equivariant moment map for this action, and $z\in Z:=J^{-1}(0)$. Let us choose a complementary subspace $\mathfrak h\subset \mathfrak g$ to the isotropy subalgebra $\mathfrak g_z\subset \mathfrak g$. We fix a basis $v_1,\dots,v_d$ for $\mathfrak g_z$ and a basis $w_1,\dots,w_e$ for $\mathfrak h$ and denote by $v^1,\dots,v^d$ and $w^1,\dots,w^e$ the corresponding dual basis for $\mathfrak g_z^*$ and  $\mathfrak h^*$, respectively. 
Then there is  a neighborhood $U\subset M$ of $z$ and a local coordinate system around $z$
\begin{align*}
\phi:=(x,y,s,t)=(x_1,\dots,x_m,y_1,\dots,y_m,s_1,\dots,s_e,t_1,\dots,t_e):U\to \mathbbm  R^{2n}, z\mapsto 0,
\end{align*}
such that the following conditions are true.
\begin{enumerate}
\item The subset $\mathcal M:=U\cap \phi^{-1}\big(\{(x,y,0,0)\in\mathbbm  R^{2n}\}\big)$ is a $G_z$-invariant symplectic submanifold of $U$, which is called \emph{the linear reduced space}. There is a representation $G_z\to U(m)$ such that the coordinate map $(x,y):\mathcal M\to \mathbbm  R^{2m}$ intertwines the $G_z$ actions. Here we think  of $\mathbbm  R^{2m}$ provided with the standard  K\"ahler structure with standard $U(m)$-action (see example \ref{standardex}). Accordingly, the moment map $J_z$ for the $G_z$-action on $\mathcal M$ is given by the formula:
\begin{align}\label{polpart}
J_z(x,y)=-\half \sum _{a=1}^d\sum_{i,j=1}^m \big(A^{ij}_a(x_iy_j-x_jy_i)+S^{ij}_a(x_iy_j+x_jy_i)\big)v^a,
\end{align}
where $A_a=(A^{ij}_a)$ and $S_a=(S^{ij}_a)$ are real antisymmetric (resp. symmetric) $m\times m$-matrices for $a=1,\dots,d$.
\item The restriction of $J$ to $U$ is given by the formula:
\begin{align}
J(x,y,s,t)=\operatorname{Ad}^\natural_{\operatorname{exp}(-\sum_{a=1}^e t^a w_a)}\big(J_z(x,y)+\sum_{a=1}^e s^a w^a\big).
\end{align} 
Here $\operatorname{exp}$ stands for the exponential map from $\mathfrak g$ to $G$, and $\operatorname{Ad}^\natural$ denotes the coadjoint action of $G$ on $\mathfrak g^*$.
\item Each of the sets $\{(x,y,s,t)\:|\: x=x_0,y=y_0, s=s_0\}\subset U$ is contained in a single $G$-orbit.
\end{enumerate}
Such a coordinate system will be called \emph{normal coordinates centered at} $z\in Z$. If $z'\in \mathcal M\cap Z$, then there exists an analytic change of normal coordinates centered at $z$ to coordinates centered at $z'$.
\end{SATZ}
\begin{BEWEIS} See \cite[p.62--65]{AGJ}.
\end{BEWEIS}
\\

It is well known that every symplectic manifold $M$ admits a compatible almost complex structure. For a proof of this fact see e.g. \cite{BatesWein}. If in addition a compact Lie group acts on $M$ then this almost complex structure can be made equivariant by averaging over the Haar measure. Note, that in the examples \ref{harmosc},\ref{standardex},\ref{drehimpuls},\ref{mpart},\ref{stratres},\ref{tzwo} and \ref{commvar} the moment map is already given in (global) normal coordinates around the fix point zero. 
\subsection{The generating hypothesis}\label{ghsubsection}
In this subsection we shall explain the tools developed in \cite{AGJ} to decide, whether the components of a given moment map $J:M\to \mathfrak g^*$ generate the vanishing ideal $I_Z\subset \mathcal C^\infty(M)$. In this case we will say for short that $J$ satisfies the \emph{generating hypothesis}.  
For completeness and better intelligibility let us mention the following criteria for the constraint set of a moment map to be first class \cite[Proposition 5.2]{AGJ}.
\begin{SATZ} \label{frstclasscrit} Let $M$ be a symplectic manifold with a Hamiltonian action of a compact Lie group $G$ and let $Z=J^{-1}(0)$ be the zero fibre of the moment map $J$. Then the following statements are equivalent
\begin{enumerate}
\item The \emph{spanning condition} $T_z Z=\ker T_z J$ holds at every $z\in Z$.
\item $Z$ is first class.
\item For every $z\in Z$ all polynomial constraints on the linear reduced space at $z$ (see Theorem \ref{normcoord}) are first class. 
\end{enumerate}
\end{SATZ}
Here the tangent space at $z\in Z$ is the linear span of tangent vectors , which are obtained by taking the derivatives of smooth curves $\gamma:[0,\epsilon[\to Z$, with $\gamma(0)=z$. Note that there are examples of noncompact group actions for which Theorem \ref{frstclasscrit} is wrong.
The spanning condition is  quite a practical tool to sort out second class examples, like the harmonic oscillator. Moreover, the theorem says that the question of $Z$ being a first class constraint set is actually a question of the local real algebraic geometry of the moment map. It is important to note that there are examples of moment maps with first class constraint set, which do not satisfy the generating hypothesis (see \cite[Example 7.13]{AGJ}). Another necessary criterion, which is even more easy to check, is the following \emph{nonpositivity condition} \cite[Proposition 6.7]{AGJ}.

\begin{PROPOSITION} \label{nonpcond} If the spanning condition holds at $z\in Z$, then $J$ fulfills the following \emph{nonpositivity condition}  at $z\in Z$: for every $\xi\in\mathfrak g$ one has either
\begin{enumerate}
\item $J(\xi)=0$ in a neighborhood $U\subset M$ of $z$, or
\item in every neighborhood $U\subset M$ of $z$ the function $J(\xi)$ takes strictly positive as well as strictly negative values.
\end{enumerate}  
\end{PROPOSITION}
If  $G$ is a torus, then, due to the following theorem (see \cite[Theorem 6.8]{AGJ}), the nonpositivity condition is also sufficient for $J$ to satisfy the generating hypothesis.
\begin{SATZ} \label{torus}Let $M$ be a symplectic manifold with a Hamiltonian action of a torus $G$, with moment map $J:M\to\mathfrak g^*$. Then $J$ satisfies the generating hypothesis $\Leftrightarrow$ the spanning condition holds for all $z\in Z=J^{-1}(0)$ $\Leftrightarrow$  the nonpositivity condition holds for all $z\in Z$.
\end{SATZ}

Now we are ready to address the torus actions from our list of examples. As we already indicated the harmonic oscillator, Example \ref{harmosc}, does not satisfy the nonpositivity condition. The same is true for Example \ref{lemon} (the `lemon') at the singular values $\pm 1$. It is easy to see, that at the singular value $0$ the nonpositivity condition is fulfilled. Furthermore, the nonpositivity condition is clearly true for zero angular momentum for $m$ particles in the plane, Example \ref{mpart}, and for the $(1,1,-1,-1)$-resonance, Example \ref{stratres}. For the $\mathbbm  T^2$-action of Example \ref{tzwo} the nonpositivity condition holds iff $\alpha<0$.

As the nonpositivity condition, in the case of nonabelian group actions, is only 
\emph{necessary} for the ideal $I_Z\subset \mathcal C^{\infty}(M)$ to be generated by $J_1,\dots ,J_\ell$, the reasoning here is usually more intricate. As a first step, one notices that it is enough to analyse the problem locally. In fact, with the assumptions of Theorem \ref{normcoord} the following statements are equivalent \cite[Corollary 4.6]{AGJ}
\begin{enumerate}
\item the components of the moment map $J$ generate the vanishing ideal $I_Z\subset \mathcal C^\infty(M)$ of $Z$ in $M$,
\item \label{locgh}for every $z\in Z$, in a normal coordinate system around $z$ the components of the moment map $J_z$ (cf. equation (\ref{polpart})) generate the ideal $I_{J_z^{-1}(0)}\subset \mathcal C^\infty(\mathcal M)$ of the zero level set $J_z^{-1}(0)$ in the linear reduced space $\mathcal M$.
\end{enumerate}
Of course,  it is tempting to view the quadratic polynomial function $J_z$ as a polynomial in the polynomial ring $\mathbbm  R[x_1,\dots,x_m,y_1,\dots,y_m]$, where $m=\half\dim M+\dim G_z-\dim G$. As a matter of fact, if $I_{pol}(J_z)\subset \mathbbm  R[x_1,\dots,x_m,y_1,\dots,y_m]$ denotes the ideal generated polynomial $J_z$ then it follows from the proof of \cite[Theorem 6.3]{AGJ} that if 
\begin{enumerate} \setcounter{enumi}{2}
\item \label{realcrit} the ideal $I_{pol}(J_z)$ is a \emph{real ideal}  in the sense of real algebraic geometry,
\end{enumerate} 
then the local generating hypothesis \ref{locgh}. above is true in $z$ ,i.e., $J_z$ generates the ideal $I_{J_z^{-1}(0)}\subset\mathcal C^\infty(\mathcal M)$.
Let us recall that an ideal $I$ in an $\mathbbm  R$-algebra $R$ is defined to be a \emph{real ideal}, if it coincides with its \emph{real radical} 
\begin{displaymath}
\begin{split}
  \sqrt[\mathbbm R]{I} := \big\{f\in R \mid f^{2i}+\sum_{s=1}^j g_s^2
  \in I\text{ for some $i,j$ and } 
  g_1,\dots , g_j\in R \big\}.
\end{split}
\end{displaymath} 
If $R=\mathbbm  R[x_1,\dots,x_k]$ the \emph{real Nullstellensatz} says that an ideal $I\subset R$  is real iff it coincides with the ideal $I_{V(I)}=\{f\in R \mid f_{|V(I)}=0\}$ of its real locus $V(I)=\{(a_1,\dots,a_k)\in\mathbbm  R^k\mid f(a_1,\dots, a_k)=0\quad \forall f\in I\}$. The real Nulstellensatz is \emph{not} valid for ideals in the ring of smooth functions of some manifold.  The full statement of \cite[Theorem 6.3]{AGJ} is that for the above argument also the converse is true: the components of the moment map $J$ generate the ideal $I_Z\in \mathcal C^\infty(M)$ $\Leftrightarrow$ they generate a real ideal in $\mathcal C^\infty(M)$ $\Leftrightarrow$ for all $z\in Z$ the ideal $I_{pol}(J_z)$ is a real ideal in $\mathbbm  R[x_1,\dots,x_m,y_1,\dots,y_m]$.

As a last step, we follow the advice of \cite{AGJ} and translate the reality criterion \ref{realcrit}. above into the more amenable language of complex algebraic geometry.
\begin{SATZ}[\cite{DubEf}]
\label{GotayThm} 
Let $I$ be an ideal in $\mathbbm R[x_1,\dots, x_k]$. Then $I$ is real, if and only 
if the following two conditions hold:
\begin{enumerate}
\item 
  $I_\mathbbm C:=I\otimes_\mathbbm R \mathbbm C$ is radical in 
  $\mathbbm C[x_1,\dots, x_k]$, and
\item  for every irreducible component $W\subset \mathbbm C^m$ of the (complex) locus of $I_\mathbbm C$ we have that  $\dim_{\mathbbm R} (W\cap \mathbbm R^m)=\dim_{\mathbbm C}(W)$.\footnote{Here, the dimension of a variety is the dimension of the smooth part. In general, we have  $\dim_{\mathbbm R} (W\cap \mathbbm R^m)\le \dim_{\mathbbm C}(W)$.}
\end{enumerate}
\end{SATZ} 
In other words, in order to check the generating hypothesis  it is necessary to  gain 
detailed insight into the complex algebraic geometry behind the scene (e.g.~knowing 
the primary decomposition of $I_\mathbbm C$). Regardless of the fact that the varieties in 
question are cones, there is no straightforward way to provide this information. A  basic and physically interesting example, for which the algebraic geometry is well-studied is the system of one  particle in dimension $n$ with zero angular momentum, Example \ref{drehimpuls}. Since the components of the moment map can (up to a sign) be written as 
the $2\times 2$-minors of the $2\times n$-matrix
\begin{align*}
\left(
\begin{array}{cccc}
q^1&q^2&\dots&q^n\\
p_1&p_2&\dots&p_n
\end{array}
\right)
\end{align*}
 the ideal $I_{\mathbbm C}$ which is generated by the components of the moment map is an instance of a determinantal ideal. By a theorem of Hochster \cite{Hochster} the ideal $I_\mathbbm C$ is prime, and the complex locus is of dimension $n+1$. It is easy to see \cite[Example 7.10]{AGJ}, that the dimension condition is true as well. In fact, the $2\times 2$-minors of the above matrix are zero if and only if the vectors $(q^1,\dots,q^n)$ and $(p_1,\dots,p_n)$ are proportional. In particular, for any $(q^1,\dots q^n)\in\mathbbm  R^n$ and $\lambda\in \mathbbm  R^\times$  the vector $(\lambda q^1,\lambda q^2,\dots,\lambda q^n,
,\lambda^{-1}q^1,\lambda^{-1}q^2,\dots,\lambda^{-1}q^n)\in T^*\mathbbm  R^n$ is in the real locus, whose dimension is thus $\ge n+1$. 
We conclude that the ideal $I$ is real.
Unfortunately, this example is not a 
complete intersection for $n\ge3$ (more on this in section \ref{kc}). 

The only class of nonabelian examples, which the 
author is aware of, where the generating and (as we will see in section \ref{kc}) the complete intersection hypothesis are true 
at the same time, is Example \ref{commvar}. 
The complex locus $Z_{\mathbbm C}$ defined by  
these $\frac{1}{2}n(n-1)$ quadratic equations is an instance of what is called a \emph{commuting variety}. In \cite{BrenPiVas} it was shown that $Z_{\mathbbm C}$ is 
irreducible of codimension $\frac{1}{2}n(n-1)$, and the ideal generated by the 
coefficients of $J$ in the complex polynomial ring is prime. Let 
$S_{\mathrm{reg}}\subset S$ be the open subset of symmetric matrices with pairwise 
distinct eigenvalues. Since the action of $SO(n)$ on $T^*S_{\mathrm{reg}}$ is locally 
free, it follows that $Z\cap T^*S_{\mathrm{reg}}$ is of codimension $\frac{1}{2}n(n-1)$ 
likewise. As a consequence of Theorem \ref{GotayThm}, the components of $J$ generate
the vanishing ideal $I_Z$ in $\mathcal C^{\infty}(T^*S)$. 

Summarizing, we have seen that the generating hypothesis is true for the Examples \ref{drehimpuls}, \ref{mpart}, \ref{lemon}, \ref{stratres}, \ref{tzwo} for $\alpha<0$ and \ref{commvar}.
We would like to close this subsection with a discussion of the elementary example of the free particle on the line, Example \ref{allwrong}, where essentially all what we have done so far \emph{goes wrong}. Clearly, $0$ is a singular value of $J(q,p)=\half p^2$. The constraint surface $Z=J^{-1}(0)=\{(q,0)\in T^* \mathbbm R\}=\mathbbm  R$ is in fact a Lagrangian submanifold, hence first class. The action of $\mathbbm  R$ on $Z$ is trivial, therefore the action of $\mathbbm  R$ on $T^*\mathbbm  R$ can not be proper (proper actions have compact isotropy groups). The orbits of the $\mathbbm  R$-action on $T^*\mathbbm  R$ are closed, but the quotient space $T^*\mathbbm  R/\mathbbm  R$ is not Hausdorff. Furthermore, the quotient space $Z/\mathbbm  R=Z$ is for obvious reasons not a stratified symplectic space, i.e., the theorem of Sjamaar and Lerman does not apply. Any function on $Z$ is invariant, but the only smooth function which can be extended to a smooth invariant functions on $T^* \mathbbm  R$ are the constants. Hence, the Poisson algebra obtained by universal reduction is just $\mathbbm  K$. In contrast, the Dirac reduced algebra is $\mathcal C^\infty(\mathbbm  R)$ with trivial Poisson structure. Finally, $p$ is a constraint which is not a  multiple of $J$. According to \'Sniatycki and Weinstein \cite{SnW} $\mathcal C^\infty(T^*\mathbbm R)/p^2\mathcal C^\infty(T^*\mathbbm R)$ is naturally a Poisson algebra, which in this example differs from the universal and Dirac Poisson reduced structure. 
 
\section{Deformation quantization}
In this section we recall some basic notions and examples from the theory of deformation quantization. We recall Fedosov's construction of star products on symplectic manifolds. We discuss the various invariance properties of star products, which will be the basic requisites for quantum phase space reduction in chapter \ref{qbrst}. 
\subsection{Formal deformations of associative algebras}
Let $A$ be an associative $\mathbbm K$-algebra with unit $\mathbbm 1$ and let us denote the multiplication map by $\mu:A\otimes_\mathbbm K A\to A$,, $a\otimes b\mapsto ab$. A \emph{formal deformation} of the algebra $A$ is a sequence of linear operations  $\mu_i:A\otimes_\mathbbm K A\to A$ for $i\ge 1$, such that the deformed multiplication $*$
\begin{align}\label{formaldef}
a*b:=ab+\sum_{i\ge 1}\nu^i\mu_i(a,b), 
\end{align}
where $a,b\in A$ defines (by $\nu$-linear extension) an associative algebra structure on $A\fnu$ with unit $\mathbbm 1$. The variable $\nu$ is called the \emph{formal parameter}. Two formal deformations $*$ and $*'$ are called \emph{equivalent} if there is a sequence of linear maps $S_i:A\to A$, for $i\ge 1$, such that $S:=\id+\sum_{i\ge 1}\nu^i S_i$ defines an isomorphism of unital $\mathbbm K\fnu$-algebras from  $(A\fnu,*)$ to $(A\fnu,*')$, i.e. \
\begin{align}  \label{equiv}
S(a*b)=S(a)*'S(b)
\end{align} 
for all $a,b\in A$. If the algebra $A$ is commutative, it is easy to show, that the semiclassical limit of a formal  deformation of $A$ makes $A$ into a Poisson algebra, i.e.
\begin{align}\{a,b\}:=\mu_1(a, b)-\mu_1(b, a),\quad \mbox{ for }a,b\in A
\end{align} 
is a Poisson bracket. Conversely, if we start with a Poisson algebra $(A,\cdot,\{\:\:\})$, then a formal deformation of $(A,\cdot)$ is said to be a deformation of this Poisson algebra, if the semiclassical limit reproduces the original Poisson structure. 

An easy example of such a formal deformation of a commutative algebra $(A,\mu)$ arises if there is a family of pairwise commuting derivations $D_i$ , $i=1,\dots,n$, of $A$. Then for any  tensor $P^{ij}\in \mathbbm K$, $i,j=1\dots,n$, 
\begin{align*}
a*b:=\mu\circ \ee^{\nu\sum_{i,j=1}^nP^{ij}D_i\otimes D_j}(a\otimes b)
\end{align*}
defines an associative deformation of $(A,\mu)$. If $A$ is the real polynomial ring $\mathbbm R[x^1,\dots,x^n]$ or $\mathcal C^\infty(\mathbbm R^n)$ and the $D_i$ are the partial derivatives $\frac{\partial}{\partial x^i}$ then the above multiplication is called the \emph{Moyal-Weyl}-multiplication. The semicassical limit of the Moyal-Weyl multiplication is the Poisson structure  corresponding to the constant Poisson bivector field $\half\sum_{ij}\Pi^{ij}\frac{\partial}{\partial x^i}\wedge\frac{\partial}{\partial x^i}$, where $\Pi^{ij}=P^{ij}-P^{ji}$. Another important example of an associative deformation arises if  $A=S\mathfrak h$ is the symmetric algebra of an real Lie algebra $\mathfrak h$ from the PBW-symmetrization map (more on this in subsection \ref{qmom}). The notion of a formal deformation of an associative algebra, straightforwardly generalizes to $\mathbbm Z_2$ or $\mathbbm Z$-graded algebras (with Koszul sign rule). The purely odd analogon of Moyal-Weyl product is Clifford multiplication (see also equation (\ref{clmult})). We leave it to the reader to fill in the details. 

   It should be mentioned that the stepwise obstructions to construct a formal deformation of the multiplication $\mu$ lie in the third Hochschild cohomology group ${HH}^3(A,A)$ of the algebra $A$. In fact, the space  $C^n(A,A)$ of Hochschild  $n$-cochains is just the space linear maps $\varphi: A^{\otimes n}=A\otimes \dots \otimes A\to A$ and the differential of Hochschild cohomology $\delta_H:C^n(A,A)\to C^{n+1}(A,A)$ is given by the formula
\begin{eqnarray}
(\delta_H\varphi)(a_1,\dots, a_{n+1})&=&a_1\varphi(a_2,\dots, a_n)+\sum_{i=1}^{n}(-1)^i\varphi(a_1,\dots, a_ia_{i+1},\dots,a_{n+1})\nonumber\\
&&\qquad\qquad\qquad\qquad+(-1)^{n+1}\varphi(a_1,\dots, a_n)\:a_{n+1}.\label{Hochdiff}
\end{eqnarray}
for $a_1,\dots, a_{n+1}\in A$.
The associativity of the product (\ref{formaldef}) in degree $\nu^k$ can be written as 
\begin{align}
(\delta_H \mu_k)(a_1,a_2,a_3)=\sum_{i=1}^{k-1} \mu_i(\mu_{k-i}(a_1,a_2),a_3)-\mu_i(a_1,\mu_{k-i}(a_2,a_3))\qquad\forall a_1,a_2,a_3\in A.
\end{align}
Therefore, if we would construct such a  product inductively, then we would have to assure that the right hand side of this equation is a Hochschild 3-coboundary at every step.
Using the differential graded Lie algebra structure on the shifted Hochschild cochain complex $C^\bullet(A,A)[1]$ which is given by the so-called Gerstenhaber bracket, the set of equivalence classes of formal deformations of $\mu$ can be described in terms of deformation functors and Maurer-Cartan equations. As we will not use this slightly more sophisticated language we refer the interested reader to the exposition \cite{KellerMini}. 
\subsection{Star products}
If we are looking for a formal deformation of the Poisson algebra $A:=\mathcal C^\infty(M)$ of smooth functions on a Poisson manifold $(M,\Pi)$, then the general opinion is that full Hochschild cochain complex $C^\bullet (A,A)$ is too big for the deformation problem to make sense. As a consequence one usually considers the subcomplex $C^\bullet_{\diffind}(A,A)$ of \emph{differential} Hochschild cochains
\begin{align*}
C^k_{\diffind}(A, A)=\{D:\underbrace{A\otimes\dots\otimes A}_{\mbox{k times}},\to A\:|\:D \mbox{ is a polydifferential operator}\}
\end{align*} 
instead. 
Accordingly, a formal deformation 
\begin{align}
f*g:=fg+\sum_{i\ge 1}\nu^i\mu_i(f\otimes g),\qquad f,g\in \mathcal C^\infty(M)
\end{align}
of the Poisson algebra $(\mathcal C^\infty(M),\{\:,\})$ is defined to be a \emph{star product}, if the bilinear operations $\mu_i$, $i=1,2,\dots$, are in fact \emph{bidifferential operators}. Moreover, two such star products $*$ and $*'$ are defined to be equivalent, if the equivalence transformation $S:=\id+\sum_{i\ge 1}\nu^i S_i$ of equation (\ref{equiv}) is in fact a series of differential operators. Essentially all known explicit constructions of star products yield bidifferential operators $\mu_i$, which are differential operators of order at most $i$ in each argument. These star products are also called \emph{natural}. It can be shown \cite{GuttRawnsley03} that an equivalence transformation $S$ between two natural star products $*$ and $*'$ is of the form $S=\operatorname{exp}(\sum_{i=1}^\infty \nu^iD_i)$, where the $D_i$ are differential operators of order at most $i+1$. Moreover, in the symplectic case every star product is equivalent to a Fedosov star product (see the next subsection) and, hence, to a natural one.  

Instead of considering differential star products one can also exploit the fact that $\mathcal C^\infty(M)$ is a nuclear Fr\'echet algebra and look at \emph{continuous} star products (defined below). Even though the meaning of continuity is somehow obscure,  this approach has the great advantage that it applies for other interesting nuclear Fr\'echet algebras such as the singular reduced algebra $\mathcal C^\infty(Z)^{\mathfrak g}$ of subsection \ref{Diracsubsection}. In fact, $\mathcal C^\infty(Z)$, being a quotient of a nuclear Fr\'echet space modulo a closed subspace, is nuclear Fr\'echet (cf. \cite[Proposition 50.1]{Treves}), as well as $\mathcal C^\infty(Z)^{\mathfrak g}$, being a closed complemented subspace of a nuclear Fr\'echet space (remember that $G$ is assumed to be compact and connected). 

If $A$ is a nuclear Fr\'echet algebra with commutative multiplication $\mu$, then the space $C^\bullet_{\contind}(A, A)$ of \emph{continuous} Hochschild cochains as follows
\begin{align*}
C^k_{\contind}(A,A)=\{D:\underbrace{A\widehat{\otimes}\dots\widehat{\otimes} A}_{\mbox{k times}},\to A\:|\:D \mbox{ is linear continuous}\}.
\end{align*} 
Here, $\widehat{\otimes}$ denotes the topological tensor product (since $A$ is a nuclear space all topological tensor products coincide). Accordingly, by a \emph{continuous star product} we mean a formal deformation as in formula (\ref{formaldef}), such that the operations $\mu_i:A\widehat{\otimes}A\to A$ are continuous. An equivalence transformation between continuous star products such as in formula (\ref{equiv}) is said to be \emph{continuous} if the operations $S_i:A \to A$ are continuous. The notion of a continuous star product has been studied before in \cite{BFGP}, together with the topological version of Gerstenhaber deformation theory. Actually, for the main application we have in mind, i.e., the singular reduced algebra $\mathcal C^\infty(Z)^\mathfrak g$, it is not at all obvious what the `correct' definition of  (multi-) differential operator should be (see \cite{PflaumHab}), in order to define a feasible notion of a differential star product. As an indication that the continuous setup is not too weak, let us mention the well known fact that the natural map $C_{\diffind}(A,A)\to C_{\contind}(A,A)$ is a quasiisomorphism if $A=\mathcal C^\infty(M)$ is the algebra of smooth functions on a smooth manifold $M$. The statement that the cohomology of $C_{\diffind}(A,A)$ (respectively $C_{\contind}(A,A)$) is isomorphic to the space of polyvector fields on $M$ is just the well known differential (resp. continuous) Hochschild-Kostant-Rosenberg theorem.

\subsection{Fedosov construction}
For the sake of completeness and because it is a nice application of the perturbation lemma \ref{BPL1}, let us recall Fedosov's construction of star products on a symplectic manifold  $(M,\omega)$. 
The point of departure of Fedosov's construction is an algebraic, fibrewise version of the de Rham complex. This complex occurs at many places in mathematics, e.g. in \cite[\S 9 No.3 exemple 1]{Bourbdix}, but viewed the opposite way as a Koszul complex of a `linear constraint'. More precisely, we consider $S\otimes \wedge \:\:T^*M\to M$, the tensor product of the symmetric and Grassmann algebra bundle of $T^*M$, the cotangent bundle of the manifold $M$. A choice of local coordinates $x^1,\dots,x^n$ for $M$ gives a local frame $dx^1,\dots,dx^n$ for $T^*M=\wedge^1 T^*M$. The corresponding frame for $T^*M=S^1 T^*M$ will be written as $y^1,\dots,y^n$. We will interprete a section of   $S^k\otimes \wedge^l \:\:T^*M\to M$ as a polynomial valued $l$-form, i.e. the local frames are written as $y^{i_1}\dots y^{i_k}dx^{j_1}\wedge\dots\wedge dx^{j_l}$ for $i_1\le i_2,\dots\le i_k$ and $j_1<j_2\dots<j_l$. We introduce the algebraic (fibrewise) de Rham differential
\[\delta:=\sum_i dx^i\wedge \frac{\partial}{\partial y^i}:\Gamma^\infty(M,S^k\otimes \wedge^l\:T^*M)\to\Gamma^\infty(M,S^{k-1}\otimes \wedge^{l+1}\:T^*M).\]
Obviously, $\big(\Gamma^\infty(M,S\otimes \wedge^\bullet\:T^*M),\wedge,\delta\big)$ is a $\mathbbm Z$-graded super-commutative differential graded algebra. Moreover, there is an algebraic Poincar\'e lemma, i.e., the differential $\delta$ is acyclic. A contracting homotopy for $\delta$ is given as follows. We introduce the Koszul differential $\delta^*=\sum_i y^i i(\partial/\partial x^i):\Gamma^\infty(M,S^k\otimes \wedge^l\:T^*M)\to\Gamma^\infty(M,S^{k+1}\otimes \wedge^{l-1}\:T^*M)$. An easy calculation yields the commutation relation $\delta\delta^*+\delta^*\delta=(k+l)\id:\Gamma^\infty(M,S^k\otimes \wedge^l\:T^*M)\to\Gamma^\infty(M,S^k\otimes \wedge^l\:T^*M)$. Renormalizing the Koszul differential $\delta ^*$, we obtain a contracting homotopy: $\delta^{-1}:=(k+l)^{-1}\delta^*:\Gamma^\infty(M,S^k\otimes \wedge^l\:T^*M)\to\Gamma^\infty(M,S^{k+1}\otimes \wedge^{l-1}\:T^*M)$ for $k+l>0$ and for $k=0=l$ we define $\delta^{-1}$ to zero. Let us introduce the canonical projection $\pi:\Gamma^\infty(M,S\otimes \wedge\:T^*M)\to \mathcal C^\infty(M)$ and the canonical injection $\iota:\mathcal C^\infty(M)\to\Gamma^\infty(M,S\otimes \wedge\:T^*M)$. then the commutation relation above can be rewritten as $\delta\delta^{-1}+\delta^{-1}\delta=\id-\iota\pi$. Traditionally $\iota\pi$ is denoted by $\sigma$. What we have done so far may be neatly subsumed in the language of appendix \ref{HPT}:
\begin{align} \label{elemcontr}
  \xymatrix{(\mathcal C^\infty(M),0)\ar@<-0.6ex>[r]_{\iota\quad\qquad}^{}&(\Gamma^\infty(M,S\otimes \wedge T^*M),\delta),\delta^{-1} \ar@<-0.6ex>[l]_{\pi\quad\qquad}^{}},
\end{align}
is a \emph{contraction} fulfilling all side conditions. Moreover, $\iota$ and $\pi$ are homomorphisms of super-commutative algebras and $\delta$ is a derivation.

It is clear, that if we replace the symmetric algebra part by the algebra of formal power series (i.e., we  take the completion with respect to the ideal generated by $\Gamma^\infty(M, S^1\otimes \wedge^0 T^*M)$), we do not spoil the contraction (\ref{elemcontr}). The same applies, if we adjoin a formal variable $\nu$. Therefore, we may replace the above de Rham complex by the so called \emph{Weylalgebra}  
\begin{align}\mathcal W\otimes \Omega^\bullet(M):=\prod_{k\ge 0}^\infty \Gamma^\infty(M,S^k\otimes \wedge^ \bullet T^*M)\fnu,
\end{align}
which is in an obvious way a $\mathbbm Z$-graded super-commutative $\mathbbm K\fnu$-algebra. Again, we have a contraction
\begin{align} \label{formelemcontr}
  \xymatrix{(\mathcal C^\infty(M)\fnu,0)\ar@<-0.6ex>[r]_{\iota\quad}^{}&\big(\mathcal W\otimes \Omega^\bullet(M),\delta\big),\delta^{-1} \ar@<-0.6ex>[l]_{\pi\quad}^{}},
\end{align}
which fulfills all side conditions. 

The next ingredient is a formal deformation $\circ$ of the super-commutative product on $\mathcal W=\prod_{k\ge 0}^\infty \Gamma^\infty(M,S^k T^*M)\fnu$ into an associative multiplication. It is given by a fibrewise Moyal-Weyl-multiplication
\begin{align}
a\circ b:= \mu \exp\Big(\frac{\nu}{2}\sum_{i,j}\Pi^{ij}\frac{\partial}{\partial y^i}\otimes\frac{\partial}{\partial y^j}\Big)(a\otimes b),
\end{align}
where $a,b\in\mathcal W$ and $\mu$ denotes the commutative multiplication. Recall that according to our sign convention $\sum_k\Pi^{ik}\omega_{kj}=\delta_j^i$. It is clear, that the  Moyal-Weyl-multiplication $\circ$ respects the form degree. But in fact, there is a second $\mathbbm Z$-grading for which $\circ$ is graded. This \emph{total degree} is given by counting the symmetric degree and twice the $\nu$ degree simultaneously. The homogeneous components of this grading are given by  the eigenspaces $\mathcal W^{(k)}=\{a\in \mathcal W\mid \operatorname{Deg}a=ka\}$ of the derivation $\operatorname{Deg}:=\sum_i y^i\frac{\partial}{\partial y^i}+ 2\nu\frac{\partial}{\partial \nu}$, which is a derivation of $\circ$. For the descending filtration induced by $\operatorname{Deg}$ we write $\mathcal W_i:=\prod_{k\ge i}^\infty \mathcal W^{(k)}$. The multiplication $\circ$ extends naturally to $\mathbbm Z$-graded multiplication $\circ$ for $\mathcal W\otimes\Omega^\bullet(M)$. The super-center of the algebra $(\mathcal W\otimes\Omega^\bullet(M);\circ)$ is just $\Omega(M)[[\nu]]$. The derivation $\delta$ can be written as an inner derivation of $\mathcal W\otimes\Omega^\bullet(M)$: $\delta=\nu^{-1}\ad(\widetilde{\omega})$, where is $\widetilde{\omega}=\sum_{i,j}\omega_{ij}y^i dx^j$ and $\ad(\widetilde \omega)$ means taking the super-commutator with $\widetilde \omega$.  
Let $\nabla$ be a torsion free symplectic connection on the tangent bundle of $M$, i.e.,$\nabla_X Y-\nabla_Y X=[X,Y]$ for all $X,Y\in \Gamma^\infty(M,TM)$ and $\nabla \omega =0$. Such a connection can always be found. In contrast to the Riemannian case a symplectic connection is not uniquely determined. Let $\hat R\in \Gamma^\infty(M, \End(TM)\otimes \wedge^2 T^*M$ denote the curvature endomorphism of the connection: $\hat R(X,Y)Z=\nabla_X\nabla_Y Z-\nabla_Y\nabla_X Z-\nabla_{[X,Y]}Z$. Because $\nabla$ is symplectic one can show that 
\[R(X,Y,Z,W):=\omega(X,\hat R(Z,W),Y)\]
is in fact symmetric in $X$ and $Y$. This means that $R$ is in fact a section in $\Gamma^\infty(M,S^2T^*M\otimes\wedge ^2 T^*M)$, i.e., an element of $\mathcal W\otimes \Omega$. We extend $\nabla$ in a standard fashion to 
\begin{eqnarray*}
\nabla=\sum_i dx^i\wedge\nabla_{\frac{\partial}{\partial x^i}}:\mathcal W\otimes \Omega ^\bullet\to \mathcal W\otimes \Omega ^{\bullet+1}.
\end{eqnarray*}   
A little calculation yields that 
\[\nabla^2=\nu^{-1}\ad_\circ (R).\] 
Since the  connection is symplectic, one can prove that $\nabla$ is an odd derivation of $\circ$ and that $\nabla\widetilde{\omega}=0$. 
The latter equation entails that $[\nabla,\delta]=\nabla\delta+\delta\nabla=0$. From $[\delta,[\nabla,\nabla]]=0$ one derives $\delta R=0$, which is also known as the first Bianchi identity. From $[\nabla,[\nabla,\nabla]]=0$ one infers 
that $\nabla R=0$, which is known as the second Bianchi identity. 

The beautiful insight of Fedosov \cite{Fed94} has been that even though the derivation $-\delta +\nabla$ is, in general, not of square zero, it can be made into a differential by adding an (almost) inner derivation $\frac{1}{\nu}\ad_\circ(r)$. This $r$ can be found recursively.

\begin{SATZ} [Fedosov] Let $\Omega=\sum_{i\ge1}^\infty\nu^i \Omega_i\in \nu Z_{dR}^2(M)[[\nu]]$ be a series of closed 2-forms on $M$ and $s\in \mathcal W_3\otimes \Omega^0(M)$ such that $\pi(s)=0$. Then there is a unique $r\in  \mathcal W_2\otimes \Omega^1(M)$ such that $\delta r=R+\nabla r+\frac{1}{\nu}r\circ r+\Omega$ and $\delta^{-1}r=s$. It follows that  the \emph{Fedosov derivation}
\begin{align}
D:=-\delta+\nabla+\frac{1}{\nu}\ad_\circ(r)
\end{align}
is a differential $\mathcal W\otimes \Omega^\bullet(M)\to\mathcal W\otimes \Omega^{\bullet+1}(M)$, i.e., $D^2=0$.
\end{SATZ}
\begin{BEWEIS} The original proof (for $s=0$) can be found in \cite{Fed}. For a more elaborate exposition the reader may, e.g.,  consult \cite[subsection 6.4.2]{Waldmannsoevre}.
\end{BEWEIS}
\\

The Fedosov derivation can be seen as a perturbation of the differential $-\delta$ of the contraction (\ref{formelemcontr}). The filtration in question is the aforementioned one, which is associated to the degree $\operatorname{Deg}$. Clearly, this perturbation complies with the premises of perturbation lemma \ref{contrEins}. As a result we obtain a contraction
\begin{align} \label{Fedcontr}
  \xymatrix{(\mathcal C^\infty(M)[[\nu]],0)\ar@<-0.6ex>[r]_{\tau\quad}^{}&(\mathcal W\otimes \Omega(M),D),D^{-1} \ar@<-0.6ex>[l]_{\pi\quad}^{}}.
\end{align}
The contracting homotopy is given by $D^{-1}=-\delta^{-1}(-D\delta^{-1}-\delta^{-1}D)^{-1}$. A little calculation yields  that the fomula for the map $I$ in Lemma \ref{contrEins} reproduces the well-known fomula for the \emph{Fedosov Taylor series} $\tau=(-D\delta^{-1}-\delta^{-1}D)^{-1}\iota$. 

The contraction (\ref{Fedcontr}) is used to transfer the associative algebra structure from the differential graded associative algebra $(\mathcal W\otimes \Omega(M),\circ,D)$ to $\mathcal C^\infty(M)[[\nu]]$, i.e., we obtain an associative product
\begin{align}
f* g:=\pi (\tau(f)\circ\tau(g)),
\end{align} 
for $f,g \in \mathcal C^\infty(M)[[\nu]]$. The associativity of such a product follows from general considerations (cf. also the computation (\ref{assoccomp}) in chapter \ref{qbrst}).   The product $*$ will be called the \emph{Fedosov star product} obtained from the data $(\nabla,\Omega,s)$. The name is justified by the following well known theorem.
\begin{SATZ} If $*$ is obtained from  the data $(\nabla,\Omega,s)$ then it is natural star product and the equivalence class of $*$ does only depend on the cohomology class $[\Omega]\in\mathrm \nu H_{dR}^2(M)[[\nu]]$. Moreover, every star product on $M$ is equivalent to some Fedosov star product.
\end{SATZ}
\begin{BEWEIS} An elementary proof for the first statments can be found in \cite{Waldmannsoevre}. The last statement follows from classification results of the thesis of N. Neumaier \cite{Neudiss}. 
\end{BEWEIS}
\\

We would also like to mention the following nice result on derivations of Fedosov star products.

\begin{SATZ} \label{nicethm} Let $M$ be a symplectic manifold with a Fedosov star product $*$, which is obtained from the data $(\nabla, \Omega,s)$. If $X$ is a symplectic vector field on $M$, then $X$ is a derivation of $*$ if and only if $\nabla$ is affine with respect to $X$ and $\mathcal L_X\Omega=0=\mathcal L_X s$. Such a vector field can be written as an (almost) inner derivation: 
\begin{align*}
X(h)=\frac{1}{\nu}ad_*(f) h \qquad\forall h\in\mathcal C^\infty(M)
\end{align*} 
for some series $f=\sum_{i\ge 0}\nu^i\,f_i\in \mathcal C^\infty(M)\fnu$ if and only if $f$ is a solution of
\begin{align*} 
df=i_X(\omega+\Omega).
\end{align*} In this case we have $X=X_{f_0}$.
\end{SATZ}
\begin{BEWEIS} A proof can be found in \cite[section 3]{Neubahns}.
\end{BEWEIS}

\subsection{Quantum moment maps and strong invariance}\label{qmom}

There are several inequivalent notions of compatibility of a star product $*$ on $M$ with action of some Lie group $G$ on $M$. For example, if the Lie group $G$ acts by automorphisms of the algebra $(\mathcal C^\infty(M)\fnu,*)$, i.e., $(\Phi_g^*f)*(\Phi_g^*h)=\Phi_g^*(f*h)$ for all $f,h\in\mathcal C^\infty(M)$ and $g\in G$ one simply says that $*$ is \emph{$G$-invariant}. The infinitisemal version of this notion is that of \emph{$\mathfrak g$-invariance} of the star product $*$. Here one requires that the fundamental vector fields act by derivations. If the group $G$ is connected both notions are clearly equivalent. 

A slightly stronger notion, which will become important for us (cf. chapter \ref{qbrst}), is the notion of a \emph{strongly invariant star product}: $*$ is said to be \emph{strongly invariant} with respect to a Hamiltonian action of a Lie group $G$ with moment map $J:M\to \mathfrak g^*$ if 
\begin{eqnarray}
J(X)*f-f*J(X)=\nu \{J(X),f\}\qquad \forall X\in \mathfrak g,f\in\mathcal C^\infty(M).  
\end{eqnarray}
It follows easily from Theorem \ref{nicethm} that for symplectic manifolds with a Hamiltonian action of a compact Lie group a strongly invariant star product can always be found. In fact, the symplectic connection can be made invariant by averaging, the same is true for $s$ (alternatively, one can assume $s=0$). Finally, one has to assure that $i_{X_{J(\xi)}}\Omega=0$ for all $\xi\in \mathfrak g$. 

Another example of a strongly invariant star product is provided by the BCH star product \cite{Gutt83} on the dual space $\mathfrak h^*$ of a $\mathbbm R$-Lie algebra $\mathfrak h$. Let $\mathfrak h_\nu:=\mathfrak h\fnu$ be the $\mathbbm R\fnu$-Lie algebra with the modified bracket $[\:,\:]_\nu:=\nu[\:,\:]$ and let $U\mathfrak h_\nu$ be the universal enveloping algebra of this $\mathbbm R\fnu$-Lie algebra. Let us denote the canonical multiplication on $U\mathfrak h_\nu$ by $\cdot_\nu$. The (rescaled) PBW symmetrization map $\sigma_\nu$ is defined on monomials $X_1X_2\dots X_k\in S^k \mathfrak h$ by the formula
\begin{align}\label{PBWsymm}
\sigma_\nu(X_1X_2\dots X_k)=\frac{\nu^k}{k!}\sum_{\tau\in S_k} X_{\tau(1)}\cdot_\nu  X_{\tau(2)}\cdot_\nu\dots\cdot_\nu X_{\tau(n)}.
\end{align}
It is well-known that $\sigma_\nu$ extends to an injective map of $\mathbbm R\fnu$-modules $\sigma _\nu:S\mathfrak h\fnu\to U\mathfrak h_\nu$, such that $\sigma_\nu\left(S\mathfrak h\fnu\right)$ is a subalgebra. It induces an associative multiplication $*_{BCH}$ on  $ S\mathfrak h\fnu$ which is uniquely defined by the formula
\begin{align}\label{BCHstar}
\sigma_\nu(f*_{BCH}g)=\sigma_\nu(f)\cdot_\nu\sigma_\nu(g).
\end{align}
In fact, the bilinear composition $*_{BCH}$ is given by  a series of bidifferential operators for the algebra $S\mathfrak h$. Viewing $S\mathfrak h$ as polynomial functions on $\mathfrak h^*$, the BCH product $*_{BCH}$ extends uniquely to a star product on $\mathcal C^\infty(\mathfrak h^*)\fnu$, which deforms the linear Poisson structure. By an elementary calculation one can show  that for any $X\in \mathfrak h$ and any monomial $X_1X_2\dots X_k\in S^k\mathfrak h$
\begin{align}
X\cdot_\nu \sigma_\nu(X_1X_2\dots X_k)-\sigma_\nu(X_1X_2\dots X_k)\cdot_\nu X=\nu \sigma_\nu\left( \{X,X_1X_2\dots X_k\}\right),
\end{align}
where $\{,\}$ is the Poisson bracket arising from the linear Poisson structure.
It follows that $*_{BCH}$ is strongly invariant for any  moment map as in Example \ref{linPoiss}.

There is yet another notion of compatibility of the $G$-action with the star product $*$, that is interesting for our purposes (cf. chapter \ref{qbrst}). This is the notion of a \emph{quantum moment map} introduced by Xu \cite{Xuqmom}, which is an deformed analog of equation (\ref{JLiemorph}). If $M$ is a symplectic manifold with a Hamiltonian action of a Lie group $G$ with moment map $J:M\to \mathfrak g^*$ then a \emph{quantum moment map} is a linear map $\qimp:\mathfrak g\to \mathcal C^\infty(M)\fnu$, such that $\qimp(\xi)=J(\xi)+\sum_{i\ge 1}\nu^i\,J_i(\xi)$ and we have 
\begin{align}\label{qmomprop} 
\qimp(X)*\qimp(Y)-\qimp(Y)*
\qimp(X)=\nu\qimp([X,Y])\quad \forall X,Y\in \mathfrak g.
\end{align}
If we have found such a $\qimp$, then we say that $*$ is \emph{quantum covariant with moment map $\qimp$}.
It follows that the linear map $\qrep$ which associates to $X\in \mathfrak g$ the operator $\qrep_X:=\frac{1}{\nu}\ad_*(\qimp(X))$ makes $\mathcal C^\infty(M)\fnu$ into a $\mathfrak g$-module. Note that we \emph{do not assume} that this representation has to coincide with the representation given by the classical moment map (the latter condition is defined to be part of the data, e.g., in \cite{Xuqmom,Neubahns}). The existence and uniqueness question for quantum moment maps has been discussed in detail so far merely for the case when these representations coincide, see \cite{Neubahns} (presumably most of the statements generalize somhow to the above situation). It is clear that one has always the freedom to add to a quantum moment map a Lie algebra 1-cocycle with values in the center of the algebra $(\mathcal C^\infty(M)\fnu,*)$. Besides, a quantum moment map gives rise to a ring homomorphism $U\mathfrak g_\nu\to \mathcal C^\infty(M)\fnu$. A general formalism to treat quantum reduction for such change of ring maps has been proposed in \cite{Sevost}.



\chapter{The classical BFV-construction}\label{BFVchapter}
In this section we collect some results and techniques, which are related to the classical BFV-construction. First we recall the basic rules of multilinear super-algebra. We recall the notion of graded Poisson/Gerstenhaber algebras and review the derived bracket construction introduced by Koszul and Kosmann-Schwarzbach \cite{KosSchw96}.  We study (generalized) graded manifolds in the setting of graded Lie-Rinehart pairs. We propose a method how to adjoin momenta for the `antighost variables'.  We generalize the construction of the Rothstein Poisson-bracket to generalized graded manifolds, and discuss thereby the finitely and infinitely generated case seperately. We give a criterion how to check the acyclicity of a Koszul complex of an analytic map over the ring of smooth functions and show that in this situation there are continuous contracting homotopies. We apply the criterion to our list of examples. We introduce the notion of a projective Koszul-Tate resolution. We generalize  the `Existence of the BRST-charge'-theorem to the vector bundle setting. We discuss the special cases of coisotropic submanifolds and of moment maps which satisfy the generating and the complete intersection hypothesis. 

\section{Multilinear super-algebra}

Let $A$ be a commutative $\mathbbm K$-algebra and $\Gamma$ be one of the abelian groups $\mathbbm Z_2=\mathbbm Z/2 \mathbbm Z$ or $\mathbbm Z$. Since $A$ is commutative there is no distinction between left and right modules. The tensor product and the space of morphisms of two modules are $A$-modules in a canoncical way. We would like to recall the abelian tensor category  $\Gmod{A}$ of $\Gamma$-graded $A$-modules. Objects in $\Gmod{A}$ are $\Gamma$-graded $A$-modules $\mathcal V=\oplus_{n\in \Gamma} \mathcal V^k$. Elements of $v\in\mathcal V$ such that $v\in V^k$ for some $k\in \Gamma$ are said to be homogeneous of degree $|v|:=k$. The sign of a homogeneous element $v$ is defined to be $(-1)^{|v|}$. Homogeneous elements with sign 1 are said to be \emph{even} and homogeneous elements with sign -1 are said to be \emph{odd}. Given two graded $A$-modules $\mathcal V$ and $\mathcal W$ the space of linear maps of degree $i\in \Gamma$ is defined to be \[\Hom^i_{\Gmod{A}}(\mathcal V,\mathcal W)=\{\varphi\in\Hom_{A}(\mathcal V, \mathcal W)\:|\: \varphi(\mathcal V^j)\subset \mathcal W^{j+i}\quad\forall j\in \Gamma\}.\]
The space of morphisms in $\Gmod{A}$ between two graded $A$-modules $\mathcal V$ and $\mathcal W$ is defined to be 
\[\Hom_{\Gmod{A}}(\mathcal V,\mathcal W):=\bigoplus_{i\in \Gamma}\Hom^i_{\Gmod{A}}(\mathcal V,\mathcal W).\]
Clearly, the Hom-sets are $A$-modules and the composition of morphisms is $A$-linear. Moreover, every morphism has a kernel and a cokernel. 
The tensor product $\mathcal V\otimes \mathcal W:=\oplus _{n\in \Gamma}(\mathcal V\otimes \mathcal W)^n$ of two objects $\mathcal V$ and $\mathcal W$ in $\Gmod{A}$ is the direct sum of its homogenous components
\[(\mathcal V\otimes \mathcal W)^n:=\bigoplus_{i,j\in \Gamma,\:i+j=n}\mathcal V^i\otimes _A \mathcal W^j.\] 
The tensor product is a biadditive bifunctor in $\Gmod{A}$.
The \emph{commutativity morphism} is given by
\begin{align*}
\tau_{\mathcal V,\mathcal W}:\mathcal V\otimes \mathcal W\to\mathcal W\otimes \mathcal V,
\quad v\otimes w\mapsto (-1)^{|v||w|}w\otimes v 
\end{align*}
for homogeneous $v$ and $w$. Together with the obvious associativity morphism $\sigma_{\mathcal U,\mathcal V,\mathcal W}:(\mathcal U\otimes \mathcal V)\otimes \mathcal W\to \mathcal U\otimes(\mathcal V\otimes\mathcal W)$ it satisfies the triangle, pentagon and the hexagon axiom (for details see e.g. \cite{MacLaneWorking}). 

Note that there is a forgetful functor of tensor categories  $\Zmod{A}\to \Zzmod{A}$ which associates to a $\mathbbm Z$-graded $A$-module $\mathcal V=\oplus_{i\in \mathbbm Z} \mathcal V^i$ the $\mathbbm Z_2$-graded $A$-module $\mathcal V_{\underline 0}\oplus\mathcal V_{\underline 1}$, where  $\mathcal V_{\underline 0}:=\oplus_{i\text{ even}} \mathcal V^i$ and $\mathcal V_{\underline 1}:=\oplus_{i\text{ odd}} \mathcal V^i$. For every $j\in \mathbbm Z$ the \emph{shift} $\Zmod{A}\to\Zmod{A}$, $\mathcal V\mapsto \mathcal V[j]$ is given by 
\begin{align*}
\mathcal V[j]^i:=\mathcal V^{i+j} 
\end{align*}
for all $i\in \mathbbm Z$. It is compatible with the tensor product in the following sense: $\mathcal V[i]\otimes \mathcal W[j]=\mathcal V\otimes \mathcal W[i+j]$. The canonical map $\mathcal V\to \mathcal V[j]$ has degree $-j$. The analog of the shift functor $[1]$ in $\Zzmod{A}$ is the parity change.

By considering a graded $A$-module as a complex with zero differential the category  $\Zmod{A}$ is a full subcategory of the category of complexes $\chain{A}$ over $A$. Here we adopt the convention that maps of chain complexes (which will be cochain complexes if not otherwise specified) may carry a nonzero degree. On the other hand, there is a forgetful functor $\chain{A}\to \Zmod{A}$. Note that $\chain{A}$ has also the structure of a tensor category. The differential $d_{X\otimes Y}$ on the tensor product $X\otimes Y$ (understood in $\Zmod{A}$) of the complexes $(Y,d_Y)$ and $(X,d_X)$ is given by the formula 
\begin{align*}d_{X\otimes Y}(x\otimes y):=d_X x \otimes y+(-1)^{|x|}x\otimes d_Y y
\end{align*}
for homogeneous $x\in X$ and $y\in Y$. The forgetful functor is a functor of tensor categories. Furthermore, the space of graded $A$-linear morphisms $\Hom_{\Zmod{A}}(X,Y)$ between two complexes $(X,d_X)$ and $(Y,d_Y)$ is in a natural way a complex with differential:
\begin{align*}
D(\varphi):=d_Y \varphi-(-1)^{|\varphi|}\varphi\: d_X
\end{align*}
for homogeneous $\varphi\in\Hom_{\Zmod{A}}(X,Y)$. The chain maps $\Hom^0_{\chain{A}}(X,Y)$ are precisely the $0$-cycles of this complex. Two chain maps are homotopic iff they are homologous. The shift functor extends to $\chain{A}$ by setting $d_{X[j]}:=(-1)^j d_X$.

By reinterpreting the structural diagrams of the basic algebraic structures, such as that of a commutative algebra, Lie algebra, module etc., these notions straightforwardly translate into the language of tensor categories. For instance, super-commutative $\mathbbm K$-algebras and $\mathbbm K$-Lie algebras are nothing but  commutative algebras and Lie algebras in the category $\Zzvect$, respectively. A differential graded associative (Lie) algebra over $\mathbbm K$ is a associative (Lie) algebra in $\chain{\mathbbm K}$.

Given an object $\mathcal V$ in $\Gmod{A}$ we can define the \emph{tensor algebra} $T_A\mathcal V:=\oplus_{k\ge 0} T_A^k \mathcal V:=\oplus_{k\ge 0} \mathcal V^{\otimes k}$ generated by $\mathcal V$. This is an associative algebra in $\Gmod{A}$, which carries an additional $\mathbbm Z$-grading given by tensor power. The \emph{symmetric algebra} $S_A\mathcal V=T_A\mathcal V/<v\otimes w-(-1)^{|v||w|}w\otimes v>$ generated by $\mathcal V$ is obtained  by dividing  out the two sided ideal generated by expressions of the form $v\otimes w-(-1)^{|v||w|}w\otimes v$ for homogeneous  $v,w \in \mathcal V$. Similarly, the \emph{Gra{\ss}mann algebra} generated by $\mathcal V$ is given by the quotient $\wedge_A\mathcal V=T_A\mathcal V/<v\otimes w+(-1)^{|v||w|}w\otimes v>$. In addition to their natural $\Gamma$-grading, which is refered to as \emph{total degree}, $S_A\mathcal V$ and $\wedge_A\mathcal V$ inherit from the tensor algebra $\mathbbm Z$-grading, which we refer to as \emph{tensor power}. It is well known that $T_A\mathcal V$, $S_A\mathcal V$ and $\wedge_A\mathcal V$ are actually bialgebras in $\Gmod{A}$. In all three cases the comultiplication $\Delta$ is uniquely determined by the requirement that $\mathcal V$ is the space of primitives, i.e.,  $\Delta(v)=v\otimes 1+1\otimes v$ for all $v\in \mathcal V$.

The symmetric group $S_n=\Aut_{Set}(\{1,\dots,n\})$ acts on the set of multiindices $\mathbbm Z^n$ from the right by $\sigma(x_1,\dots,x_n):=(x_{\sigma(1)},\dots,x_{\sigma(n)})$.
Given such a multiindex $x=(x_1,\dots,x_n)$ the \emph{Koszul sign} $\sign(\tau_{i,i+1},x)$ of the transposition $\tau_{i,i+1}\in S_n$, which interchanges $i$ and $i+1$ is defined to be $(-1)^{x_ix_{i+1}}$. It is well known that according to the rule
\begin{align}
\sign(\sigma\circ\sigma',x)=\sign(\sigma,\sigma'(x))\sign(\sigma',x)
\end{align}
the Koszul sign unambiguously extends to a map $\sign(\sigma,x):S_n\to\{\pm1\}$. Instead of the Koszul sign of $x=(x_1,\dots,x_n)$ one may also consider the Koszul sign of the shifted multiindex $x[1]:=(x_1+1,\dots,x_n+1)$.  

\section{Derived brackets, Gerstenhaber algebras, etc.}\label{derbrsection}

A \emph{left Leibniz} (or \emph{Loday}) bracket of degree $n$ on a $\mathbbm Z$-graded vector space $L=\oplus_i L^i$ is a graded linear map $[\:,\:]:L\otimes L\to L$ of degree $n$, such that the following Leibniz rule holds for all homogeneous $a,b,c \in L$:
\begin{align}\label{leibniz}
[a,[b,c]]=[[a,b],c]+(-1)^{(|a|+n)(|b|+n)}[b,[a,c]].
\end{align}
If the bracket  $[\:,\:]$ is graded antisymmetric, i.e., $[a,b]=-(-1)^{|a|+n)(|b|+n)}[b,a]$, then it is a Lie bracket of degree $n$; in other words:, $L[-n]$ is a graded  Lie algebra. If the space $L$ is, in addition, a super-commutative $\mathbbm Z$-graded algebra with multiplication $\mu(a\otimes b)=ab$, such that the following Leibniz rule holds
\begin{align}\label{leibgerst}
[a,bc]=[a,b]c+(-1)^{|b|(|a|+n)}b[a,c],
\end{align}  
then we say that $[\:,\:]$ is a Poisson bracket of degree $n$. Under these circumstances we also say: $\left(L,\mu,[\:,\:]\right)$ is an \emph{$n$-Poisson algebra}. Thus a $0$-Poisson algebra  is just a \emph{$\mathbbm Z$-graded Poisson algebra} and a $(-1)$-Poisson algebra is what is usually called a \emph{Gerstenhaber algebra}.
\begin{SATZ}[Derived brackets] \label{derKS} Let $(L,[\:,\:])$ be a $\mathbbm Z$-graded left Leibniz algebra, where the  bracket $[\:,\:]$ is of degree $n$, and let $d:L^\bullet\to L^{\bullet+m}$ a differential (i.e., $d^2=0$) which is a derivation of $[\:,\:]$, i.e.,
\begin{align}\label{mder}
d[a,b]=[da,b]+(-1)^{m(|a|+n)}[a,db]
\end{align}
for all homogeneous $a,b\in L$. 
\begin{enumerate}
\item
Then the derived bracket 
\begin{align*}
[a,b]_d:=(-1)^{m(n+|a|)+1}[da,b]
\end{align*}
is a left Leibniz bracket of degree $n+m$. Moreover, $d$ is a derivation of $[\:,\:]_d$. 
\item \label{derbr1}If, more specifically, $[\:,\:]$ is a Lie bracket of degree $n$ and $L_0$ is an abelian subalgebra of $L$, such that $[L_0,L_0]_d\subset L_0$ and $m$ is odd, then $[\:,\:]_d$ is a Lie bracket on $L_0$ of degree $n+m$.
\item \label{derbr2} There is an natural left Leibniz bracket of degree $n+m$ on $L/dL$: 
\begin{align}\label{coex}
[\bar{a},\bar{b}]_d:=(-1)^{m(n+|a|)+1}\overline{[d{a},{b}]},
\end{align} 
where the bar indicates taking classes in $L/dL$. The left Leibniz algebra $L/dL$ contains the homology space  $\mathrm H=\mathrm Z/dL$, where $\mathrm Z=\ker(d)$, naturally as a graded subalgebra. If the original bracket is a graded Lie bracket of degree $n$ and $m$ is odd, then the derived bracket $[\:,\:]_d$ is a Lie bracket of degree $n+m$. In this case $\big(L/dL,[\:,\:]_d\big)$ is also called the \emph{Lie algebra of co-exact elements} of $L$.
\end{enumerate} 
If the original bracket is a Poisson bracket of degree $n$ and $d$ is an odd derivation of this Poisson structure, then the derived bracket defined in \ref{derbr1}. is a Poisson bracket of degree $n+m$. Moreover, the derived bracket of equation (\ref{coex}) restricted to the homology subalgebra $\mathrm H=\mathrm Z/dL$ a Poisson bracket of degree $n+m$. 
\end{SATZ}
\begin{BEWEIS} We proof these statements along the lines of \cite{KosSchw96}. First of all, let us write $\ad: L\to\End_{\Zmod{\mathbbm K}}(L)$, $a\mapsto [a,\:]$ for the adjoint action of $L$ on itself. The degree of the endomorphism $\ad(a)$ is $|a|+n$. We rewrite equation (\ref{leibniz}) as  $[\ad(a),\ad(b)]=\ad\left(\ad(a)b\right)$, this time $[\:,\:]$ denotes the super-commutator of graded endomorphisms. Analogously, we rewrite equation (\ref{mder}) as $[d,\ad(a)]=\ad(da)$. For the adjoint action corresponding to the derived bracket $[\:,\:]_d$ we write $\ad_d(a):=[a,\:]_d$. The degree of the endomorphism $\ad_d(a)$ is $|a|+n+m$. In fact, using super-commutators we have $\ad_d(a)=[\ad(a),d]$. An easy calculation using the Jacobi identitiy for the super-commutator and $[d,d]=2d^2=0$ yields that $d$ is a derivation of the derived bracket:
\begin{align*}
[d,\ad_d(a)]=\big[d,[\ad(a),d]\big]=\big[[d,\ad(a)],d\big]+0=[\ad(da),d]=\ad_d(da).
\end{align*}
Let us proof that $[\:,\:]_d$ is a left Leibniz bracket:
\begin{eqnarray*}\lefteqn{\left[\ad_d(a),\ad_d(b)\right]-\ad_d\left(\ad_d(a)b\right)}\\
&=&\big[[\ad(a),d],[\ad(b),d]\big]-\big[\ad\big([\ad(a),d]b\big),d\big]\\
&=& (-1)^{m(|a|+|b|)}\big[\ad(da),\ad(db)\big]+(-1)^{m(|a|+n)}\big[\ad\big(\ad(da)b\big),d\big]\\
&=&(-1)^{m(|a|+|b|)}\ad\big(\ad(da)db\big)-(-1)^{m(|a|+n)}(-1)^{m(m+2n+|a|+|b|)}\ad\big(d\ad(da)b\big)\\
&=&(-1)^{m(|a|+|b|)}\ad\big(\ad(da)db\big)-(-1)^{m(|b|+n+1)}(-1)^{m(m+n+|a|)}\ad\big(\ad(da)db\big)=0.
\end{eqnarray*}
In order to proof item \ref{derbr1}.) we only need to show that the restriction of the derived bracket to the subalgebra is graded antisymmetric. In fact, since $L_0$ is an abelian with respect to the graded Lie bracket $[\:,\:]$ we have: \begin{align}\label{triveq}
0=d[a,b]=[da,b]+(-1)^{m(|a|+n)}[a,db]=[da,b]-(-1)^{(|a|+n)(|b|+n)}[db,a].
\end{align}
Thus we have to see that $(-1)^{m(|a|+n)+1}[da,b]=(-1)^{m(|a|+n)+1}(-1)^{(|a|+n)(|b|+n)}[db,a]$ coincides with $-(-1)^{m(|b|+n)+1}(-1)^{(|a|+n+m)(|b|+n+m)}[db,a]$. These signs match iff $m$ is odd. 

In order to proof item \ref{derbr2}.) we observe that $dL$ is a two-sided ideal in $L$: $[dL,L]_d\subset dL\supset [L,dL]_d$. We conclude that the derived bracket of equation (\ref{coex}) is well-defined on $L/dL$. Moreover, since the space of cycles $Z=\ker(d)$ is obviously a subalgebra of $(L,[\:,\:]_d)$, and hence $H=Z/dL$ is a subalgebra of $(L/dL,[\:,\:]_d)$. If $[\:,\:]$ is in fact a Lie bracket of degree $n$ and $m$ is odd we again have, due to equation (\ref{triveq}), that the derived bracket is Lie of degree $n+m$.

Finally, let us address the question what happens in the Poisson case. Let $\mu: L\otimes L\to L, a\otimes b\mapsto ab$ be a super-commutative multiplication of degree zero. Accordingly, the operator of left multiplication $\lambda: L\to\End_{\Zmod{\mathbbm K}}(L)$, $\lambda(a)b:=ab$ for $b\in L$, is of degree zero. We impose the Leibniz rules $[d,\lambda(a)]=\lambda(da)$  (i.e., $d$ is a derivation of $\mu$) and $[\ad(a),\lambda(b)]=\lambda(\ad(a)b)$, which is equivalent to (\ref{leibgerst}). The analogue of equation (\ref{leibgerst}) for the derived bracket works out as follows
\begin{eqnarray*}\big[\ad_d(a),\lambda(b)\big]&=&\big[[\ad(a),d],\lambda(b)\big]=\big[\ad(a),[d,\lambda(b)]\big]+(-1)^{m|b|}\big[[\ad(a),\lambda(b)],d\big]\\
&=&\big[\ad(a),\lambda(db)\big]+(-1)^{m|b|}\big[\lambda(\ad(a)b),d\big]\\
&=&\lambda(\ad(a)db)-(-1)^{m|b|}(-1)^{m(|a|+|b|+n)}\lambda(d\ad(a)b)\\
&=&\lambda\big((\ad(a)d-(-1)^{m(|a|+n)}d\ad(a))b\big)=\lambda(\ad_d(a)b).
\end{eqnarray*}
Consequently, if the abelian subalgebra $L_0$ of item \ref{derbr1}.) is  also a subalgebra for the multiplication $\mu$, then the derived bracket is a Poisson bracket. In general, the subspace $dL\subset L$ is not an ideal for the multiplication $\mu$, but $dL\subset \mathrm Z$ actually is an ideal! For every $x\in dL\,\mathrm Z$ can be written as $x=\sum_{i=1}^k da_i\:b_i=d\big(\sum_i^k a_ib_i\big)$ for some  $a_1,\dots a_k\in L$ and $b_1,\dots b_k\in Z$. 
\end{BEWEIS}
\\

There is yet another way to produce a derived bracket, which generalizes the construction of item \ref{derbr2}.) of the above theorem. Let $(L,\mu,[\:,\:],d)$ be a Poisson algebra with a bracket of degree $n$ and differential of degree $m$ (assumed to be $\pm 1$), such that $d$ is a derivation of the supercommutative multiplication $\mu(a\otimes b)=ab$ and of the bracket $[\:,\:]$. Under these circumstances we also we call  $(L,\mu,[\:,\:],d)$ a \emph{differential graded n-Poisson algebra}. We will call a graded subspace $K\subset L$  a \emph{coisotropic ideal} if the following conditions are true:
\begin{eqnarray*}
K\cdot L&\subset&K\\
\:[K, K]&\subset& K\\
dK&\subset&K.
\end{eqnarray*} 
Given such a coisotropic ideal $K\subset L$ the quotient space $V:=L/K$ is a differential graded commutative algebra. In particular $V$ is a complex. For the differential we write just $d$ and for the graded spaces of cycles and boundaries of this complex we will write $\mathrm Z\,V$ and $dV$, the homology is $\mathrm H\,V:=\mathrm Z\,V/dV$. Let us write $\widetilde a$ for  the image of $a\in L$ under the projection onto $V$. Note that a cycle in $V$ is a class $\widetilde a$ of an $a\in V$ such that  $da\in K$.
\begin{SATZ} [Reduced bracket]\label{redderbr} Let $a,b\in L$ represent the cycles $\widetilde a,\widetilde b\in \mathrm Z\,V$ then 
\begin{align}
[\widetilde a,\widetilde b]_{d,K}:=(-1)^{m(|a|+n)+1}\widetilde{[da,b]}+dV
\end{align}
is a well defined class in $\mathrm H\,V$. $[\:,\:]_{d,K}$ is a Poisson bracket of degree $n+m$ on the super-commutative $\mathbbm Z$-graded algebra $\mathrm H\, V$.
\end{SATZ}
\begin{BEWEIS} In order to see, that the bracket is in fact well defined, it is more comfortable to write the cycles as $\widetilde{a}=a+K$, $\widetilde{b}=b+K$. Now we get
\begin{eqnarray*}
[d(a+K),b+K]=[da,b]+[dK,K]+[da,K]+[dK, b].
\end{eqnarray*}
The second term is obviously in $K$. The third term is in $K$ since $\widetilde{a}$ is a cycle. Finally, writing the last term as $\pm d[K,b]\pm[K,db]$, it is in $K$ up to a boundary since $\widetilde{b}$ is a cycle. Note that $[da,b]$ is a cycle since $d[da,b]=\pm[da,db]\subset[K,K]$.  Again, $dV$ is an ideal in $\mathrm Z\,V$ with respect to the commutative multiplication. The remaining statements follow from the proof of Theorem \ref{derKS}.
\end{BEWEIS}
\\

In any of the cases when the derived bracket $[\:,\:]_d$ is a Poisson bracket of degree $n+m$ the \emph{opposite bracket}  $[a,b]^{\oppind}_d:=(-1)^{(|a|+n+m)(|b|+n+m)}[b,a]_d$, for homogeneous $a,b \in L$, is a  Poisson bracket of degree $n+m$ as well (see Appendix \ref{oppPoiss}). Since we have 
\begin{eqnarray*}
[a,db]&=&-(-1)^{(|b|+m+n)(|a|+n)}[db,a]=(-1)^{m(|b|+n)}(-1)^{(|b|+m+n)(|a|+n)}[b,a]_d\\
&=&(-1)^{(|a|+n+m)(|b|+n+m)}(-1)^m [b,a]_d,
\end{eqnarray*}
we conclude that the opposite derived bracket is just
\begin{align}
[a,b]^{\oppind}_d:=(-1)^m[a,db].
\end{align}

An important example of Poisson bracket of degree $-1$ is given by the Schouten-Nijenhuis bracket on the  symmetric algebra $S_A(L[-1])$ of a Lie-Rinehart pair $(A,L)$ in the category $\Zmod{\mathbbm K}$. For convenience of the reader we have collected some basic material on Lie-Rinehart pairs in Appendix \ref{LieRinehart}. We caution the reader that for certain applications that we have in mind, this notion of polyvector field is  a little too restricitive. 

\begin{SATZ}[Schouten-Nijenhuis bracket]\label{SchoutenNijenhuis} If $(A,L)$ is a Lie-Rinehart pair in the tensor category $\Zmod{\mathbbm K}$, then there is a unique Lie bracket $[\:,\:]$ on $S_A(L[-1])$ extending the structure maps on $(A,L)$ and making $S_A(L[-1])$ into a Gerstenhaber algebra. More precisely, the bracket is the unique Gerstenhaber bracket determined by the requirements 
\begin{enumerate}
\item $A\subset S_A(L[-1])$ is an abelian subalgebra, i.e., $[A,A]=0$,
\item $[X,a]=-(-1)^{(|X|+1)|a|}[a,X]=X(a)$, 
\item if $[\:,\:]'$ denotes, for the moment, the bracket in $L$ then we have $[X,Y]'=[X,Y]$
\end{enumerate}
for all homogeneous $X,Y\in L$ and $a\in A$.
Moreover, any morphism from $(A,L)$ to another Lie Rinehart pair $(A',L')$ in $\Zmod{\mathbbm K}$ extends uniquely to a morphism of Gerstenhaber algebras $S_A(L[-1])\to S_{A'}(L'[-1])$.
\end{SATZ}
\begin{BEWEIS} Presumably, there is some simple `operadic' reason for the statement. Nonetheless, we have included an elementary proof in the appendix \ref{Schoutenproof}. 
\end{BEWEIS} 
\\

We will call this Gerstenhaber algebra the \emph{algebra of polyvector fields} of the Lie-Rinehart pair $(A,L)$. We will use for it the notation $\mathfrak X(A,L):=S_A(L[-1])$. It carries two degrees: the \emph{total degree} and the \emph{tensor power}. For the tensor power we will also use the term \emph{arity}. For the subspace of arity $k$ we will write
\[\poly^k(A,L)=S^k_A(L[-1]).\]
If $a\in A$ and $X_1,\dots,X_i\in L$ are homogeneous, then the total degree of a monomial $a X_1\dots X_i$ is $|a|+|X_1|+\dots+|X_i|+i$. The space of polyvector fields of total degree $k$ will be written as 
\[\poly(A,L)^k.\]
Of course, the degree for which $\poly(A,L)$ is a Gerstenhaber algebra is the total degree. However, we also have an inclusion $[\poly^k(A,L),\poly^l(A,L)]\subset \poly^{k+l-1}(A,L)$ for all $k,l\ge 0$.

For the  important special case, where $A=\mathcal C^\infty(M)$ is the ring of smooth functions on some manifold $M$ and $L=\Gamma^\infty(M,TM)$ is the space of vector fields we will use the shorthand $\poly(M)$. If $\Pi\in\poly^2(M)$ is a Poisson tensor, i.e., $[\Pi,\Pi]=0$, then 
\begin{align}
\delta_{\Pi}:=[\Pi,\:]:\poly^\bullet(M)\to \poly^{\bullet+1}(M)
\end{align} is a codifferential. The cohomology of this differential is known as \emph{Lichnerowicz-Poisson cohomology}. The Poisson bracket corresponding to $\Pi$ is just the derived bracket of $-\delta_\Pi$, or, equivalently, the opposite derived bracket of $\delta_\Pi$. Here, in the ungraded case, we will use the standard notations and identify $\poly^\bullet(M)$ with the space of sections $\Gamma^\infty(M,\wedge^\bullet TM)$ of the Grassmann-algebra bundle of the tangent bundle. The Lichnerowicz-Poisson cohomological complex $(\poly^\bullet(M),\delta_\Pi)$ can alternatively be interpreted as the space of cochains of Lie algebroid cohomology of the Lie algebroid $(T^*M,\#_\Pi,\{\:,\:\}_{KB})$ associated to $\Pi$. The anchor $\#_\Pi:T^*M\to TM$ of this Lie algebroid is given the `musical map' $\#_\Pi(\alpha):=\alpha^\#:=i(\alpha)\Pi$, where $i$ denotes the insertation derivation. The bracket $\{\:,\:\}_{KB}$ is the \emph{Koszul-Brylinski bracket} which is given by the formula $\{\alpha,\beta\}_{KB}=L_{\alpha^\#}\beta-L_{\beta^\#}\alpha-d\Pi(\alpha,\beta)$.

Given a commutative algebra $A$ in $\Zmod{\mathbbm K}$ there  is still another Gerstenhaber algebra $\mathrm{Der}(A)=\oplus_{n\ge 0}\mathrm{Der}^n(A)$, the algebra of \emph{multiderivations}, which is related to the algebra of polyvector fields. A multiderivation $D\in \mathrm{Der}^n(A)$ of arity $n$ is by definition a graded linear map $D:A^{\otimes n}\to A$
which is graded symmetric in the following sense
\begin{align}\label{symmprop}
D(a_1,\dots,a_i,a_{i+1},\dots,a_n)=(-1)^{(|a_i|+1)(|a_{i+1}|+1)}D(a_1,\dots,a_{i+1},a_{i},\dots,a_n)
\end{align}
for $i=1,\dots,n-1$ and which is a derivation in every argument, i.e.,
\begin{align}\label{derprop}
D(a_1,\dots,a_{n-1},ab)=D(a_1,\dots,a_{n-1},a)b+ (-1)^{|a||b|}D(a_1,\dots,a_{n-1},b)a
\end{align}
for all homogeneous $a,b,a_1,\dots,a_n\in A$. If $D$ is a multiderivation of arity $n$ and degree $k$, then we define the total degree  of $D$ to be $|D|:=n+k$. If $(A,L)$ is a Lie-Rinehart pair in $\Zmod{A}$, then a polyvector field $X\in \poly^n(A,L)$ can be interpreted as a multiderivation. More precisely, we let $X$ act on  $(a_1,\dots,a_n)$ by taking iterated commutators 
\begin{align}\label{itcom} 
B_X(a_1,\dots,a_n):= X(a_1,\dots,a_n):=[\dots[[X,a_1],a_2,],\dots ,a_n].
\end{align}
Let us now define an analogue $\cup$ of the super-commutative product of polyvector fields
\begin{eqnarray}\label{cupform}
D\cup E(a_1,\dots,a_{k+l})&:=&\sum_{\sigma \in S_{k,l}}\sign(\sigma,|a|[1])(-1)^{|E|(|a_{\sigma(1)}|+|a_{\sigma(2)}|+\dots +|a_{\sigma(k)}|+k)}\nonumber\\
&&\qquad\qquad D(a_{\sigma(1)},\dots,a_{\sigma(k)})E(a_{\sigma(k+1)},\dots, a_{\sigma(k+l)}),
\end{eqnarray}
where $D$ and $E$ are homogeneous multiderivations of arity $k$ and $l$, respectively, and $S_{k,l}$ is the set of  $(k,l)$-unshuffle permutations, i.e., $\sigma(1)<\sigma(2)<\dots<\sigma(k)$ and $\sigma(k+1)<\sigma(k+2)<\dots<\sigma(k+l)$. In general, by $|a|:=(|a_1|,\dots,|a_n|)\in \mathbbm Z^{n}$ we mean the multiindex of the $n$-tuple $(a_1,\dots,a_n)$ of homogeneous elements of the graded space $A$.
Next we define a bilinear operation $\bullet:\mathrm{Der}^k(A)\times\mathrm{Der}^l(A)\to\mathrm{Der}^{k+l-1}(A)$ as follows
\begin{eqnarray*}\label{bullform}
D\bullet E(a_1,\dots, a_{k+l-1}):=\sum_{\sigma\in S_{l,k-1}}\sign(\sigma,|a|[1])D(E(a_{\sigma(1)},\dots,a_{\sigma(l)}),a_{\sigma(l+1)},\dots, a_{\sigma(k+l-1)}).
\end{eqnarray*}
The Richardson-Nijenhuis bracket \cite{RichNij} of two multiderivations is defined by
\begin{align}\label{RNformula} 
[D,E]_{RN}:=D\bullet E-(-1)^{(|D|-1)(|E|-1)}E\bullet D.
\end{align}
\begin{SATZ}\label{multiderisgerst}$\left(\Der(A),\cup,[\:,\:]_{RN}\right)$ is a Gerstenhaber algebra. The iterated bracket map $B$ (cf. equation (\ref{itcom})) is a morphism of Gerstenhaber algebras, i.e., for all $X,Y \in \poly(A,L)$ we have $B_{XY}=B_X\cup B_Y$ and $B_{[X,Y]}=[B_X,B_Y]_{RN}$.
\end{SATZ}
\begin{BEWEIS}
It is well known (see e.g. \cite[Th\'eor\`eme III.2.1]{Manchon}) that the space \[\mathscr L:=\prod_n\Sym_{\mathbbm K}^n(A[1],A[1])\]
of graded $\mathbbm K$-linear maps $A[1]^{\otimes n}\to A[1]$ which have the symmetry property (\ref{symmprop}) above is isomorphic to the graded Lie algebra $\Coder(S_{\mathbbm K}(A[1]))$ of coderivations of the free symmetric $\mathbbm K$-coalgebra cogenerated by $A[1]$. More specifically, if $D\in\Sym_{\mathbbm K}^n(A[1],A[1])$ then the corresponding coderivation $D$ of $S_{\mathbbm K}(A[1])$ is given by the formula 
\begin{align*}
D(a_1\dots a_m):=\sum_{\sigma\in S_{n,m-n}}\sign(\sigma,|a|[1])D(a_{\sigma(1)},\dots,a_{\sigma(n)})\:a_{\sigma(n+1)}\dots a_{\sigma(n+m)}.
\end{align*} 
The induced Lie bracket on $\mathscr L$ is just the bracket of equation (\ref{RNformula}). It is not difficult to see that the space of multiderivations $\Der(A)$ is in fact a Lie subalgebra of $\mathscr L$. In this way a multiderivation $D$ of degree $|D|$ is a coderivation of degree $|D|-1$. In order to show the Leibniz rule for $[\:,\:]_{RN}$ it suffices to prove 
\begin{eqnarray*}
D\bullet(E\cup F)&=&(D\bullet E)\cup F+(-1)^{(|D|-1)|E|}E\cup (E\bullet F)\qquad \mbox{ and} \\
(E\cup F)\bullet D&=&E\cup (F\bullet D)+(-1)^{(|D|-1)|F|}(E\bullet D)\cup F.
\end{eqnarray*}
These relations can be shown by an unpleasant but straightforward computation. The compatibility of $B$ with the cup product $B_{XY}=B_X\cup B_Y$ is straight forward to check. In order to prove $B_{[X,Y]}=[B_X,B_Y]_{RN}$ one can  now use induction over the arities of $X$ and $Y$. To start the induction we need only to check the formulas $B_{[X,Y]}=[B_X,B_Y]_{RN}$ and $[X,a]=[B_X,a]_{RN}$ for $X,Y\in\poly^1(A,L)$ and $a\in A$, which are obviously fulfilled.
\end{BEWEIS}
\section{Generalized graded manifolds}\label{ggmf}
There are several possible ways to define the notion of a super-manifold (see e.g. \cite{superkostant,tuynman,Eckel}). The most basic version is that of a space of sections $\Gamma^\infty(M,\wedge E)$ of a Gra{\ss}mann algebra bundle $\wedge E\to M$ of a finite dimensional vector bundle over $M$ seen as a $\mathbbm Z$-graded super-commutative $\mathcal C^\infty(M)$-algebra.   
It is known that many geometric constructions, such as tangent vectors, the cotangent bundle, the de Rham complex and vector bundles with  connections generalize to super-manifolds. 

In view of the applications we have in mind (cf. Sections \ref{RothPB}, \ref{projkostate} and \ref{BRSTchargesection}) we are forced to consider more general $\mathbbm Z$-graded algebras which serve as algebras of super-functions. We will propose a way to define the notion of cotangent space for these algebras and introduce the notion of polyvector fields in the spirit of the preceding section. Besides, we do not think that we have found a general geometric theory for this more general graded algebras (and it may well be that there are important references related to these questions, which we are not aware of). Let $M$ be a smooth manifold and let $V=\oplus_{k \in \mathbbm Z} V^k$ be the direct sum of finite rank vector bundles $V^k$ over $M$. For simplicity, let us always assume that these vector bundles admit simultaneous bundle charts. Let us write $\mathcal V$ for the space $\Gamma^\infty(M,V)$ of smooth sections of $V$. For convenience we use the shorthand notation $A:=C^\infty(M)$. According to the theorem of Serre and Swan we have that $\mathcal V=\oplus_{k\in \mathbbm Z} \mathcal V^k$ is the direct sum of the finitely generated projective $A$-modules $\mathcal V^k=\Gamma^\infty(M,V^k)$. The support of $\mathcal V$ is defined to be $\operatorname{supp}(\mathcal V):=\{k\in\mathbbm Z\mid \mathcal V^k< 0\}$. Assigning to each $\mathcal V^k$ the degree $k$ we will consider $\mathcal V$ as an object in $\Zmod{A}$. 

Let $S_A\mathcal V=\oplus_{i\ge 0} S^i_A \mathcal V$ be the symmetric algebra in $\Zmod{A}$ generated by the $A$-module $\mathcal V$. 
In fact, $S_A\mathcal V$ is a free commutative algebra in the category $\Zmod{A}$. It is easy to see, that the $A$-linear derivations $\Der_A(S_A\mathcal V)^k=\Der^1_A(S_A\mathcal V)^k$ of degree $k$ of the $\mathbbm Z$-graded algebra $S_A\mathcal V$ may be identified with the $A$-module $\prod_i(S_A\mathcal V)^{k+i} \otimes_A {\mathcal V^i}^*$. In every degree $k$ we have an exact sequence
\begin{align}\label{derseq}
0\to\Der_A(S_A\mathcal V)^k\to\Der_{\mathbbm K}(S_A\mathcal V)^k\stackrel{D\mapsto D_{|A}}{\to}\Der_{\mathbbm K}(A,S_A\mathcal V)^k\to 0.
\end{align}
Let us now choose for every vector bundle $V^i$ a connection $\nabla^{V^i}$. For the derivation of $S_A\mathcal V$ which gives for every section $v_j\in\mathcal V^j$  and every vector field $X\in \poly(M)$ the value $\nabla^{V^j}_{X}v_j$ we will write  
\begin{align}
\nabla_{X}=\sum_{j\in \mathbbm Z}\nabla^{V^j}_{X}.
\end{align} 
Using this family of connections we can identify  $\Der_{\mathbbm K}(A,S_A\mathcal V)^k=(S_A\mathcal V)^k\otimes \Der_{\mathbbm K}A=(S_A\mathcal V)^k\otimes \poly^1(M)$ with a submodule of $\Der_{\mathbbm K}(S_A\mathcal V)^k$ which is complementary to $\Der_A(S_A\mathcal V)^k$, i.e., we use the connections to split the sequence (\ref{derseq}). More precisely, this split is the $S_A\mathcal V$-linear extension of the map which associates to every vector field $X$ the derivation $\nabla_X$. This complementary subspace to the space of $A$-linear derivations will be refered  to as the space of \emph{geometric derivations}. Using the exact sequence (\ref{derseq}) one can show that $\mathrm{Der}_\mathbbm K(S_A\mathcal V)$ is a projective $S_A\mathcal V$-module. Hence, one can use the machinery of Appendix \ref{LieRinehart} to define the Lie-Rinehart cohomology of the graded Lie-Rinehart pair  $(S_A\mathcal V,\mathrm{Der}_\mathbbm K(S_A\mathcal V))$.

The reader may have noticed that, in general, the $A$-linear derivations as well as the geometric derviations appear to be infinite sums. We therefore introduce the space of \emph{finite type derivations}, which are derivations as above, such that the sums are actually finite sums. In order to give the precise definition let us write the derivations in local coordinates. Let $(x_1,\dots,x_n)$ be local coordinates for $M$ and let $\xAG{i}{1},\dots ,\xAG{i}{\ell_i}$ be  local frames for the bundles $V^i$ for all $i\in \mathbbm Z$. Then a  \emph{finite type derivations} of degree $k$ is a derivation $X$ which writes locally as
\begin{eqnarray}
{\sum_{j\in \mathbbm Z}}^\prime \sum_{a=1}^{\ell_j} X_a^j\:\ddAG{i}{a}+{\sum_{j\in \mathbbm Z}}^\prime\sum_{i=1}^n X_i\:\nabla^{V^j}_{\frac{\partial}{\partial x^i}},
\end{eqnarray}
where the $^\prime$  indicates that all except finitely many summands vanish. In the above formula $\ddAG{i}{a}$ is the unique $A$-linear derivation such that\begin{align} \label{algebraicder}
\frac{\partial \xAG{j}{b}}{\partial \xAG{i}{a}}:=\ddAG{i}{a}\xAG{j}{b}=\delta_a^b\:\delta_i^j\qquad \forall i,j\in\mathbbm Z,\:a=1,\dots \ell_i,\:b=1,\dots,\ell_j.
\end{align}
Of course, the coefficients $X_a^j$ are in $(S_A\mathcal V)^{j+k}$ and the $X_i$ have to be in $(S_A\mathcal V)^k$. We write for the graded space of finite type derivations $\poly^1(S_A\mathcal V)$. It is clear that $(S_A\mathcal V,\poly^1(S_A\mathcal V))$ form a graded Lie-Rinehart pair. Of course, if $|\operatorname{supp}(\mathcal V)|<\infty$ then every graded derivation is of finite type. Note that if $|\operatorname{supp}(\mathcal V)|=\infty$ the \emph{Euler vector field}
\begin{align*}
\sum_{i\in \mathbbm Z} \sum_{a=1}^{\ell_i} \xAG{i}{a}\:\ddAG{i}{a}
\end{align*}
is not a finite type derivation. Analogously, $\nabla_X$ is a finite type derivation if and only if $|\operatorname {supp}(\mathcal V)|<\infty$.

For the applications we have in mind  (cf. Section \ref{projkostate}) we have that the module $\mathcal V=\oplus_{k\ge 1}\mathcal V^k$ is positively graded. This entails that the graded components of $S_A\mathcal V$ are actually finitely  generated. In fact, the dimensions can be determined recursively from the Poincar\'e series
\begin{align}
\sum_{i\ge 0}\operatorname{rank}((S_A\mathcal V
)^i)\:t^i=\prod_{j\ge 1}\big(1-(-t)^j\big)^{(-1)^{j+1}\operatorname{rank}(\mathcal V^j)}.
\end{align}
The product on the right hand side converges in the $t$-adic topology of $\mathbbm Z[[t]]$.
\begin{DEFINITION} By a \emph{generalized graded manifold} with base manifold $M$, $A:=\mathcal C^\infty(M)$, we mean a Lie-Rinehart pair in the category $\Zmod{A}$of the form $\big(S_A\mathcal V,\poly^1(S_A\mathcal V)\big)$, where $\mathcal V=\oplus_{k\ge 1}\mathcal V^k$ is positively graded. We will say for short: $S_A \mathcal V$ is a graded manifold. The  generalized graded manifold is said to be \emph{finitely generated} of level $L:=\inf\{j\:|\: \mathcal V^i=0\:\forall i>j\}$ if $L<\infty$. Otherwise it is said to be \emph{infinitely generated}. In the finitely generated of level $L$ case we will sometimes write $\mathcal V^{\le L}$ instead of $\mathcal V$. If all the vector bundles $V^k$, $k=1,2,\dots$ are trivial, we will say that the generalized graded manifold is \emph{freely generated}.
\end{DEFINITION}

For the BFV-construction we also need to adjoin extra variables which are dual to the local frames $\xAG{i}{1},\dots ,\xAG{i}{\ell_i}$ for the vector bundles $ V^i$, $i=1,2,\dots$. This construction is a graded analog of taking the cotangent bundle of a manifold. Let ${\mathcal V^i}^*=\Gamma^\infty(M,{V^i}^*)$ be space of sections of the dual vector bundle of $V^i$ with dual local frames $\xG{i}{1},\dots ,\xG{i}{\ell_i}$ (these sections will be also called momenta). The naive guess for an analog of cotangent space for the graded manifold $S_A\mathcal V$  would probably be the $\mathbbm Z$-graded algebra $S_A(\mathcal V\oplus \mathcal V^*)=S_A\mathcal V^*\otimes S_A\mathcal V$, where $\mathcal V^*=\oplus_{i\ge i}{\mathcal V^i}^*$. Unfortunately, this proposal does not work for our purposes (i.e. the constructions of the BRST charge in Section \ref{BRSTchargesection} and of the Rothstein bracket in Section \ref{RothPB}). The `correct' cotangent space $\brsalg=\oplus_{k\in\mathbbm Z} \brsalg^k$, which is considerably bigger, is given by
\begin{align}
\brsalg^k=\prod_{j\ge 0} (S_A\mathcal V^*)^{k+j}\otimes (S_A\mathcal V)^j.
\end{align}
It is important to note that in this definition we have reversed the natural grading: if the variables $\xAG{i}{1},\dots ,\xAG{i}{\ell_i}$ are considered as elements of $\brsalg$ they have degree $-i$. Analogously, if the momenta  $\xG{i}{1},\dots ,\xG{i}{\ell_i}$  are considered as elements of $\brsalg$ they have dgeree $i$ (which is $-1$ times their natural degree). The multiplication for $S_A(\mathcal V\oplus \mathcal V^*)$ extends naturally to a $\mathbbm Z$-graded super-commutative multiplication on $\brsalg$. More specifically, if $\sum_{i\ge 0}\alpha_{i+n}\otimes v_i\in \brsalg^n$ and  $\sum_{j\ge 0}\beta_{j+m}\otimes w_j\in \brsalg^m$ then the product is defined to be 
\begin{align}
\left(\sum_{i\ge 0}\alpha_{i+n}\otimes v_i\right)\left(\sum_{j\ge 0}\beta_{j+m}\otimes w_j\right):=\sum_{k\ge 0}\sum_{i+j=k} (-1)^{i(j+m)} \alpha_{i+n}\beta_{j+m}\otimes  v_i w_j.
\end{align}
In calculations we will usually drop the $\otimes$-sign for convenience.

Even though we will not make use of it, let us mention the following interpretation of the algebra $\brsalg$. Recall that the symmetric algebra $S_A\mathcal V$ is actually super-commutative (let us write for the multiplication $\mu$), super-cocommutative bialgebra (for the comultiplication we write $\Delta$). Moreover, there is a canonical isomorphism between the $\mathbbm Z$-graded $A$-modules $\brsalg^k$ and $\End_{\Zmod{A}}^{-k}(S_A\mathcal V)$. The bialgebra structure on $S_A\mathcal V$ induces on $\End_{\Zmod{A}}^{-\bullet}(S_A\mathcal V)$ the structure of a super-commutative algebra in $\Zmod{A}$ with respect to the \emph{convolution product} $\varphi\diamond \psi:=\mu\circ(\varphi\otimes \psi)\circ \Delta$. We claim that the aforementioned canonical isomorphism is an isomorphism of graded algebras.

Let us now address the more subtle issue of derivations of the algebra $\brsalg$. In analogy to the derivations $\ddAG{i}{1},\dots,\ddAG{i}{\ell_i}$ for $i=1,2,\dots$ (defined in equation (\ref{algebraicder})) we now have in addition the $A$-linear derivations 
\begin{align}
 \ddG{i}{1},\dots,\ddG{i}{\ell_i}\qquad \mbox{for }i=1,2,\dots. 
\end{align}
The above $A$-linear derivations do pairwise super-commute. 
Let us stipulate that the connections $\nabla^{V^i}$ and $\nabla ^{{V^i}^*}$ are chosen in such a way that 
\begin{align} 
X<\alpha,v>=<\alpha,\nabla^{V^i}_Xv>+<\nabla ^{{V^i}^*}_X\alpha,v>
\end{align}
for all $X\in\poly^1(M)$, $\alpha\in{\mathcal V^i}^*$ and $v\in\nabla^{V^i}$. In other words, the connection $\nabla^i:=\nabla^{V^i}+\nabla ^{{V^i}^*}$ on the super-Riemannian vector bundle $({V^i}^*\oplus V^i,g:=<,>)$ is required to be \emph{metric}. Here, super-Riemannian means that $g:=<,>$ is symmetric for all odd $i$ and antisymmetric for even $i$. The above derivations are derivations of $S_A(\mathcal V\oplus \mathcal V^*)$ and $\brsalg$ as well. Unfortunately, not every graded (and infinite) linear combination over $\brsalg$ of the above derivations gives a well-defined operation. In order to define the sensible tangent space to $\brsalg$ we therefore consider again finite type derivations. A derivation $X$ of $\brsalg$ on degree $k$ is said to be of \emph{finite type} if it writes in local coordinates as follows
\begin{align}
{\sum_{j\in \mathbbm Z}}^\prime \sum_{a=1}^{\ell_j} X_a^j\:\ddAG{j}{a}+
{\sum_{j\in \mathbbm Z}}^\prime \sum_{a=1}^{\ell_j} X^a_j\:\ddG{j}{a}+
{\sum_{j\in \mathbbm Z}}^\prime\sum_{i=1}^n X_i\:\nabla^j_{\frac{\partial}{\partial x^i}}.
\end{align} 
Once again, the prime indicates that all except finitely many summands vanish. For the coefficients we have  $X_a^j\in \brsalg^{k-j}$, $X^a_j\in \brsalg^{k+j}$ and $X_i\in \brsalg^{k}$. The graded $A$-module of finite type derivations of $\brsalg$ will be denoted by $\poly^1(\brsalg)$. Of course, $(\brsalg,\poly^1(\brsalg))$ form a Lie-Rinehart pair in $\Zmod{A}$.
\begin{DEFINITION} If $(S_A\mathcal V,\poly^1(S_A\mathcal V))$ is a generalized graded manifold then its \emph{ghost-cotangent space} is defined to be the Lie-Rinehart pair $(\brsalg,\poly^1(\brsalg))$ in the category $\Zmod{A}$.
\end{DEFINITION}

In the infinitely generated case this definition has some drawbacks. In fact, the Euler vector field 
\begin{align}
\sum_{i\in Z}\sum_{a=1}^{\ell_i} \xAG{i}{a} \ddAG{i}{a}+\xG{i}{a} \ddG{i}{a}
\end{align}
is a finite type derivation if and only if $|\operatorname{supp}(\mathcal V)|<\infty$. The same is true for the (even) derivations which write  locally as follows
\begin{align}
\frac{\partial}{\partial \xi^i}:=\sum_{j\in\mathbbm Z}\nabla^j_{\frac{\partial}{\partial x^i}}\qquad i=1,\dots,n=\dim M.
\end{align}

Note that there  are several canonical  morphism of Lie-Rinehart pairs around. For example, there are the obvious inclusions 
\begin{align*}
(S_A\mathcal V, \poly^1 (S_A\mathcal V))\hookrightarrow(S_A(\mathcal V\oplus\mathcal V^*), \poly^1 (S_A(\mathcal V\oplus\mathcal V^*)))\hookrightarrow(\brsalg ,\poly^1(\brsalg)).
\end{align*} 
Moreover, we have an obvious surjection
\begin{align}
(\brsalg ,\poly^1(\brsalg))\stackrel{\epsilon^\prime}{\to} (A=\mathcal C^\infty(M),\poly^1(M)).
\end{align}
More specifically, $\epsilon^\prime_{|\brsalg}$ is just the canonical augmentation map $\brsalg \to A$. The image under $\epsilon^\prime$ of an $A$-linear derivation is defined to be zero. Finally, if $\mathcal V^j\ne 0$ then we define 
\begin{align*}
\epsilon^\prime(\nabla^j_{\frac{\partial}{\partial x^i}}):=\frac{\partial}{\partial x^i}
\end{align*} 
(otherwise the latter is defined to be zero). Note that this map is in fact compatible with the brackets, since the curvature terms are killed by the augmentation map.

Following the ideas of Section \ref{derbrsection}, we define the space of polyvector fields $\poly^\bullet(\brsalg)$ to be that which associated to the Lie-Rinehart pair $(\brsalg,\poly^1(\brsalg))$. It will become clear, that (at least for the infinitely generated case) the more correct definition would be the space of multiderivations $\Der^{\bullet}(\brsalg)$ of the $\mathbbm Z$-graded algebra $\brsalg$. Unfortunately, computations in  this object tend to be clumsy. Note that, according to Theorem \ref{SchoutenNijenhuis}, the morphism $\epsilon'$ extends to a morphism of Gerstenhaber algebras
\begin{align}
\epsilon': \poly^\bullet(\brsalg)\to \poly^\bullet(M).
\end{align}
In the finitely generated case we prefer to work with a morphism $\epsilon$ which is obtained from $\epsilon^\prime$ by a renormalization. This $\epsilon$ is defined in the same way as $\epsilon^\prime$ except that 
\begin{align}
\epsilon(\nabla^j_{\frac{\partial}{\partial x^i}})=|\operatorname{supp}(\mathcal V)|^{-\half}\frac{\partial}{\partial x^i}.
\end{align}
\section{ The Rothstein-Poisson bracket}\label{RothPB}
In \cite{Rothstein} M. Rothstein has shown, that on every super-manifold with symplectic base there is  a super-symplectic 2-form. This fact has been rediscovered by M. Bordemann \cite{SuperB} as a byproduct of the Fedosov construction for super-manifolds. Here we shall extend this result to ghost-cotangent spaces of generalized graded manifolds with a Poisson base manifold. We will treat first the case when the graded manifold is finitely generated. The reason for that is that in the finitely generated case the computations can be done in the Gerstenhaber algebra of polyvector fields defined in the preceding section. In contrast, in the infinitely generated case they make sense merely in the Gerstenhaber algebra of multiderivations. In the infinitely generated case our argument is  still incomplete.

\subsection{The finitely generated case}\label{RothPBfin}
\begin{SATZ}\label{Roth}
Let $M$ be a Poisson manifold with Poisson tensor $\Pi\in\poly^2(M)$ and let $V=V^{\le L}=\oplus_{k=1}^L V^k$ a positively graded finite dimensional vector bundle over $M$. Let us write $A:=\mathcal C^\infty(M)$ for the algebra  of smooth functions on $M$ and $\mathcal V=\oplus_{i\ge 1}^L\mathcal V^i$ for the space of smooth sections of $V$. Then there exists  a Poisson tensor $\Pi_R\in \poly^2(\ibrsalg{L})$ on the ghost cotangent space of the finitely generated graded manifold $S_A\mathcal V^{\le L}$, such that
\begin{enumerate}
\item \label{Poiss} $[\Pi_R,\Pi_R]=0$,
\item  \label{degzero}$\Pi_R$ is of total degree $2$,
\item \label{proj}the image under $\epsilon$ of $\Pi_R$ in $\poly^2(M)$ is $\Pi$.
\end{enumerate}
\end{SATZ}

The proof will be done by checking, that the explicit solution $\Pi_R$, called the \emph{Rothstein Poisson tensor}, which can be read off from \cite{Rothstein} does in fact fulfill the requirements in this more general context. 

The first thing that  we will do is to introduce a more condensed notation using the super-Riemannian metric $g$ which we already mentioned in the preceeding section. We denote by $\xi_1^{(j)},\xi_2^{(j)}\dots \xi_{\ell_j}^{(j)}$, $j=1,\dots ,L$, a local frame for $V^j$ and by  $\xi^1_{(j)},\xi^2_{(j)}\dots \xi^{\ell_j}_{(j)}$ the corresponding dual frame for $V^{j*}$. Recall that $\ell_j$ is the dimension of $V^j$ and $\ell:=\sum_j \ell_j$ is the dimension of $V=\oplus_j V^j$. On $V\oplus V^*$  metric the $g$ is given by:
\begin{align}
g\big(\xG{i}{a},\xAG{j}{b}\big):=\delta_j^i\:\delta^a_b=:(-1)^{i+1}g\big(\xAG{j}{b},\xG{i}{a}\big)\\
g\big(\xAG{i}{a},\xAG{j}{b}\big):=0=:g\big(\xG{i}{a},\xG{j}{b}\big),
\end{align}
for indices $a=1,\dots ,\ell_i$, $b=1,\dots,\ell_j$ and $i,j=1,\dots,r$. Recall that we have chosen a connection $\nabla=\sum_j \nabla^j$ on $V\oplus V^*$ which is metric with respect to $g$. Imposing the graded Leibniz rule, we let 
\begin{align}
\ddxi{i}:=\nabla_{\frac{\partial}{\partial x^i}}
\end{align} 
act on the super-commutative algebra $\ibrsalg{L}$ as even derivations. In order to simplify the computations, we reindex the above frames into a frame $\xi^1,\xi^2,\dots,\xi^{2\ell}$ for  the vector bundle $V\oplus V^*$ as follows:
\begin{align*}
\xi^A:=\xi_a^{(i)},\qquad\xi^{\ell+A}:=\xi^a_{(i)}.
\end{align*}
for $A=a+\sum_k^{i-1}\ell_k$ (here we set $\ell_0:=0$). The parity $p(A)$ of the index $A$ is defined to be the parity of $\xi^A$, which is $i(\operatorname{mod}2)$. These indices will run through the capital letters $A,B,C$, etc.

Quite importantly, we have the following commutation relations
\begin{align}
\Big[\:\ddxi{i},\ddxi{j}\:\Big]=\sum_{A,B}R^B_{Aij}\,\xi^A\,\ddxi{B}\label{commrel1}\\
\:\Big[\:\ddxi{i},\ddxi{A}\:\Big]=\sum_{B}\Gamma_{iA}^B\,\ddxi{B}\label{commrel2}
\end{align}
 for the `super-coordinate' vector fields, where $\Gamma_{iA}^B$ and $R^B_{Aij}$ are the Christoffel symbols and the components of the curvature tensor of the connection $\nabla$. Using the inverse of the metric tensor $\sum_A g_{AB}\:g^{AC}=\delta_B^C$ we define the following algebraic bivector field
\begin{align*}
\Pi_0:=\half\sum_{A,B}g^{AB}\ddxi{A}\ddxi{B}\in \poly^2(\ibrsalg{L}).
\end{align*} 
We need to define the insertation derivations $i({\xi_A})$ and $i(\frac{\partial}{\partial \xi_A})$, for $A=1,\dots,2\ell$, of the supercommutative algebra structure on $\poly(\ibrsalg{L})$, which extend the dual pairings. More precisely, we set $i({\xi_A})\xi^B:=\delta^B_A$ and $i(\frac{\partial}{\partial \xi_A})\ddxi{B}:=\delta_B^A$. Of course, the remaining generators will be killed by these derivations. Moreover, let us introduce the derivations $d:=[\Pi_0,\:]$ and  $d^*:=(-1)^{p(A)+1}\sum_{A,B}g_{AB}\:\xi^A\:
i(\frac{\partial}{\partial \xi_B})$ and the degree derivations $\degree_\xi:=\sum_A \xi^A\:i({\xi_A})$ and $\degree_{\partial\xi}:=\sum_A \ddxi{A}\:i({\frac{\partial}{\partial \xi_A}})$. Finally, we introduce the operator $d^{-1}$ as follows. If $X\in \poly(\ibrsalg{L})$ is a polyvector field such that $\degree_\xi X=m X$ and $\degree_{\partial\xi}X=n X$, then we define $d^{-1}X:=(m+n)^{-1}d^*X$. On $\operatorname{ker}(\degree_\xi)\cap\operatorname{ker}(\degree_{\partial\xi})$ we define $d^{-1}$ to be zero.
\begin{LEMMA} \label{pinull}
Since the connection $\nabla$ is metric, we have  
\begin{eqnarray}
d=[\Pi_0,\:]=\sum_{A,B}g^{AB}\ddxi{A}\:i(\xi_B)\label{adpi}\\
d^2=0.\label{Komplex}
\end{eqnarray}
Moreover, we obtain a ``Hodge identity'': $dd^*+d^* d=\degree_\xi+\degree_{\partial\xi}$, which implies that $d^{-1}$ is a contracting homotopy for $d$. More specifically, we have that 
\begin{align} \label{dcontr}
\xymatrix{(\poly^\bullet(M),0)\ar@<-0.6ex>[r]_{\iota\quad}^{}&(\poly^\bullet(\ibrsalg{L}),d),d^{-1} \ar@<-0.6ex>[l]_{\epsilon\quad}^{}}
\end{align}
is a contraction fulfilling all side conditions.
\end{LEMMA}
\begin{BEWEIS} Since $d$ is a derivation of the super-commutative product, it is sufficient to evaluate $d$ on the generators. Since the connection is metric we have $[\Pi_0,\ddxi{i}]=0$ for all $i=1,\dots,n$. Moreover, $\Pi_0$ commutes with all $\ddxi{C}$ and functions $f\in\mathcal C^\infty(M)$. Hence the only nonvanishing contribution comes from:
\begin{align}
\Big[\Pi_0,\xi^{C}\Big]=\half\Big[\sum_{A,B}g^{AB}\ddxi{A}\ddxi{B},\xi^{C}\Big]=\sum_{AB}g^{AB}\ddxi{A}\:\delta^C_B.
\end{align}
This proves equation (\ref{adpi}). Equation (\ref{Komplex}) follows immediately. The Hodge identity follows from the supercommutator:
$[\ddxi{A}\:i_{\xi_B},\xi^C\: i_{\frac{\partial}{\partial \xi_D}}]=\delta_B^C\:\ddxi{A}i_{\frac{\partial}{\partial \xi_D}}+\delta_A^D\:\xi^C\:i_{\xi_B}$ for $p(A)=p(B)$ and $p(C)=p(D)$. This is straightforward to check (however, one has to be careful with the signs). All the remaining statements are obvious.
\end{BEWEIS}
\\

Next, we  introduce the following matrix of super-functions, which has incorporated the Poisson structure and curvature as well:
\begin{align*}
\hat{R}_j^i:=-\half \sum_k\sum_{A,B,C}\:\Pi^{ik}\:R^A_{Bkj}\:g_{AC}\:\xi^B\xi^C\quad\in \ibrsalg{L}^0.
\end{align*}
One of the  reasons to work with the completion, i.e. $\brsalg$ instead of $S_A(\mathcal V^\oplus{\mathcal V}^*)$,  is that otherwise  the matrix $(\id-\hat R)^{-1}$ is well defined only under special circumstances. More specifically, the matrices  
\begin{align*}
\mathcal A^i_j:=\big((\id-\hat{R})^{-\half}\big)^i_j
     =\big(\id+\half\:\hat{R}+\frac{1\cdot 3}{2\cdot 4}\:\hat{R}^2+\frac{1\cdot 3\cdot 5}{2\cdot 4\cdot 8}\:\hat{R}^3+\dots\big)_j^i \quad\in \ibrsalg{L}^0\\
\mathcal B^i_j:=\sum_k \mathcal A^i_k \mathcal A^k_j=\big((\id-\hat{R})^{-1}\big)^i_j
=\big(\id+\hat{R}+\hat{R}^2+\hat{R}^3+\dots\big)_j^i\quad\in \ibrsalg{L}^0
\end{align*}
of superfunctions will play a vital role in what follows. Let us establish some calculation rules.

\begin{LEMMA} We have the following local formulas:
\begin{align}
\sum _j\Pi^{ij}\,\hat{R}^k_j=-\sum_j\Pi^{kj}\,\hat{R}_k^i \label{sym1},\\
\sum_{i,j}\Pi^{ij}\,\mathcal A_i^k\mathcal A_j^l=\sum_j \Pi^{kj}\mathcal B_j^l=-\sum_j \Pi^{lj}\mathcal B_j^k \label{sym2},\\
\ddxi{i}\Big(\sum_{A,B,C} \:R^A_{Bkj}\:g_{AC}\:\xi^B\xi^C\Big)=0\label{Bianchi},\\
d\hat R_j^i= - \sum_{k,A,B}\Pi^{ik}\: R_{Bkj}^A\:\xi^B\ddxi{A}\label{ortho1},\\
d\mathcal B_j^i=- \sum_{k,m,n,A,B}\Pi^{mk}\: \mathcal B^i_m\:\mathcal B^n_j\:R_{Bkn}^A\:\xi^B\ddxi{A}\label{ortho2},\\
 \ddxi{k}\Big(\sum_l \Pi^{il} \mathcal B_l^j\Big)=\sum_{m,n} (\partial_k\Pi ^{mn})\,\mathcal B^i_m\,\mathcal B^j_n.\label{fortune}
\end{align}
\end{LEMMA}
\begin{BEWEIS} Throughout the proof we use Einstein summation convention. Equation (\ref{sym1}) is an immediate consequence of the symmetry properties of the curvature tensor and (\ref{sym2}) follows easily.  Equation (\ref{Bianchi}) is equivalent to Bianchi's identity. In order to check identity (\ref{ortho1}) we note that since $\nabla$ is metric, the curvature is orthogonal with respect to $g$ (regardless of the symmetry properties of $g$):
\begin{align}
g^{AC}\:R^B_{Cij}=g(\xi^A,R(\partial_i,\partial_j)\xi^B)=-g(R(\partial_i,\partial_j)\xi^A,\xi^B)=-g^{CB} \:R^A_{Cij}.
\end{align}
By inverting $g$ we obtain $R^C_{Dij}=-g^{AC}\:g_{BD} R_{Aij}^B$. After these preparations we get 
\begin{eqnarray*}
d\hat R_j^i&=& -\half\:\Pi^{il}\:R_{Blj}^A\: g_{AC}d(\xi^B\xi^C)\\
&=&-\half\:\Pi^{il}\:R_{Blj}^A\Big(g_{AC}\:g^{DE}\:\ddxi{D}\:\delta^B_E\:\xi^C+(-1)^{p(B)}g_{AC}\:g^{DE}\:\ddxi{D}\xi^B\:\delta^C_E\Big)\\
&=&-\half\:\Pi^{il}\:R_{Blj}^A\Big((-1)^{p(B)}g_{AC}\:g^{DB}\:\xi^C+(-1)^{2p(B)}g_{AC}\:g^{DC}\xi^B \:\Big)\ddxi{D}\\
&=& - \Pi^{ik}\: R_{Bkj}^A\:\xi^B\ddxi{A}.
\end{eqnarray*}
Note the following easy to prove matrix identity
\begin{align}
d\big((\id-\hat{R})^{-1}\big)=(\id-\hat{R})^{-1}(d\hat{R})(\id-\hat{R})^{-1}\label{geomid}.
\end{align}
This and the analogous statement for $\ddxi{i}$ are the key to the remaining computations. Equation (\ref{ortho2}) follows from the equations (\ref{ortho1}) and (\ref{geomid}). Finally we show equation (\ref{fortune}):
\begin{eqnarray*}
\ddxi{k}\Big(\Pi^{il}\mathcal B_l^j\Big)&\stackrel{(\ref{geomid})}{=}&\partial_k \Pi^{il}\:\mathcal B^j_l+\Pi^{il}\mathcal B^m_l\Big(\ddxi{k}\hat R_m^n\Big)\mathcal B_n^j\\
&\stackrel{(\ref{Bianchi})}{=}&\partial_k \Pi^{il}\:\mathcal B^j_l-\half\:\Pi^{il}\mathcal B^m_l\Big(\partial_k \Pi^{nr}\: R_{Brm}^A\:g_{AC}\:\xi^B\xi^C\Big)\mathcal B_n^j\\
&\stackrel{(\ref{sym2})}{=}&\partial_k \Pi^{il}\:\mathcal B^j_l-\half\:\Pi^{lm}\mathcal B^i_l\Big(\partial_k \Pi^{nr}\: R_{Brm}^A\:g_{AC}\:\xi^B\xi^C\Big)\mathcal B_n^j\\
&=&\partial_k \Pi^{il}\:\mathcal B^j_l-\hat R ^l_r\,\mathcal B_l^i\,\mathcal B_n^j\:\partial_k \Pi^{nr}\\
&\stackrel{(*)}{=}&\partial_k \Pi^{il}\:\mathcal B^j_l-(\mathcal B^i_r-\delta_r^i)\mathcal B_n^j\:\partial_k \Pi^{nr}=-\mathcal B^i_r\,\mathcal B_n^j\:\partial_k \Pi^{nr}.
\end{eqnarray*}
At step (*) we again have used the geometric series $\hat R(\id-\hat R)^{-1}=(\id-\hat R)^{-1}-\id$.
\end{BEWEIS}
\\

Now we are ready to prove that 
\begin{align}
\widetilde{\Pi}:=\half\sum_{i,j,k,l}\Pi^{ij}\,\mathcal A^k_i\mathcal A^l_j\:\ddxi{k}\ddxi{l}=\half\sum_{j,k,l}\:\Pi^{kj}\: \mathcal B_j^l\:\ddxi{k}\ddxi{l}
\end{align}
is a solution of the Maurer-Cartan equation
\begin{LEMMA}[Maurer-Cartan equation]\label{MCeq} $d\widetilde{\Pi}+\half\,[\widetilde{\Pi},\widetilde{\Pi}]=0$.
\end{LEMMA}
\begin{BEWEIS} Again, Einstein summation convention is in force. Using equation (\ref{ortho2}) we obtain
\begin{eqnarray}
d\widetilde{\Pi}&=&\half \:\Pi^{kj}\:d\mathcal B_j^l\:\ddxi{k}\ddxi{l}=-\half\:\Pi^{kj}\:\Pi^{mn}\:\mathcal B_m^l\:\mathcal B_j^i\:R^A_{Bni}\:\xi^B\ddxi{A}\ddxi{k}\ddxi{l}\\
&=&\half\:\Pi^{kj}\:\Pi^{ln}\:\mathcal B_m^l\mathcal B_j^i\:R^A_{Bni}\:\xi^B\ddxi{A}\ddxi{k}\ddxi{l}.\label{sieheunten}
\end{eqnarray}
It remains to compute $[\widetilde{\Pi},\widetilde{\Pi}]$. Note this is a linear combination of two terms. The first is proportional to $\Pi^{k_1i_1}\mathcal B_{i_1}^{l_1}\big(\ddxi{k_1}(\Pi^{k_2i_2}\mathcal B_{i_2}^{l_2})\big)\ddxi{l_1}\ddxi{k_2}\ddxi{l_2}$ and we show using equation (\ref{fortune}) that this contribution vanishes:
\begin{align*}
\Pi^{k_1i_1}\mathcal B_{i_1}^{l_1}\Big(\ddxi{k_1}(\Pi^{k_2i_2}\mathcal B_{i_2}^{l_2})\Big)\ddxi{l_1}\ddxi{k_2}\ddxi{l_2}=\Pi^{k_1i_1}\:\partial_{k_1}\Pi^{mn}\:\mathcal B_{i_1}^{l_1}\mathcal B_{m}^{k_2}\mathcal B_{n}^{l_2}\ddxi{l_1}\ddxi{k_2}\ddxi{l_2}=0,
\end{align*}
since $\Pi^{k_1i_1}\partial_{k_1}\Pi^{mn}+cycl(i_1,m,n)=0$. Therefore, using (\ref{commrel1}) we conclude that
\begin{eqnarray*}
[\widetilde{\Pi},\widetilde{\Pi}]&=&\frac{1}{4}\:\Pi^{k_1j_1}\mathcal B_{j_1}^{l_1}\:\Pi^{k_2j_2}\mathcal B_{j_2}^{l_2}\:\Big[\ddxi{k_1}\ddxi{l_1},\ddxi{k_2}\ddxi{l_2}\Big]\\
&=&-\Pi^{k_1j_1}\Pi^{k_2j_2}\:\mathcal B_{j_1}^{l_1}\mathcal B_{j_2}^{l_2}\:R_{Bl_1l_2}^A\:\ddxi{k_1}\:\xi^B\ddxi{A}\:\ddxi{k_2}\\
&=&\Pi^{k_1j_1}\Pi^{k_2j_2}\:\mathcal B_{j_1}^{l_1}\mathcal B_{j_2}^{l_2}\:R_{Bl_1l_2}^A\:\xi^B\ddxi{A}\:\ddxi{k_1}\ddxi{k_2},
\end{eqnarray*}
which after comparison with equation (\ref{sieheunten}) yields the Maurer-Cartan equation.
\end{BEWEIS}
\\

We conclude that the \emph{Rothstein Poisson tensor} 
\begin{align}
\Pi_R:=\Pi_0+\widetilde{\Pi} 
\end{align}
has the desired properties \ref{Poiss}, \ref{degzero} and \ref{proj} of Theorem \ref{Roth}. Therefore, $\poly^\bullet(\ibrsalg{L})$ is a cochain complex with differential $\delta_{\Pi_R}:=[\Pi_R,\:]$, which is an analog of Lichnerowicz-Poisson cohomology. Since $\ibrsalg{L}$ is an abelian subalgebra of the Gerstenhaber algebra $\poly(\ibrsalg{L})$,  by the derived bracket construction \ref{derbr1}.) of Theorem \ref{derKS},  $\ibrsalg{L}$ acquires the structure of a $\mathbbm Z$-graded Poisson algebra. More precisely, we define the \emph{Rothstein bracket} $\{\:,\:\}_R$ to be the opposite derived bracket of $\delta_{\Pi_R}$, i.e., for two homogeneous elements $\alpha,\beta\in\ibrsalg{L}$ we have :
\begin{align}
  \{\alpha,\beta\}_R:=[\alpha,\beta]^{\oppind}_{\delta_{\Pi_R}}=-[\alpha,\delta_{\Pi_R}\beta]=-[\alpha,[\Pi_R,\beta]]
\end{align} 
In local coordinates the Rothstein bracket writes as follows
\begin{eqnarray}
\{\alpha,\beta\}_R&=&-[\alpha,[\widetilde{\Pi},\beta]]-[\alpha,[\Pi_0,\beta]]
\label{rothbr}\\
&=&\sum_{j,k,l}\:\Pi^{kj}\: \mathcal B_j^l\:\frac{\partial\alpha}{\partial \xi^k} \frac{\partial\beta}{\partial \xi^l} +\sum_{i=1}^L(-1)^{|\alpha|i+1}\sum_{a=1}^{\ell_i}\Big(\frac{\partial\alpha}{\partial \xAG{i}{a}}\frac{\partial\beta}{\partial \xG{i}{a}}+ (-1)^i\frac{\partial\alpha}{\partial \xG{i}{a}}\frac{\partial\beta}{\partial \xAG{i}{a}}\Big).\nonumber
\end{eqnarray}

We would like to emphasize that $\mathcal C^\infty(M)$ is a Poisson subalgebra of $\ibrsalg{L}$ if and only if the curvature terms vanish, i.e., the bundle $V=\oplus_i V^i$ is \emph{flat}. Nonetheless, the Lichnerowicz-Poisson cohomologies of the Poisson algebras $(\ibrsalg{L},\{\,,\,\}_R)$ and $(\mathcal C^\infty(M),\{\,,\,\})$ are always quasiisomorphic as Theorem \ref{Lichcontr} below shows. If the bundle $V=\oplus_i V^i$ is the \emph{trivial} bundle with the canonical flat connection, then  the $\frac{\partial}{\partial \xi^k}$ are just the ordinary partial derivatives, the curvature term $\mathcal B_j^l$ disappears and the Rothstein bracket boils down to the well known
\begin{align}\label{ordbr}
\{\alpha,\beta\}_R=\sum_{k,l}\:\Pi^{kl}\:\frac{\partial\alpha}{\partial \xi^k} \frac{\partial\beta}{\partial \xi^l} +\sum_{i=1}^L\sum_{a=1}^{\ell_i} \alpha\frac{\overleftarrow{\partial}}{\partial \xG{i}{a}}\frac{\overrightarrow{\partial}}{\partial \xAG{i}{a}}\beta-(-1)^{|\alpha||\beta|}\beta\frac{\overleftarrow{\partial}}{\partial \xG{i}{a}}\frac{\overrightarrow{\partial}}{\partial \xAG{i}{a}}\alpha.
\end{align}
Here we have used the physicist's denotation\footnote{The physicist's denotation is in this case also mnemonically the simplest.}:
\begin{eqnarray*}
\alpha\frac{\overleftarrow{\partial}}{\partial \xG{i}{a}}:=-\Big[\alpha,\frac{\partial}{\partial \xG{i}{a}}\Big]=(-1)^{(|\alpha|+1)i}\frac{\partial\alpha}{\partial \xG{i}{a}},\quad 
\frac{\overrightarrow{\partial}}{\partial \xAG{i}{a}}\alpha=\frac{\partial\alpha}{\partial \xAG{i}{a}}.
\end{eqnarray*}
In order to proof equation (\ref{rothbr}) and (\ref{ordbr}) we essentially have to compute
\begin{eqnarray*}
[\alpha,[\Pi_0,\beta]]&=&\sum_{i=1}^L\sum_{a=1}^{\ell_i} \Big[\alpha,\Big[\ddG{i}{a}\ddAG{i}{a},\beta\Big]\Big]\\
&=&\sum_{i=1}^L\sum_{a=1}^{\ell_i}\Big[\alpha,\ddG{i}{a}\Big[\ddAG{i}{a},\beta\Big]\Big]+(-1)^{(|\beta|+1)(i+1)}\Big[\alpha,\Big[\ddG{i}{a},\beta\Big]\ddAG{i}{a}\Big]\\
&=&\sum_{i=1}^L\sum_{a=1}^{\ell_i}\Big[\alpha,\ddG{i}{a}\Big]\Big[\ddAG{i}{a},\beta\Big]-(-1)^{(|\beta|+1)(i+1)}(-1)^{(|\alpha|+1)|\beta|}\Big[\ddG{i}{a},\beta\Big]\Big[\ddAG{i}{a},\alpha\Big]\\
&=&\sum_{i=1}^L\sum_{a=1}^{\ell_i}\Big[\alpha,\ddG{i}{a}\Big]\Big[\ddAG{i}{a},\beta\Big]-(-1)^{|\alpha||\beta|+i(|\beta|+1)}(-1)^{(|\beta|+1)i+1}\Big[\beta,\ddG{i}{a}\Big] \Big[\ddAG{i}{a},\alpha\Big]\\
&=&\sum_{i=1}^L\sum_{a=1}^{\ell_i}\Big[\alpha,\ddG{i}{a}\Big]\Big[\ddAG{i}{a},\beta\Big]-(-1)^{|\alpha||\beta|}\Big[\beta,\ddG{i}{a}\Big] \Big[\ddAG{i}{a},\alpha\Big].
\end{eqnarray*}
\begin{SATZ}\label{Lichcontr} There are linear maps $\mathcal J: \poly^\bullet(M)\to \poly^\bullet(\ibrsalg{L})$ and $\Sigma: \poly^\bullet(\ibrsalg{L})\to \poly^{\bullet+1}(\ibrsalg{L})$ such that 
\begin{align} 
\xymatrix{(\poly^\bullet(M),\delta_{\Pi})\ar@<-0.6ex>[r]_{\mathcal J\quad}^{}&(\poly^\bullet(\ibrsalg{L}),\delta_{{\Pi}_R}),\Sigma \ar@<-0.6ex>[l]_{\epsilon\quad}^{}}
\end{align}
is a contraction. In particular, the Poisson cohomologies of $\Pi$ and $\Pi_R$ are isomorphic.
\end{SATZ}
\begin{BEWEIS} 
Since $\epsilon$ is a morphism of Gerstenhaber algebras and $\epsilon(\Pi_R)=\Pi$, we have $\epsilon\:\delta_{\Pi_R}=\delta_{\Pi}\:\epsilon$. Hence we may apply perturbation lemma \ref{BPL1}.
\end{BEWEIS}
\subsection{The infinitely generated case}\label{RothPBinfin}
Let us now address the infinitely generated case $\mathcal V= \oplus_{i\ge 1} \mathcal V^i$, where $|\operatorname{supp}(\mathcal V)|=\infty$. In this case the argument of the previous section seem to break down, first of all, due to the fact that $\Pi_0$ is no longer a true bivector field
\begin{align*}
\Pi_0=\sum_{i=1}^\infty\sum_{a=1}^{\ell_i} \ddG{i}{a}\ddAG{i}{a}\quad \notin
 \poly^2(\brsalg).
\end{align*} 
In fact, there seems to be no way to write $\Pi_0$ as a second symmetric power.
Secondly, $\widetilde{\Pi}$ is no bivector field in the sense of our definition
\begin{align*}
\widetilde{\Pi}=\half\sum_{j,k,l}\:\Pi^{kj}\: \mathcal B_j^l\:\ddxi{k}\ddxi{l}\notin \poly^2(\brsalg).
\end{align*}
This is because the innocent looking derivations $\partial/\partial \xi^i$, $i=1,\dots,n=\dim(M)$ are \emph{not} derivations of finite type.

We conclude that, in the infinitely generated case, the Gerstenhaber algebra $\poly(\brsalg)$ of polyvector fields is \emph{too small} for our purposes and we need to look for an appropriate replacement. The only proposal that we are aware of and that seems to make sense is to take the Gerstenhaber algebra $\left(\Der(\brsalg),\cup,[\:,\:]_{RN}\right)$ of multiderivations of $\brsalg$ (cf. Theorem \ref{multiderisgerst}). In fact, one can show that the formal analog of the two-fold bracket 
\begin{align*}
``\:B_{\Pi_0}(\alpha,\beta)=[[\Pi_0,\alpha],\beta]\:``=\sum_{i=0}^\infty
(-1)^{(|\beta|+1)i+\alpha}\sum_{a=1}^{\ell_i}\Big(\frac{\partial\alpha}{\partial \xG{i}{a}}\frac{\partial\beta}{\partial \xAG{i}{a}}+ (-1)^i\frac{\partial\alpha}{\partial \xAG{i}{a}}\frac{\partial\beta}{\partial \xG{i}{a}}\Big)
\end{align*}
is a well-defined operation (this is somewhat surprising) and we conclude that $\Pi_0\in \Der^2(\brsalg)^0$. It is also not difficult to prove that $\Pi_0\bullet \Pi_0=\half[\Pi_0,\Pi_0]_{RN}=0$. The term $\widetilde{\Pi}$ is evidently a member of $\Der^2(\brsalg)^0$.

We strongly believe that all calculations which lead us to the Maurer-Cartan equation of Lemma \ref{MCeq}, and which now involve certain infinite summations, make sense in the algebra of multiderivations. (even though we do not know whether there is an analog of the contraction (\ref{dcontr})). As a consequence, we conjecture that also in the infinitely generated case the Rothstein bracket 
\begin{align}
\{\alpha,\beta\}_R=\sum_{j,k,l}\:\Pi^{kj}\: \mathcal B_j^l\:\frac{\partial\alpha}{\partial \xi^k} \frac{\partial\beta}{\partial \xi^l} +\sum_{i=1}^\infty\sum_{a=1}^{\ell_i}\alpha\frac{\overleftarrow{\partial}}{\partial \xG{i}{a}}\frac{\overrightarrow{\partial}}{\partial \xAG{i}{a}}\beta-(-1)^{|\alpha||\beta|}\beta\frac{\overleftarrow{\partial}}{\partial \xG{i}{a}}\frac{\overrightarrow{\partial}}{\partial \xAG{i}{a}}\alpha.
\end{align}
is a $\mathbbm Z$-graded super-Poisson structure, which
can be seen as the opposite derived bracket of the abelian subalgebra $\brsalg\subset \oplus_{n\ge 0} \Der^n(\brsalg)$ with respect to the differential $\delta_{\Pi_R}:=[\Pi_0+\widetilde{\Pi},\:]_{RN}$. Clearly, the virtues of the derived bracket construction (i.e. keeping the signs simple) disappear if we do calculations in the algebra of multiderivation. A better strategy to prove the above statement seems to be to use the concept of approximation (cf. Section \ref{BRSTchargesection}) and some `clever' continuity  argument. In the freely generated case with a flat connection the Jacobiidentity for the bracket $\{\:,\:\}_R$ seems to be folklore.
\section{The Koszul complex}\label{kc}
Given a smooth map $J:M\to \mathbbm R^\ell=:V^*$ we consider the Koszul holomogical 
complex of the sequence  of ring elements $J_1,\dots, J_\ell\in \mathcal C^{\infty}(M)$. In other words, we define 
the space of chains to be \[K_{i}:=K_{i}(\mathcal C^\infty(M),J):=S^i_{\mathcal C^{\infty}(M)}(V[-1]),\]
i.e., the free (super)symmetric $\mathcal C^{\infty}(M)$-algebra generated by the 
graded vector space $V[-1]$, where we consider $V$ to be concentrated in degree zero. 
$K_{\bullet}$ may also be viewed as the space of sections of the trivial vector bundle 
over $M$ with fibre $\wedge^\bullet V$. Denoting by $\xi^1,\dots, \xi^\ell$ the canonical 
bases of $V[-1]$ for the dual space $V$ of $V^*=\mathbbm R^\ell$, we define the Koszul differential 
\begin{align*}
\partial:=\sum_a J_a \ddxi{a},
\end{align*} where the $\ddxi{a}$ , $a=1,\dots,\ell$, are the derivations extending the 
dual pairing. We will say, in accordance with \cite{Bourbdix}, that the sequence of ring elements 
$J_1,\dots, J_\ell\in \mathcal C^{\infty}(M)$ is a \emph{complete intersection}, if the homology of the Koszul complex vanishes in degree $\ne 0$.

If zero is a regular value of $J$ is well known (and follows from Theorem \ref{acyckrit} below) that $J_1,\dots, J_\ell\in \mathcal C^{\infty}(M)$ is a \emph{complete intersection}. An elementary example of a noncomplete intersection is provided by the moment map for one particle of zero angular momentum, Example \ref{drehimpuls}, in dimension $n=3$. In this case, using the physicist's denotation, the Koszul complex can be rewritten as
\begin{align*}
0\leftarrow\mathcal C^\infty(T^*\mathbbm R^3)\stackrel{<J,\:>}{\leftarrow}\mathbbm R^3\otimes C^\infty(T^*\mathbbm R^3)\stackrel{J\times}{\leftarrow}\mathbbm R^3\otimes C^\infty(T^*\mathbbm R^3)\stackrel{J\cdot}{\leftarrow}\mathcal C^\infty(T^*\mathbbm R^3)\leftarrow 0.
\end{align*}
Here $q=(q^1,q^2,q^3)$ and $p=(p_1,p_2,p_3)$ are interpreted as vector valued functions on $T^*\mathbbm R^3$ and the angular momentum is $J=q\times p$, where $\times$ denote the vector product in $\mathbbm R^3$. Since $q$ and $p$ are orthogonal to $J$, i.e., the euclidian scalar products $<q,J>=0=<p,J>$ vanish, $q$ and $p$ are one-cycles. In fact, they cannot be boundaries, for if $q=J\times\text{'something'}$ it had to vanish as a function at $p=0$ and vice versa.

Since such findings are merely accidental, we would like to have a more systematic way to decide whether a moment map is a complete
intersection. One way could be to make a detour and use methods from the theory of commutative Noetherian  ring together with flatness arguments. Instead, if we already know that the generating hypothesis is true, then the following `Jacobian criterion' yields a more convenient method.
\begin{SATZ} \label{acyckrit} Let M be an analytic manifold and 
$J:M\to \mathbbm R^\ell$  an  analytic map, such that the following conditions are true
\begin{enumerate}
\item \label{gh}$(J_1,\dots,J_\ell)$ generate the vanishing ideal of $Z:=J^{-1}(0)$ in 
      $\mathcal C^\infty(M)$,
\item \label{jc}the regular stratum $Z_r:=\{z\in Z\:|\: T_zJ \mbox{ is surjective}\}$ is dense 
      in $Z:=J^{-1}(0)$.
\end{enumerate}
      Then the Koszul complex $K:=K(\mathcal C^\infty(M),J)$ is acyclic and $H_0=\mathcal C^\infty(Z)$.
\end{SATZ}
\begin{BEWEIS} We will show that the Koszul complex $K(\mathcal C^\omega_x(M),J) $ is 
acyclic for the ring $\mathcal C^\omega_x(M)$ of germs in $x$ of real analytic 
functions. Then it will follow that the Koszul complex $K(\mathcal C^{\infty}(M),J)$ 
is acyclic, since the ring of germs of smooth functions $\mathcal C^\infty_x(M)$ is (faithfully)
flat over $\mathcal C^\omega_x(M)$ (see  \cite[Corollary VI 1.12]{Malgrange}), and the sheaf of smooth 
functions on $M$ is fine. Since $\mathcal C^\omega_x(M)$ is Noetherian, Krull's 
intersection theorem says that $\cap_{r\ge0} I^r_x=0$, where $I_x$ is the ideal of 
germs of analytic functions vanishing on $Z$. According to \cite[A X.160]{Bourbdix}, 
it is therefore sufficient to show that  $H_1(\mathcal C^\omega_x(M),J)=0$. Note that 
since $J$ generates the vanishing ideal of $Z$ in $\mathcal C^\infty(M)$, it also 
generates the vanishing ideal of $Z$ in $\mathcal C^\omega_x(M)$. This can easily be 
seen using M. Artin's approximation theorem (see e.g.\cite{Ruiz}).\footnote{We have been hinted by L. Avramov and S. Iyengar that this also follows  from the faithful flatness of $\mathcal C^\omega_x(M)\to\mathcal C^\infty(M)$ \cite[Theorem 7.5]{Matsumura}.} Suppose 
$f=\sum_a f^a\:e_a\in K_1$ is a cycle, i.e. $\partial f=\sum_a J_af^a=0$. Since the 
restriction to $Z$ of the Jacobi matrix $D(\sum_a J_af^a)$ vanishes,  we conclude 
(using condition b)) that $f^a_{|Z}=0$ for all $a=1,\dots,\ell$. Since $J$ generates 
the vanishing ideal, we find an $\ell\times\ell$-matrix $F= (F^{ab})$ with smooth 
(resp. analytic) entries such that $f^a=\sum_b F^{ab}J_b$.\footnote{We have been hinted by L. Avramov and S. Iyengar that using  the theorem of Vasconcelos \cite[Theorem 19.9 and the remark that follows]{Matsumura} we are done at this point.} It remains to be shown, that 
this matrix can be choosen to be \emph{antisymmetric}. We have to distinguish two 
cases.  If $x\notin Z$, the claim is obvious, since then one can take for example 
$F^{ab}:=(\sum_a J_a^2)^{-1}(J_bf^a-J_a f^b)$.
So let us consider the case  $x \in Z$. We then introduce some formalism
to avoid tedious symmetrization arguments. 
Let $E$ denote the free $k:=\mathcal C^\omega_x(M)$-module on $\ell$ generators, 
and consider the Koszul-type complex $SE\otimes \wedge E$. Generators of the symmetric part 
will be denoted by $\mu_1,\dots,\mu_\ell$, generators of the Grassmann part by 
$e_1,\dots,e_\ell$, respectively. We have two derivations 
$\delta:=\sum_a e_a\wedge\frac{\partial}{\partial \mu_a}:S^nE\otimes \wedge^m E\to 
S^{n-1}E\otimes \wedge^{m+1}E$, and 
$\delta^*:=\sum_a \mu_a i(e^a):S^nE\otimes \wedge^m E\to S^{n+1}E\otimes \wedge^{m-1}E$. 
They satisfy the well known identities: $\delta^2=0$, $(\delta^*)^2=0$ and 
$\delta\delta^*+\delta^*\delta=(m+n)\id$. Furthermore, we introduce the two commuting 
derivations $i_J:=\sum_a J_a i(e^a)$ and $d_J=\sum_a J_a\frac{\partial}{\partial \mu_a}$. 
They obey the identities $i_J^2=0$, $[i_J,\delta]=d_J$, $[d_J,\delta^*]=i_J$ and 
$[i_J,\delta^*]=0=[d_J,\delta]$. We interprete the cycle $f$ above as being in $E\otimes k$ 
and the matrix $F$ as a member of $E\otimes E$. We already know that $d_Jf=0$ implies $f=i_JF$. 
This argument may be generalized as follows: if $a\in S^nE\otimes k$ obeys $d_J^na=0$, then 
there is an $A\in S^nE\otimes E$ such that $a=i_J A$. The proof is easily provided by taking 
all $n$-fold partial derivatives of $d_J^na=0$, evaluating the result on $Z$ and using 
conditions a) and b).
We now claim that there is a sequence of $F_{(n)}\in S^{n+1}E\otimes E$, $n\ge 0$, such that 
$F=F_{(0)}$, $\delta^*F_{(n)}=(n+2)i_JF_{(n+1)}$ and 
\begin{eqnarray}\label{nteGl}
 f=d_J^ni_JF_{(n)}+i_J\delta^*\underbrace{\Big(\sum_{i=0}^{n-1}\frac{1}{i+2}\:d_J^i\delta 
 F_{(i)}\Big)}_{=:B_{n-1}}\quad\mbox{ for all }n\ge1.
\end{eqnarray} 
We prove this by induction. Setting $B_{-1}:=0$, we may start the induction with $n=0$, where 
nothing has to be done. Suppose now, that the claim is true for $F_{(0)},\dots ,F_{(n)}$. 
We obtain $f=\frac{1}{n+2}d_J^ni_J\big(\delta \delta^*F_{(n)}+\delta^*\delta F_{(n)}\big)+
i_J\delta^*B_{n-1}=\frac{1}{n+2}d_J^{n+1}\delta^*F_{(n)}+i_J\delta^*B_n$,
where we made use of the relations $[d^n_Ji_J,\delta^*]=0$ and $[d^n_Ji_J,\delta]=d^{n+1}_J$. 
Since $0=d_Jf=d_J^{n+2}\delta^*F_{(n)}$, we find an $F_{(n+1)}$ such that 
$\frac{1}{n+2}\delta^*F_{(n)}=i_JF_{(n+1)}$, and the claim is proven. 
Finally, we would like to take the limit of equation (\ref{nteGl}) 
as $n$ goes to $\infty$. 
For this limit to make sense, we have to change the ring to the ring of formal power series. 
Let us denote this change of rings by 
$\hat{}:\mathcal C_x^\omega(M)\to\mathbbm K[[x^1,\dots,x^n]]$. 
Since by Krull's intersection theorem $\cap_{r\ge0} \hat{I}^r=0$ ($\hat I$ the ideal generated by 
$\hat J_1,\dots,\hat J_{\ell}$), we obtain a formal solution of the problem: 
$\hat f=i_{\hat J}\delta^*B_{\infty}$, where $B_{\infty}:=\sum_{i=0}^{\infty}\frac{1}{i+2}\:
d_{\hat J}^i\delta \hat F_{(i)}$ is well defined since $\hat I$ contains the maximal ideal. 
Applying M. Artin's approximation theorem yields an analytic solution, and we are done. 
\end{BEWEIS}
\\

The above reasoning can be considered to be folklore, as the subtlety of finding an 
\emph{antisymmetric} source term is often swept under the rug in semirigorous arguments. 
 We do not know, whether, if  condition \ref{gh}.) in Theorem \ref{acyckrit} holds, \ref{jc}.) is also sufficient for the acyclicity of the Koszul complex. Nonetheless, let us, as a plausibility check, reconsider the system of one particle of zero angular momentum in dimension $n\ge 3$, Example \ref{drehimpuls}. Here, the Jacobi matrix $T_z J$ is a ${n\choose 2}\times 2n$-matrix. Since the zero fibre $Z$ is the set of points where $p$ and $q$ are proportional, it follows easily that $T_z J$ has the same rank as a certain ${n\choose 2}\times n$ submatrix, which in fact occurs twice in $T_z J$. With a little more effort, one may proof that the rank of this submatrix is in fact $\le n-1$. Since ${n\choose 2}>n-1$ for $n\ge 3$ the regular stratum is \emph{empty} here.

For the Examples \ref{mpart}, \ref{lemon}, \ref{stratres}, \ref{tzwo} (for $\alpha<0$) and \ref{commvar} we already know from subsection \ref{ghsubsection} that condition \ref{gh}.) of Therorem \ref{acyckrit}, i.e., the generating hypothesis, is fulfilled. Let us now check that condition \ref{jc}.) of Theorem \ref{acyckrit} holds for these examples as well. For Example \ref{mpart}, i.e., the system of $m\ge 1$ particles of zero angular momentum in the plane, and for the $(1,1,-1,-1)$-resonance, Example \ref{stratres}, the Jacobi matrix $T_x J$ is clearly not onto iff $x=0$, i.e., the singularity is isolated in $M$. Clearly, in both cases $Z\ne \{0\}$ and hence  condition \ref{jc}.) is true. In the case of the `lemon', Example \ref{lemon}, an easy calculation using affine coordinates yields that the singular points of the moment map are precisely the fix points of the $S^1$-action: $((1:0),(1:0))$, $((0:1),(0:1))$, $((1:0),(0:1))$  and $((0:1),(1:0))$. The latter two are isolated points in the zero level $Z=J^{-1}(0)$ (which is obviously nondiscrete), and hence the requisites of Theroem \ref{acyckrit} are fulfilled here. 

For the $\mathbbm T^2$-action of Example \ref{tzwo}, $\alpha<0$, the Jacobi matrix of $J$ works out as follows 
\begin{align*}
\left(
\begin{array}{cccccccc}
\alpha \bar{z}_1&\alpha z_1  &  0      &0  &\bar{z}_3&z_3&0        &0\\
\beta  \bar{z}_1&\beta z_1   &\bar{z}_2&z_2&    0    & 0 &\bar{z}_4&z_4\\
\end{array}
\right).
\end{align*}
The set of points where this matrix has not the full rank is just the union $\cup _{i=1}^4 L_i$ of the coordinate lines $L_1:=\{(z,0,0,0)\:|\:z\in \mathbbm C\}$, $L_2:=\{(0,z,0,0)\:|\:z\in \mathbbm C\}$,  $L_3:=\{(0,0,z,0)\:|\:z\in \mathbbm C\}$ and  $L_4:=\{(0,0,0,z)\:|\:z\in \mathbbm C\}$. Clearly, $Z\cap\cup_i L_i=\{0\}\ne Z$ and hence condition \ref{jc}.) is true. Finally, let us address Example \ref{commvar}, i.e., the commuting variety. Let $(Q,P)\in S\times S=T^*S$ be a pair of symmetric $n\times n$-matrices. Then the Jacobi-matrix of $J$ applied on the tangent vector  $(V,W)\in S\oplus S$ yields
\begin{align*}
TJ_{(Q,P)}(V,W)=\frac{d}{dt}[Q+tV,P+tW]_{|t=0}=[V,P]+[Q,W].
\end{align*}
It follows easily that $TJ_{(Q,P)}$ is surjective for $(Q,P)$ from an open dense subset of $S\times S$. In fact, let $Q\in S_{\operatorname{reg}}$ have pairwise distinct Eigenvalues $q_i$, $q_i\neq q_j$ for $i\ne j$. After an orthogonal change of the basis we may assume that $Q$ is diagonal: $Q=\sum_i q_i E_{ii}$. Here, $\left(E_{ij}\right)_{kl}=\delta_{ik}\delta_{jl}$ is the matrix with the only nonzero entry $1$ at the $i$th row and $j$th column. For any $P\in S$ we have
$TJ_{(Q,P)}(E_{ij}+E_{ji},0)=[Q,E_{ij}+E_{ji}]=(q_i-q_j)(E_{ij}+E_{ji})$.
Since $(E_{ij}+E_{ji})_{ij}$ constitute is a basis for $S$, it follows that $TJ_{|S_{\operatorname{reg}}\times S}$ is surjective. Since $(Q,P)\in Z\cap S_{\operatorname{reg}}\times S$ $\Leftrightarrow$ $Q\in S_{\operatorname{reg}}$ and $P\in S$ are simultaneously diagonalizable, the claim follows.

The next theorem is a consequence of rather 
deep analytic results.
The problem of splitting the Koszul resolution in the context of
Fr\'echet spaces has been also addressed in \cite{DomJac} from a different perspective.

\begin{SATZ} 
  Let $M$ be a smooth manifold, $J:M\to \mathbbm R^\ell$ be a smooth map such that around every 
  $m\in M$ there is a local chart in which $J$ is real analytic. Moreover, assume that the 
  Koszul complex $K=K(M,J)$ is a resolution of $\mathcal C^\infty(Z)$, $Z=J^{-1}(0)$. Then 
  there are a prolongation map 
  $\prol: \mathcal C^{\infty}(Z)\to \mathcal C^{\infty}(M)$ 
  and contracting homotopies $h_i: K_i\to K_{i+1}$, $i\ge 0$, 
  which are continuous in the respective Fr\'echet topologies, such that
\begin{eqnarray}
\label{KosContr}
  \big(\mathcal C^{\infty}(Z),0\big)\:
  \begin{array}{c}{\res}\\\leftrightarrows\\{\prol}
  \end{array}\:\big(K,\partial\big),h
\end{eqnarray}
  is a contraction, i.e., $\res$ and $\prol$ are chain maps and $\res\:\prol=\id$ and 
  $\id-\prol\:\res=\partial h+h\partial$. If necessary, these can be adjusted in 
  such a way, that the side conditions (see Appendix \ref{HPT})
  $h_0\:\prol=0$ and 
  $h_{i+1}\:h_i=0$
  are fulfilled. If, moreover, a compact Lie group $G$ acts smoothly
  on $M$, $G$ is represented on $\mathbbm R^\ell$ and $J:M\to \mathbbm R^\ell$ is equivariant, 
  then $\prol$ and $h$ can additionally be chosen to be equivariant.
\end{SATZ}
\begin{BEWEIS} A closed subset $X\subset\mathbbm R^n$ is defined to have the \emph{extension property}, if there is a continuous linear map $\lambda:\mathcal C^\infty(X)\to \mathcal C^\infty(\mathbbm R^n)$, such that $\res\: \lambda=\id$. 
The extension theorem of E. Bierstone and G. W. Schwarz, \cite[Theorem 0.2.1]{Schwarzbier} says that Nash subanalytic sets (and hence closed analytic sets) have the extension property. 
Using a partition of unity, we get a continuous linear map $\lambda:\mathcal C^\infty(Z)\to \mathcal C^\infty(M)$, such that $\res\: \lambda=\id$. 
In the same reference, one finds a ``division theorem'' (Theorem 0.1.3.), which says that for a matrix $\varphi\in \mathcal C^\omega(\mathbbm R^n)^{r,s}$ of analytic functions the image of $\varphi:\mathcal C^\infty(\mathbbm R^n)^s\to C^\infty(\mathbbm R^n)^r$ is closed, and there is a continuous split $\sigma:\mathrm{im}\:\varphi\to\mathcal C^\infty(\mathbbm R^n)^s$ such that $\varphi\:\sigma=\id$. 
Using a partition of unity, we conclude that there are linear continuous splits $\sigma_i:\mathrm{im}\:\partial_{i+1}\to K_{i+1}$ for the Koszul differentials $\partial_{i+1}:K_{i+1}\to K_i$ for $i\ge 0$, i.e., $\partial_{i+1}\:\sigma_i=\id$. We observe that $\mathrm{im}\:\lambda\oplus\mathrm{im}\:\partial_1=K_0$, since for every $x\in K_0$ the difference $x-\lambda\:\res x$ is a boundary due to exactness and the sum is apparantly direct. 
Similarly, we show that $\mathrm{im}\:\sigma_i\oplus\mathrm{im}\:\partial_{i+2}=K_{i+1}$ for $i\ge 0$. The next step is to show that $\mathrm{im}\:\sigma_i$ is a \emph{closed} subspace of $K_0$. Therefor we assume that $(x_n)_{n\in\mathbbm N}$ is a sequence in $\mathrm{im}\:\partial_{i+1}$ such that $\sigma_i(x_n)$ converges to $y\in K_{i+1}$. Then $x_n=\partial_{i+1}\sigma_i(x_n)$ converges to $\partial_{i+1} y$, since $\partial_{i+1}$ is continuous. Since $\partial_{i+1} y$ is in the domain of $\sigma_i$, we obtain that $\sigma_i(x_n)$ converges to $\sigma_i\partial_{i+1} y=y\in \mathrm{im}\:\sigma_i$. Similarly, we have that $\mathrm{im}\:\lambda$ is a closed subspace of $K_0$. Altogether, it is feasible to extend $\sigma_i$ to a linear continuous map $K_i\to K_{i+1}$ (cf. \cite[ p.133]{Rudin}). If necessary, $\lambda$ and $\sigma_i$ can be made equivariant by averaging over $G$, since $\res$ and $\partial$ are equivariant. We observe that we have 
$\lambda \res_{|\mathrm{im}\lambda}=\id$ and 
$\lambda \res_{|\mathrm{im}\partial_1}=0$ and analogous equations in higher degrees. 
We now replace $\lambda$ by $\prol:=\lambda-\partial_1\sigma_0 \lambda$ and $\sigma_i$ 
by $h_i:=\sigma_i-\partial_{i+2}\sigma_{i+1}\sigma_i$ for $i\ge0$. These maps share 
all of the above mentioned properties with $\lambda$ and $\sigma_i$. Additionally, we 
have ${\partial_1 h_0}_{|\mathrm{im}(\prol)}=0$ and 
${\partial_{i+2}h_{i+1}}_{|\mathrm{im}(h_i)}=0$ for $i\ge0$. This concludes the 
construction of (\ref{KosContr}). The side conditions can be achieved by algebraic 
manipulations (see Appendix \ref{HPT}).  Note that these modifications do not ruin 
the equivariance.
\end{BEWEIS}
\\

If the constraint surface $Z$ is singular it seems to be hopeless to find explicit formulas for $\prol$ and $h$. However, if $Z$ is a closed submanifold there is  a general recipe \cite{BHW} to produce such formulas. We sketch this procedure for the case of linear Poisson structure, Example \ref{linPoiss}, which essentially contains already the whole idea. Recall that, if $x_1,\dots,x_n$ are linear coordinates for $M=\mathfrak h^*$, the moment map is given by the projection $J:M\to\mathfrak g^*$, $J(x_1,\dots,x_n)=(x_1,\dots ,x_\ell)$. The zero fibre $Z=J^{-1}(0)$ is the linear subspace on which the first $\ell$ coordinates vanish. The prolongation map $\prol:\mathcal C^\infty(Z)\to\mathcal C^\infty(M)$ is just
\begin{align}\label{prolformula}
(\prol f)(x_1,\dots,x_n):=f(x_{\ell+1},\dots,x_n).
\end{align}
If $f\in\mathcal C^\infty(M)$ and $v\in S^k_\mathbbm K(\mathfrak g[1])$ then the contracting homotopy is given by the formula
\begin{align}\label{hformula} 
h(fv):=\sum_{a=1}^\ell\left(\int_0^1 dt\:t^k\frac{\partial f}{\partial x_a}(tx_1\dots,tx_\ell,x_{\ell+1},\dots,x_n)\right)\:\xi_av.
\end{align}
By linear extension this defines the contracting homotopy $h:K_k\to K_{k+1}$. An easy calulation yields that in this case (\ref{KosContr}) is a contraction fulfilling the side conditions $h_0\prol=0$ and $h^2=0$.

For completeness, we conclude this section by mentioning two consequences of the Koszul resolution, which are of course well known to commutative algebraists. 
\begin{SATZ} \label{freebasis}Let $J=(J_1,\dots,J_\ell): M\to \mathbbm R^\ell=:V^*$ such that the  Koszul complex $K(A,J)$ is a resolution of $\mathcal C^\infty(Z)=\mathcal C^\infty(M)/I$. Then the image of $J_1,\dots,J_\ell$ under the projection to $I/I^2$ is a free system of generators for the  $A/I$-module $I/I^2$.  
\end{SATZ}
\begin{BEWEIS}
Let $\sum_a f^a [J_a]=0$, where $f^1\dots, f^\ell\in C^\infty(Z)$ and $[J_a]$ are the representatives of the $J_a$'s in $I/I^2$. By choosing $F^1, \dots F^\ell\in C^\infty(M)$ such that $f_a={F_a}_{|Z}$ for $a=1\dots\ell$, we may rewrite this as $\sum_a F_a J_a\in I^2$. Defining $F:=\sum_a F^a e_a \in K_1(A,J)$ we obtain $\partial F=\sum_{a,b}G^{ab}J_aJ_b\in I^2$, for some $G^{ab}\in \mathcal C^\infty(M)$. So $F-\sum_{ab}G^{ab}J_ae_b$ is a 1-cycle. Since the first homology of $K(A,J)$ vanishes, there is an $H:=\sum H^{ab} e_ae_b\in K_2(A,J)$ such that $F=\sum_{ab}G^{ab}J_ae_b+\partial H$, and thus all $F^1,\dots ,F^\ell$ are in $I$. 
\end{BEWEIS}

\begin{KOROLLAR} \label{exttor} Let $J=(J_1,\dots,J_\ell)$ be as above. Then there are isomorphisms of $\mathbbm K$-vector spaces
\begin{eqnarray*}
\Tor^A_i(A/I,A/I)&\cong&\wedge^i_{A/I}I/I^2\quad\mbox{and}\\
\Ext^i_A(A/I,A/I)&\cong& \Hom_{A/I}(\wedge^i_{A/I}I/I^2,A/I)
\end{eqnarray*}
for all $i\ge 0$.
\end{KOROLLAR}
\begin{BEWEIS} In order to compute $\Tor^A_i(A/I,A/I)$, we have to compute the homology of the complex $(K(A,J)\otimes _A A/I,\partial \otimes 1)$. But this is evidently the homology of the Koszul complex $K(A/I,(0,\dots,0))$. The claim follows from Theorem \ref{freebasis}. It is easy to see, that $\Ext^i_A(A/I,A/I)$ is the cohomology of the complex $\End_A(K(A,J))$ with differential $D$ equal to the super-commutator with $\partial$. We may identify $\End_{\Zmod{A}}(K(A,J))=\End_{\Zmod{A}}(\wedge v\otimes A)=\End_{\Zmod{\mathbbm K}}(\wedge V)\otimes A$ with the Clifford algebra bundle over $M$ with fibre $\wedge V^*\otimes \wedge V$. The latter is identified with $\End_{\mathbbm K}(\wedge V)$ using the standard representation: $x\in V$ acts on $\wedge V$ by left multiplication, and $\alpha\in V^*$ acts on $\wedge V$ by $-1$ times the insertation derivation $i(\alpha)$. The analog of composition in $\End_{\mathbbm K}(\wedge V)$ is the Clifford multiplication $\mu\circ(-\sum_a i(e^a)\otimes i(e_a))$, where $e_1,\dots,e_\ell$ is a basis for $V$ and $e^1,\dots,e^\ell$ is the corresponding dual basis. Now the analog of $D$ in the Clifford algebra is just the commutator with $\sum_a J_a e^a$. The latter actually is the total differential of a double complex with acyclic rows and trivial columns.
\end{BEWEIS}
\\

Actually, one can  make the above statement more precise. It is well known \cite[\S 7]{Bourbdix}, that  $\Ext^\bullet_A(A/I,A/I)$ is $\mathbbm Z$-graded $\mathbbm K$-algebra and $\Tor^A_\bullet(A/I,A/I)$ is a graded module for this algebra. The isomorphisms of corollary \ref{exttor} reflect this structure. The ``produit de composition'' is given up to a sign by the wedge product and module structure by insertation, respectively. In particular, the above isomorphisms are isomorphisms of $A/I$-modules.

\section{The projective Koszul resolution of a closed submanifold}\label{projkos}
Even in the regular case not every constraint surface admits a Koszul resolution. For a closed codimension $\ell$ submanifold $C$ of the manifold $M$, we know from Theorem \ref{freebasis}, that if there is a Koszul complex, which is a resolution of $\mathcal C^\infty(C)$ then the conormal bundle of $C$ in $M$ is trivial. Since from the homological point of view projective modules are as good as free modules, one may ask whether the situation improves if one also accepts projective Koszul resolutions. We will see that after restricting to an appropriate open neighborhood $U$ of $C$ in $M$, there is a $\mathcal C^\infty(U)$-projective resolution for every closed submanifold $C$.

Let $TC$ be the tangent bundle of $C$ and let $TM_{|C}$ and $TM^*_{|C}$ be the restrictions of the tangent bundle and the cotangent bundle of $M$ to $C$, respectively. The annihilator bundle $TC^\annind$ of $C$ is the subbundle of $TM^*_{|C}$ consisting of all 1-forms vanishing on $TC$. The dual bundle ${TC^\annind}^*$ to the annihilator bundle is canonically isomorphic to $TM_{|C}/TC$. It is well known that the space of sections of the annihilator bundle has a nice algebraic description.
\begin{LEMMA} \label{conormallemma} The the map $f\mapsto \dext f_c$ which associates to a smooth function $f\in \mathcal C^\infty(M)$ its differential evaluated at the point $c\in C$ induces an isomorphism of $\mathcal C^\infty(C)$-modules of conormal module $I/I^2$ of the vanishing  ideal $I$ of $C$ and  the space of sections $\Gamma^\infty(C,TC^\annind)$ of the annihilator bundle $TC^\annind$ of the submanifold $C$.
\end{LEMMA}
\begin{BEWEIS} For a detailed proof which uses the Koszul complex constructed below see \cite{opusmagnum}.
\end{BEWEIS}
\\

Let us recall the tubular neighborhood theorem for the submanifold $C\subset M$. There exists a subfibre bundle  $N_C$ of $TM_{|C}$ which is complementary to $TC$
\begin{align*}
TM_{|C}=TC\oplus N_C,
\end{align*}
an open neighborhood $U'$ of the zero section of $N_C$ which is a disk bundle over $C$ and a diffeomorphism $\varphi$
\begin{align*}
\xymatrix
{N_C\supset U'\qquad\ar[rr]^{\operatorname{\varphi}}&&\qquad U\subset M\\
&\quad C\quad \ar@{^{(}->}[ul]\ar@{^{(}->}[ur]&}
\end{align*}
into an open neighborhood $U\subset M$ whose restriction to $C$ is the identity. In fact, $\varphi$ may be given as a restriction of the exponential map corresponding to an auxiliary Riemannian metric on $M$. In this way we obtain a (noncanonical) identification of $N_C$ and ${TC^\annind}^*$. Moreover, we have a submersion
\begin{align*}
\tau:U\to C.
\end{align*}
Since $U'$ is fibrewise convex it makes sense to take a convex linear combination of any set of points of $U$ which lie in the same fibre of $\tau$. In particular for any $t\in [0,1]$ we have a fibre preserving shrinking map $\Phi_t:U\to U$, which is obtained from the convex combination $tu+(1-t)u_0$, where $u_0\in C\subset U$ and $u$ are in the same fibre.

Let $V^*:= \kernel T\tau$ be the vertical subbundle of $TU$ of the projection  $\tau:U\to C$. There is a section $J\in\Gamma^\infty(U,V^*)$, which we will call the \emph{tautological section}, such that the submanifold $C$ is the zero locus of $J$. It is given by the image under $T\varphi$ of the Euler vector field $\in\Gamma^\infty(U',TN_C)$ of the vector bundle $N_C$. Note that the restriction $V^*_{|C}$ to $C$ of the bundle $V^*$ is isomorphic to ${TC^{\annind}}^*$.

The \emph{projective (homological)  Koszul complex} $(K_{\bullet}(U,J),\partial)$ on the tautological section $J$ is defined as follows. The space of  chains is $K_{\bullet}(U,J):=S_{\mathcal C^\infty(U)}(V[-1])$. In more conventional terminology, it also may be  viewed as $\Gamma^\infty(U,\wedge^\bullet V)$. The differential is the unique $\mathcal C^\infty(U)$-linear superderivation, such that $\partial(f)=<J,f>$ for all $f\in K_1(U,J)$, where $<\:,\:>$ denotes the dual pairing. By augmenting this complex with the restriction map we obtain a sequence
\begin{align*}
0\longleftarrow\mathcal C^\infty(C)\stackrel{\res}{\longleftarrow}\mathcal C^\infty(U)=K_0(U,J)\stackrel{\partial}{\longleftarrow}K_1(U,J)\stackrel{\partial}{\longleftarrow}K_2(U,J)\longleftarrow\dots
\end{align*}
In fact, this sequence is \emph{exact}, which can of course  be proven locally without difficulty. We will exhibit explicit contracting homotopies similar to the equations (\ref{prolformula}) and (\ref{hformula}) for the linear case. For the prolongation map 
\[\prol:=\tau^*:\mathcal C^\infty(C)\to\mathcal C^\infty(U)\]
we simply take the pull back of the projection $\tau:U\to C$. The contracting homotopies $h_i:K_i(U,J)\to K_{i+1}(U,J)$ for $i=0,1,\dots,\ell$  are defined by the formula (which is due to M. Bordemann \cite{opusmagnum})
\[h_i(a)(u):=\int_0^1dt\: (\Phi_t^*(da))(u)\]
for $a\in K_i(U,J)$ and $u\in U$.
In order to understand the right hand side  of this equation note that the space $K_i(U,J)$ can be identified with the space of cochains of Lie algebroid cohomology with coefficients in $\mathcal C^\infty(U)$ of the Lie algebroid $V^*:= \kernel T\pi\to U$, the differential is denoted by $d:K_i(U,J)\to K_{i+1}(U,J)$. Since $\Phi_t$ is fibre preserving, the pullback $\Phi_t^*$  with respect to $\Phi_t$ for differential forms descends to a well defined map $\Phi_t^*:K_i(U,J)\to K_i(U,J)$ for all $i$.
\begin{SATZ} The maps $\prol$ and $h$ defined as above are continuous $\mathbbm K$-linear maps, such that 
\begin{align} \label{prKosContr}
  \big(\mathcal C^{\infty}(C),0\big)\:
  \begin{array}{c}{\res}\\\leftrightarrows\\{\prol}
  \end{array}\:\big(K_\bullet(U,J),\partial\big),h
\end{align}
is a contraction.
\end{SATZ}
\begin{BEWEIS} The proof can be reduced to a  local computation, which is well-known.
\end{BEWEIS}
\\

Similar to Corollary \ref{exttor} we can use the projective Koszul resolution to compute 
\begin{eqnarray*}
\Tor_\bullet^{\:\mathcal C^\infty(U)}\big(\mathcal C^\infty(C),\mathcal C^\infty(C)\big)&=&\Gamma^\infty(C,\wedge^\bullet\:TC^\annind)\\
\Ext^\bullet_{\:\mathcal C^\infty(U)}\big(\mathcal C^\infty(C),\mathcal C^\infty(C)\big)&=&\Gamma^\infty(C,\wedge^\bullet\:{TC^\annind}^*).
\end{eqnarray*}

\section{Projective Koszul-Tate complexes}\label{projkostate}
For the moment, let $A$ be an arbitrary commutative $\mathbbm K$-algebra, think of it as the algebra of smooth functions on a manifold or the algebra of germs of real analytic functions etc., and $I$ be an ideal in $A$. In the spirit of the preceding subsection we will need the notion of a  \emph{projective presentation} of an  $A$-module 
$\mathcal W$,  that is, a short exact sequence of $A$-modules 
\[\mathcal V\stackrel{J}{\to}\mathcal W\to 0,\]
where $\mathcal V$ is a projective $A$-module. It is called \emph{finite} if $\mathcal V$ is finitely generated. We will call $J$ a \emph{system of projective generators} for $\mathcal W$. A typical example is the tautological section of the preceding subsection, which is a system of projective generators for the vanishing ideal of the submanifold $C$. In the following we will be exclusively interested in the situation, where either 
\begin{enumerate}
\item  $A$ is arbitrary, but $\mathcal V$ is a finitely generated free $A$-module, or 
\item \label{vectorbundle}  $A=\mathcal C^\infty(M)$ and  $\mathcal V$ is a finitely generated projective $A$-module, i.e., the space of sections of a vector bundle $V$ over $M$.
\end{enumerate}
In any case, there is a good notion of a basis for $\mathcal V$. In case \ref{vectorbundle}. this will be a local frame for $V$.

Generalizing the notion of Koszul resolution for modules over a (commutative) Noetherian ring, Tate \cite{Tate} introduced what is nowadays called a Tate resolution. For an exository article on the theory of Tate-resolutions for Noetherian (local) rings we refer to \cite[section 6 and 7]{Avramov}. Since we are concerned with the case $\mathbbm Q\subset \mathbbm K$ we can use a slightly simplified version of this construction (we replace the algebra of divided powers by the symmetric algebra). On the other hand, the ring of primary interest for us, $A=\mathcal C^\infty(M)$, is not a Noetherian local ring. Therefore, it makes sense to slightly generalize the construction by considering projective Tate generators. The potential of Tate resolutions for being useful in phase space reduction has already been recognized by the theoretical physicists Batalin, Fradkin and Vilkovisky in the early 80s \cite{BF,BV1,BV2,BV3}. In \cite{Stashbull} Stasheff recognized that their ``ghost for ghost'' procedure and that of Tate \cite{Tate} essentially coincide. Traditionally in the mathematical physics literature Tate resolutions are called Koszul-Tate resolutions.

A projective Koszul-Tate resolution of the $A$-module $A/I$ is a direct limit of a directed system of projective Koszul-Tate complexes $(KT_\bullet^{\le i},\partial^{\le i})$ of level $i\ge 1$. These are defined inductively as follows. 
 The (nonunique) \emph{projective Koszul-Tate complex}  of level $i\ge 1$ over the  module $A/I$ for an ideal $I$ of the commutative ring $A$ is a complex $(KT_\bullet^{\le i},\partial^{\le i})$
\begin{eqnarray*}
0 \leftarrow A=KT^{\le i}_0\stackrel{\partial^{\le i}}{\leftarrow}KT^{\le i}_1\stackrel{\partial^{\le i}}{\leftarrow}KT^{\le i}_2\leftarrow\dots \leftarrow KT^{\le i}_i\leftarrow KT^{\le i}_{i+1}\leftarrow\dots,
\end{eqnarray*}
such that 
\begin{enumerate}
\item $KT^{\le i}_k$ is the degree $k$ part of a graded symmetric algebra $S_A(\mathcal V^{\le i})$ generated by an $\mathbbm N$-graded projective $A$-module $\mathcal V^{\le i}=\oplus_{j= 1}^i \mathcal V^j$.
\item $\partial^{\le i}$ is an $A$-linear derivation of $KT^{\le i}=S_A(\mathcal V^{\le i})$.
\item The restriction map $\res:A\to A/I$ gives an isomorphism $\operatorname H_0 KT^{\le i}\cong A/I$.
\item $KT^{\le i}_k$ is acyclic up to degree $i-1$, that is $\operatorname H_k KT^{\le i}=0$ for $k=1,\dots ,i-1$.
\item For all $i\le j$ the obvious injections $KT^{\le i}\to KT^{\le j}$ are maps of supercommutative differential graded $\mathbbm K$-algebras. 
\item For $i=1$ the restriction $J_{(1)}$ of $\partial ^{\le 1}$ to $\mathcal V^1=KT^{\le 1}_1$ is a projective presentation of the $A$-module $I$. For $i\ge2$ the restriction $J_{(i)}$ of $\partial^{\le i}$ to $\mathcal V^i\subset KT^{\le i}_i$ is a system of projective generators for the $A$-module $\operatorname H_i KT^{\le i}$.  More precisely, we lift a projective presentation $\mathcal V^i$ of the $A$-module $\operatorname H_{i-1} KT^{\le i-1}$ to the space of cycles
\begin{align}\label{ghostforghost}
\xymatrix{&&&\mathcal V^i\ar[dl]_{J_{(i)}}\ar[d]&\\
0\ar[r]&\operatorname B_{i-1} KT^{\le i-1}\ar[r]&\operatorname Z_{i-1} KT^{\le i-1}\ar[r]&\operatorname H_{i-1} KT^{\le i-1}\ar[r]\ar[d]&0.\\
&&&0&}
\end{align}
\end{enumerate}  
The elements of $\mathcal V^j$ are called \emph{Tate generators of level $j$}. In physics they are called  \emph{antighosts of level $j$}. Note, that $(KT^{\le 1}_\bullet,\partial^{\le1})$ coincides with (a projective version of) the usual Koszul complex for $A/I$. The process described by diagram (\ref{ghostforghost}), which is called by the physicists the ``ghost for ghost''-procedure, is named by the commutative algebraists ``killing cycles by adding variables''. 

Given the data $(KT^{\le i}_\bullet,\partial^{\le i})$ for $i=1,2,\dots$ there is a unique differential $\partial$ on $KT_\bullet:=S_A(\oplus_{i\ge 1}\mathcal V^i)$ such that the obvious injections $(KT^{\le i}_\bullet,\partial^{\le i})\to (KT_\bullet,\partial)$ are morphisms of supercommutative differential graded algebras. By construction it is a projective resolution of $A/I$ and will be called the \emph{Koszul-Tate resolution} henceforth. Sometimes we would like to stress the dependence on the ring and the presentation $J=J_{(1)}$ of the ideal $I$ and will write in this case $KT_\bullet=KT_\bullet(A,J)$. 

Even though the above definition applies quite generally, there seems to be no  \emph{a priori} reason for the Koszul-Tate complexes to be complexes. We still have to convince ourselves that the operators $\partial^{\le i}$ and $\partial$ are of square zero. At this point we use a (local) basis $\xi_1^{(i)},\dots,\xi_{\ell_i}^{(i)}$ for the $A$-module\footnote{Let us assume for convenience that it is finitely generated.} $\mathcal V^i$ for $i=1,2,\dots$. The Koszul-Tate differentials $\partial^{\le i}$ and $\partial$ are given by the (local) formulas
\begin{align}
\partial^{\le i}=\sum_{j=1}^i\sum _{a=1}^{\ell_j} J^{(j)}_a\ddAG{j}{a},\\
\partial=\sum_{j=1}^\infty\sum _{a=1}^{\ell_j} J^{(j)}_a\ddAG{j}{a},
\end{align}
where $J^{(j)}_1,\dots, J^{(j)}_{\ell_j}\in KT_{j-1}^{\le j-1}$ are the components of the map $J_{(j)}$ in diagram (\ref{ghostforghost}). We proof by induction that $\partial^{\le i}$ is of square zero for $i=1,2,\dots$. For $i=1$ this is obvious. For $i\ge 2$ we obtain
\begin{align}\label{dsquareiszero}
(\partial^{\le i})^2=(\partial^{\le i-1})^2+\sum_{a=1}^{\ell_{i}}\big(\partial^{\le i}J^{(i)}_a\big) \ddAG{i}{a}+\sum_{a,b=1}^{\ell_{i}}J^{(i)}_a\frac{\partial J_b^{(i)}}{\partial\xi_a^{(i)}}\:\ddAG{i}{b}=0.
\end{align}
In the above formula the first term  on the right hand side vanishes by induction. The second term vanishes due to diagram (\ref{ghostforghost}). Finally, the third term vanishes since $J_1^{(i)},\dots ,J_{\ell_i}^{(i)}\in KT_{i-1}^{\le i-1}$ clearly do not depend on $\xi_1^{(i)},\dots,\xi_{\ell_i}^{(i)}$. Since for every homogeneuos $x\in KT_k$ there is an $i\in\mathbbm N$ such that $x$ is in the subcomplex $KT_k^{\le i}$ we conclude that $\partial ^2=0$. It seems to be an interesting question, whether equation (\ref{dsquareiszero}) holds in more complicated situations.

Obviously, the Koszul-Tate resolution of $A/I$ is nonunique. However, if $(A,\mathfrak m)$ is a Noetherian local ring\footnote{In this case all projective modules are free.}(e.g. the ring of germs of real analytic functions), then there is a distinguished Tate resolution, the so-called \emph{minimal model} of $A/I$, which is uniquely determined up to isomorphism of differential graded algebras \cite[Proposition 7.2.4]{Avramov}. Moreover, the number of Tate generators is bounded from below by invariants of the ring $(A,\mathfrak m)$, the so-called \emph{deviations}. The lower bounds are realized by the minimal model \cite[Proposition 7.2.5]{Avramov}. The deviations can be read off from the Poincar\'e series of $A/I$. It is known \cite{AvramovPrivate} for the important example  \ref{drehimpuls} of one particle in dimension $n\ge 3$ with angular momentum zero that the number of Tate generators grows exponentially with the level (in the polynomial setup one can show that $A/I$ is a Golod ring). In fact, due to the rigidity theorem of S. Halperin \cite{Halperin} we know that if $A/I$ is not a complete intersection, then \emph{none} of the deviations vanish.

Before we show that in the case of moment maps there are reasonable Koszul-Tate resolutions let us introduce some terminology.
If there is a smallest integer $L$ such that $\mathcal V^i=0$ for all $i>L$, then $KT_\bullet$ is said to be \emph{finitely generated of level $L$}. Otherwise it is said to be \emph{infinitely generated}. It is clear from what is said above that, in general, there is no reason to expect that the Koszul-Tate resolution of a singular moment map is finitely generated. More reasonable is the following property. The Koszul-Tate resolution is said to be \emph{locally finite} if all the modules $\mathcal V_i$ are finitely generated projective $A$-modules. In this case the rank of the projective module $KT_i$ can be recursively determined from the ranks of the projective modules $\mathcal V^j$, $1\le j\le i$ according to the product formula
\begin{align}
\sum_{i\ge 0}\operatorname{rank}(KT_i)\:t^i=\prod_{j\ge 1}\big(1-(-t)^j\big)^{(-1)^{j+1}\operatorname{rank}(\mathcal V^j)}.
\end{align}
The product on the right hand side converges in the $t$-adic topology of $\mathbbm Z[[t]]$.
If all the modules $\mathcal V_i$ are free $A$-modules,  we say that the Koszul-Tate resolution is \emph{free}.

\begin{PROPOSITION} \label{existsfiniteTate}
Let $J:M\to \mathfrak g^*$ be the moment map of a linear Hamiltonian action of a compact Lie group $G$ on the real symplectic vector space $M$. Moreover, suppose that $J$ generates the vanishing ideal $I_Z \subset \mathcal C^\infty(M)$ of the zero fibre $Z=J^{-1}(0)$. Then there is a locally finite free Koszul-Tate resolution $\left(KT_\bullet=KT_\bullet(C^\infty(M),J),\partial\right)$ of the ring of smooth functions $C^\infty(Z)$ in $Z$ and a continuous contracting homotopy $h:KT_\bullet\to KT_{\bullet+1}$, a continuous prolongation map $\prol:\mathcal C^\infty(Z)\to KT_0$, such that 
\begin{align}
\label{KTContr}
  \big(\mathcal C^{\infty}(Z),0\big)\:
  \begin{array}{c}{\res}\\\leftrightarrows\\{\prol}
  \end{array}\:\big(KT,\partial\big),h
\end{align}
is a contraction.
\end{PROPOSITION}
\begin{BEWEIS} First of all, let us identify $M$ with $\mathbbm R^{2n}$ and note that the components of the moment map are quadratic polynomial functions with respect to the canonical coordinates. Since $\mathbbm R[x_1,\dots,x_{2n}]$ is a Noetherian ring, there is a free locally finite Koszul-Tate resolution $KT_\bullet (\mathbbm R[x_1,\dots,x_{2n}]),J)$ of the $\mathbbm R[x_1,\dots,x_{2n}]$-algebra $\mathbbm R[x_1,\dots,x_{2n}]/<J_1,\dots,J_\ell>$ (see \cite{Avramov}). According to subsection \ref{ghsubsection} the latter is nothing but the coordinate ring of the real-variety determined by $J$. The complex $KT_\bullet(C^\infty(M),J)$ we are looking for is obtained by tensoring the complex $KT_\bullet (\mathbbm R[x_1,\dots,x_{2n}]),J)$ with the $\mathbbm R[x_1,\dots,x_{2n}]$-module $\mathcal C^\infty(M)$, and we have to show that this complex is  still a resolution of the $\mathcal C^\infty(M)$-module $\mathcal C^\infty(Z)$. 
By a standard result \cite[Theorem 7.2]{Eisenbud}  for every $x\in M$ the ring of formal power series $\mathcal F_x$ around $x$, being the completion of the real polynomial ring with respect to the maximal ideal corresponding to $x$, is a flat $\mathbbm R[x_1,\dots,x_{2n}]$-module. 
Next, $\mathcal F_x$ is a faithfully flat $\mathcal C^\omega_x(M)$-module \cite[Proposition III 4.10]{Malgrange}. 
We conclude that $KT_i(C^\omega_x(M),J)$ is exact for $i\ge 1$.
Using a partition of unity and the fact that $\mathcal C^\infty_x(M)$ is a (faithfully) flat $\mathcal C^\omega_x(M)$-module \cite[Corollary VI 1.12]{Malgrange}, it follows that  $KT_i(C^\infty(M),J)$ is exact for all $i\ge 1$. It remains to proof the existence of the continuous prolongation map and the contracting homotopies. But this follows from the results of Bierstone and Schwarz \cite{Schwarzbier} precisely along the lines of the proof of Theorem \ref{KosContr}.
\end{BEWEIS}
\\

We expect that the conclusion of the above proposition holds for essentially all moment maps of compact Hamiltonian group actions on symplectic manifolds.
\section{The BRST-charge}\label{BRSTchargesection}
We are now ready to introduce the BRST-algebra, which will turn out to be a differential graded Poisson algebra $\big(\brsalg,\{\:,\:\}_R,\brsop\big)$. To this end, we view the space of chains $KT=S_A\mathcal V$ of Koszul-Tate resolution of the preceding section as a generalized graded manifold with base manifold $M$. Note that in the infinitely generated case the Koszul-Tate differential is clearly not a derivation of finite type. We will therefore not use the module of finite type derivations in this section. The $\mathbbm Z$-graded algebra $\brsalg=\oplus_{i\in\mathbbm Z}\brsalg^i$ underlying the ghost-cotangent space of $KT=S_A\mathcal V$ is given by
\begin{align}
\brsalg^i=\prod_{j\ge 0}KT^{i+j}\otimes KT_j.
\end{align}
Here we have used the shorthand notation
\begin{align*}
KT^j:=(S_A \mathcal V^*)^j
\end{align*}
for the degree $j$ part of the symmetric algebra $S_A \mathcal V^*$ over the (positively graded) module $\mathcal V^*=\oplus_{i\ge 1} {\mathcal V^i}^*$, which is dual to the module of Tate-generators.
We will use local frames $\xG{i}{1},\dots,\xG{i}{\ell_i}$ for the modules ${\mathcal V^i}^*$, the so-called \emph{momenta of level $i$} or \emph{ghosts of level i}, which are dual to the Tate generators $\xAG{i}{1},\dots,\xAG{i}{\ell_i}$ of level $i$.

The Koszul-Tate differential of the preceeding section naturally extends to a derivation
\begin{align}
\partial=\sum_{j=0}^\infty \sum_{a=1} ^{\ell_j}J^{(j)}_a\ddAG{j}{a},\qquad J^{(j)}_a\in KT^{\le j-1}_{j-1} \text{ for } a=1,\dots \ell_j
\end{align}
of the algebra $\brsalg$.
Slightly abusing the language, we introduce the \emph{filtration by ghost degree} $F^k\brsalg=\oplus_{i\in \mathbbm Z}F^k\brsalg^i$, where
\begin{align}F^k\brsalg^i=\prod_{i+j\ge k}KT^{i+j}\otimes KT_j
\end{align}
$F^k\brsalg$ can be identified with the space of endomorphism which annihilate  $\oplus_{i=0}^{k-1}KT_{i}$. The spaces $F^k\brsalg$ form a descending Hausdorff filtration
\begin{align}
 F^0\brsalg=\brsalg\supset F^1\brsalg\supset\dots\supset F^k\brsalg\supset F^{k+1}\brsalg\supset\dots,
\end{align}
which is preserved by the supercommutative multiplication: $F^k\brsalg\:F^l\brsalg \subset F^{k+l}\brsalg$. A $\mathbbm K$-linear map $\varphi:\brsalg\to\brsalg$ is said to be of \emph{filtration degree} $l$ if $\varphi(F^k\brsalg)\subset F^{k+l}\brsalg$. 

In Section \ref{RothPB} we have introduced the Rothstein bracket $\{\:,\:\}_R$. In the finitely generated case we have seen that $\{\:,\:\}_R$ defines a $\mathbb Z$-graded super-Poisson bracket on $\brsalg$. In the infinitely generated case a full proof of the Jacobiidentity has yet to be given. In the following we postulate that the  Jacobiidentity fulfilled. (Curiously, for the construction of the BRST-charge we will merely use the Jacobiidentity for the approximating brackets defined below).
The filtration by ghost degree and the Rothstein Poisson bracket $\{\:,\:\}_R$  are, in general, \emph{not compatible} in the following sense: there is no $r\in \mathbb N$ such that $\{F ^k\brsalg,F ^l\brsalg\}_R\subset F^{k+l-r}\brsalg$ for all $k,l\in\mathbbm Z$. This phenomenon is due to the algebraic part of the bracket: a ghost of arbitrarily high level can be killed if it is paired with an antighost of the same level. The geometric part of the Rothstein Poisson bracket, however, is filtered: $[[\widetilde{\Pi},F ^k\brsalg],F^l\brsalg]\subset F^{k+l}\brsalg$ (here $[\:,\:]$ denotes the Schouten-Nijenhuis bracket). 

In order to prove the main result of this section,  Theorem \ref{chargeexists} below, we need a more refined analysis, which we will explain next.  For free Tate resolutions Theorem \ref{chargeexists} is a well known result, which goes back to \cite{BF,BV1,BV2,BV3}. The  first rigorous proof seems to be due to J. Stasheff and is sketched in \cite{Stashbull}. We have also benefited from the more elaborate exposition in \cite{Kimura}. The refinement consists in examining successively the level $i$ BRST algebras $\ibrsalg{i}^\bullet$ together with the level $i$ Koszul-Tate differentials 
\begin{align}
\partial ^{\le i}=\sum_{j=0}^i \sum_{a=1} ^{\ell_j}J^{(j)}_a\ddAG{j}{a}
\end{align}
and their level $i$ Rothstein Poisson brackets, which will be denoted by 
$\{\:,\:\}_{\le i}$ for $i=1,2,\dots$. In this manner we obtain a family of differential graded Poisson algebras 
\begin{align} 
\big(\ibrsalg{i}^\bullet,\partial^{\le i},\{\:,\:\}_{\le i}\big)_{i=1,2,\dots}. 
\end{align}
For $i<j$ the inclusion $\ibrsalg{i}\hookrightarrow\ibrsalg{j}$ is a map of differential graded commutative algebras, which is compatible with the filtration by ghost degree. We emphasize that, due to the presence of curvature, $\ibrsalg{i}\hookrightarrow\ibrsalg{j}$ is, in general, \emph{not} a Poisson subalgebra! The same remarks apply to the inclusions $\ibrsalg{i}\hookrightarrow\brsalg$ for $i=1,2,\dots$.

\begin{SATZ}[Existence of the BRST-charge] \label{chargeexists} Let $Z\subset M$ be a first class constraint set, and let $I:=I_Z\subset \mathcal C^\infty (M)=:A$ be the vanishing ideal of $Z$. Furthermore, let $(KT_\bullet(A,J),\partial)$ be a locally finite projective  Koszul-Tate resolution of the $A$-module $A/I$. Let $\brsalg$ be the corresponding BRST-algebra together with the Poisson bracket $\{\:,\:\}_R$. Then there exists an element $\theta\in\brsalg^1$ such that
\begin{enumerate}
\item $\{\theta,\:\}_R=\partial + \mbox{higher order terms}$,
\item $\{\theta,\theta\}_R=0$.
\end{enumerate}
Here, ``higher order terms'' stands for a $\mathbbm K$-linear derivation of filtration degree 1. 
\end{SATZ}
\begin{BEWEIS} We will construct $\theta$ by induction. To this end 
we make the following Ansatz $\theta_{\le i}:=\sum_{j=1}^i \theta_j$, where
\begin{align}\label{Ansatz}
\theta_j=\sum_{a=1}^{\ell_j} J^{(j)}_a\:\xi^a_{(j)}+Q_{(j)}\in \ibrsalg{j}^1,
\end{align}
such that $Q_{(j)}\in {KT_{\le j}^{j+1}}\otimes KT^{\le j}_j$. Note that $Q_{(j)}$ is at least quadratic in the momenta. We will show that the $Q_{(i)}$, $i=1,2,\dots$, may be successively chosen, such that
\begin{align}\label{filtprop}
\{\theta_{\le i},\theta_{\le i}\}_{\le i}\in F^{i+1}\ibrsalg{i}^2.
\end{align}
for all $i=1,2,\dots$. From the Ansatz (\ref{Ansatz}) it is clear that the $\theta_i$ add up to a well defined $\theta=\sum_{j=1}^\infty \theta_j\in \brsalg^1$. By an argument similar to that in equation (\ref{rothcompare}) below we conclude that 
\begin{align}
\{\theta_{\le i},\theta_{\le i}\}_R-\{\theta_{\le i},\theta_{\le i}\}_{\le i}\in F^{i+1}\brsalg^2.
\end{align} 
Since the filtration by ghost degree is Hausdorff we obtain $\{\theta,\theta\}_R=0$.

First of all, let us take a closer look at the derivations
\[\brsop^{\le i}:=\{\theta_{\le i},\:\}_{\le i}:\ibrsalg{i}^\bullet\to\ibrsalg{i}^{\bullet+1}\]
for $i=1,2,\dots$. In fact, the operators $\brsop^{\le i}$ are  (as we will see in a moment) of filtration degree $0$. We decompose them into its homogeneous components $\brsop^{\le i}=\sum_{j=0}^\infty \brsop^{\le i}_j$ (see also diagram (\ref{brsdecomp}) below), where ${\brsop_j^{\le i}}\big({KT^k_{\le i}}\otimes KT_l^{\le i}\big)\subset {KT^{k+j}_{\le i}}\otimes KT_{l+j-1}^{\le i}$ for all $k,l\ge 0$. A crucial observation is that the term of lowest order in this decomposition is exactly the level $i$ Koszul-Tate differential
\begin{align}\label{filtrationisgood}
\partial^{\le i}=\brsop^{\le i}_0.
\end{align}
This can be proven as follows. As we have already indicated, the part of $\brsop^{\le i}$ which originates from the geometric part of the Rothstein bracket, which  equals $-[\theta_{\le i},[\widetilde{\Pi}_{\le i},\:]]$, is of filtration degree 1. Therefore, we have for every $\alpha \in F^k\ibrsalg{i}^n$ 
\begin{eqnarray}
\brsop^{\le i}(\alpha)=\{\theta_{\le i},\alpha\}_{\le i}\in&-&\sum_{j=1}^i\sum_{a=1}^{\ell_j}\Big[\theta_{\le i},\Big[\ddG{j}{a}\ddAG{j}{a},\alpha\Big]\Big]+ F^{k+1}\ibrsalg{i}^{n+1}\nonumber\\
&=&\sum_{j=1}^i\sum_{a=1}^{\ell_j} \Big(\Big(J^{(j)}_a+\sum_{l=1}^i\frac{\partial Q_{(l)}}{\partial\xG{j}{a}}\Big)\frac{\partial \alpha}{\partial \xAG{j}{a}}\nonumber\\
&&\qquad\qquad+(-1)^{(n+1)(i+1)}\frac{\partial \alpha}{\partial \xG{j}{a}}\frac{\partial \theta_{\le i}}{\partial \xAG{j}{a}}\Big) +F^{k+1}\ibrsalg{i}^{n+1}\nonumber\\
&\stackrel{(*)}{=}&\partial^{\le i}(\alpha)+\sum_{l=1}^i\sum_{j=1}^{l-1}\sum_{a=1}^{\ell_j}(-1)^{(n+1)(i+1)}\frac{\partial \alpha}{\partial \xG{j}{a}}\frac{\partial J_a^{(l)}}{\partial \xAG{j}{a}}\xG{l}{a}\label{remainingterms}\\
&&+\sum_{l=1}^i\sum_{j=1}^{l}\sum_{a=1}^{\ell_j}(-1)^{(n+1)(i+1)}\frac{\partial \alpha}{\partial \xG{j}{a}}\frac{\partial Q_{(l)}}{\partial \xAG{j}{a}}+F^{k+1}\ibrsalg{i}^{n+1}.\nonumber
\end{eqnarray}
At step $(*)$ we have used the fact that since $Q_{(l)}$ is at least quadratic in the momenta $\partial Q_{(l)}/\partial \xG{j}{a}\in F^1\ibrsalg{l}^{1-j}$. We claim that the remaining terms in equation (\ref{remainingterms}) are also in $F^{k+1}\ibrsalg{i}^{n+1}$. In fact for $j+1\le l\le i$ we have
\begin{align}
\frac{\partial \alpha}{\partial \xG{j}{a}}\frac{\partial J_a^{(l)}}{\partial \xAG{j}{a}}\xG{l}{a}\in (F^{k-j}\ibrsalg{i}^{n-j})(\ibrsalg{l-1}^{n+1-l+j})(F^{l}\ibrsalg{l}^l)\subset F^{k+1}\ibrsalg{i}^{n+1}.
\end{align}
Finally, for $j\le l\le i$ we have 
\begin{align}
\frac{\partial \alpha}{\partial \xG{j}{a}}\frac{\partial Q_{(l)}}{\partial \xAG{j}{a}}\in (F^{k-j}\ibrsalg{i}^{n-j})(F^{l+1}\ibrsalg{l}^{1+j})\subset F^{k+1}\ibrsalg{i}^{n+1},
\end{align}
and the proof of equation (\ref{filtrationisgood}) is finished.

In order to start the induction we look at the projective presentation
\begin{align*}
\mathcal V^1=KT_1\stackrel{J=J_{(1)}}{\longrightarrow}I\to 0.
\end{align*}
Since the ideal $I$ is first class  the term $r_1$ of lowest order in $\half\{J_{(1)},J_{(1)}\}_{\le 1}$ vanishes when restricted to $Z$. By definition of the Koszul-Tate resolution there is a $Q_{(1)}\in KT^2_{\le 1}\otimes KT_1^{\le 1}$ such that $r_1=-\partial^{\le 1}Q_{(1)}$. Setting $\theta_{\le 1}=\theta_1=J_{(1)}+Q_{(1)}$, we obtain
\begin{eqnarray*}
\{\theta_1,\theta_1\}_{\le 1}&=&\{J_{(1)},J_{(1)}\}_{\le 1}+2\:\{J_{(1)},Q_{(1)}\}_{\le 1}+\{Q_{(1)},Q_{(1)}\}_{\le 1}\\
&\in&\{J_{(1)},J_{(1)}\}_{\le 1}+2\partial^{\le 1}q_1+F^2\ibrsalg{1}\\
&\in&2r_1+2\partial^{\le 1}Q_{(1)}+F^2\ibrsalg{1}
\end{eqnarray*}
and we are done. Note that we may write the system of projective generators in local coordinates $J_{(1)}=\sum_{a=1}^{\ell_1}J^{(1)}_a\xi^a_{(1)}$. Since the  ideal is first class, we have $\{J^{(1)}_a,J^{(1)}_b\}=\sum_{c=1}^{\ell_1} f_{ab}^c\:J^{(1)}_c$ for some (nonunique)  functions $f_{ab}^c$, the so-called structure functions. Locally, we can  write $Q_{(1)}$ in terms of the structure functions: $Q_{(1)}=-\half \sum_{a,b,c=1}^{\ell_1} f_{ab}^c\: \xi_{(1)}^a\xi_{(1)}^b\xi^{(1)}_c$.

Let us now assume that the Ansatz (\ref{Ansatz}) fulfills equation (\ref{filtprop}) for $i$. We are looking for a $Q_{(i+1)}$ such that equation (\ref{filtprop}) is true for $i\to i+1$. Taking advantage of the Jacobi identity for the bracket $\{\:,\:\}_{\le i}$
\begin{align*}
0=\{\theta_{\le i},\{\theta_{\le i},\theta_{\le i}\}_{\le i}\}_{\le i}=\brsop^{\le i}\{\theta_{\le i},\theta_{\le i}\}_{\le i},
\end{align*}
we deduce
\begin{align}\label{icycle}
\partial^{\le i}\{\theta_{\le i},\theta_{\le i}\}_{\le i}=-\sum_{j=1}^{\infty}\brsop^{\le i}_j\{\theta_{\le i},\theta_{\le i}\}_{\le i}\qquad\in F^{i+2}\ibrsalg{i}^3.
\end{align}
Let $r_i\in {KT_{\le i}^{i+2}}\otimes KT^{\le i}_i$ be the term of lowest degree in $\half\{\theta_{\le i},\theta_{\le i}\}_{\le i}$. It follows from equation (\ref{icycle}) that 
\[\partial^{\le i} r_i=0.\]
By construction of the Koszul-Tate resolution there is a $Q_{(i+1)}\in {KT_{\le i}^{i+2}}\otimes KT^{\le i+1}_{i+1}$, such that 
\[r_i=-\partial^{\le i+1} Q_{(i+1)}.\]
Setting $\theta_{i+1}:=\sum_{a=1}^{\ell_{i+1}}J_a^{(i+1)}\xi^a_{(i+1)}+Q_{(i+1)}$ we have to make sure that $\{\theta_{\le i+1},\theta_{\le i+1}\}_{\le i+1}\in F^{i+2}\ibrsalg{i+1}^2$. First of all, note that 
\[\partial^{\le i+1} \theta_{i+1}=\sum_{a=1}^{\ell_{i+1}}\big(\partial^{\le i+1}J^{(i+1)}_a\big)\:\xi_{(i+1)}^a+\partial^{\le i+1}Q_{(i+1)}=-r_i.\]
We conclude that $\partial^{\le i+1}\theta_{i+1}+\half\{\theta_{\le i},\theta_{\le i}\}_{\le i}\in F^{i+2}\ibrsalg{i}^2\subset F^{i+2}\ibrsalg{i+1}^2$. It is important to note that from the very definition of the Rothstein Poisson bracket we have  
\begin{align}\label{rothcompare}\{\theta_{\le i},\theta_{\le i}\}_{\le i+1}-\{\theta_{\le i},\theta_{\le i}\}_{\le i}\:\in  F^{i+2}\ibrsalg{i+1}^2.
\end{align} 
More precisely, this difference does merely originate from the geometric part of the Rothstein Poisson bracket. The lowest order contributions involve the curvature of the vector bundle $V^{i+1}$ of level $i+1$ Tate generators and are  generated by $\xi^a_{(i+1)}\xi_b^{(i+1)}\xi^c_{(1)}$ for $a,b=1,\dots ,\ell_{i+1}$ and $c= 1,\dots, \ell_1$.
As a result we obtain
\begin{eqnarray*} 
\{\theta_{\le i+1},\theta_{\le i+1}\}_{\le i+1}&=&\{\theta_{\le i},\theta_{\le i}\}_{\le i+1}+2\:\{\theta_{\le i},\theta_{i+1}\}_{\le i+1}+\{\theta_{i+1},\theta_{i+1}\}_{\le i+1}\\
&=&\{\theta_{\le i},\theta_{\le i}\}_{\le i+1}+2\:\brsop^{\le i+1}(\theta_{i+1})-\{\theta_{i+1},\theta_{i+1}\}_{\le i+1}\\
&\in&\{\theta_{\le i},\theta_{\le i}\}_{\le i}+2\:\brsop^{\le i+1}(\theta_{i+1})-\{\theta_{i+1},\theta_{i+1}\}_{\le i+1}+F^{i+2}\ibrsalg{i+1}^2\\
&\in& 2\sum_{j=1}^\infty\brsop^{\le i+1}_j(\theta_{i+1})-\{\theta_{i+1},\theta_{i+1}\}_{\le i+1}+F^{i+2}\ibrsalg{i+1}^2.
\end{eqnarray*}
From the very definition (\ref{Ansatz}) we see that  $\theta_{i+1}\in F^{i+1}\ibrsalg{i+1}^1$. Since  $\sum_{j=1}^\infty\brsop^{\le i+1}_j$ is of filtration degree $1$ we obtain $\sum_{j=1}^\infty\brsop^{\le i+1}_j(\theta_{i+1})\in F^{i+2}\ibrsalg{i+1}^2$.
\end{BEWEIS}
\\

The element $\theta$ is called the \emph{BRST charge}. 
Note that in the case of a moment map $J=J_{(1)}$ the first two terms of the series for the charge read
\begin{eqnarray*}
\theta_1&=&\sum_{a=1}^{\ell_1} J_a\:\xG{1}{a}-\half\sum_{a,b,c=1}^{\ell_1}f_{ab}^c\: \xG{1}{a}\xG{1}{b}\xAG{1}{c},\\
\theta_2&=&\sum_{a=1}^{\ell_2} J^{(2)}_a\:\xG{2}{a}\:,
\end{eqnarray*}
where $f_{ab}^c$ are the structure constants of the Lie algebra $\mathfrak g$. In other words, the quadratic term $Q_{(2)}$ \emph{vanishes} in this case. On the other hand, in the case of a projective Koszul-resolution (cf. section \ref{projkos}) the $\theta_i$ for $2\le i\le \ell$ consist merely of the quadratic term (here we have  of course $\theta_i=0$ for $i>\ell$). 

Associated to the charge there is a homological Hamiltonian vector field 
\begin{align*}
\brsop:=\{\theta,\:\}_R
\end{align*}
which is called the \emph{BRST differential}.
Although the homogeneous components of $\brsalg^\bullet$ are  not direct sums, there is a unique decomposition $\brsop=\sum_{i=0}^\infty \brsop_i=\partial+\sum_{i=1}^\infty \brsop_i$ of the BRST differential, such that ${\brsop_i}\big(KT^j\otimes KT_k\big)\subset KT^{j+i}\otimes KT_{k+i-1}$:
\begin{align}\label{brsdecomp}
\xymatrix{\cdot&\cdot&\cdot&\cdot&\\
          \cdot&\cdot&\cdot&\cdot&\\
          \cdot&\cdot&\cdot&\cdot&\\
          \cdot&\cdot\ar[l]^{\brsop_0}\ar[u]^{\brsop_1}\ar[uur]^{\brsop_2}\ar[uuurr]_{\dots}&\cdot&\cdot&}
\end{align}
The equation $\brsop^2=0$ translates into a sequence of relations, starting with
\begin{eqnarray}
\partial^2&=&0\\
\partial\: \brsop_1+ \brsop_1\:\partial&=&0\\
\partial\: \brsop_2+ \brsop_2\:\partial+\brsop_1^2&=&0.\label{Dsquared2}
\end{eqnarray}
The restriction map $\res:\mathcal C^\infty(M)\to \mathcal C^\infty(Z)$, which is an augmentation for the Koszul-Tate resolution extends to a degree zero map 
\begin{align}
\res: \brsalg^\bullet\to KT^\bullet_{|Z}.
\end{align} 
Here, $KT^\bullet_{|Z}$ is the space of smooth sections of the restriction of the bundle of ghosts to $Z$ and $\res$ is defined  to vanish on antighosts and acts on ghost by restricting the coefficients. By definition the smooth sections of this bundle are those which arise by restriction. Hence the restriction map $\res$ is \emph{onto}. Therefore, the formula $d\res=\res\brsop_1$ defines uniquely a $\mathbbm K$-linear degree 1 map $KT^\bullet_{|Z}\to KT^{\bullet+1}_{|Z}$. If $\brsop_1^2$ is evaluated on an element containing no antighosts, the result will be, due to equation (\ref{Dsquared2}), in the kernel of $\res$. It follows that $d^2=0$. We will call the cochain complex 
\[(KT^\bullet_{|Z},d)\]
the \emph{vertical complex}. A priori the vertical complex depends upon choices made for the space of generators $\mathcal V$, for the connection on $\mathcal V$, for the differential $\partial$, for $\theta$ and for $\prol$. It is conjectured that the homotopy class of the vertical complex -- in a sense yet to be specified -- does not depend on these choices.
The cochain complex $(KT^\bullet_{|Z},d)$ is well understood for certain special cases (see section \ref{cois} and section \ref{classred}). It may be shown by a spectral sequence argument, that $\res$ is in fact a quasiisomorphism of complexes (see also Figure 3.1). 
\begin{figure}[h]\label{quasiisom}
\begin{center}
\small
\psfrag{res}{$\res$}
\psfrag{Anull}{$\brsalg^0$}
\psfrag{del}{$\partial$}
\psfrag{D1}{$\brsop_1$}
\psfrag{CZ}{$\mathcal C^\infty(Z)$}
\psfrag{CM}{$\mathcal C^\infty(M)$}
\psfrag{dots}{$\dots$}
\psfrag{d}{$d$}
\psfrag{KT1}{$KT_1$}
\psfrag{KT2}{$KT_2$}
\psfrag{1KT}{$KT^1$}
\psfrag{2KT}{$KT^2$}
\psfrag{1KTZ}{${KT^1}_{|Z}$}
\psfrag{2KTZ}{${KT^2}_{|Z}$}
\psfrag{1KT1}{$KT^1\otimes KT_1$}
\psfrag{1KT2}{$KT^1\otimes KT_2$}
\psfrag{2KT1}{$KT^2\otimes KT_1$}
\psfrag{2KT2}{$KT^2\otimes KT_2$}
\normalsize
\includegraphics{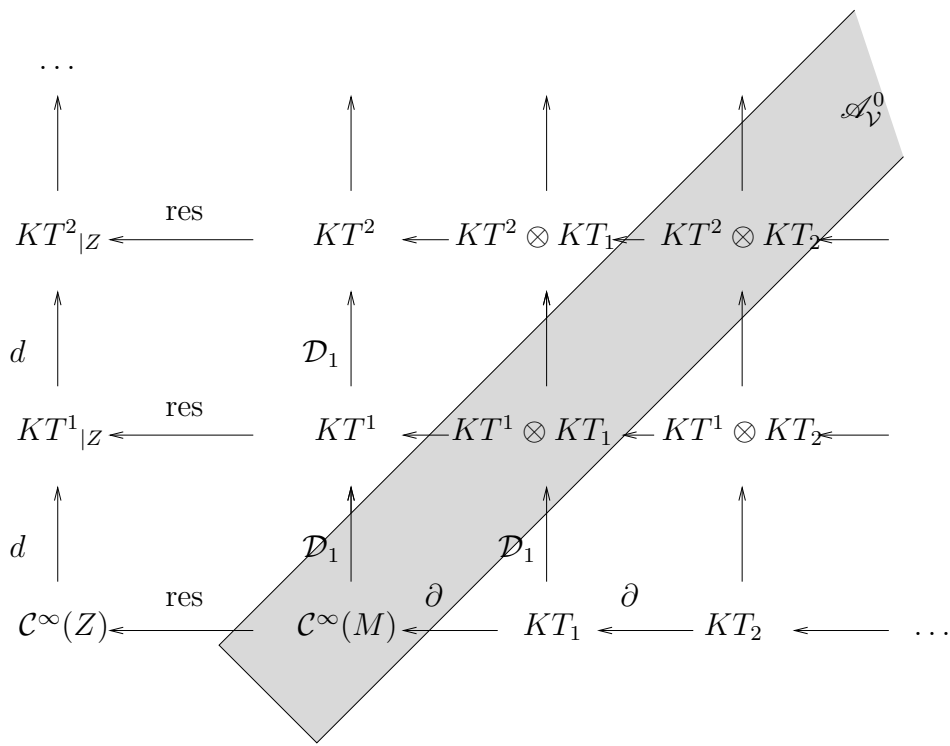}
\caption{The map $\res$  is a continuous quasiisomorphism from the BRST complex  to the vertical complex, which has a continuous split. Hence the  zeroth cohomology algebra of the BRST complex is toplogically isomorphic to the space of invariant functions $\mathcal C^\infty(Z)^I$, i.e., the Dirac reduced algebra of the first class constraint. We emphasize that the BRST complex is, in general, not a double complex.}
\end{center} 
\end{figure}
In fact, if we are given a continuous prolongation map $\prol:\mathcal C^\infty(Z)\to \mathcal C^\infty(M)$ and continuous contracting homotopies $h_i:KT_i\to KT_{i+1}$ for the Koszul-Tate resolution we can be more specific. Note that these maps naturally extend to maps $\prol:KT^\bullet_{|Z}\to KT^\bullet$ and $h_i:KT^\bullet\otimes KT_i\to KT^\bullet\otimes KT_{i+1}$.
\begin{KOROLLAR} Assume that  the premises for Theorem \ref{chargeexists} are true and, in addition, we have a contraction (\ref{KTContr}) as in Proposition \ref{existsfiniteTate}. Then there are continuous $\mathbbm K$-linear maps $\Phi:KT_{|Z}^\bullet\to \brsalg^\bullet$ and $H:\brsalg^{\bullet}\to \brsalg^{\bullet +1}$ such that
\begin{align}
\label{KTbrstcontr}
\big(KT^\bullet_{|Z},d\big)\:
\begin{array}{c} {\res}\\\leftrightarrows\\ {\Phi}
\end{array}
\:(\brsalg^\bullet,\brsop),H
\end{align}
is a contraction. If $[a],[b]\in \mathrm H\: KT_{|Z}$ are the cohomology classes of the cocycles $a,b \in\mathrm Z\: KT_{|Z}$, then 
\begin{align}\label{HPTredpoiss}
\{[a],[b]\}:=[\res\{\Phi(a),\Phi(b)\}_R]
\end{align}
 defines a $\mathbbm Z$-graded Poisson algebra structure on $\mathrm H^\bullet KT_{|Z}$. In degree zero $\mathrm H^0KT_{|Z}=\mathcal C^\infty(M)^I$ this Poisson structure coincides with the Dirac reduced Poisson structure.
\end{KOROLLAR}
\begin{BEWEIS} We apply perturbation lemma \ref{BPL1} to the contraction (\ref{KTContr}) and obtain the contraction (\ref{KTbrstcontr}). It is a straightforward matter to check that the bracket is well defined. The Jacobi identity follows immediatly from the  computation
\begin{eqnarray*}
\{[a],\{[b],[c]\}\}&=&\Big[\res\big\{\Phi(a),\Phi\big(\res\{\Phi(b),\Phi(c)\}_R\big)\big\}_R\Big]\\
&=&\Big[\res\big\{\Phi(a),\{\Phi(b),\Phi(c)\}_R\big\}_R\Big]-\Big[\res\big\{\Phi(a),\brsop H\{\Phi(b),\Phi(c)\}_R\big\}_R\Big]\\
&&-\Big[\res\big\{\Phi(a),H\{\Phi(db),\Phi(c)\}_R+H\{\Phi(b),\Phi(dc)\}_R\big\}_R\Big]\\
&=&\Big[\res\big\{\Phi(a),\{\Phi(b),\Phi(c)\}_R\big\}_R\Big]-\Big[d \res\big\{\Phi(a), H\{\Phi(b),\Phi(c)\}_R\big\}_R\Big]\\
&=&\Big[\res\big\{\Phi(a),\{\Phi(b),\Phi(c)\}_R\big\}_R\Big]\\
\end{eqnarray*}
for cocycles $a,b,c\in \mathrm Z\:KT_{|Z}$. 
The Leibniz rule is a consequence of the following consideration. Given two cocycles $b,c \in \mathrm Z\: KT_{|Z}$ the difference $\Delta:=\Phi(bc)-\Phi(b)\Phi(c)$ is, in general, nonzero since $\Phi$ (as well as $\prol$) is not multiplicative. Nonetheless, because $\res$ is multiplicative and $\res\:\Phi=\id$, we have $\Delta\in \kernel(\res)$. On the other hand, it is closed, since $\Phi$ is a chain map. It follows that $\Delta= (\brsop H+ H \brsop+\Phi\res)\Delta=\brsop (H( \Delta))$ is a coboundary, and consequently
\begin{eqnarray*}
\{[a],[b][c]\}=\big[\res\{\Phi(a),\Phi(bc)\}_R\big]=\big[\res\{\Phi(a),\Phi(b)\Phi(c)\}_R\big]+\underbrace{\big[d\res\{\Phi(a),H(\Delta\}_R\big]}_{=0}.
\end{eqnarray*}
 
Since in degree zero the differential $d:\mathcal C^\infty(Z)\to KT_{|Z}^1$ is given by the formula $d(a)=\res\{J_{(1)},\prol(a)\}_R$, it is clear that $\mathrm H^0 KT_{|Z}$ is  the space of $I$ invariant smooth functions on $Z$. In order to see that the induced bracket on  $\mathrm H^0 KT_{|Z}$ is the Dirac reduced bracket, note that the image $\Phi(a)\subset \prol(a)+F^1\brsalg^0$ of a function $a\in\mathcal C^\infty(Z)$ starts with $\prol(a)$. Since $\{F^1\brsalg^0,F^1\brsalg^0\}_R\subset F^1\brsalg^0\subset\kernel(\res)$  and $\{\mathcal C^\infty(M),F^1\brsalg^0\}_R\subset F^1\brsalg^0$, the higher order terms do not contribute, and the reduced bracket of two invariant functions $a,b\in \mathcal C^\infty(Z)^I$ is given by the formula $\{a,b\}=\res\{\prol a,\prol b\}$. 
\end{BEWEIS}
\\

We would like to stress that the algebraic properties of the contraction (\ref{KTbrstcontr}) which enable us to transfer the Poisson structure are: 1.) the right hand side is a differential graded Poisson algebra and 2.) $\res$ is a map of graded commutative agebras.
\section{Coisotropic submanifolds} \label{cois}
First class constraint sets, which are closed submanifolds are, per definition, \emph{coisotropic submanifolds}. It can be easily seen that a closed submanifold $C$ of a Poisson manifold $(M,\Pi)$ is coisotropic iff
\begin{align}\label{coiscond}
\Pi(\alpha,\beta)(c)=0 \quad \forall c\in C, \quad\forall \alpha,\beta\in TC^\annind.
\end{align}
Another equivalent characterization of coisotropy of the closed submanifold $C$ is that 
\begin{align}\label{anchor}
\#_\Pi(TC^\annind)\subset TC,
\end{align}
where $\#_\Pi$ is the restriction to $C$ of the anchor map $\#_\Pi:T^*M \to TM$, $\#_\Pi(\alpha)=i(\alpha)\Pi$ of the Lie algebroid associated to $\Pi$. It is well-known that $\#_\Pi$ in equation (\ref{anchor}) is an anchor map of a Lie algebroid over $C$. The bracket on $TC^\annind$ (cf. Lemma \ref{conormallemma}) is given by the formula $[\mathrm df,\mathrm dg]_{|c}:=\mathrm d\{f,g\}_{|c}$, where $c\in C$ and $f,g\in \mathcal C^\infty(M)$ are functions vanishing on $C$. 

Coisotropic submanifolds arise naturally in many geometric situations. For example, the graph of a Poisson map $\varphi:M\to N$ is a coisotropic submanifold in $M\times\bar{N}$, where $\bar{N}$ means $N$ understood with the opposite Poisson structure. Furthermore, if $J: M\to \mathfrak g^*$ is the moment map of a Hamiltonian $G$-action and if $J$ intersects the coadjoint orbit $\mathcal O\subset \mathfrak g^*$ cleanly, then $J^{-1}(\mathcal O)\subset M$ is a coisotropic submanifold.

Given a coisotropic submanifold $C\subset M$ there is an important subspace $\poly(C,M)$ of the Gerstenhaber algebra of polyvector fields $\poly(M)$, which is defined as follows:
\begin{eqnarray*}
\poly^0(C,M)&:=&I_C=\{f\in\mathcal C^\infty(M)\mid f_{|C}=0\}, \mbox{ and}\\
\poly ^k(C,M)&:=&\{\:X\in \poly^k(M)\mid X_c(\alpha_1,\dots,\alpha_k)=0 \:\:\forall c\in C\mbox{ and }\alpha_1,\dots,\alpha_k\in T_cC^{\annind}\:\}
\end{eqnarray*}
for $k\ge 1$. Moreover, we consider the canonical  map $\Psi:\Gamma^\infty(M,TM)\to \Gamma^\infty(C,TM_{|C}/TC)$, which restricts a vector field to $C$ and takes the residue class  of the result modulo $TC$. After identifying $TM_{|C}/TC$ with the dual of the annihilator bundle we obtain a map 
\[\wedge ^k\Psi:\poly^k(M)\to\Gamma^\infty(C,\wedge^k{TC^\annind}^*).\]
The latter space is just the space of cochains of Lie algebroid cohomology of $TC^\annind$. In fact, it follows from the next proposition that $\poly(C,M)$ is a \emph{coisotropic ideal}  in the differential Gerstenhaber algebra $(\poly^\bullet(M),\wedge,[\:,\:],\delta_\Pi)$, in the sense of the definition which was given in subsection \ref{derbrsection}. It is a curious fact that in this picture the Poisson tensor is a `first class constraint' by itself.
\begin{PROPOSITION} $\poly(C,M)$ has the following properties:
\begin{enumerate}
\item \label{wedgeideal}$\poly(C,M)^i\wedge \poly^j(M)\subset\poly^{i+j}(C,M)$,
\item \label{firstclassprop}$[\poly^i(C,M),\poly^j(C,M)]\subset\poly^{i+j-1}(C,M)$,
\item \label{subcomplex}$\delta_\Pi\poly^i(C,M)\subset\poly^{i+1}(C,M)$
\end{enumerate}
for all $i,j\ge0$. Moreover, $\wedge ^k\Psi$ is onto and the kernel of $\wedge ^k \Psi$ is just $\poly^k(C,M)$ for all $k\ge 0$. Thus, we can identify the quotient space $\poly^k(M)/\poly^k(C,M)$ with $\Gamma^\infty(C,\wedge^k{TC^\annind}^*)$ for all $k\ge 0$. The induced differential on $\Gamma^\infty(C,\wedge^k{TC^\annind}^*)$ coincides with the differential of Lie algebroid cohomology.
\end{PROPOSITION}
\begin{BEWEIS} \ref{wedgeideal}.) and \ref{firstclassprop}.) follow from the formulas (\ref{cupform}), (\ref{bullform}) and Theorem \ref{multiderisgerst} since for the Lie Rinehart pair $\left(\mathcal C^\infty(M),\Gamma^\infty(M,TM)\right)$ the Gerstenhaber algebras of polyvector fields and multiderivations coincide. Note that for \ref{firstclassprop}.) one also makes use of Lemma \ref{conormallemma}. 
Statement \ref{subcomplex}.) is a consequence of \ref{firstclassprop}.). The claim that $\wedge^k\Psi$ is onto can easily be proven by a partition of unity argument. The last statement follows straightforwardly from the fact that the Koszul-Brylinsky bracket of two exact forms $\mathrm df$ and $\mathrm dg$ is given by $\mathrm d\{f,g\}$.
\end{BEWEIS}
\begin{KOROLLAR}\label{derbrPoissonstructure} There is a natural $\mathbbm Z$-graded super-Poisson structure on the cohomology of the Lie algebroid $TC^\annind$. The induced Poisson strucure in degree zero coincides with the Dirac reduced bracket. 
\end{KOROLLAR}

\begin{BEWEIS} Use the derived bracket of Theorem \ref{redderbr}.
\end{BEWEIS}
\\

Alternatively, we could use a projective Koszul resolution (\ref{prKosContr}) of $\mathcal C^\infty(C)$ for some tubular neighborhood $U$ of $C$ in $M$ and the construction of Section \ref{BRSTchargesection} to aquire the cohomology of the vertical complex $KT_{|C}$ with a $\mathbbm Z$-graded super-Poisson structure according to formula (\ref{HPTredpoiss}). It is not difficult to prove that the vertical complex coincides with the cochain complex of Lie algebroid cohomology of the Lie algebroid $TC^{\annind}$ as above. However, it is not clear to the author whether the $\mathbbm Z$-graded Poisson structure of formula (\ref{HPTredpoiss}) coincides with that of Corollary \ref{derbrPoissonstructure}. Moreover, these Poisson structures still have  to be compared with the $P_\infty$-algebra structure on the vertical complex which has been introduced in \cite{relform}.
\section{Classical BRST-algebra for Hamiltonian group actions} \label{classred}

In the remainder of this work we will exclusively be concerned with the following important, and most simple, special case of the BFV-construction. We will consider a Hamiltonian $G$-space $M$, $G$ compact and connected, with a moment map $J:M\to \mathfrak g^*$ for which the Koszul complex
\begin{align}
K_{-\bullet}=K_{-\bullet}(C^\infty(M),J)=S^\bullet_{C^\infty(M)}\big(\mathfrak g[1]\big),\qquad\partial=\sum_{a=1}^\ell J_a\uddxi{a},
\end{align} 
is a free resolution of the $A:=\mathcal C^\infty(M)$-module $\mathcal C^\infty(Z)$ of smooth functions on the zero fibre $Z=J^{-1}(0)$. Here $\xi_1,\dots,\xi_\ell$ denote a basis for $\mathfrak g[1]$ and $\uddxi{1},\dots,\uddxi{\ell}$ the corresponding algebraic vector fields on the graded manifold $S^\bullet_{A}\big(\mathfrak g[1]\big)$. We have readily seen in Subsection \ref{ghsubsection} and Section \ref{kc} that the Examples \ref{mpart}, \ref{lemon}, \ref{stratres}, \ref{tzwo} for $\alpha<0$ and \ref{commvar} provide in fact examples of such moment maps.

Since we do not need to mention the level here, we will simply call the elements of $\mathfrak g^*[-1]$ and $\mathfrak g[1]$  \emph{ghosts} and \emph{antighosts}, respectively. Dually to the basis $\xi_1,\dots ,\xi_\ell$ for $\mathfrak g[1]$, we also will need a basis $\xi^1,\dots ,\xi^\ell$ for $\mathfrak g^*[-1]$. The respective indices will run over  latin letters: $a,b,\dots$. Since all ghost variables are purely odd, we will identify the BRST algebra $\brsalgo:=\brsalgo_{\mathcal V}$, where $\mathcal V=\Gamma^\infty\left(M,\mathfrak g[-1]\times M\right)$ is the space of sections of the trivial vector bundle with fibre $\mathfrak g[-1]$, with the space of polynomials in the ghosts and antighosts:
\begin{align*}
\brsalgo=S_A\left(\mathfrak g[1]\oplus \mathfrak g^*[-1]\right).
\end{align*}
Alternatively, we could identify $\brsalgo$ with the space of sections of the trivial vector bundle over $M$ with fibre $\wedge(\mathfrak g\oplus \mathfrak g^*)=\wedge\mathfrak g^*\otimes \wedge\mathfrak g$ understood with the appropriate grading. There is an even graded Poisson bracket on $S_{\mathbbm K}\left(\mathfrak g[1]\oplus \mathfrak g^*[-1]\right)$, which is defined by the formula
\begin{align}\label{oddpoiss}
\{v,w\}=-2\mu\circ \left(\sum_{a=1}^\ell\uddxi{a}\otimes \ddxi{a}+\ddxi{a}\otimes\uddxi{a}\right)(v\otimes w).
\end{align}
Here, $\mu$ denotes the super-commutative multiplication in $S_{\mathbbm K}\left(\mathfrak g[1]\oplus \mathfrak g^*[-1]\right)$ and $\otimes$ is the graded tensor product. This Poisson bracket is the unique even super-Poisson bracket such that $\{\xi_a,\xi^b\}=2\delta_a^b$ and $\{\xi_a,\xi_b\}=0=\{\xi^a,\xi^b\}$ for all $a,b=1,\dots,\ell$. Note the slight \emph{change of convention}: the above bracket differs from Rothstein bracket by a factor of $2$! If $f,g\in\mathcal C^\infty(M)$ and $v,w\in S_{\mathbbm K}\left(\mathfrak g[1]\oplus \mathfrak g^*[-1]\right)$, then the formula
\begin{align}
\{fv,gw\}:=\{f,g\}vw+fg\{v,w\}
\end{align}
defines a super-Poisson bracket on $\brsalgo$. 

From  the Lie bracket  of $\mathfrak g$ and the moment map 
we build the so-called \emph{BRST-charge}
\begin{align}
\theta:= -\frac{1}{4}\sum_{a,b,c}\:f_{ab}^c\:\xi^a\xi^b\xi_c+\sum_a J_a \xi^a\in\mathscr A^1,
\end{align}
 where the $f_{ab}^c$
are the structure constants of $\mathfrak g$. An easy calculation yields 
$\{\theta,\theta\}=0$. In other words, the recursion of Theorem \ref{chargeexists} breaks off after one step here. Hence, the \emph{classical BRST-differential} 
\[\brsop:=\{\theta,\:\}\] is in fact of square zero. 
Summing up, we obtain a differential graded Poisson algebra 
\[(\brsalgo,\{,\},\brsop=\{\theta,\:\}),\] 
which is called the \emph{classical BRST algebra} henceforth. 

Closer examination shows that $(\brsalgo,\brsop)$ is the total complex of  a certain double complex. In fact, there is a canonical identification $\brsalgo=S_A(\mathfrak g^*[-1])\otimes S_A(\mathfrak g[1])=C^\bullet\left(\mathfrak g,S_A(\mathfrak g[1])\right)$, where the latter denotes the space of Lie algebra cochains with values in the $\mathfrak g$-module $S_A\mathfrak g[1]$. This representation will be denoted henceforth by $\mathrm L$. In this way $\brsalgo=\oplus_{i,j=0}^\ell\brsalgo ^{i,j}$ acquires a $\mathbbm Z\times\mathbbm Z$ -grading
\[\brsalgo ^{i,j}:=S^i_A(\mathfrak g^*[-1])\otimes S^j_A(\mathfrak g[1])\]
The differential of the Lie algebra cohomology $\delta:\brsalgo ^{i,j}\to\brsalgo ^{i+1,j} $  corresponding to the $\mathfrak g$-module 
$S_{A}(\mathfrak g[1])$ can be writen as a super-differential operator
\begin{align*}
\delta=-\frac{1}{2}\sum_{a,b,c}f_{ab}^c\:\:\xi^a\xi^b\:\uddxi{c} +\sum_{a,b,c}f_{ab}^c\:\xi^a\xi_c\:\ddxi{b}+\sum_a\:\xi^a\{J_a,\:\}.
\end{align*} 
The Koszul-differential $\partial=\sum_{a=1}^\ell J_a\ddxi{a}$ extends
naturally to a differential $\partial:\brsalgo^{i,j}\to\brsalgo^{i,j-1}$. Because $J$ is equivariant, these two differentials super-commute: $\partial \delta+\delta\partial=0$. An easy calculation yields that
\begin{align*}
\brsop=2\partial+\delta.
\end{align*}
We will view $\brsop$ as a perturbation 
(see Appendix \ref{HPT}) of the acyclic differential $2\partial$.
 
We extend the restriction map $\res$ to a map 
$\res:\mathscr A\to S_{\mathcal C^\infty(Z)}(\mathfrak g^*[-1])$ by setting it zero for all 
terms containing antighosts and restricting the coefficients. In the same fashion, we extend 
$\prol$ to a map $S_{\mathcal C^\infty(Z)}(\mathfrak g^*[-1])\to\mathscr A$ extending the 
coefficients. 

Since the moment map $J$ is $G$-equivariant, $G$  acts on $Z=J^{-1}(0)$. Hence 
$\mathcal C^\infty(Z)$ is a $\mathfrak g$-module, this representation will be denoted by 
$\mathrm{L}^z$. Note that $\mathrm{L}^z_X=\res\;\mathrm L_X\;\prol$ for all $X\in\mathfrak g$. 
We identify  $S_{\mathcal C^\infty(Z)}(\mathfrak g^*[-1])$  with the space of cochains of Lie 
algebra cohomology $C^\bullet\big(\mathfrak g,\mathcal C^\infty(Z)\big)$.
Let us denote 
$d:C^\bullet\big(\mathfrak g,\mathcal C^\infty(Z)\big)\to C^{\bullet+1}
\big(\mathfrak g,\mathcal C^\infty(Z)\big)$ the codifferential of Lie algebra cohomology 
coresponding to $\mathrm{L}^z$. Since $\res$ is a morphism of $\mathfrak g$-modules we obtain 
$d\;\res=\res\;\delta$.

\begin{SATZ}There are $\mathbb K$-linear maps 
$\Phi:C^\bullet\big(\mathfrak g, \mathcal C^{\infty}(Z)\big)\to \mathscr A^\bullet$
and $H: \mathscr{A}^\bullet\to \mathscr{A}^{\bullet-1}$ which are continuous in the 
respective Fr\'echet topologies such that
\begin{eqnarray}\label{cbrstSDR}
\Big(C^\bullet\big(\mathfrak g, \mathcal C^{\infty}(Z)\big),d\Big)\:
\begin{array}{c} {\res}\\\leftrightarrows\\ {\Phi}
\end{array}
\:(\mathscr A^\bullet,\brsop),H
\end{eqnarray}
is a contraction.
\end{SATZ}
\begin{BEWEIS} Apply lemma \ref{BPL1} to the perturbation $\brsop$ of $2\partial$. Explicitly, 
we get 
\begin{align*}
H=\frac{1} {2} h\sum_{j=0}^{\ell}\Big(-\frac{1}{2}\Big)^j(h\delta+\delta h)^j,\\
\Phi=\prol-H(\delta\:\prol-\prol\:d),
\end{align*}
which are obviously Fr\'echet continuous. 
Note that from $h\prol=0$ and $h^2=0$ it follows that $H\Phi=0$ and $H^2=0$. If $\prol$ is 
chosen to be equivariant, then the expression for $\Phi$ simplifies to $\Phi=\prol$. In the 
same way one gets $H=\frac{1}{2}h$, if $h$ is equivariant.
\end{BEWEIS}

\begin{KOROLLAR} There is a  graded Poisson structure on 
$\mathrm H^\bullet\big(\mathfrak g,\mathcal C^\infty(Z)\big)$. If $[a],[b]$ are the cohomology 
classes of $a,b\in C^\bullet\big(\mathfrak g, \mathcal C^{\infty}(Z)\big)$, then the bracket is 
 given by 
\begin{align*}\{[a],[b]\}:=[\res\{\Phi(a),\Phi(b)\}].\end{align*}
The restriction of this bracket to 
$\mathrm H^0\big(\mathfrak g,\mathcal C^\infty(Z)\big)=\mathcal C^\infty (Z)^{\mathfrak g}$ 
coincides with the Dirac reduced Poisson structure. 
\end{KOROLLAR}

\chapter{Quantum BRST reduction}\label{qbrst}
In this chapter we will construct a deformation quantization of the classical BRST algebra for the situation, when the moment map satisfies the generating and complete intersection hypothesis, cf. section \ref{classred}. In order to define the quantum BRST algebra it is sufficient to assume that there is some quantum moment map, which deforms the original moment map. It has been observed in \cite{BHW}, that it is most convenient to use for the ghost variables the so-called standard ordered Clifford multiplication instead of the Weyl-ordered multiplication, since this renders the quantum BRST-complex a double complex. We will see that the quantum Koszul differential can also be found using standard homological algebra. This entails that  the quantum BRST-algebra can essentially be viewed as an Ext-algebra. We are able to compute the BRST-cohomology if 1.) the star product is assumed to be strongly invariant or 2.) the group is compact and semisimple. In both cases we find deformation quantizations for the classical reduced algebra. In these cases continuous star products, which deform the reduced Poisson algebra, can be found.
\section{The quantum BRST algebra}\label{brstalgconstr}
In this section we will introduce the quantum BRST algebra, which is 
$\mathbbm K[[\nu]]$-differential graded associative algebra $(\mathscr{A}^\bullet[[\nu]],*,\qbrst)$ deforming the differential graded Poisson algebra 
$(\mathscr A^\bullet,\{,\},\mathscr D)$. In order to define a graded product $*$ on $\mathscr A[[\nu]]$, 
we use on the one hand a formal Clifford multiplication 
\begin{eqnarray}
v\cdot w&:=&\mu\big(\operatorname{e}^{-2\nu\: \sum_a\ddxi{a}\otimes \uddxi{a}}(v\otimes w)\big)\nonumber\\
&=&\sum_{i=0}^\infty \frac{(-2\nu)^i}{i!}\sum_{a_1,\dots,a_i}(-1)^{i|x|}\frac{\partial ^i v}{\partial \xi^{a_1}\dots\partial \xi^{a_i}}\frac{\partial ^i w}{\partial \xi_{a_1}\dots\partial \xi_{a_i}}\label{clmult}
\end{eqnarray} 
for homogeneous $v,w\in S_{\mathbbm K}(\mathfrak g[1]\oplus \mathfrak g^*[-1])$.
Here $\mu$ denotes the super-commutative multiplication and $\otimes$ is the graded tensor product. The product $\cdot$ satisfies the Clifford relation $\xi_a\cdot \xi^b+\xi^b\cdot\xi_a=2\nu \delta_a^b$ for all $a,b=1\dots \ell$. There are of course other (equivalent) ways to define the Clifford multiplication, e.g. by symmetrisation. The above product, which is sometimes called the antistandard ordered product, has the advantage to render the quantum BRST-complex a double complex (see Theorem \ref{doppelcompl}). 

On the other hand, we will need a \emph{quantum  covariant} star product $\star$ on $M$ with \emph{quantum moment map} $\qimp$ (cf. Subsection \ref{qmom}).  Here the quantum moment map will be viewed as an element of the BRST algebra: $\qimp=J+\sum_{i\ge1}\nu^i J_i \in \mathfrak g^*[-1]\otimes \mathcal C^\infty(M)\fnu\subset \mathscr A^1\fnu$.
We will frequently refer to a basis $\xi_1,\dots,\xi_\ell$ of $\mathfrak g[1]$ and write for short $<\qimp,\xi_a>=\qimp_a$ for $a=1,\dots,\ell$. The representation property (\ref{qmomprop}) can be rewritten as 
\begin{eqnarray}\label{qmomentmap}
\qimp_a\star \qimp_b-\qimp_b\star 
\qimp_a=\nu\sum_c f_{ab}^c\:\qimp_c\qquad \mbox{ for }a,b=1,\dots,\ell,
\end{eqnarray}
where $f_{ab}^c$ are the structure constants of the Lie algebra $\mathfrak g$.

For $f,g\in \mathcal C^{\infty}(M)$ and 
$v,w \in S\big(\mathfrak g[1]\oplus \mathfrak g^*[-1]\big)$ we define 
\[(fv)*(gw):=(f\star g)\: (v\cdot w).\] 
Note, that $*$ defines a $\mathbbm Z$-graded associative product, which is a formal deformation of the super-Poisson structure of Section \ref{classred}. 

The next step is to quantize the BRST-charge. 
It was observed by Kostant and Sternberg \cite{KostStern} that
\[\qch:= -\frac{1}{4}\sum_{a,b,c}\:f_{ab}^c\:\xi^a\xi^b\xi_c+
\mathbbm \sum_a\qimp_a\:\xi^a+\frac{\nu}{2}\sum_{a,b} f^b_{ab}\:\xi^a \:\in\mathscr A^1[[\nu]],\]
is a good guess. Here $f_{ab}^c$ denote the structure constants of the Lie algebra. Note, that the trace term is a side effect of the operator ordering. It may be absorbed by redefining the quantum moment map $\qimp':=\qimp+ \frac{\nu}{2}\sum_{a,b} f^b_{ab}\:\xi^a$. In fact, since the trace of the adjoint representation vanishes on commutators, this modification does not spoil the representation property (\ref{qmomentmap}). To start with $\qimp'$ from scratch is considered as slightly incorrect from the point of view of representation theory of deformed algebras. Nevertheless, we will sometimes do it and will indicate that by the prime.
\begin{SATZ} \label{qsqzero} $\qch*\qch=0$.
\end{SATZ} 
\begin{BEWEIS} Let us write for short $\qch=Q+\qimp'$, where $Q =-\frac{1}{4}\sum_{a,b,c}f_{ab}^c\:\xi_c\xi^a\xi^b$ and $\qimp'=\sum_a\qimp_a\xi^a+\half \sum_{a,b}f_{ab}^b\xi^a$. We have
\begin{align*}
Q*Q=Q^2+2\nu\:\frac{1}{4^2}\:2 \sum_{a_1,b_1,c_1,a_2,b_2,c_2}f_{a_1b_1}^{c_1}\:f_{c_1b_2}^{c_2}\:\xi_{c_2}\xi^{a_1}\xi^{b_1}\xi^{b_2}=0,
\end{align*}
since $Q$ is odd and the bracket satisfies the Jacobi identity. The terms of higher order in $\nu$ vanish due to degree reasons. Furthermore, we have 
\begin{align*}
\qimp'*\qimp'=\sum_{a,b}\qimp'_a\star\qimp'_b\:\xi^a\xi^b=\half\sum_{a,b}[\qimp'_a,\qimp'_b]_{\star}\:\xi^a\xi^b=\frac{\nu}{2}\sum_{a,b}f_{ab}^c\:\qimp'_c\:\xi^a\xi^b.
\end{align*}
Hence, it remains to compute
\begin{align*}
Q*\qimp'+\qimp'*Q=Q\qimp'- 2\nu \frac{1}{4}\sum_{a,b,c} f_{ab}^c\:\qimp'_c \:\xi^a\xi^b+\qimp'Q=-\frac{\nu}{2}\sum_{a,b}f_{ab}^c\:\qimp'_c\:\xi^a\xi^b.
\end{align*}
It follows that $\qch*\qch=Q*Q+Q*\qimp'+\qimp'*Q+\qimp'*\qimp'=0$.
\end{BEWEIS}
\\

Now we are ready to define the \emph{quantum BRST differential} to be 
\[\qbrst:=\frac{1}{\nu}\ad_*(\qch).\]
 Before we take a closer look at $\qbrst$, let us introduce some terminology. We define the super-differential operators $\qce,\mathscr R,q,u:\mathscr A^\bullet\to\mathscr A^{\bullet+1},$
\begin{eqnarray*}
\qce(fv)&:=&-\frac{1}{2}\sum_{a,b,c}f_{ab}^c\:f\:\xi^a\xi^b\:\uddxi{c} v+\sum_{a,b,c}f_{ab}^c\:f\:\xi^a\xi_c\:\uddxi{b}v+\sum_a\:\frac{ 1}{\nu}[\qimp_a,f]_*\:\xi^av, \\      
\mathscr R(fv)&:=&\sum_a f * \qimp_a\:\ddxi{a}v,\quad\quad\quad\quad\quad\quad\quad\mbox{``right multiplication''}\\
q(fv)&:=&-\frac{1}{2}\sum_{a,b,c}f\:f_{ab}^c\:\xi_c\:\frac{\partial^2 v}{\partial\xi^a\partial\xi^b},\quad\quad\quad\quad\quad\quad\mbox {``quadratic ...''} \\
u(fv)&:=&\sum_{a,b}f\:f_{ab}^b\:\ddxi{a}v,\quad\quad\quad\quad\quad\quad \quad\quad\mbox{``unimodular term''}
\end{eqnarray*} 
for $f\in \mathcal C^\infty(M)$ and $v\in S_{\mathbbm K}(\mathfrak g[1]\oplus\mathfrak g^*[-1])$. Note that $\qce$ is the coboundary operator of Lie algebra cohomology corresponding to the representation 
\begin{eqnarray}
\qrep_X:S_{\mathcal C^\infty(M)}(\mathfrak g[1])[[\nu]]&\to& S_{\mathcal C^\infty(M)}(\mathfrak g[1])[[\nu]],\nonumber \\
fv&\mapsto& f(\ad_X(v)) +\nu^{-1}(\qimp(X)\star f- f\star \qimp(X))v,
\end{eqnarray} 
where $X\in \mathfrak g$, $v\in S_{\mathbbm K}(\mathfrak g[1])$ and $f\in \mathcal C^\infty(M)[[\nu]]$. Finally, we set
\begin{align}\label{defkos}
\qkos:=\mathscr R+\nu\Big(\frac{1}{2}u-q\Big).
\end{align}
This operator will be called the \emph{deformed} or \emph{quantum Koszul differential}. Clearly, adding to $\qimp$ a scalar multiple of the trace form does only have an effect on $\qkos$, the Lie algebra differential $\qce$ stays unchanged. If we would have started with $\qimp'$, then the unimodular term would not occur in formula (\ref{defkos}). It will become clear later, that the unimodular term does not have an effect on the homology of $\qkos$.
\begin{SATZ} \label{doppelcompl}The quantum BRST differential 
\begin{align}\label{total}
\qbrst=\qce+2\qkos
\end{align}
is a linear combination of two super-commuting differentials $\qce$ and $\qkos$. 
\end{SATZ}
\begin{BEWEIS} With the shorthand notation as in the proof of Theorem \ref{qsqzero} we compute
\begin{eqnarray*}
[\qimp',fv]&=&\qimp'*fv-(-1)^{|v|}fv*\qimp'\\
&=&\sum_a (\qimp'_a\star f-f\star\qimp'_a)\xi^av+2\nu\sum_a(f\star \qimp'_a)\ddxi{a}v,
\end{eqnarray*}
where $f\in \mathcal C^\infty(M)$ and $v\in S_{\mathbbm K}(\mathfrak g[1]\oplus\mathfrak g^*[-1])$. One also has to compute the Clifford commutator with the element $Q=-\frac{1}{4}\sum_{a,b,c}f_{ab}^c\:\xi_c\xi^a\xi^b$ with a homogeneous element $v\in S_{\mathbbm K}(\mathfrak g[1]\oplus\mathfrak g^*[-1])$:
\begin{eqnarray*}
[Q,v]_\cdot&=&Q\cdot v-(-1)^{|v|}v\cdot Q\\
&=&Q v+2\nu\sum_{a,b,c}\frac{-1}{4} f_{ab}^c\:\xi^a\xi^b\uddxi{c}v\\
&&-(-1)^{|v|}\Big(vQ+2\nu(-1)^{|v|}\sum_{a,b,c}\frac{-1}{4}\:f_{ab}^c\:2\xi_c\xi^b \ddxi{a}v+\frac{4\nu^2}{2}\sum_{a,b,c}\frac{-1}{4}\:f_{ab}^c\:2\frac{\partial^2 v}{\partial \xi^b\partial \xi^a}\xi_c\Big)\\
&=&-\frac{\nu}{2}\sum_{a,b,c}f_{ab}^c\:\xi^a\xi^b\uddxi{c}v+\nu\sum_{a,b,c}f_{ab}^c\:f\:\xi^a\xi_c\:\uddxi{b}v+2\nu^2 q.
\end{eqnarray*}
Collecting the terms, equation (\ref{total}) follows easily.
Since each of the three terms in $0=\qbrst^2=4\qkos^2+2(\qkos\qce+\qce\qkos)+\qce^2$ lives in different degrees (see below) the claim follows.
\end{BEWEIS}

\begin{KOROLLAR}  For all $X\in\mathfrak g$ we have 
\begin{align}\label{qkosequiv}
\qrep_X\:\qkos-\qkos\:\qrep_X=[\qrep_X,\qkos]=0.
\end{align}
\end{KOROLLAR}
\begin{BEWEIS} Let us write $X=\sum_a X^a\xi_a\in\mathfrak g$ with respect to a basis $\xi_1,\dots,\xi_\ell\in\mathfrak g$. For the insertation derivation $i_X:=\sum_a X^a\uddxi{a}$ we have the well known Cartan homotopy formula $\qrep_X=i_X\:\qce+\qce\:i_X=[i_X,\qce]$. Since $\qkos$ obviously commutes with $i_X$ the claim follows.
\end{BEWEIS}
\\

We conclude that the deformed Koszul complex $(K_\bullet[[\nu]],\qkos)$ is in fact a complex of $\mathfrak g$-modules. The BRST-cochains may be identified with the Lie algebra cochains of this representation. More precisely, we have $\mathscr A^n[[\nu]]=\oplus_{i,j,n=i-j}\mathscr A^{i,j}[[\nu]]$, where $\mathscr A^{i,j}[[\nu]]=C^i\big(\mathfrak g,K_j[[\nu]]\big)$. The quantum BRST-differential is (up to a trivial  factor of 2) the total differential of the double complex formed be the deformed Lie algebra cohomology differential $\qce:\mathscr A^{i,j}[[\nu]]\to \mathscr A^{i+1,j}[[\nu]]$ and the deformed Koszul differential $\qkos: \mathscr A^{i,j}[[\nu]]\to \mathscr A^{i,j-1}[[\nu]]$.
\section{Quantum BRST as an Ext-algebra}\label{extalg}

In this section we will give a conceptual explanation for the quantum BRST algebra similar to \cite{Sevost}. As the material is not needed in the following, the reader may skip this section. We will work over the field $\mathbbm K\launu$ of formal Laurent series, mainly because the standard reprensentation $\rho$ of Lemma \ref{strep} is not onto for formal power series in $\nu$. This has the drawback, that the classical limit makes no sense. Nonetheless, for the cases under consideration the classical limit is already at hand.

From the quantum symmetry point of view it is more natural to replace the Lie bracket $[\,,\,]$ on $\mathfrak g$ by $\nu[\,,\,]$. More precisely, we consider the $\mathbbm K\fnu$-Lie algebra $\mathfrak g\fnu$ with bracket $[\,,\,]_\nu:=\nu[\,,\,]$ and, accordingly, the universal envelopping algebra $U\mathfrak g_\nu=T\mathfrak g\fnu/<x\otimes y-y\otimes x-\nu [x,y]>$. This is an augmented $\mathbbm K\fnu$-algebra, the augementation map $\epsilon:U\mathfrak g_\nu\to\mathbbm K \fnu$ is induced by the obvious augmentation of $T\mathfrak g\fnu$. We consider the complex $X^\nu_\bullet=U\mathfrak g_\nu\otimes_{\mathbbm K} \wedge^\bullet \mathfrak g[[\nu]]$ with differential 
\begin{eqnarray*}
d(u\otimes x_1\wedge\cdots \wedge x_n)&:=&\sum_{i=1}^n(-1)^{i+1}ux_i\otimes x_1\wedge\cdots\widehat{x_i}\cdots\wedge x_n\label{cliffmult}\\
&&+\sum_{1\le i<j\le n}(-1)^{i+j}u\otimes [x_i,x_j]_\nu\wedge x_1\wedge\cdots\widehat{x_i}\cdots\widehat{x_j}\cdots\wedge x_n.
\end{eqnarray*}
In the literature this complex is frequently called Koszul resolution. For obvious reasons we refrain from using this terminology.
\begin{SATZ} $(X^\nu_\bullet,d)$ is a free resolution of the $U\mathfrak g_\nu$-module $\mathbbm K[[\nu]]$. 
\end{SATZ}
\begin{BEWEIS} The proof (see e.g. \cite[p.279--282]{CE}) relies on the PBW-theorem, which applies since $\mathfrak g[[\nu]]$ is a free $\mathbbm K[[\nu]]$-module. 
\end{BEWEIS}
\\

It is tempting to interprete the Clifford algebra as an algebra of super-differential operators. However, in the formal situation not every super-differential operator arises in this way, since the partial derivative is decorated with the formal parameter $\nu$. We solve this problem by brute force by formally inverting $\nu$. Secondly, we have to take opposite $\clopp$ of the Clifford algebra multiplication (\ref{cliffmult}). Here opposite is understood in the graded sense, i.e., $v\clopp w=(-1)^{|v||w|}w\cdot v$.  Let us identify $S_{\mathbbm K}(\mathfrak g[1]\oplus \mathfrak g^*[-1])$ with  $\wedge(\mathfrak g\oplus \mathfrak g^*)$, which is understood with the induced grading. We will use the symbols $e_1,\dots,e_{\ell}$  to denote a basis of $\mathfrak g$ and  $e^1,\dots,e^{\ell}$
will denote the corresponding dual basis. We will write $i(\alpha)$ for the super-derivation of $\wedge \mathfrak g$ which extends dual pairing with $\alpha\in\mathfrak g^*$. The Clifford multiplication now writes $\mu \exp\big(-2\nu \sum_a i(e^a)\otimes i(e_a)\big)$.

\begin{LEMMA} \label{strep}The so called standard representation
\begin{eqnarray*}\rho :\wedge(\mathfrak g\oplus \mathfrak g^*)\launu&\to& \End_\mathbbm K(\wedge \mathfrak g)\launu\\
(x_1\wedge\dots\wedge x_n\otimes \alpha_1\wedge\dots\wedge\alpha_m&\mapsto& \big (v\mapsto(-2\nu)^m x_1\wedge\dots\wedge x_n\wedge i(\alpha_1)\circ\dots\circ i(\alpha_m)v\big),
\end{eqnarray*}
where $x_1,\dots,x_n\in\mathfrak g$ and $\alpha_1,\dots,\alpha_m\in\mathfrak g$, is an isomorphism of algebras for the reversed Clifford multiplication $\clopp$, i.e $\rho(a\clopp b)=\rho(a)\circ\rho(b)$ for all $a,b\in\wedge(\mathfrak g\oplus \mathfrak g^*)\launu$. In fact, if we reverse the canonical $\mathbbm Z$-grading on $\End_\mathbbm K(\wedge \mathfrak g)\launu$ then $\rho$ is $\mathbbm Z$-graded.
\end{LEMMA}

A quantum moment map  $\qimp'$ gives rise to an algebra morphism $U\mathfrak g\fnu\to\mathcal C^\infty(M)\launu$. Applying the functor $\mathcal C^\infty(M)\launu\otimes_{U\mathfrak g_\nu}-$ on the complex $X^\nu_\bullet$ yields the complex
\[\mathcal C^\infty(M)\launu\otimes_{U\mathfrak g_\nu}X^\nu_\bullet=\mathcal C^\infty(M)\otimes \wedge^\bullet \mathfrak g\launu=K_\bullet\launu\]
with differential
\begin{align} \label{dtilde}
\widetilde d (f\otimes x)=\sum_a f\star \qimp'(e_a)\otimes i(e^a)x-\frac{\nu}{2}\sum_{a,b,c}f_{ab}^c\:f\otimes e_c\wedge i(e^a)i(e^b)x,
\end{align}
where $f\in\mathcal C^\infty(M)$ and $x\in\wedge\mathfrak g$. In fact, this differential essentially coincides with the quantum Koszul differential $\qkos$.
\begin{KOROLLAR} There is an isomorphism $\rho$ (extending the standard representation) of $\mathbbm Z$-graded algebras between the BRST algebra $(\mathscr A\launu,\stopp)$ with the reversed multiplication and the algebra of endomorphisms of the Koszul complex $\big(\End_{\mathcal C^\infty(M)\launu}(K^\bullet\launu),\circ\big)$, such that $\rho(\qch)=-2\nu\:\widetilde{d}$. 
\end{KOROLLAR}
\begin{BEWEIS} Let us write for short $A:=C^\infty(M)$. There is a canonical isomorphism of algebras
\[\End_{A\launu}(K^\bullet\launu)=\End_{A\launu}(A\otimes\wedge \mathfrak g\launu)\cong A\launu^{\oppind}\otimes\End_{\mathbbm K\launu}(\wedge \mathfrak g\launu).\]
The left hand side is an algebra with respect to composition and the right hand side is isomorphic to $(\brsalg\fnu,*^{\oppind})$. If we view the algebra on the left hand side as a $\mathbbm Z$-graded algebra with the reverse of the canonical grading,  then this isomorphism is actually an isomorphism of $\mathbbm Z$-graded algebras. The formula $\rho(\qch)=-2\nu\:\widetilde{d}$ follows by inspection.
\end{BEWEIS}

\begin{KOROLLAR} If the zeroth order term $J$ of the quantum moment map $\qimp'$ is a moment map satisfying the generating and complete intersection hypothesis, then the opposite of the BRST-cohomology algebra $\mathrm H^\bullet\mathscr A\launu$ is isomorphic to $\Ext^\bullet_{\mathcal C^\infty(M)\launu}(\mathcal B,\mathcal B)$ with composition product, where $\mathcal B$ is the left $\mathcal C^\infty(M)\launu$-module which is obtained by dividing out the left ideal which is generated by $\qimp_1',\dots,\qimp_{\ell}'$.
\end{KOROLLAR}
\begin{BEWEIS} It is clear (see Proposition \ref{qkoscontr} below) that the complex of equation (\ref{dtilde}) is a resolution of $\mathcal B$. Therefore, the cohomology algebra of the differential graded algebra 
\[\left(\End^{-\bullet}_{\mathcal C^\infty(M)\launu}\big(K_\bullet\launu,\big),\circ,[\widetilde{d},\:]\right)\]
is just the $\mathbbm Z$-graded algebra $\Ext^\bullet_{\mathcal C^\infty(M)\launu}(\mathcal B,\mathcal B)$ (cf. \cite[\S 7]{Bourbdix}).
\end{BEWEIS}

\section{Computation of the quantum BRST-Cohomology}\label{brstcomp}

The main idea which we follow  in order to compute the quantum BRST cohomology (i.e., the 
cohomology of $(\mathscr A[[\nu]],\qbrst)$), is to provide a deformed version of the 
contraction (\ref{cbrstSDR}). This will be done by applying Lemma \ref{BPL2} to the 
contraction (\ref{KosContr}) for the perturbation $\qkos$ of $\partial$ and then applying 
Lemma \ref{BPL1} for the perturbation $\qbrst$ of $2\qkos$. We will also need to examine a 
deformed version of the representation $\mathrm{L}^z$ of $\mathfrak g$ on $\mathcal C^\infty(Z)$. 

\begin{PROPOSITION} \label{qkoscontr} If we choose $h_0$ such that $h_0\prol=0$, then there are deformations of 
the restriction map 
$\qres=\res+\sum_{i\ge 1}\nu^i\;\res_{i}:\mathcal C^{\infty}(M)\to 
\mathcal C^{\infty}(Z)[[\nu]]$ and of the contracting homotopies 
$\qhop=h+\sum_{j\ge 1}\nu^j\;h_{(j)}:K_\bullet[[\nu]]\to K_{\bullet+1}[[\nu]]$,
which are a formal power series of Fr\'echet continuous maps and such that 
\begin{eqnarray}\label{qKosContr}
\big(\mathcal C^{\infty}(Z)[[\nu]],0\big)\:\begin{array}{c} {\qres}\\\leftrightarrows\\ {\prol}
\end{array}\:\big(K[[\nu]],\qkos\big),\qhop
\end{eqnarray}
is a contraction with $\qhop_0\;\prol=0$. Explicitly, we have 
\[\qres:=\res\;(\id+(\qkos_1 -\partial_1)h_0)^{-1}.\]
If we choose $h$ to be $\mathfrak g$-equivariant, the same is true for $\qhop$.
\end{PROPOSITION}
\begin{BEWEIS} Apply lemma \ref{BPL2} to the perturbation $\qkos$ of $\partial$.
\end{BEWEIS}
\\

We are now ready to define the quantized representation $\qrep^z$ of $\mathfrak g$ on  $\mathcal C^\infty (Z)[[\nu]]$ by setting 
\[\qrep^z_X:=\qres\;\qrep_X \;\prol \quad\mbox{ for $X\in \mathfrak g$.}\]
\begin{LEMMA} We have $\qrep^z_X\:\qrep^z_Y -\qrep^z_X\:\qrep^z_Y=\qrep^z_{[X,Y]}$ for all $X,Y\in\mathfrak g$.
\end{LEMMA}
\begin{BEWEIS} The claim follows from
\begin{eqnarray*}
\qrep^z_X\:\qrep^z_Y&=&\qres\:\qrep_X\:\prol\:\qres\:\qrep_X\:\prol\\
&\stackrel{(\ref{qKosContr})}{=}&\qres\:\qrep_X(\id-\qkos\qhop)\qrep_Y\:\prol\:\stackrel{(\ref{qkosequiv})}{=}\qres(\id-\qkos\qhop)\qrep_X\:\qrep_Y\:\prol
\stackrel{(\ref{qKosContr})}{=}\qres\:\qrep_X\:\qrep_Y\:\prol
\end{eqnarray*}
and the fact that $\qrep$ is a representation.
\end{BEWEIS}
\\

In the same fashion 
as in Section \ref{classred}, we define 
$\qdz:C^\bullet(\mathfrak g,\mathcal C^\infty (Z)[[\nu]])\to C^{\bullet+1}
(\mathfrak g,\mathcal C^\infty (Z)[[\nu]])$ to be the differential of Lie algebra cohomology 
of the representation $\qrep^z$, i.e., $\qdz\;\qres=\qres\;\qce$. 
In the same manner, we extend $\qres$ and $\qhop$  as in Section \ref{classred} to maps 
$\qres:\mathscr A\to C(\mathfrak g,\mathcal C^{\infty}(Z)[[\nu]]\big)$ and 
$\qhop:\mathscr A^\bullet[[\nu]]\to\mathscr A^{\bullet-1}[[\nu]]$. 

\begin{SATZ} 
\label{Hauptsatz}\label{maincontr}
  There are $\mathbbm K[[\nu]]$-linear maps 
  $\qPhi:C^\bullet\big(\mathfrak g, \mathcal C^{\infty}(Z)\big)\to \mathscr A^\bullet[[\nu]]$ 
  and $\qHop: \mathscr{A}^\bullet\to \mathscr{A}^{\bullet-1}[[\nu]]$,
  which are series of  Fr\'echet continuous maps such that
\begin{eqnarray}\label{maincontrformula}
  \Big(C^\bullet\big(\mathfrak g, \mathcal C^{\infty}(Z)[[\nu]]\big),\qdz\Big)\:
\begin{array}{c} {\qres}\\\leftrightarrows\\ {\qPhi}
\end{array}
  \:(\mathscr A^\bullet[[\nu]],\qbrst),\qHop
\end{eqnarray}
is a contraction.
\end{SATZ}
\begin{BEWEIS} 
Since the requisite condition $\qres\;\qhop=0$ is obviously fulfilled, we apply 
Lemma \ref{BPL1} to the perturbation $\qbrst$ of $2\qkos$. Explicitly, this means that 
$\qHop:=\frac{1} {2} \qhop\sum_{j=0}^\ell(-\frac{1}{2})^j(\qhop\qce+\qce \qhop)^j $ and 
$\qPhi=\prol-\qHop(\qce\;\prol-\prol\:\qdz)$, which are obviously series of Fr\'echet 
continuous maps. Note that from $h_0 \prol=0$ and $h^2=0$, we get $\qHop\qPhi=0$ and 
$\qHop^2=0$. If $\prol$ is chosen to be equivariant, then the expression for $\qPhi$ simplifies 
to $\qPhi=\prol$. If $h$ and (hence $\qhop$) is equivariant, then it follows that 
$\qHop=\frac{1}{2}\qhop$.
\end{BEWEIS}
\\

For better intelligibility, let us mention that the above argument is a strengthening of the well known tic-tac-toe lemma \cite[p.135]{BottTu}: The homology of the total differential of a double complex with acyclic rows is isomorphic to the homology of the differential, which is induced by the action of the vertical differential on the horizontal homology. The double complex in question is of course that of Theorem \ref{doppelcompl}. It is depicted in the following diagram from the second column onwards.
\begin{align*}
\xymatrix{
\cdots&\cdots&\cdots&\\
C^1\big(\mathfrak g, \mathcal C^{\infty}(Z)[[\nu]]\big)\ar[u]^\qdz&\mathscr A^{1,0}[[\nu]]\ar[l]^{\qquad\qres}\ar[u]^\qce&\mathscr A^{1,1}[[\nu]]\ar[l]^{\:\:\:\:2\qkos}\ar[u]&\qquad\cdots\qquad\ar[l]\\
C^0\big(\mathfrak g, \mathcal C^{\infty}(Z)[[\nu]]\big)\ar[u]^\qdz&K_0(J,M)[[\nu]]\ar[l]^{\quad\qres}\ar[u]^\qce&K_1(J,M)[[\nu]]\ar[l]^{\:\:\:\:2\qkos}\ar[u]^\qce&K_2(J,M)[[\nu]]\cdots\ar[l]^{2 \qkos}}
\end{align*}
The first column is the Lie algebra cochain complex of the $\mathfrak g$-module
$\mathcal C^\infty(Z)[[\nu]]$, which is quasiisomorphic to the BRST complex via the quasiisomorphism $\qres$.
  
We use the contraction (\ref{maincontrformula}) to transfer the associative algebra structure from $\mathscr A[[\nu]]$ to the Lie algebra cohomology $\mathrm H^\bullet\big(\mathfrak g,\mathcal C^\infty(Z)[[\nu]]\big)$ of the representation $\qrep^z$
by setting 
\begin{eqnarray}\label{Formel1}
  [a]*[b]:=[\qres\big(\qPhi(a)*\qPhi(b)\big)],
\end{eqnarray} 
where $[a],[b]$ denote the cohomology classes of the cocyles 
$a,b\in Z^\bullet\big(\mathfrak g, \mathcal C^{\infty}(Z)[[\nu]]\big)$. The associativity of this operation follows from
\begin{eqnarray}
[a]*([b]*[c])&=&\Big[\qres\Big(\qPhi(a)*\qPhi\big(\qres\big(\qPhi(b)*\qPhi(c)\big)\big)\Big)\Big]\nonumber\\
&=&\Big[\qres\big(\qPhi(a)*\qPhi(b)*\qPhi(c)\big)\Big]-\Big[\qres\Big(\qPhi(a)*\big(\qbrst\qHop\big(\qPhi(b)*\qPhi(c)\big)\big)\Big)\Big]\nonumber\\
&&-\Big[\qres\Big(\qPhi(a)*\big(\qHop\qbrst\big(\qPhi(b)*\qPhi(c)\big)\big)\Big)\Big]\nonumber\\
&=&\Big[\qres\big(\qPhi(a)*\qPhi(b)*\qPhi(c)\big)\Big]-\Big[\qdz\:\qres\Big(\qPhi(a)*\big(\qHop\big(\qPhi(b)*\qPhi(c)\big)\big)\Big)\Big]\nonumber\\
&&+\Big[\qres\Big(\qPhi(\qdz a)*\big(\qHop\big(\qPhi(b)*\qPhi(c)\big)\big)\Big)\Big]\nonumber\\
&&-\Big[\qres\Big(\qPhi(a)*\big(\qHop\big(\qPhi(\qdz b)*\qPhi(c)+\qPhi(b)*\qPhi(\qdz c)\big)\big)\Big)\Big]\nonumber\\
&=&\Big[\qres\big(\qPhi(a)*\qPhi(b)*\qPhi(c)\big)\Big]\label{assoccomp},
\end{eqnarray}
which coincides, as a result of a similar calculation, with $([a]*[b])*[c]$.

However, this is 
\emph{not exactly}, what we want to accomplish.
The primary obstacle on the way to the main result, Corollary \ref{redprod}, is that, in general, 
we have 
\begin{align}\label{anomaly}\mathrm H^0\big(\mathfrak g,\mathcal C^\infty(Z)[[\nu]]\big)\neq \mathrm H^0\big(\mathfrak g,\mathcal C^\infty(Z)\big)[[\nu]],
\end{align} 
since there is no \emph{a priori} reason that the representations $\qrep^z$ and $\mathrm L^z$ have the same space of invariants.
An example where this phenomenon does in fact occur has been 
given in \cite[section 7]{BHW}. One way out is to sharpen the compatibility condition (\ref{qmomentmap}).  
We require, that $\qimp=J$ and
\[J(X)\star f- f \star J(X)=\nu\{J(X),f\}\quad\mbox{ for all }X\in \mathfrak g, 
  f\in \mathcal C^{\infty}(M).\]
This property, which has been discussed in Subsection \ref{qmom}, is referred to as \emph{strong invariance} of the star product $\star$ with 
respect to the Lie algebra action. It can always be achieved for the cases under consideration. Of course, now the representations $\qrep$ and 
$\mathrm L$ coincide and we get $\delta=\qce$. But with some mild restrictions on the 
contracting homotopy $h$ of the Koszul resolution we also have the following.

\begin{LEMMA} 
  If $h_0$ is $\mathfrak g$-equivariant and $h_0\prol=0$, then 
  $\qrep^z=\mathrm L^z$.
\end{LEMMA}
\begin{BEWEIS} For $X\in\mathfrak g$ we have 
  $\qrep^z_X=\qres\;\mathrm L_X \;\prol=\res\;(\id+(\qkos_1 -\partial_1)h_0)^{-1}
  \mathrm L_X\;\prol$. Since $\mathrm L_X$ commutes with $\qkos_1$, $\partial_1$ and $h_0$, the 
  last expression can be written as 
  $\res\;\mathrm L_X(\id+(\qkos_1 -\partial_1)h_0)^{-1}\prol=\res\;\mathrm L_X \;\prol$. 
\end{BEWEIS}

\begin{KOROLLAR} 
  \label{redprod} With the assumptions made above, the product defined by equation 
  (\ref{Formel1}) makes $\mathrm H^\bullet\big(\mathfrak g,\mathcal C^\infty(Z)\big)[[\nu]]$ into a $\mathbbm Z$-graded associative algebra. For the subalgebra  of invariants $\mathrm H^{0}\big(\mathfrak g, \mathcal C^{\infty}(Z)\big)[[\nu]]= \big(\mathcal C^{\infty}(Z)\big)^{\mathfrak g}[[\nu]]$ this formula simplifies to 
\begin{eqnarray}\label{Formel2}
f*g:=\qres\big(\prol(f)*\prol(g)\big)\quad\mbox{ for }f,g\in\big(\mathcal C^{\infty}(Z)\big)^{\mathfrak g}.
\end{eqnarray}
Since $\big(\mathcal C^{\infty}(Z)\big)^{\mathfrak g}[[\nu]]$ is
$\mathbbm K[[\nu]]$-linearly isomorphic to the algebra of smooth functions 
on the symplectic stratified space $M_{\mathsf{red}}$, we obtain an associative
product on $\mathcal C^\infty(M_{\mathsf{red}})[[\nu]]$ which gives rise to a \emph{continuous} Hochschild cochain. 
\end{KOROLLAR}

There is another strategy to attack problem (\ref{anomaly}).
If $\mathrm H^1(\mathfrak g, C^\infty(Z))$ vanishes, it is possible to find a topologically linear isomorphism between the spaces of invariants for the classical and the deformed representation. 
\begin{PROPOSITION} Let $G$ be a compact, connected semisimple Lie group acting on the Poisson manifold $M$ in a Hamiltonian fashion. Assume that the equivariant moment map $J$ satisfies the generating and complete intersection hypothesis. Then for any star product $*$ on $M$ with quantum moment map $\qimp$ there is a invertible sequence of continuous maps 
\[S=\sum_{i\ge 0}\nu^i\:S_i :\mathrm H^0(\mathfrak g, \mathcal C^\infty(Z))[[\nu]]=\mathcal C^\infty(Z)^\mathfrak g[[\nu]]\to \mathrm H^0(\mathfrak g,\mathcal C^\infty(Z)[[\nu]])\] 
such that the formula 
\begin{align*}
f* g:= S^{-1}\big(S(f)*S(g)\big)=S^{-1}\Big(\qres\big(\qPhi(S(f))*\qPhi(S(g))\big)\Big)   
\end{align*}
defines a continuous formal deformation of the Poisson algebra $\mathcal C^\infty(Z)^\mathfrak g$ into an associative algebra. 
\end{PROPOSITION}
\begin{BEWEIS} According to Viktor L. Ginzburg (see \cite[Theorem 2.13]{Ginz98}) we have for any compact connected Lie group $G$ with a smooth representation on a Fr\'echet space $W$ an isomorphism 
\begin{align*} \mathrm H^\bullet(\mathfrak g,W)\cong \mathrm H^\bullet(\mathfrak g,\mathbbm K)\otimes W^\mathfrak g.
\end{align*}
In particular, this implies that if $\mathfrak g$ is semisimple the first and the second cohomology groups of the $\mathfrak g$-module $\mathcal C^\infty(Z)$ vanish (for the so called Whitehead lemmata see e.g. \cite{Hilton}). Note that
the Fr\'echet subspace of invariant functions $\mathcal C^\infty(Z)^\mathfrak g\subset \mathcal C^\infty(Z)$ has a closed complementary subspace $V$
\[\mathcal C^\infty(Z)=\mathcal C^\infty(Z)^\mathfrak g\oplus V.\] This can be achieved by taking $V$ to be the kernel of the averaging projection $\pi:\mathcal C^\infty(Z)\to \mathcal C^\infty(Z)^G$, $\pi(f)(x):=\operatorname{vol}(G)^{-1}\int_G f(gx)dg$.
Hence, the restriction of the Lie algebra cohomology differential $d$ to the closed complementary subspace $V$ is a bijection onto the closed supspace $Z^1(\mathfrak g,\mathcal C^\infty(Z))\subset C^1(\mathfrak g,\mathcal C^\infty(Z))$. Since every continuous linear bijection of Fr\'echet spaces has a continuous inverse (see \cite[corollary 2.12]{Rudin}), we have a continuous inverse map, which we call $d^{-1}$.
\begin{align*}
\xymatrix{\mathcal C^\infty(Z)\ar[r]^{d\quad}&\mathfrak g^*\otimes \mathcal C^\infty(Z)\\
V\ar@{^{(}->}[u]\ar[r]^{\cong\qquad}&Z^1(\mathfrak g,\mathcal C^\infty(Z))\ar@{^{(}->}[u]\ar@/^2pc/[l]^{d^{-1}}
}
\end{align*}

Let $\varphi_0\in\mathcal C^\infty(Z)^\mathfrak g$. We will inductively construct an element $\varphi=\sum_i\nu^i\varphi_i$ which is invariant for the deformed representation, i.e., $\qdz\varphi=0$.
Let us assume that we have found $\varphi_0,\varphi_1,\dots,\varphi_n\in \mathcal C^\infty(Z)$, such that
\begin{align}\label{IV}
\sum_{j=0}^i d_j\varphi_{i-j}=0\quad\forall i=0,\dots,n.
\end{align}
We are looking for an element $\varphi_{n+1}$, such that $\sum_{j=0}^{n+1}d_j \varphi_{n-j}=0$. Rewriting the equation $\qdz^2=(\sum_{i=0}^\infty \nu^i d_i)^2$ order by order in the powers of $\nu$  we obtain $d_0d_j=-\sum_{i=0}^{j-1}d_{j-i}d_i$. Now, an easy calculation yields that $\sum_{i=0}^nd_{i+1}\varphi_{n-i}\in Z^1(\mathfrak g,\mathcal C^\infty(Z))$:
\begin{align*}
d_0\big(\sum_{i=0}^nd_{i+1}\varphi_{n-i}\big)=-\sum_{i=0}^n\sum_{k=0}^id_{i-k}d_k\varphi_{n-i}=-\sum_{r=0}^n d_{n-r}\sum_{s=0}^{n-r}d_s\varphi_{r-s}\stackrel{(\ref{IV})}{=}0.
\end{align*}
Setting $\varphi_{n+1}:=-d^{-1}(\sum_{i=0}^nd_{i+1}\varphi_{n-i})$, we are done. Obviously, $\varphi_1=-d^{-1}d_1\varphi_0=:S_1(\varphi_0)$, $\varphi_2=-d^{-1}(d_2\varphi_0-d^{-1}d_1\varphi_0)=:S_2(\varphi_0)$, etc. arise by successive continuous operations acting on $\varphi_0$ and we acquire the desired sequence $S=\sum_{i=0}^\infty\nu^i S_i$. 

Conversely, let $\varphi=\sum_i\nu^i \varphi_i\in \mathrm H^0(\mathfrak g,\mathcal C^\infty(Z)[[\nu]])$. Rewriting $\qdz\varphi=0$ order by order in powers of $\nu$ we get $\sum_{j=0}^n d_j\varphi_{n-j}=0$ for all $n\ge0$. In particular we have $d_0(\sum_{j=1}^nd_j\varphi_{n-j})=0$. Setting $\psi_n:=\varphi+d^{-1}(\sum_{j=1}^nd_j\varphi_{n-j})$, we obtain a series $\psi:=\sum_{i=0}^\infty\nu^i\psi_i\in \mathcal C^\infty(Z)^\mathfrak g[[\nu]]$. It is clear that $S\psi=\varphi$.
\end{BEWEIS}

\begin{appendix}

\chapter{Auxiliary material}
\section{Two perturbation lemmata} \label{HPT}
We consider (cochain) complexes in an additive $\mathbb K$-linear category $\mathscr C$ (e.g. the category of Fr\'echet spaces). A \emph{contraction} in $\mathscr C$ consists of the following data
\begin{align} \label{contr}
\xymatrix{(X,d_X)\ar@<-0.6ex>[r]_{i\:\:\:\:}^{}&(Y,d_Y),h_Y, \ar@<-0.6ex>[l]_{p\:\:\:\:}^{}}
\end{align}
where $i$ and $p$ are chain maps between the chain complexes $(X,d_X)$ and $(Y,d_Y)$, $h_Y:Y\to Y[-1]$ is a morphism, and we have
$pi=\id_X$, $d_Yh_Y+h_Yd_Y=\id_Y-ip$. The contraction is said to satisfy the \emph{side conditions} (sc1--3), if moreover, $h_Y^2=0$ , $h_Yi=0$ and $ph_Y=0$ are true. It was observed in \cite{LambeSt}, that in order to fulfill (sc2) and (sc3), one can replace $h_Y$ by $h'_Y:=(d_Yh_Y+h_Yd_Y) \:h_Y\:(d_Y h_Y+h_Y d_Y)$. If one wants to have in addition (sc1) to be satisfied, one may replace $h_Y'$ by $h_Y'':=h_Y'd_Yh'_Y$. 

Let $C:=Cone(p)$ be the mapping cone of $p$, i.e., $C=X[1]\oplus Y$ is the complex with differential $d_C(x,y):=(d_Xx+(-1)^{|y|}py,d_Y y)$. The homology of  $C$ is trivial, because $h_C(x,y):=(0,h_Yy+(-1)^{|x|}ix)$ is a contracting homotopy, i.e., $d_Ch_C+h_Cd_C=\id_C$, if (sc3) is true. In fact, we calculate
\begin{eqnarray*}
(d_Ch_C+h_Cd_C)(x,y)&=&d_C\big(0,h_Yy+(-1)^{|x|}ix, d_Y y\big)+h_C\big(d_Xx+(-1)^{|y|}py,d_Yy\big)\\
&=&\big((-1)^{|y|+1}ph_Yy+pix,d_Yh_Yy+(-1)^{|x|}d_Yix\big)\\
&&\qquad\qquad\qquad\qquad\qquad\qquad+\big(0,h_Yd_Yy+(-1)^{|x|+1}i\:d_Xx+ipy\big)\\
&\stackrel{(sc3)}{=}&\big(pix,(h_Yd_Y+d_Yh_Y)y+ipy\big)=(x,y).
\end{eqnarray*}

Let us now assume that the objects $X$ and $Y$ carry  descending 
Hausdorff filtrations and the structure maps are filtration preserving. Moreover, pretend that we have found a \emph{perturbation} $D_Y=d_Y+t_Y$ of $d_Y$, i.e., $D_Y^2=0$ and $t_Y:Y\to Y[1]$, called the \emph{initiator}, having the property that $t_Yh_Y+h_Yt_Y$ raises the filtration. Since, in general, $t_X:=pt_Yi$ does not need to be a perturbation of $d_X$, we impose that as an extra condition: we \emph{assume} that $D_X=d_X+t_X$ is a differential. Setting $t_C:=(t_X,t_Y)$, we will get a perturbation $D_C:=d_C+t_C$ of $d_C$, if we have  in addition $t_Xp=pt_Y$ (this will imply that $(d_X+t_X)^2=0$). Then an easy calculation yields that $H_C:=h_C(D_Ch_C+h_CD_C)^{-1}=h_C(\id_C+t_Ch_C+h_Ct_C)^{-1}$ is well defined and satisfies $D_CH_C+H_CD_C=\id_C$. In fact, we have
\begin{eqnarray*}
D_CH_C+H_CD_C&=&\big(D_Ch_C+H_CD_C(D_Ch_C+h_CD_C)\big)(D_Ch_C+h_CD_C)^{-1}\\
&=&\big(D_Ch_C+H_C(D_Ch_C+h_CD_C)D_C\big)(D_Ch_C+h_CD_C)^{-1}\\
&=&\big(D_Ch_C+h_CD_C\big)(D_Ch_C+h_CD_C)^{-1}=\id_C.
\end{eqnarray*} 
Defining the morphism $I:X\to Y$, $H_C(x,0)=:(0,(-1)^{|x|}Ix)$ and the homotopy $H_Y:Y\to Y[-1]$, $H_C(0,y)=:(0,H_Yy)$ we obtain the following 
\begin{LEMMA}[\emph{Perturbation Lemma -- Version 1}] \label{BPL1}If the contraction (\ref{contr}) satisfies (sc3) and  $D_Y=d_Y+t_Y$ is a perturbation of $d_Y$ such that $t_Xp=pt_Y$, then 
\begin{align} \label{contrEins}
  \xymatrix{(X,D_X)\ar@<-0.6ex>[r]_{I\:\:\:\:}^{}&(Y,D_Y),H_Y \ar@<-0.6ex>[l]_{p\:\:\:\:}^{}}
\end{align}
is a contraction fulfilling (sc3). Moreover, we have $H_Y= h_Y(\id_Y+t_Yh_Y+h_Yt_Y)^{-1}$ and $Ix=ix-H_Y(t_Yix-it_Xx)$. If all side conditions are true for (\ref{contr}), then they are for (\ref{contrEins}), too.
\end{LEMMA}
\begin{BEWEIS} For all homogenuous $x\in X$, $y\in Y$ we have
\begin{eqnarray*}
(x,y)&=&(D_CH_C+H_CD_C)(x,y)\\
&=&D_C\big(0,(-1)^{|x|}Ix+H_Yy\big)+H_C\big(D_Xx+(-1)^{|y|}py,D_Yy\big)\\
&=&\big(pIx+(-1)^{|y|+1}pH_Yy,(-1)^{|x|}D_YIx+D_YH_Yy\big)\\
&&\qquad\qquad\qquad\qquad+\big(0,(-1)^{|x|-1}ID_Xx+Ipy+H_YD_Yy\big),
\end{eqnarray*} 
and we conclude that $pI=\id_X$, $D_YI=ID_X$, $pH_Y=0$ and $D_YH_Y+H_YD_Y+Ip=\id_Y$. Let us verify the formula for $H_Y$:
$(0,H_Yy)=H_C(0,y)=h_C\sum_{i\ge 0}(-1)^i(t_Ch_C+h_Ct_C)^i(0,y)
=\big(0,h_Y\sum_{i\ge 0}(-1)^i(t_Yh_Y+h_Yt_Y)^iy\big).$
Note that $(t_Ch_C+h_Ct_C)(x,0)=(-1)^{|x|}(0,t_Y ix-it_Xx)$.
Using this result, it is straight forward to check:
$(0,Ix)=(0,ix)-\big(0,h_Y\sum_{i\ge0}(-1)^i(t_Yh_Y+h_Yt_Y)^i(t_Y ix-it_Xx)
=\big(0,ix-H_Y(t_Y ix-it_Xx)\big).$
Finally, let us address the side conditions. With the above formula for $H_Y$ the condition $i h_Y=0$ entails $iH_Y=0$. Note that $h_Y^2=0$ implies $h_Y(D_Yh_Y+h_YD_Y)^{-1}=(D_Yh_Y+h_YD_Y)^{-1}h_Y$, and we conlude that $H_Y^2=0$.
\end{BEWEIS}
\\
\\

Starting with the mapping cone $K=Cone(i)$, i.e., the complex $K=Y[1]\oplus X$ with the differential $d_K(y,x)=(d_Yy+(-1)^{|x|}ix,d_Xx)$, we may give a version of the above argument arriving at a contraction with all data perturbed except $i$. More precisely, we have a homotopy $h_K(y,x):=(h_Yy,(-1)^{|y|}py)$, for which $d_Kh_K+h_Kd_K=\id_K$ follows from (sc2). In fact, we calculate
\begin{eqnarray*}
(d_Kh_K+h_Kd_K)(y,x)&=&d_K\big(h_Y,(-1)^{|y|}py\big)+h_K\big(d_Yy+(-1)^{|x|}ix,d_X x\big)\\
&=&\big(d_Yh_Yy+ipy,(-1)^{|y|}d_Xpy)\big)\\
&&\qquad\qquad\qquad\qquad+\big(h_Yd_Y+(-1)^{|x|}h_Yix,(-1)^{|y|+1}pd_Yy+pix\big)\\
&\stackrel{(sc2)}{=}&\big(d_Yh_Yy+h_Yd_Y+ipy,pix\big)=(y,x).
\end{eqnarray*}
Mimicking the above argument for $C$, we get a differential $D_K:=d_K+t_K$ with $t_K:=(t_Y,t_X)$, if $t_Yi=it_X$ (this will imply $D_X^2=0$). Assuming (\ref{contr}) to satisfy (sc2), $H_K:=h_K(D_Kh_K+h_KD_K)^{-1}$ will become a contracting homotopy $D_KH_K+H_KD_K=\id_K$. Defining $P:Y\to X$ and $H'_Y:Y\to Y[-1]$ by  $H_K(y,0)=H_K(y,x)=:(H'_Yy,(-1)^{|y|}Py)$ we get the following 
\begin{LEMMA} [\emph{Perturbation Lemma -- Version 2}] \label{BPL2}If the contraction (\ref{contr}) satisfies (sc2) and $D_Y=d_Y+t_Y$ is a perturbation of $d_Y$ such that $t_Yi=it_X$, then 
\begin{align} \label{contrZwei}
  \xymatrix{(X,D_X)\ar@<-0.6ex>[r]_{i\:\:\:\:}^{}&(Y,D_Y),H'_Y \ar@<-0.6ex>[l]_{P\:\:\:\:}^{}}
\end{align}
is a contraction fulfilling (sc2). Moreover, we have $H'_Y= h_Y(\id_Y+t_Yh_Y+h_Yt_Y)^{-1}$ and $P=p(\id+t_Yh_Y+h_Yt_Y)^{-1}$. If all side conditions are true for (\ref{contr}), then they are for (\ref{contrZwei}), too.
\end{LEMMA}
\begin{BEWEIS} For all homogenuous $x\in X$, $y\in Y$ we have
\begin{eqnarray*}
(y,x)&=&(D_KH_K+H_KD_K)(y,x)\\
&=&D_K\big(H'_Yy,(-1)^{|y|}Py\big)+H_K\big(D_Y+(-1)^{|x|}ix,D_X x\big)\\
&=&\big(D_Y H'_Yy+iPy,(-1)^{|y|}D_XPy\big)\\
&&\qquad\qquad\qquad\qquad+\big(H'_YD_Yy+(-1)^{|x|}H'_Yix,(-1)^{|y|+1}PD_Yy+Pix\big),
\end{eqnarray*}
and we conclude that $D_YH'_Y+H'_YD_Y+iP=\id_Y$, $H'_Yi=0$, $D_XP=PD_Y$ and $Pi=\id_X$. Let us check the formulas for $H_Y'$ and $P$:
\begin{eqnarray*}
(H_Y'y,(-1)^{|y|}Py)&=&H_K(y,0)=\big(h_K(\sum_{i\ge0}(-1)^i(t_Yh_Y+h_Yt_Y)^iy,0\big)\\
&=&\big(h_Y\sum_{i\ge0}(-1)^i(t_Yh_Y+h_Yt_Y)^iy,(-1)^{|y|}p\sum_{i\ge0}(-1)^i(t_Yh_Y+h_Yt_Y)^iy\big). 
\end{eqnarray*}
Finally, let us address the side conditions. If $h_Y^2=0$ then $h_Y(1+t_Yh_Y+h_yt_Y)^{-1}=(1+t_Yh_Y+h_yt_Y)^{-1}h_Y$, which entails ${H'}^2_Y=0$. If in addition $ph_Y=0$, then we conclude that $PH'_Y=0$.
\end{BEWEIS}

\section{Graded Lie-Rinehart pairs}\label{LieRinehart}

The notion of a Lie-Rinehart pair is the algebraic counterpart of the notion of a Lie algebroid. It admits a more or less obvious translation to the graded situation. $\mathbbm Z_2$-graded Lie-Rinehart have been studied for example in \cite{Chemla}. In the following graded objects, morphisms etc. are  understood to be in one of the abelian tensor categories $\Zzvect$ and $\Zvect$, the category of $\mathbbm Z_2$-graded and $\mathbb Z$-graded vector spaces, respectively. Rinehart \cite{Rinehart} introduced the notion of an universal enveloping algebra of a Lie-Rinehart pair generalizing the universal enveloping algebra of a Lie algebra. In \cite{HuebPoiss} H\"ubschmann gave an alternative construction for it, which we will translate to the graded situation. 
\begin{DEFINITION} A  \emph{graded Lie Rinehart pair} $(A,L)$ is a graded commutative $\mathbbm K$-algebra $A$ and a graded $\mathbbm K$-Lie algebra $L$ such that
\begin{enumerate} 
\item $L$ is a graded left $A$-module $A\otimes L\to L$, $(a,X)\mapsto aX$,
\item $L$ acts on $A$ by graded left derivations $L\otimes A\to A$, $(X,a)\mapsto X(a)$,
\item the actions are compatible in the following sense 
\begin{eqnarray*}
aX(b)&=&(aX)(b),\\
{[X,aY]}&=&X(a)Y+(-1)^{|a||X|}a[X,Y] \mbox{ for all homogeneous } X,Y\in L\mbox{ and }a,b\in A.
\end{eqnarray*}
\end{enumerate}
Sometimes $L$ will also be called a $(\mathbbm K,A)$-\emph{Lie algebra}. An $(A,L)$-\emph{module} is a graded $\mathbbm K$-vector space $V$ which is at the same time an $A$-and an $L$-module such that the actions are compatible in the following sense
 \begin{eqnarray*}
aX(v)&=&(aX)(v),\\
X(av)&=&X(a)v+(-1)^{|a||X|}aX(v) \mbox{ for all homogeneous }X\in L,\:a\in A\mbox{ and }v\in V.
\end{eqnarray*}
\end{DEFINITION}
\paragraph{Examples}
\begin{enumerate}
\item If $A$ is a graded commutative $\mathbbm K$-algebra then $(A,\mathrm{Der}_{\mathbbm K}A)$ is a graded Lie-Rinehart pair. 
\item If $A$ is a graded commutative algebra, $\mathfrak g$ a graded $\mathbbm K$-Lie algebra and $\rho:\mathfrak g\to \mathrm{Der}_\mathbbm K A$ a morphism of graded Lie algebras, then $A\otimes_\mathbbm K \mathfrak g$ is a graded Lie-Rinehart pair with bracket given by \[[a\otimes X,b\otimes Y]:=(-1)^{|b||X|}ab\otimes [X,Y]+a\rho(X)b \otimes Y+(-1)^{(|a|+|b|+|X|)|Y|} (\rho(Y)a)b\otimes X. \]
\item The above example can be generalized to the notion of a \emph{Lie algebroid}, which we will use merely in the ordinary (even) manifold setup. By a Lie algebroid we mean the data $(E,[\:,\:],\rho)$, where $E\to M$ is a vector bundle over a manifold $M$, $[\:,\:]$ is a $\mathbbm K$-linear Lie bracket on the space of sections $\Gamma^\infty(M,E)$ of $E$ and $\rho:E\to TM$ is a vector bundle homomorphism such that $[a,fb]=f[a,b]+\rho(a)(f)\:b$ for all $a,b\in \Gamma^\infty(M,E)$ and $f\in\mathcal C^\infty(M)$. It is clear that $(\mathcal C^\infty(M),\Gamma^\infty(M,E))$ is a Lie-Rinehart pair.
\end{enumerate}
 
Associated to the graded Lie-Rinehart pair $(A,L)$ there is an \emph{universal enveloping algebra} $(U(A,L),i_A,i_L)$, which is analog to the enveloping algebra of a Lie algebra and of the algebra of differential operators on a manifold, respectively. More precisely $U(A,L)$ is a graded $\mathbbm K$-algebra, $i_A:A \to U(A;L)$ is a morphism of graded $\mathbbm K$-algebras and $i_L:L\to U(A,L)$ is a morphism of graded Lie algebras having the following properties: $i_A(a)i_L(X)=i_L(aX)$ and $i_L(X)i_A(a)-(-1)^{|X||a|}i_A(a)i_L(X)=i_A(X(a))$ for all homogeneous $a\in A$ and $X\in L$.  $(U(A,L),i_A,i_L)$ is initial among the triples $(U,j_A,j_L)$ having these properties.

In order to construct it, we follow the approach of H\"ubschmann \cite{HuebPoiss}. Let $U(\mathbbm K,L)$ be the universal enveloping algebra of the graded Lie algebra $L$, $i_{\mathbbm K}:\mathbbm K\to U(\mathbbm K,L)$ and $i_L:L\to U(\mathbbm K,L)$ the canonical embeddings and $\Delta:U(\mathbbm K,L)\to U(\mathbbm K,L)\otimes_{\mathbbm K} U(\mathbbm K,L)$ the standard comultiplication. Remember that $\Delta$ is the unique comultiplication such that $U(\mathbbm K,L)$ is a graded bialgebra and the image of $i_L$ is the space of primitives. As an intermediate step, let us define the algebra 
\[A\odot U(\mathbbm K,L)=(A\otimes_{\mathbbm K} U(\mathbbm K,L),\mu),\] where the multiplication $\mu$ is given by
\[(a\otimes u)(b\otimes v):=(-1)^{|u||b|}ab\otimes uv+\sum (-1)^{|u''||b|}au'(b)\otimes u''v\] 
for homogeneous $a,b \in A$ and $u,v \in U(\mathbbm K,L)$. Here we used Sweedler's notation $\Delta(u):=\sum u'\otimes u''$. One easily proves that $\mu$ is associative. Let $J$ be the right ideal generated by elements of the form $ab\otimes X-a\otimes b X$ for $a,b \in A$ and $X\in L$ (here we write $X$ and $bX$ for the respective images under $i_L$). As a consequence of the identity
\begin{eqnarray*}(c\otimes Y)(ab\otimes X-a\otimes b X)&=&(-1)^{|Y|(|a|+|b|)}(cab\otimes XY-c\otimes ab YX)\\
&&\qquad+(-1)^{|Y||a|}(caY(b)\otimes X-c\otimes aY(b)X)
\end{eqnarray*}
for all homogeneous $a,b,c\in A$ and $X,Y\in L$, we have that $J$ is in fact a two sided ideal. The universal enveloping algebra is defined to be the quotient
\[U(A,L):=A\odot U(\mathbbm K, L)/J,\]
and $i_A$, $i_L$ are  the obvious morphisms.

The universal enveloping algebra $U(A,L)$ is in a natural manner a filtered algebra \[U(A,L)=U_0(A,L)\supset U_1(A,L)\supset\dots  U_{n-1}(A,L)\supset U_n(A,L)\supset\dots\quad.\]
Here $U_n(A,L)$ is the left $A$-module generated by at most $n$ products of the images of $L$ in $U(A,L)$. Clearly the induced left and right $A$-algebra structures on the associated graded algebra $\mathrm{gr}\,U(A,L)=\oplus_{n\ge 0}\mathrm{gr}_n \,U(A,L)=\oplus_{n\ge 0} U_n(A,L)/U_{n-1}(A,L)$ coincide (here we set $U_{-1}(A,L):=0$). Note that $\mathrm{gr}\,U(A,L)$ is a graded commutative algebra.

The best understood Lie-Rinehart pairs are those, for which the $A$-module $L$ is \emph{projective}. One reason for this is, that there is an analog of the Poincar\'e-Birkhoff-Witt theorem.  

\begin{SATZ}[\emph{Poincar\'e-Birkhoff-Witt theorem}] The canonical $A$-module epimorphism $S_AL\to\mathrm{gr}\,U(A,L)$ is an isomorphism of graded commutative algebras.
\end{SATZ} 
\begin{BEWEIS} Adapt the proof of Rinehart \cite[p.199--200]{Rinehart} to the graded situation.
\end{BEWEIS}
\\

A \emph{morphism} of graded Lie-Rinehart pairs $(A,L)\to (A',L')$ is a morphism of graded commutative algebras $A\to A'$, $a\mapsto a'$ and a morphism of graded Lie algebras $L\to L'$, $X\mapsto X'$ such that $(aX)'=a'X'$ and $(X(a))'=X'(a')$ for all $a\in A$ and $X\in L$. A morphism of Lie Rinehart pairs extends uniquely to a ring homorphism  of the universal envelopping algebras $U(A,L)\to U(A',L')$. There is also a natural notion of a \emph{module} for a Lie-Rinehart pair $(A,L)$ (for details see \cite{Rinehart}). Equivalently, one may think of such a module as a module for the algebra $U(A,L)$. For example, the action of $L$ on $A$ extends natually to an $U(A,L)$-module structure on $A$.

It has been shown by Rinehart, that if $L$ is a projective $A$-module, there is a Koszul resolution of the $U(A,L)$-module $A$, which generalizes the Koszul resolution of the ground field for Lie algebras. The space of chains of this complex is $K_\bullet(A,L)=U(A,L)\otimes_A \wedge ^\bullet_A L$. The differential  $\partial:K_\bullet(A,L)\to K_{\bullet-1}(A,L)$ is given by the formula
\begin{eqnarray*}
\partial(u\otimes X_1\wedge\cdots\wedge X_n)&=&\sum_{i=1}^n (-1)^{i+1}\sign(\tau_i,|X|)\:(u\: i_L(X_i))\otimes X_1\wedge \cdots \widehat{X_i}\cdots\wedge X_n\\
&+&\sum_{1\le i<j\le n} (-1)^{i+j}\sign(\tau_{i,j},|X|)\:u\otimes [X_i,X_j]\wedge X_1\wedge \cdots \widehat{X_i}\cdots\widehat{X_j}\cdots\wedge X_n.
\end{eqnarray*}
Here, $\sign(\tau_i,|X|)$ and $\sign(\tau_{i,j},|X|)$ are the Koszul signs of the permutations $\tau_i=(1\,2\dots i-1\,i)$ and 
\begin{align}
\tau_{i,j}=\left(
\begin{array}{cccccccccccc}
1&2&3&4&\dots&   &   &      &   &   &     &n    \\
i&j&1&2&\dots&i-1&i+1&\dots &j-1&j+1&\dots&n
\end{array}
\right),
\end{align}
respectively, for the multiindex $|X|=(|X_1|,\dots,|X_n|)$. Analog to \cite[Section 4] {Rinehart} one can show that $\partial$ is in fact well-defined and that $(K_\bullet(A,L),\partial)$ is a projective resolution of the $U(A,L)$-module $A$. We can therefore use this resolution to compute derived functors such as $\Ext_{U(A,L)}(A,M)$ for any $U(A,L)$-module $M$. More precisely, $\Ext_{U(A,L)}(A,M)$ is the cohomology of the complex
\begin{align}\label{LieRinehartComplex}
\Hom_A(\wedge_A^\bullet L,M)=\Alt^{\bullet}_A(L,M),
\end{align}
where the differential $d:=\Alt^n(L,M)\to \Alt^{n+1}(L,M)$ is defined for homogeneous $f\in \Alt^n(L,M)$ by the formula
\begin{eqnarray*}
(df)(X_1,\dots,X_{n+1})&=&\sum_{i=1}^{n+1} (-1)^{i+1+|f||X_i|}\sign(\tau_i,|X|)\:X_i(f(X_1,\dots\widehat{X_i}\dots,X_{n+1}))\\
&&+\sum_{1\le i<j\le n+1}(-1)^{i+j}\sign(\tau_{i,j},|X|)\:f([X_i,X_j],X_1,\dots\widehat{X_i}\dots\widehat{X_j}\dots,X_{n+1}).
\end{eqnarray*}
In the case of a Lie algebroid $\Ext_{U(A,L)}(A,A)$ is just the ordinary Lie algebroid cohomology. In particular, for $A=\mathcal C^\infty(M)$ the ring of functions on a smooth manifold $M$ and $L=\Gamma^\infty(M,TM)$ this boils down to the ordinary de Rham cohomology. Also, in the $\mathbbm Z_2$-graded setting this construction reproduces the Chevalley-Eilenberg cochain complex of Lie super-algebras (see e.g. \cite{Scheunert} and references therein) and the de Rham complex for super-manifolds. In the general $\mathbbm Z$-graded case the complex (\ref{LieRinehartComplex}) may become rather huge (note that we have to use $\Hom_A$ and not $\Hom_{\Zmod{A}}$), and there are perhaps better ways to define Lie-Rinehart cohomology for $\mathbbm Z$-graded Lie-Rinehart pairs. 
\section{The opposite of a $n$-Poisson algebra is $n$-Poisson}\label{oppPoiss}
\begin{LEMMA}
If $(L,\cdot,[\:,\:])$ is a $n$-Poisson algebra then the opposite bracket 
\begin{align}
[a,b]^{\oppind}:=(-1)^{(|a|+n)(|b|+n)}[b,a]
\end{align}
for homogeneous $a,b \in L$, is a Poisson bracket of degree $n$.
We say that $(L,\cdot,[\:,\:]^{\oppind})$ is the  \emph{opposite} of the  $n$-Poisson algebra $(L,\cdot,[\:,\:])$.
\end{LEMMA}
\begin{BEWEIS} In order to prove the graded Leibniz rule note that
\begin{eqnarray*}
[a,bc]^{\oppind}&=&(-1)^{(|b|+|c|+n)(|a|+n)}[bc,a]\\
&=&(-1)^{|a||b|+|a||c|+n(|a|+|b|+|c|+1)}\left(b[c,a]+(-1)^{(|a|+n)}[b,a]c\right)\\
&=&(-1)^{|a||b|+|a||c|+n(|a|+|b|+|c|+1)}\Big(\underbrace{(-1)^{(|a|+n)(|c|+n)}}_{=(-1)^{|a||c|+n(a|+|c|+1)}} b[a,c]^{\oppind}+\\
&&\qquad\qquad\qquad\qquad\qquad\qquad\qquad\underbrace{(-1)^{(|a|+n)|c|}(-1)^{(|a|+n)(|b|+n)}}_{=(-1)^{|a||b|+|a||c|+n(|a|+|b|+|c|+1)}}[a,b]^{\oppind}c\Big)\\
&=&(-1)^{(|a|+n)|b|}b[a,c]^{\oppind}+[a,b]^{\oppind}c.
\end{eqnarray*}
Similarly we compute:
\begin{eqnarray*}
[a,[b,c]^{\oppind}]^{\oppind}&=&(-1)^{(|b|+|c|)(|a|+n)}(-1)^{(|c|+n)(|b|+n)}[[c,b],a]\\
&=& (-1)^{|a||b|+|b||c|+|a||c|+n}\left([c,[b,a]]+(-1)^{(|a|+n)(|b|+n)}[[c,a],b]\right)\\
&=& (-1)^{|a||b|+|b||c|+|a||c|+n}\Big((-1)^{(|a|+|b|)(|c|+n)}(-1)^{(|a|+n)(|b|+n)}[[a,b]^{\oppind},c]^{\oppind}\\
&&+(-1)^{|a||b|+n(|a|+|b|+1)}(-1)^{(|a|+|c|)(|b|+n)}(-1)^{(|a|+n)(|c|+n)}[b,[a,c]^{\oppind}]^{\oppind}\Big)\\
&=& (-1)^{|a||b|+|b||c|+|a||c|+n}\Big((-1)^{|a||c|+|b||c|+n(|a|+|b|)}(-1)^{|a||b|+n(|a+|b|+1)}[[a,b]^{\oppind},c]^{\oppind}\\
&&+(-1)^{|a||b|+n(|a|+|b|+1)}(-1)^{|a||b|+|b||c|+|a||c|+n}[b,[a,c]^{\oppind}]^{\oppind}\Big)\\
&=&[[a,b]^{\oppind},c]^{\oppind}+(-1)^{(|a|+n)(|b|+n)}[b,[a,c]^{\oppind}]^{\oppind}.
\end{eqnarray*}
\end{BEWEIS}

\section{Proof of Theorem \ref{SchoutenNijenhuis}}\label{Schoutenproof}
First of all let us introduce some notation. Let $X=X_1X_2\dots X_r\in S_A(L[-1])$ be a monomial, such that $X_i\in L^{x_i}$ for $i=1,\dots,r$. We will call such monomials \emph{homogeneous}. The degree $|X|$ of $X$ in  $S_A(L[-1])$ is given by $x_1+\dots +x_r + r$. Moreover, we will use the shorthand 
\begin{eqnarray*}
X_{>i}&:=&X_{i+1}X_{i+2}\dots X_r,\qquad  |X_{>i}|=x_{i+1}+\dots+ x_r+r-i\\
X^{<i}&:=&X_1X_2\dots X_{i-1},\qquad |X^{<i}|=x_1+\dots+x_{i-1}+i-1.
\end{eqnarray*}
It is not difficult to prove that, as a consequence of the Leibniz rule, the Schouten-Nijenhuis bracket with another such monomial $Y=Y_1Y_2\dots Y_s$ is given by:
\begin{eqnarray}\label{Schoutenformel}
[X,Y]=\sum_{i=1}^r\sum_{j=1}^s (-1)^{|X_i||X_{>i}|+|Y_j||Y^{<j}|}X^{<i}X_{>i}[X_i,Y_j]Y^{<j}Y_{>j}
\end{eqnarray}
Moreover, it is easy to prove that the result of this formula does not depend of the choice of the decomposition.
Note that the sign in this formula appears naturally as the Koszul sign of the permutation which moves $X_i$ to the right of $X$ and $Y_j$ to the left of $Y$.
Conversely, the Leibniz rule follows effortlessly from equation (\ref{Schoutenformel}). In fact, if $Z=Z_1Z_2\dots Z_t$ is a third homogeneous monomial we have
\begin{eqnarray*}
[X,YZ]&=&[X,Y]Z+\sum_{i=1}^r\sum_{j=1}^t (-1)^{|X_i||X_{>i}|+|Z_j|(|Y|+|Z^{<j}|)}X^{<i}X_{>i}[X_i,Z_j]YZ^{<j}Z_{>j}\\
&=&[X,Y]Z+(-1)^{|Y|(|X|-1)}Y[X,Z].
\end{eqnarray*}
In order to see that this bracket is graded antisymmetric note that $|X^{<i}X_{>i}|=|X|-|X_i|$ and that the degree of $[X_i,Y_j]=-(-1)^{(|X_i|-1)(|Y_j|-1)}[Y_j,X_i]$ is $|X_i|+|Y_j|-1$. Moreover, the Koszul sign of the operation which replaces $X^{<i}X_{>i}[X_i,Y_j]Y^{<j}Y_{>j}$ by $Y^{<j}Y_{>j}[X_i,Y_j]X^{<i}X_{>i}$ is given by
\begin{eqnarray*}
&&(-1)^{(|X|-|X_i|)(|X_i|+|Y_j|-1)}(-1)^{(|X|-|X_i|)(|Y|-|Y_j||)}(-1)^{(|X_i|+|Y_j|-1)(|Y|-|Y_j|)}\\
&=&(-1)^{|X||X_i|+|X||Y_j|+|X_i||Y_j|+|X|}(-1)^{|X||Y|+|X||Y_j|+|X_i||Y|+|X_i||Y_j|}(-1)^{|X_i||Y|+|Y||Y_j|+|X_i||Y_j|+|Y|}\\
&=&(-1)^{|X||Y|+|X||X_i|+|Y||Y_j|+|X_i||Y_j|+|X|+|Y|}.
\end{eqnarray*}
It follows that
\begin{eqnarray*}
[X,Y]&=&-\sum_{i=1}^r\sum_{j=1}^s (-1)^{|X_i||X_{>i}|+|Y_j||Y^{<j}|}(-1)^{|X||Y|+|X||X_i|+|Y||Y_j|+|X_i||Y_j|+|X|+|Y|}\\
&&\qquad \qquad \qquad \qquad \qquad \qquad (-1)^{(|X_i|-1)(|Y_j|-1)}\:Y^{<j}Y_{>j}[Y_j,X_i]X^{<i}X_{>i}
\end{eqnarray*}
and we have to analyse the sign in the above formula. To this end note that $|X_{>i}|=|X|-|X^{<i}|-|X_i|$ and $|Y^{<j}|=|Y|-|Y_{>j}|-|Y_j|$. Some bookkeeping yields that the sign works out correctly
\begin{eqnarray*}
(-1)^{|Y_j||Y_{>j}|+|X_i||X^{<i}|}(-1)^{|X||Y|+|X|+|Y|+1},
\end{eqnarray*}
and we have thus proved $[X,Y]=-(-1)^{(|X|-1)(|Y|-1)}[Y,X]$.

In order to prove the Jacobi identity we introduce some further notation. For $1\le i<j\le r$ and a monomial $X=X_1X_2\dots X_r$ we introduce 
\begin{align*}
X_{>i}^{<j}:=X_{i+1}\dots X_{j-1}\qquad |X_{>i}^{<j}|=x_{i+1}+\dots+ x_{j-1}+j-i-1.
\end{align*}
Let $Y=Y_1Y_2\dots Y_s\in S^s_A(L[-1])$ and $Z=Z_1Z_2\dots Z_t\in S^t_A(L[-1])$ be two other homogeneous monomials. In order to proof 
\begin{eqnarray}\label{Jacobiidentitaet}
[X,[Y,Z]]=[[X,Y],Z]+(-1)^{(|X|-1)(|Y|-1)}[Y,[X,Z]]
\end{eqnarray}
we expand the three terms in the above equation according to equation (\ref{Schoutenformel}), identify the various contributions and `mow down' the signs. 
\begin{eqnarray*}
[X,[Y,Z]]&=&\sum_{i=1}^r\sum_{k=1}^s\sum_{j=1}^{k-1}\sum_{l=1}^t(-1)^{n_1} X^{<i}X_{>i}[X_i,Y_j]Y^{<j}Y_{>j}^{<k}Y_{>k}[Y_k,Z_l]Z^{<l}Z_{>l}\\
&&+\sum_{i=1}^r\sum_{k=1}^s\sum_{j=k+1}^{s}\sum_{l=1}^t(-1)^{n_2} X^{<i}X_{>i}[X_i,Y_j]Y^{<k}Y_{>k}^{<j}Y_{>j}[Y_k,Z_l]Z^{<l}Z_{>l}\\
&&+\sum_{i=1}^r\sum_{j=1}^s\sum_{l=1}^t (-1)^{n_3}X^{<i}X_{>i}[X_i,[Y_j,Z_l]]Y^{<j}Y_{>j} Z^{<l}Z_{>l}\\
&&+\sum_{i=1}^r \sum_{k=1}^s\sum_{l=1}^t\sum_{j=1}^{l-1}(-1)^{n_4}X^{<i}X_{>i}[X_i,Z_j]Y^{<k}Y_{>k}[Y_k,Z_l]Z^{<j}Z_{>j}^{<l}Z_{>l}\\
&&+\sum_{i=1}^r\sum_{k=1}^s\sum_{l=1}^t\sum_{j=l+1}^{t}(-1)^{n_5}X^{<i}X_{>i}[X_i,Z_j]Y^{<k}Y_{>k}[Y_k,Z_l]Z^{<l}Z_{>l}^{<j}Z_{>j}\\
&=:&A_1+A_2+A_3+A_4+A_5.
\end{eqnarray*}
The exponents $n_1,\dots,n_5$ in the above equation will be understood modulo $2$. They are given as follows
\begin{eqnarray*}
n_1&=&|X_i||X_{>i}|+|Y_j||Y^{<j}|+|Y_k||Y_{>k}|+|Z_l||Z^{<l}|\\
n_2&=&|X_i||X_{>i}|+|Y_j|(|Y^{<j}|-|Y_k|)+|Y_k||Y_{>k}|+|Z_l||Z^{<l}|\\
n_3&=&|X_i||X_{>i}|+(|Y|-|Y_j|)(|Y_j|+|Z_l|-1)+|Y_j||Y_{>j}|+|Z_l||Z^{<l}|\\
n_4&=&|X_i||X_{>i}|+|Z_j|(|Z^{<j}|+|Y_k|+|Z_l|-1+|Y|-|Y_k|) +|Y_k||Y_{>k}|+|Z_l||Z^{<l}|\\
&=&|X_i||X_{>i}|+|Z_j|(|Z^{<j}|+|Z_l|+|Y|+1)+|Y_k||Y_{>k}|+|Z_l||Z^{<l}|\\
n_5&=&|X_i||X_{>i}|+|Z_j|(|Z^{<j}|-|Z_l|+|Z_l|+|Y|-1)+|Y_k||Y_{>k}|+|Z_l||Z^{<l}|\\
&=&|X_i||X_{>i}|+|Z_j|(|Z^{<j}|+|Y|+1)+|Y_k||Y_{>k}|+|Z_l||Z^{<l}|.
\end{eqnarray*}
Next, we have
\begin{eqnarray*}
[[X,Y],Z]&=&\sum_{i=1}^r\sum_{k=1}^s\sum_{j=1}^{k-1}\sum_{l=1}^t(-1)^{m_1} X^{<i}X_{>i}[X_i,Y_j]Y^{<j}Y_{>j}^{<k}Y_{>k}[Y_k,Z_l]Z^{<l}Z_{>l}\\
&&+\sum_{i=1}^r\sum_{k=1}^s\sum_{j=k+1}^{s}\sum_{l=1}^t(-1)^{m_2} X^{<i}X_{>i}[X_i,Y_j]Y^{<k}Y_{>k}^{<j}Y_{>j}[Y_k,Z_l]Z^{<l}Z_{>l}\\
&&+\sum_{i=1}^r\sum_{j=1}^s\sum_{l=1}^t (-1)^{m_3}X^{<i}X_{>i}Y^{<j}Y_{>j}[[X_i,Y_j],Z_l] Z^{<l}Z_{>l}\\
&&+\sum_{i=1}^r \sum_{j=1}^s\sum_{l=1}^t\sum_{k=i+1}^{r}(-1)^{m_4}X^{<i}X_{>i}^{<k}X_{>k}[X_i,Y_j]Y^{<j}Y_{>j}[X_k,Z_l] Z^{<l}Z_{>l}\\
&&+\sum_{i=1}^r\sum_{j=1}^s\sum_{l=1}^t\sum_{k=1}^{i-1}(-1)^{m_5}X^{<k}X^{<i}_{>k} X_{>i}[X_i,Y_j]Y^{<j}Y_{>j}[X_k,Z_l]Z^{<l}Z_{>l}\\
&=:&B_1+B_2+B_3+B_4+B_5.
\end{eqnarray*}
The exponents $m_1,\dots,m_5$ in the above equation are given as follows
\begin{eqnarray*}
m_1&=&|X_i||X_{>i}|+|Y_j||Y^{<j}|+|Y_k||Y_{>k}|+|Z_l||Z^{<l}|\\
m_2&=&|X_i||X_{>i}|+|Y_j||Y^{<j}|+|Y_k|(|Y_{>k}|-|Y_j|)+|Z_l||Z^{<l}|\\
m_3&=&|X_i||X_{>i}|+|Y_j||Y^{<j}|+(|Y|-|Y_j|)(|X_i|+|Y_j|-1)+|Z_l||Z^{<l}|\\
m_4&=&|X_i||X_{>i}|+|Y_j||Y^{<j}|+|X_k|(|X_{>k}|+|Y|-|Y_j|+(|X_i|+|Y_j|-1)) +|Z_l||Z^{<l}|\\
&=&|X_i||X_{>i}|+|Y_j||Y^{<j}|+|X_k|(|X_{>k}| +|Y|+|X_i|+1)+|Z_l||Z^{<l}|\\
m_5&=&|X_i||X_{>i}|+|Y_j||Y^{<j}|+|X_k|(|X_{>k}|-|X_i|+|Y|-|Y_j|+(|X_i|+|Y_j|-1)+|Z_l||Z^{<l}|\\
&=&|X_i||X_{>i}|+|Y_j||Y^{<j}|+|X_k|(|X_{>k}|+|Y|+1)+|Z_l||Z^{<l}|.
\end{eqnarray*}
Finally, we have
\begin{eqnarray*}
[Y,[X,Z]]&=&\sum_{j=1}^s\sum_{k=1}^r\sum_{i=1}^{k-1}\sum_{l=1}^t(-1)^{k_1} Y^{<j}Y_{>j}[Y_j,X_i]X^{<i}X_{>i}^{<k}X_{>k}[X_k,Z_l]Z^{<l}Z_{>l}\\
&&+\sum_{j=1}^s\sum_{k=1}^r\sum_{i=k+1}^{s}\sum_{l=1}^t(-1)^{k_2} Y^{<j}Y_{>j}[Y_j,X_i]X^{<k}X_{>k}^{<i}X_{>i}[X_k,Z_l]Z^{<l}Z_{>l}\\
&&+\sum_{j=1}^s\sum_{i=1}^r\sum_{l=1}^t (-1)^{k_3}Y^{<j}Y_{>j}[Y_j,[X_i,Z_l]]X^{<i}X_{>i} Z^{<l}Z_{>l}\\
&&+\sum_{k=1}^s \sum_{i=1}^r\sum_{j=1}^t\sum_{l=1}^{j-1}(-1)^{k_4}Y^{<k}Y_{>k}[Y_k,Z_l]X^{<i}X_{>i}[X_i,Z_j]Z^{<l}Z_{>l}^{<j}Z_{>j}\\
&&+\sum_{i=1}^s\sum_{k=1}^r\sum_{j=1}^t\sum_{l=j+1}^{t}(-1)^{k_5}Y^{<k}Y_{>k}[Y_k,Z_l]X^{<i}X_{>i}[X_i,Z_j]Z^{<j}Z_{>j}^{<l}Z_{>l}\\
&=:&C_1+C_2+C_3+C_4+C_5.
\end{eqnarray*}
The exponents $k_1,\dots,k_5$ in the above equation are given as follows
\begin{eqnarray*}
k_1&=&|Y_j||Y_{>j}|+|X_i||X^{<i}|+|X_k||X_{>k}|+|Z_l||Z^{<l}|\\
k_2&=&|Y_j||Y_{>j}|+|X_i|(|X^{<i}|-|X_k|)+|X_k||X_{>k}|+|Z_l||Z^{<l}|\\
k_3&=&|Y_j||Y_{>j}|+(|X|-|X_i|)(|X_i|+|Z_l|-1)+|X_i||X_{>i}|+|Z_l||Z^{<l}|\\
k_4&=&|Y_k||Y_{>k}|+|Z_l|(|Z^{<l}|+|X_i|+|Z_j|-1+|X|-|X_i|) +|X_i||X_{>i}|+|Z_j||Z^{<j}|\\
&=&|Y_k||Y_{>k}|+|Z_l|(|Z^{<l}|+|Z_j|+|X|+1)+|X_i||X_{>i}|+|Z_j||Z^{<j}|\\
k_5&=&|Y_k||Y_{>k}|+|Z_l|(|Z^{<l}|-|Z_j|+|Z_j|+|X_i|-1+|X|-|X_i|)+|X_i||X_{>i}|+|Z_j||Z^{<j}|\\
&=&|Y_k||Y_{>k}|+|Z_l|(|Z^{<l}|+|X|+1)+|X_i||X_{>i}|+|Z_j||Z^{<j}|.
\end{eqnarray*}
First of all we note that 
\begin{eqnarray}\label{A12id}
A_1=B_1, \qquad A_2=B_2
\end{eqnarray}
and we have to check some similar relations among the $A,B,C$'s. The only one which involves the Jacobi identity for $L$ is
\begin{align}\label{Jac}
A_3=B_3+(-1)^{(|X|-1)(|Y|-1)}C_3,
\end{align}
which we would like to check right now. Let us compare $A_3$ and $B_3$:
\begin{eqnarray*}
A_3&=&=\sum_{i=1}^r\sum_{j=1}^s\sum_{l=1}^t(-1)^{n_3} X^{<i}X_{>i}[X_i,[Y_j,Z_l]]Y^{<j}Y_{>j} Z^{<l}Z_{>l}\\
B_3&=&\sum_{i=1}^r\sum_{j=1}^s\sum_{l=1}^t (-1)^{m_3}X^{<i}X_{>i}Y^{<j}Y_{>j}[[X_i,Y_j],Z_l] Z^{<l}Z_{>l}\\
&=&\sum_{i=1}^r\sum_{j=1}^s\sum_{l=1}^t (-1)^{m_3+l_1}X^{<i}X_{>i}[[X_i,Y_j],Z_l]Y^{<j}Y_{>j}Z^{<l}Z_{>l},
\end{eqnarray*}
where $l_1$ is the Koszul sign of the permutation involved. Hence we need to check whether 
\begin{eqnarray*}
n_3&=&|X_i||X_{>i}|+(|Y|-|Y_j|)(|Y_j|+|Z_l|-1)+|Y_j||Y_{>j}|+|Z_l||Z^{<l}|\\
&=&|X_i||X_{>i}|+|Y||Y_j|+|Y||Z_l|+|Y|+|Y_j||Z_l|+|Y_j||Y_{>j}|+|Z_l||Z^{<l}|
\end{eqnarray*} 
coincides with $l_1+m_3$ 
\begin{eqnarray*}
l_1&=&(|Y|-|Y_j|)(|X_i|+|Y_j|+|Z_l|-2)\\
&=&|X_i||Y|+|Y||Z_l|+|X_i||Y_j|+|Y_j||Z_l|+|Y_j|\\
m_3&=&|X_i||X_{>i}|+|Y_j||Y^{<j}|+(|Y|-|Y_j|)(|X_i|+|Y_j|-1)+|Z_l||Z^{<l}|\\
&=&|X_i||X_{>i}|+|Y_j||Y^{<j}|+|X_i||Y|+|X_i||Y_j|+|Y||Y_j|+|Y|+|Z_l||Z^{<l}|.
\end{eqnarray*} 
By striking out all twice occuring terms, and using the fact that $|Y_j||Y^{<j}|+|Y_j||Y_{>j}|+|Y_j|=0\mod 2$, it follows that $n_3=l_1+m_3\mod 2$. On the other hand, we have 
\begin{eqnarray*}
C_3&=&\sum_{j=1}^s\sum_{i=1}^r\sum_{l=1}^t (-1)^{k_3}Y^{<j}Y_{>j}[Y_j,[X_i,Z_l]]X^{<i}X_{>i} Z^{<l}Z_{>l}\\
&=&\sum_{j=1}^s\sum_{i=1}^r\sum_{l=1}^t(-1)^{k_3+l_2}X^{<i}X_{>i}[Y_j,[X_i,Z_l]]Y^{<j}Y_{>j} Z^{<l}Z_{>l}.
\end{eqnarray*}
Again we have to check whether $k_3+l_2+|X||Y|+|X|+|Y|+|X_i||Y_j|+|X_i|+|Y_j|=n_3\mod 2$.
\begin{eqnarray*}
l_2&=&(|Y|-|Y_j|+|X|-|X_i|)(|X_i|+|Y_j|+|Z_l|-2)+(|Y|-|Y_j|)(|X|-|X_i|)\\
&=&|X||Y|+|X_i||Y_j|+|X||X_i|+|Y||Y_j|+|X_i|+|Y_j|+(|X|+|Y|+|X_i|+|Y_j|)|Z_l|\\
k_3&=&|Y_j||Y_{>j}|+(|X|-|X_i|)(|X_i|+|Z_l|-1)+|X_i||X_{>i}|+|Z_l||Z^{<l}|\\
&=&|Y_j||Y_{>j}|+|X||X_i|+|X||Z_l|+|X_i||Z_l|+|X|+|X_i||X_{>i}|+|Z_l||Z^{<l}|.
\end{eqnarray*}
The reader may convince himself that the sign works out correctly, and we have thus proved equation (\ref{Jac}).

The next identity we would like to prove is
\begin{eqnarray}\label{A4id}
A_4=(-1)^{(|X|-1)(|Y|-1)}C_5.
\end{eqnarray}
To this end we need to compare
\begin{eqnarray*}
A_4&=&\sum_{i=1}^r \sum_{k=1}^s\sum_{l=1}^t\sum_{j=1}^{l-1}(-1)^{n_4}X^{<i}X_{>i}[X_i,Z_j]Y^{<k}Y_{>k}[Y_k,Z_l]Z^{<j}Z_{>j}^{<l}Z_{>l}\qquad\mbox{with}\\
C_5&=&\sum_{i=1}^s\sum_{k=1}^r\sum_{j=1}^t\sum_{l=j+1}^{t}(-1)^{k_5}Y^{<k}Y_{>k}[Y_k,Z_l]X^{<i}X_{>i}[X_i,Z_j]Z^{<j}Z_{>j}^{<l}Z_{>l}\\
&=&\sum_{i=1}^r \sum_{k=1}^s\sum_{l=1}^t\sum_{j=1}^{l-1}(-1)^{k_5+l_3}X^{<i}X_{>i}[X_i,Z_j]Y^{<k}Y_{>k}[Y_k,Z_l]Z^{<j}Z_{>j}^{<l}Z_{>l}.
\end{eqnarray*}
The exponent $l_3$  works out as follows
\begin{eqnarray*}
l_3&=&(|X|-|X_i|+|X_i|+|Z_j|-1)(|Y|-|Y_k|+|Y_k|+|Z_l|-1)\\
&=&(|X|+|Z_j|-1)(|Y|+|Z_l|-1)\\
&=&|X||Y|+|X||Z_l|+ |Y||Z_j|+ |Z_j||Z_l|+|Z_j|+|Z_l| +|X|+|Y|+1.
\end{eqnarray*}
The reader may convince himself that $n_4=k_5+l_3+(|X|-1)(|Y|-1)\mod 2$
\begin{eqnarray*}
n_4&=&|X_i||X_{>i}|+|Z_j|(|Z^{<j}|+|Z_l|+|Y|+1)+|Y_k||Y_{>k}|+|Z_l||Z^{<l}|\\
k_5&=&|Y_k||Y_{>k}|+|Z_l|(|Z^{<l}|+|X|+1)+|X_i||X_{>i}|+|Z_j||Z^{<j}|
\end{eqnarray*}
and we have proved identity (\ref{A4id}).

Next we claim
\begin{eqnarray}\label{A5id}
A_5=(-1)^{(|X|-1)(|Y|-1)}C_4,
\end{eqnarray}
that is, we have to compare 
\begin{eqnarray*}
A_5&=&\sum_{i=1}^r\sum_{k=1}^s\sum_{l=1}^t\sum_{j=l+1}^{t}(-1)^{n_5}X^{<i}X_{>i}[X_i,Z_j]Y^{<k}Y_{>k}[Y_k,Z_l]Z^{<l}Z_{>l}^{<j}Z_{>j}\qquad\mbox{and}\\
C_4&=&\sum_{k=1}^s \sum_{i=1}^r\sum_{j=1}^t\sum_{l=1}^{j-1}(-1)^{k_4}Y^{<k}Y_{>k}[Y_k,Z_l]X^{<i}X_{>i}[X_i,Z_j]Z^{<l}Z_{>l}^{<j}Z_{>j}\\
&=&\sum_{i=1}^r\sum_{k=1}^s\sum_{l=1}^t\sum_{j=l+1}^{t}(-1)^{k_4+l_4}X^{<i}X_{>i}[X_i,Z_j]Y^{<k}Y_{>k}[Y_k,Z_l]Z^{<l}Z_{>l}^{<j}Z_{>j}.
\end{eqnarray*}
Here, $l_4$ coincides with $l_3$. The reader may convince himself that $n_5=k_4+l_3+(|X|-1)(|Y|-1)\mod 2$ usiing
\begin{eqnarray*}
n_5&=&|X_i||X_{>i}|+|Z_j|(|Z^{<j}|+|Y|+1)+|Y_k||Y_{>k}|+|Z_l||Z^{<l}|\\
k_4&=&|Y_k||Y_{>k}|+|Z_l|(|Z^{<l}|+|Z_j|+|X|+1) +|X_i||X_{>i}|+|Z_j||Z^{<j}|.
\end{eqnarray*}

Now let us check
\begin{eqnarray}\label{B4id}
B_4=-(-1)^{(|X|-1)(|Y|-1)}C_1.
\end{eqnarray}
We have to compare
\begin{eqnarray*}
B_4&=&\sum_{i=1}^r \sum_{j=1}^s\sum_{l=1}^t\sum_{k=i+1}^{r}(-1)^{m_4}X^{<i}X_{>i}^{<k}X_{>k}[X_i,Y_j]Y^{<j}Y_{>j}[X_k,Z_l] Z^{<l}Z_{>l}\qquad\mbox{with}\\
C_1&=&\sum_{j=1}^s\sum_{k=1}^r\sum_{i=1}^{k-1}\sum_{l=1}^t(-1)^{k_1} Y^{<j}Y_{>j}[Y_j,X_i]X^{<i}X_{>i}^{<k}X_{>k}[X_k,Z_l]Z^{<l}Z_{>l}\\
&=&\sum_{i=1}^r \sum_{j=1}^s\sum_{l=1}^t\sum_{k=i+1}^{r}(-1)^{k_1+l_5}X^{<i}X_{>i}^{<k}X_{>k}[X_i,Y_j]Y^{<j}Y_{>j}[X_k,Z_l] Z^{<l}Z_{>l},
\end{eqnarray*}
where 
\begin{eqnarray*}
l_5&=&(|X_i|-1)(|Y_j|-1)+1+(|X|-|X_i|-|X_k|)(|Y|-|Y_j|+|X_i|+|Y_j|-1)\\
&&\qquad\qquad\qquad\qquad\qquad\qquad\qquad\qquad+(|Y|-|Y_j|)(|X_i|+|Y_j|-1)\\
&=&|X||Y|+|X|+|Y|+|X||X_i|+|Y||Y_j|+|X_i||X_k|+|X_k||Y|+|X_k|+|X_i|+|Y_j|
\end{eqnarray*}
Using the identity $|X_i|(|X|+|X_{>i}|+|X^{<i}|+1)=0=|Y_j|(|Y|+|Y_{>j}|+|Y^{<j}|+1)\mod 2$, the reader may check that $m_4=k_1+l_5+|X||Y|+|X|+|Y|\mod 2$ using
\begin{eqnarray*}
m_4&=&|X_i||X_{>i}|+|Y_j||Y^{<j}|+|X_k|(|X_{>k}|+|Y|+|X_i|+1) +|Z_l||Z^{<l}|\\
k_1&=&|Y_j||Y_{>j}|+|X_i||X^{<i}|+|X_k||X_{>k}|+|Z_l||Z^{<l}|.
\end{eqnarray*}
and we have proved equation (\ref{B4id}).

Finally, let us check
\begin{eqnarray}\label{B5id}
B_5=-(-1)^{(|X|-1)(|Y|-1)}C_2.
\end{eqnarray}
To this end we need to compare 
\begin{eqnarray*}
B_5&=&\sum_{i=1}^r\sum_{j=1}^s\sum_{l=1}^t\sum_{k=1}^{i-1}(-1)^{m_5}X^{<k}X^{<i}_{>k} X_{>i}[X_i,Y_j]Y^{<j}Y_{>j}[X_k,Z_l]Z^{<l}Z_{>l}\qquad\mbox{and}\\
C_2&=&\sum_{j=1}^s\sum_{k=1}^r\sum_{i=k+1}^{s}\sum_{l=1}^t(-1)^{k_2} Y^{<j}Y_{>j}[Y_j,X_i]X^{<k}X_{>k}^{<i}X_{>i}[X_k,Z_l]Z^{<l}Z_{>l}\\
&=&\sum_{i=1}^r\sum_{j=1}^s\sum_{l=1}^t\sum_{k=1}^{i-1}(-1)^{k_2+l_6}X^{<k}X^{<i}_{>k} X_{>i}[X_i,Y_j]Y^{<j}Y_{>j}[X_k,Z_l]Z^{<l}Z_{>l}.
\end{eqnarray*}
Again, $l_6$ coincides with $l_5$. The reader may check that $m_5=k_2+l_5+|X||Y|+|X|+|Y|\mod 2$ using
\begin{eqnarray*}
m_5&=&|X_i||X_{>i}|+|Y_j||Y^{<j}|+|X_k|(|X_{>k}|+|Y|+1)+|Z_l||Z^{<l}|\\
k_2&=&|Y_j||Y_{>j}|+|X_i|(|X^{<i}|-|X_k|)+|X_k||X_{>k}|+|Z_l||Z^{<l}|
\end{eqnarray*}
The Jacoibi identity (\ref{Jacobiidentitaet}) now follows from the identities (\ref{A12id})--(\ref{B5id}). 


\end{appendix}
\bibliography{mainpromo,stratbrst,neumaier}
\bibliographystyle{wde_eng}
\end{document}